\documentclass{amsart}
\usepackage{amssymb}
\usepackage{amsmath}
\usepackage{amsfonts}

\newtheorem{theorem}{Theorem}
\theoremstyle{plain}

\newtheorem{corollary}{Corollary}

\newtheorem{example}{Example}

\newtheorem{lemma}{Lemma}

\newtheorem{remark}{Remark}

\numberwithin{equation}{section}
\input{tcilatex}

\begin{document}
\title[Quantum Chaotic States and Stochastic Integration]{Chaotic States and
Stochastic Integration in Quantum Systems }
\author{V\ P\ Belavkin}
\urladdr{}
\thanks{Published in: Uspekhi Mat. Nauk \textbf{47}:1 (1992), 47--106\\
(Russian Math. Surveys \textbf{47}:1 (1992), 53--116\ \\
This translation includes corrections submitted by the author.}
\date{Received 27 May 1991}
\subjclass{}
\keywords{Noncommutative It\^{o} algebras, Quantum L\'{e}vy processes,
Nonadapted stochastic integration, Unified It\^{o} formula, Quantum
stochastic equations}
\dedicatory{}
\thanks{Moscow Institute of Electronics and Mathematics, Moscow 109028, USSR}

\begin{abstract}
Quantum chaotic states over a noncommutative monoid, a unitalization of a
noncommutative Ito algebra parametrizing a quantum stochastic Levy process,
are described in terms of their infinitely divisible generating functionals
over the simple monoid-valued fields on an atomless `space-time' set. A
canonical decomposition of the logarithmic conditionally posive-definite
generating functional is constructed in a pseudo-Euclidean space, given by a
quadruple defining the monoid triangular operator representation and a
cyclic zero pseudo-norm state in this space.

It is shown that the exponential representation in the corresponding
pseudo-Fock space yields the infinitely-divisible generating functional with
respect to the exponential state vector, and its compression to the Fock
space defines the cyclic infinitly-divisible representation associated with
the Fock vacuum state. The structure of states on an arbitrary It\^{o}
algebra is studied with two canonical examples of quantum Wiener and Poisson
states.

A generalized quantum stochastic nonadapted multiple integral is explicitly
defined in Fock scale, its continuity and quantum stochastic
differentiability is proved. A unified non-adapted and functional quantum It%
\^{o} formula is discovered and established both in weak and strong sense,
and the multiplication formula on the exponential It\^{o} algebra is found
for the relatively bounded kernel-operators in Fock scale. The unitarity and
projectivity properties of nonadapted quantum stochastic linear differential
equations are studied, and their solution is constructed for the locally
bounded nonadapted generators in terms of the chronological products in the
underlying kernel algebra canonically represented by triangular operators in
the pseudo-Fock space.
\end{abstract}

\maketitle
\tableofcontents

\section*{Introduction. Non-commutative It\^{o} algebra}

Non-commutative stochastic analysis and calculus appeared in the eighties as
a result of the mathematical justification of the notions of quantum white
noise and the corresponding `Langevin equations' discussed by physicists
from the sixties onwards in connection with stochastic models of quantum
optics and radio-physics \cite{22}, \cite{24}, \cite{34}. The first rigorous
results in quantum stochastic calculus are due to Hudson and Parthasarathy 
\cite{34}, who in 1983 described a quantum It\^{o} formula for
operator-valued integrals with respect to non-commutative canonical
martingales of \emph{creation} $A^{+}(t)$, \emph{annihilation} $A_{-}(t)$,
and \emph{gage} (or vacuum quanta number) $N(t)$. Represented in the
symmetric Fock space $\Gamma (\mathcal{K})$ over $\mathcal{K}=L^{2}(\mathbb{R%
}_{+})$ by noncommuting operators but commuting with their increments at
each $t$, they determine three linear-independent self-adjoint combinations 
\begin{equation}
M_{1}=A_{-}+A^{+},\;\;\;\;M_{2}=\mathrm{i}\left( A_{-}-A^{+}\right) \
,\;\;\;\;\ M_{3}=N,  \label{1}
\end{equation}%
as the `classical' martingales with respect to the vacuum state. Each $%
M_{i}\left( t\right) $ can be represented as a real-valued
independent-increment classical martingale $m_{i}\left( t,\omega _{i}\right) 
$, however due to mutual noncommutativity $\left[ M_{i},M_{k}\right] \neq
0,i\neq k$ they cannot be \emph{jointly represented} as a vector-valued
stochastic process $m_{\bullet }\left( \omega ,t\right) =\left(
m_{1},m_{2},m_{3}\right) \left( \omega ,t\right) $ in any Kolmogorovian
probability space $\left( \Omega ,\mathcal{F},\mathrm{P}\right) $. They are 
\emph{quantum martingales} with respect to the conditional expectations $%
\mathrm{E}_{t}:\mathcal{A}\rightarrow \mathcal{A}_{t}$ on an operator
algebra $\mathcal{A}=\mathcal{A}\left( \Gamma \right) $ of multiple \emph{%
quantum stochastic integrals} $X$ with $\mathcal{A}_{t}=\mathcal{A}\left(
\Gamma _{t}\right) $ corresponding to the natural filtration $\left\{ \Gamma
_{t}=\Gamma (\mathcal{K}_{t}):t\in \mathbb{R}_{+}\right\} $ of the Fock
space defined by the subspaces $\mathcal{K}_{t}\subset \mathcal{K}$ of the
functions with the support in $\left[ 0,t\right] $ and the unit state $%
1_{\emptyset \text{ }}\in $ $\cap _{t>0}\Gamma _{t}$ of the vacuum state $%
\mathrm{E}_{0}[X]=\left\langle 1_{\emptyset }\mid X1_{\emptyset
}\right\rangle $. The triple $(\Gamma ,$ $\mathcal{A},$ $\mathrm{E})$ is
said to be a `quantum probability space' \cite{2}, and in general it
consists of a Hilbert space $\Gamma $, a unital algebra $\mathcal{A}$ of
operators in $\Gamma $ with involution, Hermitian conjugation $X\longmapsto
X^{\ast }$ $\in \mathcal{A}$, and the functional of mathematical expectation 
$\mathrm{E}:\mathcal{A}$ $\rightarrow $ $\mathbb{C}$ defined by the scalar
product $\left\langle 1\mid x\right\rangle $ of a unit vector $1$ $\in $ $%
\Gamma $ and the vector $x=X1$. To any `classical' probability space $\left(
\Omega ,\mathcal{F},P\right) $ there corresponds a canonical `quantum' one
consisting of the Hilbert space $\Gamma =L^{2}(\Omega )$ with the scalar
product 
\begin{equation*}
\left\langle f\mid h\right\rangle =\int f\left( \omega \right) ^{\ast
}h\left( \omega \right) P\left( \mathrm{d}\omega \right) ,
\end{equation*}%
the commutative algebra of bounded `diagonal' operators $\left( Xf\right)
(\omega )=x(\omega )f(\omega )$ given by multiplications by complex
essentially bounded $\mathcal{F}$-measurable random variables $x$ : $\Omega
\rightarrow \mathbb{C}$, and the functional 
\begin{equation}
\mathrm{E}(X)=\int x(\omega )P(\mathrm{d}\omega )=\left\langle 1\mid
x\right\rangle ,  \label{2}
\end{equation}%
defined by the probability vector $1(\omega )=1$ for all $\omega $ $\in $ $%
\Omega $, see for example \cite{32}. The converse is true only in the case
of commutative $C^{\ast }$-algebra $\mathcal{A}$ when all operators have a
joint spectrum $\Omega $ \cite{20}. This proves considerably greater
generality of the non-commutative probability theory, also covering the
purely quantum case which corresponds to a simple, or irreducible algebra,
the algebra $\mathcal{A}=\mathcal{L}\left( \Gamma \right) $ of all linear
continuous operators in a Hilbert space $\Gamma $.

Using this analogy, Hudson and Parthasarathy introduced the notion of \emph{%
adapted} operator-valued process as a family $\left\{ X\left( t\right) :t\in 
\mathbb{R}_{+}\right\} $ of operators in $\Gamma (L^{2}(\mathbb{R}_{+}))$,
each affiliated to the subalgebra $\mathcal{A}_{t}$ generated by the
canonical operators $\left\{ M_{\bullet }(s):s\leq t\right\} $. Due to the
continual tensor-product structure $\Gamma _{t+\Delta }=\Gamma _{t}\otimes
\Gamma _{\Delta }^{t}$ of $\Gamma _{t}$ with $\Gamma _{\Delta }^{t}=\Gamma
\left( \mathcal{K}_{\Delta }^{t}\right) $ for the subspaces $\mathcal{K}%
_{\Delta }^{t}$ of square-integrable functions with the support in $%
[t,t+\Delta t)$ the forward increments $\Delta M_{i}(t)=M_{i}(t+\Delta
t)-M_{i}(t)$ turn out to commute with adapted $D_{i}\left( t\right) $, which
allowed to introduce quantum stochastic integrals $X_{t}=\sum_{i}%
\int_{0}^{t}D_{i}\left( s\right) \mathrm{d}M_{i}(s)$ as the limits of
integral It\^{o} sums $\sum_{t\epsilon \tau }D_{i}\left( t\right) \Delta
M_{i}(t)$, where $\tau =\{t_{1}<\cdots <t_{N}\}$, $\Delta
t_{n}=t_{n+1}-t_{n}\rightarrow 0$ as $N\rightarrow \infty $. Building on
this approach, a quantum evolution was constructed in [26] as a solution of
the linear stochastic differential equation $\mathrm{d}U_{t}=U_{t}L_{j}%
\mathrm{d}\Lambda _{t}^{j}$, $U_{0}=I$, with constant bounded
operator-valued coefficients and non-commutative increments $\mathrm{d}%
\Lambda _{t}^{j}=\mathrm{d}M_{j}\left( t\right) $, $j=1,2,3$, and $\mathrm{d}%
\Lambda _{t}^{0}=\mathrm{d}t$ (Here and it what follows we employ Einstein
summation convention $L_{j}\Lambda ^{j}=\sum_{j\geq 0}L_{j}\Lambda ^{j}$).

The unitarity condition $U_{t}^{\ast }=U_{t}^{-1}$ was studied using the 
\emph{quantum It\^{o} formula} 
\begin{eqnarray}
\mathrm{d}(X_{t}^{\ast }X_{t}) &=&\mathrm{d}X_{t}^{\ast }X_{t}+X_{t}^{\ast }%
\mathrm{d}X_{t}+\mathrm{d}X_{t}^{\ast }\mathrm{d}X_{t},\text{ }  \notag \\
\mathrm{d}X_{t}^{\ast }\mathrm{d}X_{t} &=&D_{i}^{\ast }c_{0}^{ik}D_{k}%
\mathrm{d}t+\sum_{j\geq 1}D_{i}^{\ast }c_{j}^{ik}D_{k}\mathrm{d}M_{j}\left(
t\right) =D_{i}^{\ast }c_{j}^{ik}D_{k}\mathrm{d}\Lambda _{t}^{j},  \label{3}
\end{eqnarray}%
where $c_{j}^{ik}\in \mathbb{C}$ are the structural coefficients defining
the product of quantum-stochastic differentials $\mathrm{d}X_{t}=D_{j}%
\mathrm{d}\Lambda _{t}^{j}$ and $\mathrm{d}X_{t}^{\ast }=D_{j}^{\ast }%
\mathrm{d}\Lambda _{t}^{j}$ corresponding to the Hudson-Parthasarathy (HP)
multiplication table 
\begin{equation*}
\mathrm{d}N\mathrm{d}N=\mathrm{d}N,\text{ \ }\mathrm{d}N\mathrm{d}A^{+}=%
\mathrm{d}A^{+},\text{ \ }\mathrm{d}A_{-}\mathrm{d}N=\mathrm{d}A_{-},\text{
\ }\mathrm{d}A_{-}\mathrm{d}A^{+}=\mathrm{d}t
\end{equation*}%
(other combinations are equal to zero). It follows from this table that $%
c_{j}^{i0}=0=c_{j}^{0k}$ for all $i,j,k=0,1,2,3$ corresponding to the
completely degenerate adjoint representation of $\mathrm{d}\Lambda _{t}^{0}$%
, with $c_{0}^{\bullet 3}=0=c_{0}^{3\bullet }$, $\left( c_{0}^{ik}\right)
^{i,k=1,2}=\left( 
\begin{array}{cc}
1 & -\mathrm{i} \\ 
+\mathrm{i} & 1%
\end{array}%
\right) $ and three Hermitian $3\times 3$-matrices $c_{j}^{\bullet \bullet
}=[c_{j}^{ik}]$, 
\begin{equation*}
\text{ }c_{1}^{\bullet \bullet }=\frac{1}{2}\left[ 
\begin{array}{ccc}
0 & 0 & 1 \\ 
0 & 0 & +\mathrm{i} \\ 
1 & -\mathrm{i} & 0%
\end{array}%
\right] ,\text{ \ }c_{2}^{\bullet \bullet }=\frac{1}{2}\left[ 
\begin{array}{ccc}
0 & 0 & -\mathrm{i} \\ 
0 & 0 & 1 \\ 
+\mathrm{i} & 1 & 0%
\end{array}%
\right] ,\;\;c_{3}^{\bullet \bullet }=\left[ 
\begin{array}{ccc}
0 & 0 & 0 \\ 
0 & 0 & 0 \\ 
0 & 0 & 1%
\end{array}%
\right] ,
\end{equation*}
indexed by $i,k=1,2,3$, define the adjoint representations $c_{\bullet
}^{j\bullet },c_{\bullet }^{\bullet j}$ of the martingale differential
algebra $\mathrm{d}\Lambda ^{j}=\mathrm{d}M_{j}$, $j=1,2,3$.

It can be directly verified that the three-dimensional subspace $\mathfrak{a}%
_{\bullet }$ of complex four-vectors $a_{\bullet }=\left( 0,\alpha _{\bullet
}\right) $, given by the rows $\alpha _{\bullet }=\left( \alpha _{1},\alpha
_{2},\alpha _{3}\right) \in \mathbb{C}^{3}$, is an associative $\#$-algebra
with respect to the complex conjugation $\alpha _{\bullet }^{\#}=\left(
\alpha _{1}^{\ast },\alpha _{2}^{\ast },\alpha _{3}^{\ast }\right) $ as an
involution (not to be mixed up with Hermitian conjugation $\alpha _{\bullet
}^{\ast }=\left[ \alpha _{j}^{\ast }\right] $ defining the adjoint column to 
$\alpha _{\bullet }$) and the composition $\alpha _{\bullet }\beta _{\bullet
}\equiv \alpha _{\bullet }\#\beta _{\bullet }^{\#}$ given by polarization of
the Hermitian 3-vector-form 
\begin{equation*}
\alpha _{\bullet }\#\alpha _{\bullet }=\left( \alpha _{i}c_{1}^{ik}\alpha
_{k}^{\ast },\alpha _{i}c_{2}^{ik}\alpha _{k}^{\ast },\alpha
_{i}c_{3}^{ik}\alpha _{k}^{\ast }\right) .
\end{equation*}%
Moreover, since four fundamental differentials $\mathrm{d}\Lambda ^{j}$ form
an associative algebra,%
\begin{equation*}
\left( \alpha _{\bullet }\mid \beta _{\bullet }\gamma _{\bullet }\right)
=\left( \alpha _{\bullet }\gamma _{\bullet }\mid \beta _{\bullet }\right) 
\end{equation*}%
for any triple $\alpha _{\bullet },\beta _{\bullet },\gamma _{\bullet }\in 
\mathfrak{a}_{\bullet }$ with respect to the semi-positive scalar product%
\begin{equation*}
\left( \alpha _{\bullet }\mid \gamma _{\bullet }\right) =\left( \alpha _{1}+%
\mathrm{i}\alpha _{2}\right) \left( \gamma _{1}+\mathrm{i}\gamma _{2}\right)
^{\ast }=\alpha _{i}c_{0}^{ik}\gamma _{k}^{\ast }.
\end{equation*}%
Thanks to this property one can combine the composition and inner product in 
$\mathfrak{a}_{\bullet }$ into a four-dimensional composition in the $\star $%
-algebra $\mathfrak{a}=\mathbb{C}\oplus \mathfrak{a}_{\bullet }\equiv 
\mathfrak{a}_{\bullet }+\mathbb{C}d_{t}$ of the quadruples $a=\left( \alpha
_{0},\alpha _{\bullet }\right) $ with involution $a^{\star }=\left( \alpha
_{0}^{\ast },\alpha _{\bullet }^{\#}{}\right) $ and self-adjoint nilpotent
element $d_{t}=\left( 1,0^{\bullet }\right) =d_{t}^{\star }$ with respect to
the Hermitian sesquilinear composition%
\begin{equation*}
a\star a=a_{\bullet }\#a_{\bullet }+\left( a_{\bullet }\mid a_{\bullet
}\right) \text{$d_{t}$}=\left( \alpha _{i}c_{j}^{ik}\alpha _{k}^{\ast
}\right) _{j=0,\ldots ,4}
\end{equation*}%
defining the associative multiplication 
\begin{equation}
a\cdot b=\left( \left\langle \alpha _{\bullet },\beta _{\bullet
}\right\rangle ,\alpha _{\bullet }\beta _{\bullet }\right) =a\star b^{\star
}.  \label{4}
\end{equation}%
Here $a_{\bullet }=a-l\left( a\right) d_{t}$ is given be the linear
functional $l\left( a\right) =\alpha _{0}$ for $a=\left( \alpha _{0},\alpha
_{\bullet }\right) \in \mathfrak{a}$, and%
\begin{equation*}
\left\langle a_{\bullet },b_{\bullet }\right\rangle :=\left( \alpha
_{\bullet }\mid \beta _{\bullet }^{\#}\right) \equiv \left\langle \alpha
_{\bullet },\beta _{\bullet }\right\rangle 
\end{equation*}%
is a bilinear form defining the semiscalar product $\left\langle a\mid
b\right\rangle =\left\langle a_{\bullet }^{\star },b_{\bullet }\right\rangle 
$ in $\mathfrak{a}$. We shall call this four-dimensional $\star $-algebra
the Hudson-Prthasarathy quantum It\^{o} algebra (HP-algebra) $\mathfrak{b}%
\left( \mathsf{k}\right) $ of the "Hilbert" space $\mathsf{k}=\mathbb{C}$,
or simply vacuum It\^{o} algebra with respect to the "state" $l$. Note that
this $\mathfrak{a}=\mathfrak{b}\left( \mathbb{C}\right) $ has no identity
but killing nilpotent element $d_{t}$ normalizing the linear functional $l$
as $l\left( \text{$d_{t}$}\right) =1$, which is positive with respect to the
multiplication in the usual sense $l\left( a\star a\right) \geq 0,$%
satisfying $\star $-property $l\left( a^{\star }\right) =l\left( a\right)
^{\ast }$. One can easily see that $\mathbb{C}d_{t}$ is the ideal of $%
\mathfrak{a}$ corresponding to $ad_{t}=0=d_{t}a$ such that $\mathfrak{a}%
_{\bullet }=\left\{ a\in \mathfrak{a}:l\left( a\right) =0\right\} $ is
identified with the factor-algebra $\mathfrak{a}/\mathbb{C}d_{t}$. Moreover,
the two-sided ideal 
\begin{equation*}
\mathfrak{i}=\left\{ b\in \mathfrak{a}:l\left( b\right) =l\left( a\cdot
b\right) =l\left( b\cdot c\right) =l\left( a\cdot b\cdot c\right)
=0,\,\forall a,\,c\in \mathfrak{a}\right\} ,
\end{equation*}%
which, obviously doesn't contain $d_{t}$, is trivial in the vacuum algebra $%
\left( \mathfrak{a},l\right) $: $\mathfrak{i}=\left\{ 0\right\} $.

We take these all properties as the definition of an (abstract
noncommutative) \emph{quantum It\^{o} algebra} $\left( \mathfrak{a}%
,\,l\right) $, and in this capacity we can consider any associative
involutory algebra $\mathfrak{a}=\mathfrak{a}_{\bullet }+\mathbb{C}d_{t}$ by
adding killing element $d_{t}$ to a $\star $-algebra $\mathfrak{a}_{\bullet }
$ equipped with a semi-positive scalar product such that $\mathfrak{i}%
=\left\{ 0\right\} $, taking to $l\left( a\right) =a_{0}$.\ As for $%
\mathfrak{a}_{\bullet }$ one can take any $\star $-algebra with
semi-positive scalar product $\left\langle a_{\bullet }\mid c_{\bullet
}\right\rangle \equiv \left\langle a_{\bullet }^{\star },c_{\bullet
}\right\rangle $ given by a bilinear form satisfying $\left\langle
a_{\bullet },b_{\bullet }c_{\bullet }\right\rangle =\left\langle a_{\bullet
}b_{\bullet },c_{\bullet }\right\rangle $ and factorize it with respect to
the ideal%
\begin{equation*}
\mathfrak{i}=\left\{ b\in \mathfrak{a}_{\bullet }:\left\langle a_{\bullet
},b_{\bullet }\right\rangle =\left\langle b_{\bullet },c_{\bullet
}\right\rangle =\left\langle a_{\bullet },b_{\bullet }c_{\bullet
}\right\rangle =\left\langle a_{\bullet }b_{\bullet },c_{\bullet
}\right\rangle =0,\,\forall a,\,c\in \mathfrak{a}_{\bullet }\right\} 
\end{equation*}%
if $\mathfrak{i}\neq \left\{ 0\right\} $. In this general case one can also
write $\left\langle a_{\bullet },b_{\bullet }\right\rangle =\left( a\cdot
b\right) _{0}=\left\langle a,b\right\rangle $ and implement the $\star $%
-composition notation $a\star b=a\cdot b^{\star }$ which should be
distinguished from the $\#$-composition%
\begin{equation*}
a\#b:=a\star b-l\left( a\star b\right) \text{$d_{t}$}=\left( 0,\alpha
_{\bullet }\#\beta _{\bullet }\right) 
\end{equation*}%
with the values in $\mathfrak{a}_{\bullet }$ representing the composition in
the factor-algebra $\mathfrak{a}/\mathbb{C}d_{t}$. Choosing a selfadjoint
basis $\left\{ e_{j}=e_{j}^{\star }:j=0,\,1,\,\ldots \right\} $ of $\left( 
\mathfrak{a},l\right) $ in such a way that $l\left( a\right) =\alpha _{0}$
if $a=\sum \alpha _{j}e_{j}$, one can describe every finite-dimensional It%
\^{o} algebra as above by the Hermitian structure coefficients 
\begin{equation}
c_{j}^{ik}=\left( c_{j}^{ki}\right) ^{\ast
},\,\;\;\;\;c_{i}^{nj}c_{j}^{km}=c_{j}^{nk}c_{i}^{jm},\;%
\;c_{j}^{0k}=0=c_{j}^{i0},  \label{5}
\end{equation}%
defining a multiplication table $\mathrm{d}\Lambda _{t}^{i}\mathrm{d}\Lambda
_{t}^{k}=c_{j}^{ik}\mathrm{d}\Lambda _{t}^{j}$ of basis quantum stochastic
differentials with $\mathrm{d}\Lambda _{t}^{0}=\mathrm{d}t$.

Note that for the Abelian It\^{o} algebras all structure matrices $%
c_{j}^{\bullet \bullet }$ are real and symmetric with strictly
positive-definite $c_{0}^{\bullet \bullet }$, as it is always so in the case
of one-dimensional $\mathfrak{a}_{\bullet }$ . For example, the standard
Poisson calculus given by It\^{o} multiplication rule%
\begin{equation*}
\mathrm{d}m_{t}\mathrm{d}m_{t}=\lambda \mathrm{d}t+\mathrm{d}m_{t}
\end{equation*}%
for the compensated Poisson increments of the intensity $\lambda $ is
associated with one-dimensional algebra $\mathfrak{a}_{\bullet }\sim \mathbb{%
C}$ of 
\begin{equation*}
a_{\bullet }\sim \alpha ,a_{\bullet }^{\star }\sim \alpha ^{\#},\,a_{\bullet
}\star a_{\bullet }\sim \left\vert \alpha \right\vert ^{2},\;\left\langle
a_{\bullet }\mid a_{\bullet }\right\rangle =\lambda \left\vert \alpha
\right\vert ^{2}
\end{equation*}%
containing the unit $1\in \mathfrak{a}_{\bullet }$ such that $\mathrm{d}%
m_{t}=\mathrm{d}m_{t}^{\ast }$ can be identified with the real two-vector $%
e=\left( 0,1\right) $ in $\mathfrak{a}=\mathfrak{a}_{\bullet }+\mathbb{C}%
d_{t}$ and $\left( \mathrm{d}m_{t}\right) ^{2}$ with $e\star e=e^{2}=\left(
\lambda ,1\right) $. The standard Wiener calculus%
\begin{equation*}
\mathrm{d}w_{t}\mathrm{d}w_{t}=\mathrm{d}t,\,\mathrm{d}w_{t}\mathrm{d}t=%
\mathrm{d}t\mathrm{d}t=\mathrm{d}t\mathrm{d}w_{t}=0
\end{equation*}%
is also associated with one-dimensional but nilpotent algebra $\mathfrak{a}%
_{\bullet }\sim \mathbb{C}$, $\,a_{\bullet }\ast a_{\bullet }\sim 0$ without
unit such that $\mathrm{d}w_{t}=\mathrm{d}w_{t}^{\ast }$ is identified with
the element $e=$ $\left( 0,1\right) =e^{\star }$ of second order nilpotent
algebra $\mathfrak{a}=\mathfrak{a}_{\bullet }+\mathbb{C}d_{t}$ with respect
to the multiplication $e^{2}=\left( 1,0\right) \equiv d_{t}$ defined by the
semi-scalar product $\;\left\langle a\mid a\right\rangle =\left\vert \alpha
\right\vert ^{2}$ for $a=\left( \alpha _{0},\alpha \right) $.

It is well known \cite{31} that the Poisson calculus, as well as the Wiener
one, can be realized as a sub-calculus of the quantum stochastic calculus in
the Fock space with respect to the vacuum state $1_{\emptyset }$ putting,
for example, 
\begin{equation*}
w\left( t\right) =A_{\_}\left( l\right) +A^{+}(t),\,m(t)=\sqrt{\lambda }%
A_{-}\left( t\right) +\sqrt{\lambda }A^{+}\left( t\right) +N\left( t\right) .
\end{equation*}%
A natural question arises as to whether we can realize in this way any
(non-commutative) calculus corresponding to an (abstract quantum) It\^{o}
algebra $(\mathfrak{a},l)$ as defined above. To be more precise, the
question concerns a non-commutative calculus of stochastic integrals with
respect to operator representations of the processes $\Lambda _{t}\left(
a\right) =\alpha _{j}\Lambda _{t}^{j}$ with given expectations $\mathrm{E}%
\left[ \Lambda _{t}\left( a\right) \right] =\alpha _{0}t$, with independent
increments $\mathrm{d}\Lambda _{t}\left( a\right) =\Lambda _{t+\mathrm{d}%
t}\left( \,a\right) -\Lambda _{t}\left( a\right) ,a\in \mathfrak{a}$, and
realizing the multiplication table $\mathrm{d}\Lambda _{t}^{i}\mathrm{d}%
\Lambda _{t}^{k}=\sum_{j\geq 0}c_{j}^{ik}\mathrm{d}\Lambda _{t}^{j}$: 
\begin{equation}
\mathrm{d}\Lambda _{t}\left( a\right) \mathrm{d}\Lambda _{t}\left( a^{\star
}\right) =\alpha _{i}c_{j}^{ik}\alpha _{k}^{\ast }\mathrm{d}\Lambda _{t}^{j}=%
\mathrm{d}\Lambda _{t}\left( \,a\star a\right) .  \label{6}
\end{equation}%
We shall give a positive answer to this question, reducing it to the
construction of canonical representations of infinitely divisible generating
functions 
\begin{equation}
\varphi _{t}\left( b\right) =\mathrm{E}\,[\pi _{t}\left( b\right) ]=\exp
\left\{ tl\,\left( b\right) \right\} ,  \label{7}
\end{equation}%
defined by vacuum expectation of adapted `exponential' operators $\pi
_{t}\left( b\right) \,$\ representing in Fock space a $\star $-monoid $%
\mathfrak{b}$ as a unitalization $b=u+a$ of It\^{o} $\star $-algebra $%
\mathfrak{a}$ with $l$ trivially extended on the unit $u=u^{\star }$ as $%
l\left( u\right) =0$. These representations are constructed as solutions of
quantum stochastic differential equations 
\begin{equation}
\mathrm{d}\pi _{t}\left( b\right) =\pi _{t}\left( b\right) \mathrm{d}\Lambda
_{t}\left( a\right) ,\;\;\;\;\ \ \,\pi _{0}\left( b\right) =I  \label{8}
\end{equation}%
with $a=b-u\in \mathfrak{a}$ such that $b\star b=u+a+a^{\star }+a\star a$.
Note that one can always identify $\mathfrak{b}$ with $\mathfrak{a}$ by
taking $u=0$ such that $b\star b\equiv a+a^{\star }+a\star a\equiv a\bigstar
a$.

In Chapter I we define such functions as solutions $\varphi _{t}\left(
b\right) =\exp \left\{ tl\left( b\right) \right\} $ of the equation 
\begin{equation*}
\mathrm{d}\varphi _{t}\left( b\right) =\varphi _{t}\left( b\right) l\left(
b\right) \mathrm{d}t,\,\;\;\;\varphi _{0}\left( b\right) =1
\end{equation*}%
obtained by the averaging $\mathrm{E}$ of (\ref{8}), taking into account the
independence of the increments $\mathrm{d}\Lambda \left( t,\,a\right) $ and $%
\pi _{t}\left( b\right) $, and that $l\left( b\right) =l\left( a\right) $.

Application of the It\^{o} formula 
\begin{eqnarray*}
\mathrm{d}\left( \pi _{t}\left( b\right) \pi _{t}(b)^{\ast }\right) &=&%
\mathrm{d}\pi _{t}\left( b\right) \mathrm{d}\pi _{t}\left( b\right) ^{\ast }+%
\mathrm{d}\pi _{t}\left( b\right) \pi _{t}\left( b\right) ^{\ast }+\pi
_{t}\left( b\right) \mathrm{d}\pi _{t}\left( b\right) ^{\ast } \\
&=&\pi _{t}\left( b\right) \pi _{t}\left( b\right) ^{\ast }\mathrm{d}\Lambda
_{t}(\,a^{\star }+a\cdot a^{\star }+a)=\pi _{t}\left( b\right) \pi
_{t}\left( b\right) ^{\ast }\mathrm{d}\Lambda _{t}\left( a\bigstar a\right)
\end{eqnarray*}%
gives the multiplication rule $\pi _{t}\left( b\right) ^{\ast }\pi
_{t}\left( b\right) =\pi _{t}\left( b^{\star }b\right) $. Hence we have
positive definiteness $\sum_{a,c}\varphi _{t}\left( a\bigstar c\right)
\lambda _{a}\lambda _{c}^{\ast }\geq 0$ and normalization $\varphi
_{t}\left( 0\right) =1$ of $\varphi _{t}$ defined on $\mathfrak{b}=\mathfrak{%
a}$ as the monoid for each $t$ with respect to this new $\star $-semigroup
composition $\bigstar $ and unit $u=0$. This results from positivity $%
\mathrm{E}\left[ X^{\ast }X\right] \geq 0$ and normalization $\mathrm{E}%
\left[ I\right] =1$ of the vacuum (and any) expectation on the operator
algebra generated by linear combinations $X=\sum \lambda _{b}\pi _{t}\left(
b\right) $. Any such function $\varphi _{t}$ that is included into a
continuous one-parameter semigroup $\left\{ \varphi _{r}:r\in \mathbb{R}%
_{+}\right\} $,%
\begin{equation*}
\varphi _{r}\left( a\right) \varphi _{s}\left( a\right) =\varphi
_{r+s}\left( a\right) ,\;\;\ \;\varphi _{0}\left( a\right) =1
\end{equation*}%
of generating functionals on It\^{o} $\star $-algebra $\mathfrak{a}$ as the
monoid $\mathfrak{b}$ is called infinitely divisible law \cite{16}.

In Chapter 2 we fulfil the It\^{o} programme for quantum stochastic calculus
in a dimension-free form, proving continuity of quantum stochastic integrals
in Fock scales and constructing a noncommutative theory of multiple adapted
and nonadapted quantum stochastic integrals which give solutions to linear
quantum stochastic differential equations in the Wick form of time-ordered
exponentials. We shall use the approach based upon explicit definition of
these integrals in Fock representation, which allows to extend them to
nonadapted operator-functions. We will also obtain a functional quantum It%
\^{o} formula for a quantum stochastic "curve" $X$ with adapted, or even
nonadapted operator values $X_{t}\in \mathcal{A}$, having noncommuting
quantum stochastic increments $\boldsymbol{D}=\left( D_{j}\right) $ with
values in the tensor product $\mathbb{A=}\mathcal{A}\otimes \mathfrak{a}$ of
an operator algebra $\mathcal{A}$ with the It\^{o} algebra $\mathfrak{a}$.
For a `nice' function $f$ the adapted It\^{o} formula with respect to a
filtration $\left( \mathcal{A}_{t}\right) _{t>0}$ generated by an initial
algebra $\mathcal{A}_{0}$ and $\left( \Lambda _{t}^{\bullet }\right) _{t>0}$
can be written for $D_{j}\in \mathcal{A}_{t}$ in the Pseudo-Poisson form 
\cite{18} as%
\begin{equation}
\mathrm{d}f\left( X_{t}\right) =\left( f\left( \mathbf{X}_{t}+\mathbf{D}%
_{t}\right) -f\left( \mathbf{X}_{t}\right) \right) _{j}\mathrm{d}\Lambda
_{t}^{j}.  \label{9}
\end{equation}%
Here $\mathbf{X}$ and $\mathbf{D}$ are canonical images $X\oplus \boldsymbol{%
O}$ and $O\oplus \boldsymbol{D}$ of $X_{t}\in \mathcal{A}_{t}$ and $%
\boldsymbol{D\in }\mathcal{A}_{t}\otimes \mathfrak{a}\equiv \mathbb{A}_{t}$
in the formal sums $\mathbf{X}+\mathbf{D}:=X\oplus \boldsymbol{D}$ as the
elements of the algebra $\mathbb{B}_{t}=\mathcal{A}_{t}\otimes \mathfrak{b}=%
\mathcal{A}_{t}\oplus \mathbb{A}_{t}$ equipped with the involution $\left(
X\oplus \mathbf{D}\right) ^{\dagger }=X^{\ast }\oplus \mathbf{D}^{\star }$
and the product%
\begin{equation}
\left( \mathbf{X}+\mathbf{D}\right) \left( \mathbf{X}+\mathbf{D}\right)
^{\dagger }=XX^{\ast }\oplus \left( X\boldsymbol{D}^{\star }+\boldsymbol{D}%
X^{\ast }+\boldsymbol{D\cdot D}^{\star }\right) ,  \label{10}
\end{equation}%
where $X\mathbf{D}^{\star }=\left( XD_{j}^{\ast }\right) $, $\mathbf{D}%
X^{\ast }=\left( D_{j}X^{\ast }\right) $ and $\boldsymbol{D\cdot D}^{\star
}=\left( D_{i}c_{j}^{ik}D_{k}^{\ast }\right) $. Since $f\left( \mathbf{X}%
\right) =f\left( X\right) \oplus \boldsymbol{O}$, the whole problem is
reduced to computing the operator function $f\left( \mathbf{X}+\mathbf{D}%
\right) $ using the product in $\mathbb{B}_{t}$. Thus, in the case $f\left(
X\right) =X^{m}$, where%
\begin{equation*}
\left( \left( \mathbf{X}+\mathbf{D}\right) ^{m}-\mathbf{X}^{m}\right)
_{j}=D_{j}^{\left( m\right) }
\end{equation*}
with $D_{j}^{\left( 0\right) }=0$, and 
\begin{equation*}
D_{j}^{\left( n+1\right) }=XD_{j}^{\left( n\right)
}+D_{i}c_{j}^{ik}D_{k}^{\left( n\right) }.
\end{equation*}

Note that in the nonstochastic case this new formula also gives an
interesting difference form of the non-commutative chain rule $\mathrm{d}%
f\left( X_{t}\right) =B_{t}\mathrm{d}t$ for a smooth curve $X_{t}$ in an
initial algebra $\mathcal{A}_{0}$ with non-commuting derivative $D_{t}\in 
\mathcal{A}_{0}$. In this case the algebra $\mathfrak{a}_{\bullet }$ is
zero-dimensional, $\mathfrak{a}=\mathbb{C}d_{t}$, and $\mathbb{A}_{0}=%
\mathcal{A}_{0}\otimes d_{t}$ is nilpotent algebra of first order, $A\cdot
A^{\ast }=0$, coinciding as the linear space with $\mathcal{A}_{0}$ such that%
\begin{equation*}
\left( \mathbf{X}+\mathbf{D}\right) ^{\star }\left( \mathbf{X}+\mathbf{D}%
\right) =X^{\ast }X\oplus \left( X^{\ast }D+D^{\ast }X\right) .
\end{equation*}%
In particular, for any polynomial, $f\left( X\right) =X^{m}$ say, one
immediately obtains%
\begin{equation*}
\mathrm{d}X_{t}^{m}=\left( \left( \mathbf{X}_{t}+\mathbf{D}_{t}\right) ^{m}-%
\mathbf{X}_{t}^{m}\right) \mathrm{d}t=%
\sum_{n=1}^{m}X_{t}^{m-n}D_{t}X_{t}^{n-1}\mathrm{d}t
\end{equation*}%
as a particular case of (\ref{9}). Here $X=(X,0),\,D=\left( 0,\,D\right) $
and we took into account that 
\begin{equation*}
\left( \mathbf{X}+\mathbf{D}\right) ^{m}=\mathbf{X}^{m}+\sum_{n=1}^{m}%
\mathbf{X}^{m-n}\mathbf{DX}^{n-1}=\left(
X^{m},\sum_{n=1}^{m}X^{m-n}DX^{n-1}\right) ,
\end{equation*}%
since $\mathbf{DX}^{n}\mathbf{D}=\mathbf{0}$ for $\mathrm{d}\Lambda _{t}^{0}%
\mathrm{d}\Lambda _{t}^{0}=0$ corresponding to $\mathrm{d}\Lambda ^{0}=%
\mathrm{d}t$.

In the nonadapted case the formula (\ref{9}) also remains valid, with $%
\mathbf{X}_{t}=X_{t}\oplus \mathbf{\nabla }X_{t}$ given by quantum
stochastic derivatives $\mathbf{\nabla }_{t}X_{t}=\left( \bigtriangledown
_{t,j}X_{t}\right) \in \mathbb{A}=\mathcal{A}\otimes \mathfrak{a}$, the
noncommutative analog of Malliavin derivative with respect to the canonical
integrators $\mathbf{\Lambda }_{t}=\left( \Lambda _{t}^{j}\right) $. As to
author's knowledge, this general formula is not known even in the classical
(commutative) case.

The author expresses his gratitude to R.L. Hudson, Ya.G. Sinai, and A.S.
Kholevo for discussion on the article and helpful remarks.

\part{Infinitely divisible positive-definite functions and their
representations}

\section{Introduction}

In this paper we study two types of representations associated with a
positive infinitely divisible state on an arbitrary $\star $-semigroup $%
\mathfrak{b}$ \cite{6} with a unit $u\in \mathfrak{b}$. The first,
`differential' type, is connected with an indefinite metric space
representation of conditionally positive functions $\mathfrak{b}\rightarrow 
\mathbb{C}$ in pseudo-Euclidean Minkowski space constructed in \cite{14}. In
the case when $\mathfrak{b}$ is a group, this representation was obtained by
simple generalization \cite{50} of the Gelfand-Naimark-Segal (GNS)
construction from positive definite to conditionally positive definite
functions on $\mathfrak{b}$. However our main interest will be the case when 
$\mathfrak{b}$ is obtained by a unitalization of a noncommutative It\^{o}
algebra $\mathfrak{a}$ as a parametrizing algebra for the quantum stochastic
differentials of a quantum Levy process as operator-valued processes with
independent increments in a quite general noncommutative sense.

In our construction the Hilbert space of the GNS representation is replaced
by a pseudo-Hilbert (Minkowski) space which can be decomposed into a direct
integral sum of a pre-Hilbert space and a one-dimensional complex space in
accordance with the fact that the conditional positiveness (\ref{one d}) has
co-dimension one. In the first section we show that this representation can
be realized by block-triangular matrices of the form 
\begin{equation}
\mathbf{B}=\left[ 
\begin{array}{lll}
1 & b^{-} & \beta \\ 
0 & B & b_{+} \\ 
0 & 0 & 1%
\end{array}%
\right] ,\quad \mathbf{B}^{\dagger }=\left[ 
\begin{array}{lll}
1 & b_{+}^{\ast } & \beta ^{\ast } \\ 
0 & B^{\ast } & b^{-\ast } \\ 
0 & 0 & 1%
\end{array}%
\right]  \label{zero a}
\end{equation}%
with pseudo-Hermitian conjugation $(\boldsymbol{k}\mathbf{B}^{\dagger }\mid 
\boldsymbol{k})=(\boldsymbol{k}\mid \boldsymbol{k}\mathbf{B})$ defined by
the indefinite scalar product 
\begin{equation}
(\boldsymbol{k}\mid \boldsymbol{k}^{\prime })=k_{-}^{\ast }k_{+}^{\prime
}+(k_{\circ }\mid k_{\circ }^{\prime })+k_{+}^{\ast }k_{-}^{\prime },
\label{zero b}
\end{equation}%
on the rows $\boldsymbol{k}=\left( k_{-},k_{\circ },k_{+}\right) $, where $%
k_{+}\in \mathbb{C}\ni k_{-}$ and $k_{\circ }$ is a vector-row from a
complex Euclidean space $\mathcal{K}$. The algebra of the triangular
matrices $\mathbf{A=B-I}$ realizes the non-matrix multiplication table 
\begin{equation}
\left( 
\begin{array}{ll}
\alpha ^{\ast } & a_{+}^{\ast } \\ 
a^{-\ast } & A^{\ast }%
\end{array}%
\right) \cdot \left( 
\begin{array}{ll}
\alpha & a^{-} \\ 
a_{+} & A%
\end{array}%
\right) =\left( 
\begin{array}{ll}
a_{+}^{\ast }a_{+} & a_{+}^{\ast }A \\ 
A^{\ast }a_{+} & A^{\ast }A%
\end{array}%
\right)  \label{zero c}
\end{equation}%
in terms of the $2\times 2$ block-matrices (which are not matrices but
tables)%
\begin{equation*}
\boldsymbol{A}=\left( 
\begin{array}{ll}
\alpha & a^{-} \\ 
a_{+} & A%
\end{array}%
\right) ,\;\boldsymbol{A}^{\ast }=\left( 
\begin{array}{ll}
\alpha ^{\ast } & a_{+}^{\ast } \\ 
a^{-\ast } & A^{\ast }%
\end{array}%
\right) ,
\end{equation*}%
defining the stochastic It\^{o} differentials in the Hudson and
Parthasarathy \cite{26}, \cite{29} quantum calculus. Here $%
a^{-}=b^{-},\,a_{+}=b_{+},\,A=B-I,\,\alpha =\beta $, with involution $%
\mathbf{A}^{\dagger }=\mathbf{B}^{\dagger }-\mathbf{I}$ defined in (\ref%
{zero a}) by the usual Hermitian conjugation $\boldsymbol{A}^{\ast }$ of the
tables $\boldsymbol{A}$ in terms of $A^{\ast }=B^{\ast }-I$ in $\mathcal{K}$%
, where $I$ is the unit operator in $\mathcal{K}$.

This observation, which lays the foundation of a new formulation \cite{7}, 
\cite{9} of quantum stochastic calculus, allows us to extend it to arbitrary
algebras with infinitely divisible state $\varphi $. We mention two
particular algebras of classical stochastic differentials in the case of
one-dimensional $\mathcal{K}=\mathbb{C}$:

\begin{enumerate}
\item the Wiener case: $A=0,\,a^{-}=a_{+}^{\ast },\,\alpha \in \mathbb{C}$,

\item the Poisson case: $A\neq 0,\,a^{-}=a_{+}^{\ast }=0=\alpha $.
\end{enumerate}

If we consider $A$ as the coefficient $A_{\circ }^{\circ }$ at the standard
Poisson differential $\mathrm{d}n=\mathrm{d}\Lambda _{\circ }^{\circ
},\,a^{-}=a_{+}^{\ast }$ as the coefficient $A_{\circ }^{-}=A_{+}^{\circ
\ast }$ at the Wiener standard differential $\mathrm{d}w=\mathrm{d}\Lambda
_{-}^{\circ }+\mathrm{d}\Lambda _{\circ }^{+}$, and $\alpha $ as the
coefficient $A_{+}^{-}$ at $\mathrm{d}t=\mathrm{d}\Lambda _{-}^{+}$, then in
both cases we obtain the realization of the classical It\^{o} formula for
stochastic differential $\mathrm{d}x=\sum_{\mu ,\nu }A_{\nu }^{\mu }\mathrm{d%
}\Lambda _{\mu }^{\nu }\equiv \left\langle \mathbf{A},\mathrm{d}\Lambda
\right\rangle $ in the form 
\begin{equation*}
\mathrm{d}(x^{\ast }x)=x^{\ast }\mathrm{d}x+\mathrm{d}x^{\ast }x+\mathrm{d}%
x^{\ast }\mathrm{d}x=\left\langle x^{\ast }\mathbf{A}+\mathbf{A}^{\mathbf{%
\dagger }}x+\mathbf{A}^{\mathbf{\dagger }}\mathbf{A},\mathrm{d}\Lambda
\right\rangle
\end{equation*}%
of difference multiplication $\mathbf{Y}^{\mathbf{\dagger }}\mathbf{Y}%
-x^{\ast }x\mathbf{I=}x^{\ast }\mathbf{A}+\mathbf{A}^{\mathbf{\dagger }}x+%
\mathbf{A}^{\mathbf{\dagger }}\mathbf{A}$ of the triangular matrices $%
\mathbf{Y}=x\mathbf{I}+\mathbf{A}$, $\mathbf{Y}^{\mathbf{\dagger }}=x^{\ast }%
\mathbf{I}+\mathbf{A}^{\mathbf{\dagger }}$, where $\mathbf{I}$ is the unit $%
3\times 3$ matrix and $\mathbf{A}^{\mathbf{\dagger }}\mathbf{A}$ is defined
by the multiplication table (\ref{zero c}).

In the second section we construct a second `integral' type of
representation of an infinitely divisible chaotic state on $\mathfrak{b}$ by
means of exponential indefinite metric representation and we establish its
relation with the calculus of Maassen-Meyer kernels \cite{35}, \cite{38}, 
\cite{40}, which define chaotic distribution of quantum random variables and
processes.

The algebra of these kernels turns out to be isomorphic to the group algebra
of the exponential representation of $\mathfrak{b}$ in a pseudo-Fock space,
and its Fock projection defines an associated infinitely divisible
representation of $\mathfrak{b}$ generating the corresponding quantum
stochastic calculus in an appropriate Hilbert scale \cite{15}. We note that
this leads in a natural way to the Araki-Woods construction \cite{4}
associated with an infinitely divisible state in the case when $\mathfrak{b}$
is a group.

Finally, in the third section, we study the structure and consider examples
of pseudo-Poisson chaotic states characterized by the linearity of
conditionally positive functions $l\left( b\right) =\ln f\left( b\right) $
on a $\star $-algebra $\mathfrak{b}$. To this type belong the quantum Wiener
states of Heisenberg commutation relations, as well as quantum Poisson
states on noncommutative $C^{\ast }$-algebras $\mathfrak{b}$, studied in 
\cite{10}. Unitary representations connected with infinite divisibility of
states and their applications to the quantum probability theory were studied
in \cite{23}, \cite{27}, \cite{28}, \cite{43}, \cite{49} on groups and in 
\cite{46} on bi-algebras.

\section{Representations of conditionally positive functionals on $\star $%
-semigroups}

Let $(X,\mathfrak{F},\mu )$ be a measurable space $X$ with a $\sigma $%
-algebra $\mathfrak{F}$ and a positive $\sigma $-finite atomless measure $%
\mu :\mathfrak{F}\ni \Delta \mapsto \mu _{\Delta },\,\mu _{\mathrm{d}%
x}\equiv \mathrm{d}x:=\mathrm{d}\mu (x)$ , and let $\mathfrak{b}$ be a
semigroup with involution%
\begin{equation*}
b\mapsto b^{\star },\;\;\;(a\cdot c)^{\star }=c^{\star }\cdot a^{\star },
\end{equation*}%
and with neutral element (unit) $u=u^{\star }$, $u\cdot b=b=b\cdot u$ for
any $b\in \mathfrak{b}$. Typically $\mathfrak{b}$ will be a unitalization of
a noncommutative It\^{o} $\star $-algebra $\mathfrak{a}$, in which case%
\begin{equation*}
a\cdot c=a+c+ac\equiv a\bullet c
\end{equation*}%
if $u$ is identified with zero, or simply write $a\cdot c=ac$ if $\mathfrak{a%
}$ is realized as a $\star $-subalgebra of a unital algebra by taking $u$ $%
=1 $. However in what follows one can take any group with $u=1$ and $%
b^{\star }=b^{-1}$ or any $\star $-submonoid of an operator algebra $%
\mathcal{B}$, a unit ball of a unital C*-algebra say, or even a filter (i.e.
a submonoid) of an idempotent, Boolean say, algebra $\mathfrak{B}$ with
trivial involution $b^{\star }=b$.

Denote $\mathfrak{m}$ the monoid of integrable step-maps $g:X\rightarrow 
\mathfrak{b}$, that is $\mathfrak{b}$-valued functions $x\mapsto g(x)$
having finite images $g(X)=\{g(x):x\in X\}\subseteq \mathfrak{b}$, $%
\left\vert g\left( X\right) \right\vert <\infty $ and integrable co-images $%
\Delta (b)=\{x\in X:g(x)=b\}\in \mathfrak{F}$ in the sense$\,\mu _{\Delta
(b)}<\infty $ for all $b\in \mathfrak{b}$ except $b=u$. We define on $%
\mathfrak{m}$ an inductive structure of a $\star $-monoid with pointwise
defined operations $g^{\star }(x)=g(x)^{\star }$, $(f\cdot h)(x)=f(x)\cdot
h(x)$ and unit $e(x)=u$ for all $x\in X$, considering $\mathfrak{m}$ as the
union $\cup \mathfrak{m}_{\Delta }$ of subsemigroups $\mathfrak{m}_{\Delta
},\,\mu _{\Delta }<\infty $ of step measurable functions $g:X\rightarrow 
\mathfrak{b}$ having integrable supports%
\begin{equation*}
\Delta =\mathrm{supp}g=\{x\in X:g(x)\neq u\}.
\end{equation*}

It is convenient to describe the $\star $-monoid $\mathfrak{b}$ by means of
a single Hermitian operation $a\star c=a\cdot c^{\star }$ satisfying the
relations 
\begin{equation*}
b\star u=b,\quad u\star (u\star b)=b\quad \forall b\in \mathfrak{b}
\end{equation*}%
defining $u=u^{\star }$ as right unit for the composition $\star $, $%
b^{\star }$ as $u\star b$, and 
\begin{equation*}
u\star ((c\star b)\star a)=a\star ((u\star b)\star c)
\end{equation*}%
corresponding to $\left( a\cdot c^{\star }\right) ^{\star }=(a\star
c)^{\star }=c\star a=c\cdot a^{\star }$ and associativity of the semigroup
operation $a\cdot c$. This allows one to define both the product and
involution in a $\star $-monoid $\mathfrak{m}$ by a single Hermitian binary
operation $f\star h=g$,$\,g(x)=f(x)\star h(x)$ with left unit $e\in 
\mathfrak{m}$ which recovers the involution by $g^{\star }(x)=e\star g$ and
the associative product by $f\cdot h=f\star (e\star h)$ for all $f,\,h\in 
\mathfrak{m}$.

Following \cite{6} we say that a generating state functional over the monoid 
$\mathfrak{m}$, or briefly a state over\ $\mathfrak{m}$, is a mapping $%
\varphi :\mathfrak{m}\rightarrow \mathbb{C}$ satisfying the condition $%
\varphi (e)=1$ and positive definiteness 
\begin{equation}
\sum_{f,h\in \mathfrak{m}}\kappa _{f}\varphi (f\star h)\kappa _{h}^{\ast
}\geq 0,\;\;\;\;\,\forall \kappa _{g}\in \mathbb{C}:|\mathrm{supp}\kappa
|<\infty ,  \label{one a}
\end{equation}%
where $|\cdot |$ denotes the cardinality of the set $\mathrm{supp}\kappa
=\{g\in \mathfrak{m}:\kappa _{g}\neq 0\}$.

We introduce on $\mathfrak{m}$ a partial operation $f\sqcup \,h=f\cdot h$
for any functions $f,h\in \mathfrak{m}$ with disjoint supports $\mathrm{supp}%
\,f\cap \mathrm{supp}\,h=\emptyset $ with respect to which the monoid $%
\mathfrak{m}$ turns into a $\star $-semiring in the sense of \cite{6} with
zero $0=e$ and $\Sigma g_{n}=\sqcup g_{n}$ ($\sqcup g_{n}(x)=g_{m}(x)$ for
any $m$ and $x\in \mathrm{supp}\,g_{m}$, otherwise $\sqcup g_{n}(x)=u$). We
call a state $\varphi $ over $\mathfrak{m}$ \textit{chaotic} if 
\begin{equation*}
\varphi \left( \bigsqcup_{n=1}^{\infty }g_{n}\right) =\prod_{n=1}^{\infty
}\varphi (g_{n}),
\end{equation*}%
where $\prod_{n=1}^{\infty }\varphi (g_{n})=\lim_{N\rightarrow \infty
}\prod_{n=1}^{N}\varphi (g_{n})$ for any functions $g_{n}\in \mathfrak{m}$
with pairwise disjoint supports: $\mathrm{supp}\,g_{n}\cap \mathrm{supp}%
\,g_{m}=\emptyset $ for all $n\neq m$.

This condition is fulfilled for $\varphi $ of the exponential form $\varphi
(g)=e^{\lambda \left( g\right) }$ with 
\begin{equation}
\lambda \left( g\right) =\int l(x,g)\mathrm{d}x,\quad l(x,g)=l_{x}(g(x)),
\label{one b}
\end{equation}%
which corresponds to absolute continuity (for all $\Delta \in \mathfrak{F}$
we have $\mu _{\Delta }=0\Rightarrow \lambda _{\Delta }(b)=0$) of the
measure $\lambda _{\Delta }(b):=\langle b_{\Delta }\rangle $ for each $b\in 
\mathfrak{b}$, where $b_{\Delta }(x)=b$ for all $x\in \Delta $ and $%
b_{\Delta }(x)=u$ for $x\notin \Delta $. Here $b_{\Delta }\in \mathfrak{m}$
is an `elementary $\mathfrak{b}$-valued function' called the $b$-indicator
of the subset $\Delta \subseteq X$ if $b\neq u$. In addition the function $%
\varphi _{\Delta }:\mathfrak{b}\rightarrow \mathbb{C}$ given by 
\begin{equation}
\varphi _{\Delta }(b)=\exp \left\{ \int_{\Delta }l_{x}(b)\mathrm{d}x\right\}
=\varphi (b_{\Delta }),  \label{one c}
\end{equation}%
defines an \emph{infinitely divisible state} over the monoid $\mathfrak{b}$
in the sense of the equality $\varphi _{\Delta }(b)=\prod \varphi _{\Delta
_{l}}(b)$ also in the limit of any integral decomposition $\Delta =\Sigma
\Delta _{i}$, $\mu _{\Delta _{i}}\searrow 0$, where $\varphi _{\Delta
_{i}}(b)\rightarrow 1$ for any $b\in \mathfrak{b}$, and in the sense of
positive definiteness of the functions $\varphi _{\Delta
}^{t}(b)=e^{t\lambda _{\Delta }(b)}$ forming a continuous Abelian semigroup 
\begin{equation*}
\{\varphi _{\Delta }^{t}:t\in \mathbb{R}^{+}\},\quad \varphi _{\Delta
}^{0}(b)=1,\quad \left[ \varphi _{\Delta }^{r}\cdot \varphi _{\Delta }^{s}%
\right] (b)=\varphi _{\Delta }^{r+s}(b)
\end{equation*}%
with respect to the pointwise multiplication of $\varphi _{\Delta }^{t}$.
Necessary and sufficient conditions for the function (\ref{one b})
corresponding to the infinitely divisible state (\ref{one c}) are given by
the following theorem, where we assume that $X$ admits a net of
decompositions of the Vitali system in which $\mu _{\Delta }\searrow
0,\,x\in \Delta $, as $\Delta \searrow \{x\}$.

\begin{theorem}
\label{T1} In our notation the following conditions are equivalent:

\begin{enumerate}
\item[\textup{(i)}] For any set $\Delta \in \mathfrak{F}$ of finite measure $%
\mu _{\Delta }<\infty $ the function $\varphi _{\Delta }:b\rightarrow
\varphi \left( b_{\Delta }\right) $ defined by a functional $\varphi :%
\mathfrak{m}\rightarrow \mathbb{C}$ on $b$-indicator $b_{\Delta }$ is a
generating function of an infinitely divisible state over $\mathfrak{b}$,
and for any $b\in \mathfrak{b}$ the limit%
\begin{equation}
l_{x}(b)=\lim_{\Delta \downarrow \{x\}}\frac{1}{\mu _{\Delta }}(\varphi
_{\Delta }(b)-1)  \label{one c'}
\end{equation}%
exists almost everywhere in the Lebesgue-Vitali sense \cite{47}; in addition 
$\mu _{\Delta }=0\Rightarrow \varphi _{\Delta }(b)=1$ for all $\Delta \in
A,\,B\in \mathfrak{b}$.

\item[(ii)] $\varphi (g)=\exp \{\lambda \left( g\right) \}$, where $\lambda
\left( g\right) =\lambda _{\Delta }\left( b\right) $ for $g=b_{\Delta }$ is
an absolutely continuous complex measure of $\Delta \in \mathfrak{F}$, and
for any integrable set $\Delta \subseteq X$ the function $b\mapsto \lambda
_{\Delta }(b)$ is conditionally positive definite 
\begin{equation}
\sum_{a,c\in \mathfrak{b}}\kappa _{a}\lambda _{\Delta }(a\star c)\kappa
_{c}^{\ast }\geq 0,\,\forall \kappa :|\mathrm{supp}\kappa |<\infty
,\,\sum_{b\in \mathfrak{b}}\kappa _{b}=0,  \label{one d}
\end{equation}%
where $\lambda _{\Delta }(u)=0$ and $\lambda _{\Delta }(b^{\star })=\lambda
_{\Delta }(b)^{\ast }$ for any $b\in \mathfrak{b}$.

\item[(iii)] There exists: \newline
1) an integral $\star $-functional $\lambda \left( g\right) :=\int l(x,g)\,%
\mathrm{d}x$ with complex density $l:\mathfrak{m}\rightarrow L^{1}(X)$ such
that $l(g)^{\ast }=l(g^{\star })$ and whose values $l(x,g)=0$ for all $%
g(x)=u $ and $l(x,b_{\Delta })=l_{x}(b)$ with $x\in \Delta $ are independent
of $\Delta $;\linebreak 2) a vector map $\mathrm{k}:g\mapsto \int^{\oplus }%
\mathrm{k}(x,g)\,\mathrm{d}x$ to the subspace $\mathrm{K}\subseteq
\int^{\oplus }\mathrm{K}_{x}\mathrm{d}x$ of square integrable functions $%
\mathrm{k}:x\mapsto \mathrm{k}(x)\in \mathrm{K}_{x}$, $\,\left\Vert \mathrm{k%
}\right\Vert ^{2}=\int \Vert \mathrm{k}(x)\Vert _{x}^{2}\mathrm{d}x<\infty $
with respect to scalar products $\left\langle \mathrm{k}_{x}\mid \mathrm{k}%
_{x}^{\prime }\right\rangle \equiv \mathrm{k}_{x}^{\ast }\mathrm{k}%
_{x}^{\prime }$ in the pre-Hilbert spaces $\mathrm{K}_{x}$, with values $%
\mathrm{k}(x,b_{\Delta })=\mathrm{k}_{x}(b)\in \mathrm{K}_{x}$ independent
of $\Delta \ni x$ and $\mathrm{k}(x,b_{\Delta })=0$ if $x\notin \Delta $
such that $\mathrm{k}(x,g)=0$ if $g(x)=u$; the map $\mathrm{k}$, together
with the adjoint functions $\mathrm{k}^{\star }\left( x,g\right) =\mathrm{k}%
\left( x,g^{\star }\right) ^{\ast }$ as the linear functionals $\mathrm{k}%
^{\star }\left( g\right) =\int^{\oplus }\mathrm{k}^{\star }(x,g)\mathrm{d}%
x\in \mathrm{K}^{\ast }$, satisfies the condition 
\begin{equation}
\mathrm{k}^{\star }(f)\mathrm{k}(h)=\lambda \left( f\cdot h\right) -\lambda
\left( f\right) -\lambda \left( h\right) ,\quad \forall f,h\in \mathfrak{m};
\label{one e}
\end{equation}%
3) a unital $\ast $-representation $j:g\mapsto G:=\int^{\oplus }j(x,g)%
\mathrm{d}x$, $j\left( g\right) ^{\ast }=j\left( g^{\star }\right) $ 
\begin{equation*}
\,j(x,f)j(x,h)=j(x,f\cdot h),\quad j(x,g)=I_{x}\;\forall x:g(x)=u
\end{equation*}%
of a $\star $-semiring $\mathfrak{m}$ in the $\ast $-algebra of decomposable
operators $G:\mathrm{K}\ni \mathrm{k}\mapsto \int^{\oplus }j(x,g)\mathrm{k}%
(x)\mathrm{d}x$ with $j(x,b_{\Delta })=j_{x}(b)$ independent of $\Delta \ni
x $ and $j(x,b_{\Delta })=I_{x}$ if $x\notin \Delta $, which satisfy the
cocycle property 
\begin{equation}
j(g)\mathrm{k}\left( h\right) =\mathrm{k}\left( g\cdot h\right) -\mathrm{k}%
\left( g\right) ,\;\;\,\mathrm{k}^{\star }\left( f\right) j(g)=\mathrm{k}%
^{\star }\left( f\cdot g\right) -\mathrm{k}^{\star }\left( g\right) ,\quad
\forall f,g,h\in \mathfrak{m}  \label{one f}
\end{equation}%
and are continuous in $\mathrm{K}$ with respect to the poly-norm 
\begin{equation}
\Vert \mathrm{k}\Vert ^{h}=\left( \int \Vert j(x,h)\mathrm{k}(x)\Vert
_{x}^{2}\mathrm{d}x\right) ^{1/2},\;\ \ \;h\in \mathfrak{m}.  \label{one g}
\end{equation}

\item[(iv)] For almost all $x\in X$ there exists a pseudo-Hilbert space $%
\mathbb{K}_{x}$, a unital $\dagger $-representation 
\begin{equation*}
\mathbf{j}_{x}(a\cdot b)=\mathbf{j}_{x}(x,b)\mathbf{j}_{x}(x,b),\;\;\mathbf{j%
}_{x}(x,b^{\star })=\mathbf{j}_{x}(x,b)^{\dagger },\mathbf{\;j}_{x}(u)=%
\mathbf{I}_{x}
\end{equation*}%
in the algebra of linear operators $\mathcal{L}(\mathbb{K}_{x})=\{\mathbf{L}:%
\mathbb{K}_{x}\rightarrow \mathbb{K}_{x}:\mathbf{L}^{\dagger }\mathbb{K}_{x}%
\mathbb{\subseteq K}_{x}\}$, where $\mathbf{L}^{\dagger }$ is
pseudo-Hermitian conjugation $\left( \mathbf{k}^{\prime }\mid \mathbf{L}%
^{\dagger }\mathbf{k}\right) =\left( \mathbf{Lk}^{\prime }\mid \mathbf{k}%
\right) $, $\mathbf{k},\mathbf{k}^{\prime }\in \mathbb{K}_{x}$, and a vector 
$\mathbf{e}_{x}\in \mathbb{K}_{x}$ of \ zero pseudo-norm $\left( \mathbf{e}%
_{x}\mid \mathbf{e}_{x}\right) =0$ such that the function 
\begin{equation}
l_{x}(b)=\left( \mathbf{e}_{x}\mid \mathbf{j}_{x}(b^{\star })\mathbf{e}%
_{x}\right) =\left( \mathbf{j}_{x}(b)\mathbf{e}_{x}\mid \mathbf{e}_{x}\right)
\label{one h}
\end{equation}%
is integrable for each $b\in \mathfrak{b}$ on any $\Delta \subseteq X$ with $%
\mu _{\Delta }<\infty $ and $\int_{\Delta }l_{x}(b)\mathrm{d}x=\ln \varphi
_{\Delta }(b)$. Moreover, each $\mathbb{K}_{x}$ can be chosen in the complex
Minkowski form $\mathbb{K}_{x}=\mathbb{C}\oplus \mathcal{K}_{x}^{\circ
}\oplus \mathbb{C}\equiv \mathcal{K}_{x}^{\cdot }$ with $\mathbf{e}_{x}\in 
\mathbb{K}_{x}$ given as pseudo-adjoint $\mathbf{e}_{x}=\boldsymbol{e}%
_{x}^{\dagger }\equiv e^{\cdot }$ to the row $\boldsymbol{e}_{x}=\left(
1,0,0\right) \equiv e_{\cdot }$ of the dual space $\mathbb{K}_{x}^{\dagger }=%
\mathcal{K}_{\cdot }$ of triples $k_{\cdot }=(k_{-},k_{\circ },\,k_{+})$
with the canonical pairing $\left\langle k_{\cdot },h^{\cdot }\right\rangle
_{x}=k_{\iota }h^{\iota }\equiv k_{\cdot }h^{\cdot }$ and an antilinear
embedding $\mathbf{k}\mapsto \mathbf{k}^{\dagger }$ of $\mathbf{k}\equiv
k^{\cdot }\in \mathcal{K}_{x}^{\cdot }$ into $\mathcal{K}_{\cdot }$ such
that $k^{\pm }=k_{\mp }^{\ast }\in \mathbb{C}$,$\;k^{\circ }=k_{\circ
}^{\ast }\in \mathcal{K}_{x}^{\circ }$, defining the Minkowski scalar
product 
\begin{equation}
\left( \mathbf{k}\mid \mathbf{k}\right) _{x}:=k_{-}^{\ast
}k_{+}+\left\langle \mathrm{k}\mid \mathrm{k}\right\rangle _{x}+k_{+}^{\ast
}k_{-}\equiv \left\langle \mathbf{k}^{\dagger },\mathbf{k}\right\rangle \;
\label{one i}
\end{equation}%
on $\mathcal{K}_{x}^{\cdot }$ in terms of the Euclidean scalar product $%
k_{\circ }k^{\circ }=\left\langle \mathrm{k}\mid \mathrm{k}\right\rangle
_{x} $ for $k_{\circ }=\mathrm{k}^{\ast }$ and $k^{\circ }=\mathrm{k}\in 
\mathcal{K}_{x}^{\circ }$. The representation $\mathbf{j}_{x}$ is chosen
then in the triangular form 
\begin{equation*}
j_{\cdot }^{\cdot }(x,b)=%
\begin{bmatrix}
1 & j_{\circ }^{-}(x,b) & j_{+}^{-}(x,b) \\ 
0 & j_{\circ }^{\circ }(x,b) & j_{+}^{\circ }(x,b) \\ 
0 & 0 & 1%
\end{bmatrix}%
,\;\;j_{\cdot }^{\cdot }(x,b^{\star })=%
\begin{bmatrix}
1 & j_{+}^{\circ }(b)^{\ast } & j_{+}^{-}(b)^{\ast } \\ 
0 & j_{\circ }^{\circ }(b)^{\ast } & j_{\circ }^{-}(b)^{\ast } \\ 
0 & 0 & 1%
\end{bmatrix}%
,
\end{equation*}%
defining its dual action $\mathbf{j}_{x}\left( b^{\star }\right) ^{\dagger }$
on $\mathcal{K}_{\cdot }$ as right multiplication by this operator-matrix $%
j_{\cdot }^{\cdot }(x,b)$: 
\begin{equation*}
\mathbf{j}(b):\left( k_{-},k_{\circ },k_{+}\right) \mapsto
(k_{-},\;k_{-}j_{\circ }^{-}(b)+k_{\circ }j_{\circ }^{\circ
}(b),\;k_{-}j_{+}^{-}(b)+k_{\circ }j_{+}^{\circ }(b)+k_{+})\equiv k_{\cdot
}j_{\cdot }^{\cdot }(b).
\end{equation*}%
The Hermitian conjugation $\mathbf{L}^{\dagger }=\mathbf{gL}^{\ast }\mathbf{g%
}$ of the block-matrix operators $\mathbf{L}=[L_{\nu }^{\mu }]$ with respect
to the indefinite form \textup{(\ref{one i}) is given by the metric tensor }$%
\mathbf{g=}\left[ \delta _{-\nu }^{\mu }\right] $ corresponding to the
inversion $-(-,\circ ,+)=(+,\circ ,-)$ of the ordered set $\{-<\circ <+\}$
of the indices $\mu ,\nu =-,\circ ,+$.
\end{enumerate}
\end{theorem}

\begin{proof}
We first establish the simple implications (iv) $\Rightarrow $ (iii) $%
\Rightarrow $ (ii) $\Rightarrow $ (i), and then we prove (i) $\Rightarrow $
(iv) constructing, similarly to the Gelfand-Naimark-Segal construction, a
concrete pseudo-Euclidean representation of the logarithmic derivative of
the generating functional $\varphi _{\Delta }$ of infinitely divisible state
over $\mathfrak{b}$ with respect to $\lambda _{\Delta }$.

(iv) $\Rightarrow $ (iii). If $\mathbf{e}_{x}=\left( 1,\mathrm{e}%
,\varepsilon \right) _{x}^{\dagger }\equiv e_{x}^{\cdot }$ is a zero
pseudo-norm vector-column with the components $e_{x}^{-}=\varepsilon
_{x}^{\ast }\in \mathbb{C}$, $e_{x}^{\circ }=\mathrm{e}_{x}^{\ast }\in 
\mathcal{K}_{x}^{\circ }$ and $e_{x}^{+}=1$ such that $\left\Vert \mathrm{e}%
_{x}\right\Vert ^{2}=2\func{Re}\varepsilon _{x}$, defining (\ref{one h}) in
the triangular matrix representation as $l_{x}(b)=e_{\cdot }^{x}j_{\cdot
}^{\cdot }(x,b)e_{x}^{\cdot }$, where $e_{-}^{x}=1$, $e_{\circ }^{x}=\mathrm{%
e}_{x},e_{+}^{x}=\varepsilon _{x}$, then we can take $\lambda \left(
g\right) =\int l\left( x,g\right) \mathrm{d}x$, where $l\left( x,g\right)
=l_{x}\left( g\left( x\right) \right) $ has obviously properties%
\begin{equation*}
l\left( x,g^{\star }\right) =e_{\cdot }^{x}j_{\cdot }^{\cdot }(x,g)^{\dagger
}e_{x}^{\cdot }=l\left( x,g\right) ^{\ast },
\end{equation*}%
and $l(x,b_{\Delta })=e_{\cdot }^{x}j_{\cdot }^{\cdot }(x,b)e_{x}^{\cdot }$
does not depend on $\Delta $ if $x\in \Delta $, otherwise 
\begin{equation*}
l\left( x,b_{\Delta }\right) =e_{\cdot }^{x}j_{\cdot }^{\cdot
}(x,u)e_{x}^{\cdot }=e_{\cdot }^{x}e_{x}^{\cdot }=\left\Vert \mathrm{e}%
_{x}\right\Vert ^{2}-2\func{Re}\varepsilon _{x}=0.
\end{equation*}

We denote by $\mathrm{K}_{x}$ the completion of the pre-Hilbert space $%
\mathcal{K}_{x}^{\circ }$ by the Cauchy `kets' $\mathrm{k}$ with respect to
the poly-norm%
\begin{equation*}
\Vert \mathrm{k}\Vert _{x}^{\bullet }=\left\{ \Vert \mathrm{k}\Vert
_{x}^{c}=\Vert j_{\circ }^{\circ }(x,c)\mathrm{k}\Vert ,c\in \mathfrak{b}%
\right\}
\end{equation*}%
as fundamental sequences $\{k_{n}^{\circ }\}$ in $\mathcal{K}_{x}^{\circ }$
which do not have limits in $\mathcal{K}_{x}^{\circ }$ simultaneously with
respect to all seminorms $\Vert k^{\circ }\Vert _{x}^{c}$, $c\in \mathfrak{b}
$. For any $g\in \mathfrak{m}$ we denote by $\mathrm{k}\left( g\right) $ the
vector function 
\begin{equation*}
\mathrm{k}\left( x,g\right) =(j_{\circ }^{\circ }(x,g)-1)e_{x}^{\circ
}+j_{+}^{\circ }(x,g)=j_{\mu }^{\circ }(x,g)e_{x}^{\mu }-e_{x}^{\circ }
\end{equation*}%
with values $\mathrm{k}(x,g)\in \mathrm{K}_{x}$, with the adjoint `bras' $%
\mathrm{k}^{\star }(x,g)=e_{\mu }^{x}j_{\circ }^{\mu }(x,g^{\star })\in 
\mathrm{K}_{x}^{\ast }$, where $e_{-\mu }^{x\ast }=e_{x}^{\mu }$ and for
short we use the notation $j_{\nu }^{\mu }\left( x,g\right) =j_{\nu }^{\mu
}\left( x,g\left( x\right) \right) $. This function is square integrable
since%
\begin{equation*}
\mathrm{k}\left( g\right) ^{\ast }\mathrm{k}\left( g\right) =\lambda \left(
g^{\star }\cdot g\right) -\lambda \left( g\right) -\lambda \left( g^{\star
}\right) =\left\Vert \mathrm{k}\left( g\right) \right\Vert ^{2}<\infty
\end{equation*}%
due to the condition (\ref{one e}) which is verified straightforward,%
\begin{align*}
\mathrm{k}^{\star }\left( f\right) \mathrm{k}\left( h\right) & =\int (e_{\mu
}^{x}j_{\circ }^{\mu }(x,f)-e_{\circ }^{x})(j_{\nu }^{\circ }(x,h)e_{x}^{\nu
}-e_{x}^{\circ })\mathrm{d}x \\
& =\int \{\Vert e_{x}^{\circ }\Vert _{x}^{2}+e_{\mu }^{x}\left[ j_{\lambda
}^{\mu }(f)j_{\nu }^{\lambda }(h)-j_{-}^{\mu }(f)j_{\nu }^{-}(h)-j_{+}^{\mu
}(f)j_{\nu }^{+}(h)\right] (x)e_{x}^{\nu } \\
& -e_{\mu }^{x}\left[ j_{\nu }^{\mu }(f)e^{\nu }-j_{-}^{\mu
}(f)e^{-}-j_{+}^{\mu }(f)\right] \left( x\right) \\
& -\left[ e_{\mu }j_{\nu }^{\mu }(h)-e_{-}j_{\nu }^{-}(h)-e_{+}j_{\nu
}^{+}(h)\right] (x)e_{x}^{\nu })\}\mathrm{d}x \\
& =\int [e_{\mu }^{x}j_{\nu }^{\mu }(x,f\cdot h)e_{x}^{\nu }-e_{\mu
}^{x}j_{\nu }^{\mu }(x,h)e_{x}^{\nu }-e_{\mu }^{x}j_{\nu }^{\mu
}(x,f)e_{x}^{\nu }]\mathrm{d}x \\
& =\lambda \left( f\cdot h\right) -\lambda \left( f\right) -\lambda \left(
h\right)
\end{align*}%
for any $f,h\in \mathfrak{m}$, where $e_{\mu }^{x}j_{\nu }^{\mu
}(x,g)e_{x}^{\nu }=l(x,g)$, $\int l(x,g)\mathrm{d}x<\infty $, and we have
employed the condition $e_{\mu }^{x}e_{x}^{\mu }=0$.

Let the subspace $\mathrm{K}\subseteq \dprod \mathrm{K}_{x}$ be chosen as
also the completion of the linear hull of square-integrable functions $%
\left\{ \mathrm{k}\left( g\right) :g\in \mathfrak{m}\right\} $ with respect
to all seminorms $\Vert \mathrm{k}\Vert ^{h}=(\int \Vert j(x,h)\mathrm{k}%
(x)\Vert _{x}^{2}\mathrm{d}x)^{1/2}$, $h\in \mathfrak{m}$ given by
operator-functions $j\left( x,h\right) =j_{\circ }^{\circ }(x,h)$. For any $%
g\in \mathfrak{m}$ we denote by $G=\int^{\oplus }j(x,g)\mathrm{d}x$ a linear
decomposable operator in $\mathrm{K}=\int^{\oplus }\mathrm{K}_{x}\mathrm{d}x$
with $G^{\ast }$ defined pointwise as 
\begin{equation*}
(G^{\ast }\mathrm{k})(x)=j_{\circ }^{\circ }(x,g(x)^{\star })\,\mathrm{k}%
(x)=G(x)^{\ast }\mathrm{k}(x),\,\;\;\;\;\mathrm{k}\in \mathrm{K}.
\end{equation*}%
This definition is correct since for almost all $x\in X$ and all $f,h\in 
\mathfrak{m}$ we have $j(h\cdot f)=j(h)j(f)$ pointwise, and any sequence of
functions $\{\mathrm{k}_{n}\},\,\mathrm{k}_{n}(x)\in \mathrm{K}_{x}$,
fundamental with respect to all seminorms $\Vert \cdot \Vert ^{h}$ is mapped
by the operator $j(g)$ into a sequence $\{\mathrm{k}_{n}^{g}\}$, $\mathrm{k}%
_{n}^{g}(x)=j(x,g)\mathrm{k}_{n}(x)\in \mathrm{K}_{x}$ with the same
fundamental property: 
\begin{equation*}
\Vert \mathrm{k}_{m}^{g}-\mathrm{k}_{n}^{g}\Vert ^{h}=\Vert j(g)j(h)(\mathrm{%
k}_{m}-\mathrm{k}_{n})\Vert =\Vert \mathrm{k}_{m}-\mathrm{k}_{n}\Vert
^{g\cdot h}\longrightarrow 0.
\end{equation*}%
This yields a decomposable non-degenerate representation $G\mathrm{k}%
=\int^{\oplus }G(x)\mathrm{k}(x)\mathrm{d}x$ of the $\star $-semiring $%
\mathfrak{m}$ in the poly-Hilbert space $\mathrm{K}$: 
\begin{equation*}
e\mapsto I=j(e),\;\;\,\,f\star h\mapsto FH^{\ast },\,\,F=j(f),\,\,H=j(h).
\end{equation*}%
This representation is closed in the sense of the completeness of $\mathrm{K}
$ with respect to simultaneous convergence in all seminorms $\Vert \mathrm{k}%
\Vert ^{h}=\Vert H\mathrm{k}\Vert ,\,h\in \mathfrak{m}$ (which is equivalent
to the convergence in the Hilbert norm $\Vert \mathrm{k}\Vert $ only in the
case when the operator function $G(x)=j(x,g)$ is essentially bounded for
every $g\in \mathfrak{m}$, in which case $\mathrm{K}=\int^{\oplus }\mathrm{K}%
_{x}\mathrm{d}x$ is called Hilbert integral). The map $\mathfrak{m}\ni
g\mapsto \mathrm{k}\left( g\right) $ we have constructed, as well as $%
\mathrm{k}^{\star }$, is an additive cocycle in the sense (\ref{one f})
since the derivation property 
\begin{eqnarray*}
\mathrm{k}\left( g\cdot h\right) =j_{\mu }^{\circ }(g\cdot h)e^{\mu
}-e^{\circ } &=&j_{\mu }^{\circ }(g)j_{\nu }^{\mu }(h)e^{\nu }-e^{\circ } \\
&=&j_{\circ }^{\circ }(g)j_{\nu }^{\circ }(h)e^{\nu }+j_{+}^{\circ
}(g)-e^{\circ }=j(g)\mathrm{k}\left( h\right) +\mathrm{k}\left( g\right)
\end{eqnarray*}%
with respect to the representation $j(g)\mathrm{k}\left( h\right) =j_{\circ
}^{\circ }(g)\mathrm{k}\left( h\right) $ of the monoid in $\mathrm{K}$ and
the trivial representation $1(h)=1$ of $\mathfrak{m}$ in $\mathbb{C}$.

(iii) $\Rightarrow $ (ii). It is obvious that the absolutely continuous
measure $\lambda _{\Delta }(b)=\int_{\Delta }l(x,b)\mathrm{d}x$ defined by
the functional $\lambda \left( g\right) =\int e_{\mu }^{x}j_{\nu }^{\mu
}(x,g)e_{\nu }^{x}\mathrm{d}x$ satisfies the conditions $\lambda _{\Delta
}(b^{\star })=\lambda _{\Delta }(b)^{\ast }$ and $\lambda _{\Delta }(u)=0$,
since the functional $l(x,b)$ satisfies these conditions almost everywhere
on $X$. The conditional positivity (\ref{one d}) follows from the positive
definiteness $[\mathrm{k}\left( f\right) ^{\ast }\mathrm{k}\left( h\right)
]\geq 0$ of the scalar product $\mathrm{k}^{\ast }\mathrm{k}^{\prime
}=\left\langle \mathrm{k}\mid \mathrm{k}^{\prime }\right\rangle $ which
guarantees the conditional positivity of the form $\lambda \left( g\right) $%
: 
\begin{align*}
\sum_{f,h\in \mathfrak{m}}\kappa _{f}\langle f\star h\rangle \kappa
_{h}^{\ast }& =\sum_{f,h\in \mathfrak{m}}\kappa _{f}(\lambda \left( f\star
h\right) +\lambda \left( f\right) +\lambda \left( h\right) )\kappa
_{h}^{\ast } \\
& =\sum_{f,h\in \mathfrak{m}}\kappa _{f}\lambda \left( f\star h\right)
\kappa _{h}^{\ast }+\sum_{f\in \mathfrak{m}}\kappa _{f}\sum_{h\in \mathfrak{m%
}}\lambda \left( h\right) ^{\ast }\kappa _{h}^{\ast }+\sum_{f\in \mathfrak{m}%
}\kappa _{f}\lambda \left( f\right) \sum_{h\in \mathfrak{m}}\kappa
_{h}^{\ast } \\
& =\sum_{f,h\in \mathfrak{m}}\kappa _{f}\langle \mathrm{k}\left( f^{\star
}\right) |\mathrm{k}\left( h^{\star }\right) \rangle \kappa _{h}\geq 0
\end{align*}%
for any function $\kappa =\{\kappa _{g}\}$ with finite support and
satisfying $\sum \kappa _{g}=0$.

(ii) $\Rightarrow $ (i). If the function $\lambda _{\Delta }(b)$ is a
(complex) absolutely continuous measure, then $\varphi _{\Delta }(b)=\exp
\{\lambda _{\Delta }(b)\}$ has the property $\varphi _{\sqcup \Delta
_{l}}(b)=\prod \varphi _{\Delta _{l}}(b)$ of infinite divisibility. Moreover
the limit (\ref{one c'}) exists, and by virtue of $\varphi _{\Delta
}(b)\rightarrow 1$ as $\Delta \downarrow \{x\}$ it coincides with the
Radon-Nikod\'{y}m derivative $l_{x}(b)=\mathrm{d}\ln \varphi (b)/\mathrm{d}x$
as the limit of the quotient $\lambda _{\Delta }(b)/\mu _{\Delta }$ over a
net of subsets $\Delta \ni x$ of the system of Vitali decompositions of the
measurable space $X$. For any integrable $\Delta $ the function $b\mapsto
\varphi _{\Delta }(b)$ is positive in the sense of (\ref{one a}). Indeed,
for any complex function $b\mapsto \kappa _{b}$ with finite support we have
due to (\ref{one d}) 
\begin{equation*}
\sum_{a,c\in \mathfrak{b}}\kappa _{a}(\lambda _{\Delta }(a\star c)-\lambda
_{\Delta }(a)-\lambda _{\Delta }(c^{\star }))\kappa _{c}^{\ast }=\sum \kappa
_{a}\langle \mathrm{k}^{\star }\left( a_{\Delta }\right) \mathrm{k}\left(
c_{\Delta }^{\star }\right) \rangle \kappa _{c}^{\ast }\geq 0
\end{equation*}%
since $\langle \mathrm{k}^{\star }\left( a_{\Delta }\right) \mathrm{k}\left(
c_{\Delta }^{\star }\right) \rangle =\sum_{a,c\in \mathfrak{b}}\kappa
_{a}^{\circ }\lambda _{\Delta }(a\star c)\kappa _{c}^{\circ \ast }$ with $%
\kappa _{b}^{\circ }=\kappa _{b},\,b\neq u$, and $\kappa _{\mu }^{\circ
}=\kappa _{\mu }-\sum_{b\in \mathfrak{b}}\kappa _{b}$ is a positive-definite
kernel in $a$ and $c$ as $\sum_{b\in \mathfrak{b}}\kappa _{b}^{\circ }=0$,
and we have taken into account the fact that $\lambda _{\Delta }(u)=0$.
Since the exponent of any positive-definite kernel is a positive definite
kernel, we have for any $\Delta $ 
\begin{equation*}
\sum_{a,c\in \mathfrak{b}}\kappa _{a}^{\ast }\exp \{\lambda _{\Delta
}(a\star c)\}\kappa _{c}=\sum_{a,c\in \mathfrak{b}}\kappa _{\Delta }^{a\ast
}\exp \{\langle \mathrm{k}^{\star }\left( a_{\Delta }\right) \mathrm{k}%
\left( c_{\Delta }^{\star }\right) \rangle \}\kappa _{\Delta }^{c}\geq 0,
\end{equation*}%
where $\kappa _{\Delta }^{b}=\kappa _{b}\exp \{\lambda _{\Delta }(b)\}$ and
we have taken into account (\ref{one e}) and $\lambda _{\Delta }(b^{\star
})=\lambda _{\Delta }(b)^{\ast }$.

(i) $\Rightarrow $ (iv). Since $\varphi _{\Delta }$ is an infinitely
divisible state on $\mathfrak{b}$ and $\varphi _{\Delta }(b)\rightarrow 1$
for all $b$ as $\mu _{\Delta }\rightarrow 0$, the limit $l_{x}(b)$ is
defined as the logarithmic derivative $\mu _{\mathrm{d}x}^{-1}\ln \varphi _{%
\mathrm{d}x}(b)$ of the measure $\lambda _{\Delta }(b)=\ln \varphi _{\Delta
}(b)$ in the Radon-Nikodym sense. Consequently, the function $x\mapsto
l_{x}(b)$ is integrable and almost everywhere satisfies the conditions $%
l_{x}(a\star c)^{\ast }=l_{x}(c\star a)$, $\,l_{x}(u)=0$ and 
\begin{equation*}
\sum_{b\in \mathfrak{b}}\kappa _{b}=0\;\Rightarrow (\kappa ^{\prime }\mid
\kappa )_{x}:=\sum_{a,c\in \mathfrak{b}}\kappa _{a}l_{x}(a\star c)\kappa
_{c}^{\ast }\geq 0
\end{equation*}%
for all $\kappa $ such that $|\mathrm{supp}\,\kappa |<\infty $, which can
easily be verified directly for the difference derivative $l_{\Delta
}(b)=(\varphi _{\Delta }(b)-1)/\mu _{\Delta }$ and next we can pass to the
limit $\Delta \downarrow \{x\}$. In addition $\int_{\Delta }l_{x}(b)\mathrm{d%
}x=\ln \varphi _{\Delta }(b)$ by absolute continuity.

We consider the space $\mathfrak{B}$ of complex functions $\kappa =\left(
\kappa _{b}\right) _{b\in \mathfrak{b}}$ on $\mathfrak{b}$ with finite
supports $\{b\in \mathfrak{b}:\kappa _{b}\neq 0\}$ as a unital $\star $%
-algebra with respect to the product $\kappa ^{\prime }\cdot \kappa $
defined as $\kappa ^{\prime }\star \kappa ^{\star }$ by the Hermitian
convolution%
\begin{equation*}
(\kappa ^{\prime }\star \kappa )_{b}=\sum_{a\star c=b}\kappa _{a}^{\prime
}\kappa _{c}^{\ast },\;\;\;\;\;\,\delta _{u}\star \kappa =\kappa ^{\star
}~,\,\;\;\;\;\kappa \star \delta _{u}=\kappa ~.
\end{equation*}%
with right identity $\delta _{u}$. Here $\delta _{a}=\left( \delta
_{a,b}\right) _{b\in \mathfrak{b}}$ is the Kronecker delta and it defines a $%
\star $-representation

$a\mapsto \delta _{a}$ of the monoid $\mathfrak{b}$ in $\mathfrak{B}$, 
\begin{equation*}
\delta _{a}\star \delta _{c}=\delta _{a\star c},\quad \delta _{u}\star
\delta _{b}=\delta _{b},\quad \delta _{b}\star \delta _{u}=\delta _{b^{\star
}},
\end{equation*}%
with respect to the involution $\kappa ^{\star }=\left( \kappa _{b^{\star
}}^{\ast }\right) _{b\in \mathfrak{b}}$. The linear subspace $\mathfrak{A}%
\subset \mathfrak{B}$ of distributions $\kappa $ such that the sum $\kappa
_{-}:=\sum_{b\in \mathfrak{b}}\kappa _{b}$equals zero, is a $\star $-ideal
since 
\begin{equation*}
\sum_{b\in \mathfrak{b}}(\kappa ^{\prime }\star \kappa )_{b}=\sum_{b\in 
\mathfrak{b}}\sum_{a\star c=b}\kappa _{a}^{\prime }\kappa _{c}^{\ast
}=\sum_{a\in \mathfrak{b}}\kappa _{a}^{\prime }\sum_{c\in \mathfrak{b}%
}\kappa _{c}^{\ast }=0,
\end{equation*}

Let us equip $\mathfrak{B}$ for every $x\in X$ with the\ Hermitian form $%
(\kappa ^{\prime }\mid \kappa )_{x}$ of the kernel $l_{x}\left( a\star
c\right) $ which is positive on $\mathfrak{A}$ and can be written in terms
of the kernel $\left\langle \delta _{a},\delta _{c}^{\star }\right\rangle
_{x}^{\circ }=l_{x}(a\star c)-l_{x}(a)-l_{x}(c^{\star })$ as 
\begin{equation*}
(\kappa ^{\prime }\mid \kappa )_{x}=\kappa _{-}^{\prime }\kappa _{+}^{\ast
}+\left\langle \kappa ^{\prime },\kappa ^{\star }\right\rangle _{x}^{\circ
}+\kappa _{+}^{\prime }\kappa _{-}^{\ast },
\end{equation*}%
where $\kappa _{+}:=\sum_{b}\kappa _{b}l_{x}\left( b\right) $. We notice
that the Hermitian form%
\begin{equation*}
\left\langle \kappa ^{\prime \star }\mid \kappa ^{\star }\right\rangle
_{x}^{\circ }:=\sum_{a.c\in \mathfrak{b}}\kappa _{a}^{\prime }\left\langle
\delta _{a},\delta _{c}\right\rangle _{x}^{\circ }\kappa _{c}^{\star }\equiv
\left\langle \kappa ^{\prime },\kappa ^{\star }\right\rangle _{x}^{\circ }
\end{equation*}%
is non-negative if $\kappa _{-}=0$ or $\kappa _{-}^{\prime }=0$ as $%
\left\langle \kappa ,\kappa ^{\star }\right\rangle _{x}^{\circ }=\sum \kappa
_{a}\left\langle \delta _{a},\delta _{c}^{\star }\right\rangle _{x}^{\circ
}\kappa _{c}^{\ast }\geq 0$, coinciding with $(\kappa ^{\prime }\mid \kappa
)_{x}$. Since $(\kappa ^{\prime }\mid \kappa )_{x}=\sum_{b}\left( \kappa
^{\prime }\star \kappa \right) _{b}l_{x}\left( b\right) $, the form $(\kappa
^{\prime }\mid \kappa )_{x}$ has right associativity property 
\begin{equation*}
(\kappa ^{\prime }\cdot \kappa \mid \kappa )_{x}=(\kappa ^{\prime }\mid
\kappa \star \kappa )_{x}=(\kappa ^{\prime }\mid \kappa \cdot \kappa ^{\star
})_{x},
\end{equation*}%
for all $\kappa ,\kappa ^{\prime }\in \mathfrak{B}$, and therefore its
kernel $\mathfrak{R}_{x}=\left\{ \kappa :(\kappa ^{\prime }\mid \kappa
)_{x}=0\;\forall \kappa ^{\prime }\right\} $ is the right ideal 
\begin{equation*}
\mathfrak{R}_{x}=\{\kappa ^{\prime }\in \mathfrak{B}:(\kappa ^{\prime }\cdot
\kappa \mid \kappa )_{x}=0,\,\forall \kappa \in \mathfrak{B}\}
\end{equation*}%
belonging to $\mathfrak{A}$. We factorize $\mathfrak{B}$ by this right
putting $\kappa \approx 0$ if $\kappa \in \mathfrak{R}_{x}^{\star }:=\left\{
\kappa ^{\star }:\kappa \in \mathfrak{R}_{x}\right\} $ and denoting the
equivalence classes of the left factor-space $\mathcal{K}_{x}^{\cdot }=%
\mathfrak{B}/\mathfrak{R}_{x}^{\star }$ as the ket-vectors $|\kappa
)=\{\kappa ^{\prime }:\kappa ^{\prime }-\kappa ^{\star }\in \mathfrak{R}%
_{x}^{\star }\}$. The condition $\kappa \in \mathfrak{R}_{x}$ means in
particular that $\kappa _{x}^{-}:=(\delta _{u}\mid \kappa )_{x}=0$, and
therefore 
\begin{equation*}
(\kappa \mid \kappa )_{x}=\sum_{a,c\in \mathfrak{b}}\kappa _{a}\langle
\delta _{a},\delta _{c}^{\star }\rangle _{x}\kappa _{c}^{\ast }=\left\langle
\kappa ^{\circ }\mid \kappa ^{\circ }\right\rangle _{x}=0,
\end{equation*}%
where $\kappa ^{\circ }=\left( \kappa _{b}^{\circ }\right) _{b\in \mathfrak{b%
}}$ denotes an element of $\mathfrak{A}$ obtained as $\kappa _{b}^{\circ
}=\kappa _{b}^{\star }$ for all $\,b\neq u$ and $\,\kappa _{u}^{\circ
}=\kappa _{u}^{\star }-\sum_{b\in \mathfrak{b}}\kappa _{b}^{\star }$ such
that $\left\langle \kappa ^{\circ }\mid \kappa ^{\circ }\right\rangle
_{x}=\left\langle \kappa ,\kappa ^{\star }\right\rangle _{x}^{\circ }$.
Therefore it follows also that $\kappa ^{+}:=\sum \kappa _{b}^{\star }$ is
also zero for any $\kappa \in \mathfrak{R}_{x}$ since%
\begin{equation*}
0=(\kappa ^{\prime }\mid \kappa )_{x}=\kappa _{-}^{\prime }\kappa _{+}^{\ast
}+\left\langle \kappa ^{\prime },\kappa ^{\star }\right\rangle _{x}^{\circ
}+\kappa _{+}^{\prime }\kappa _{-}^{\ast }=\kappa _{-}^{\ast }=\kappa ^{+}
\end{equation*}%
for any $\kappa ^{\prime }\in \mathfrak{B}$ with $\kappa _{+}^{\prime }=1$
by virtue of $\kappa _{+}^{\ast }=\kappa ^{-}=0$ and also due to the
Schwartz inequality $(\kappa ^{\prime }\mid \kappa )=\left\langle \kappa
^{\prime },\kappa ^{\star }\right\rangle ^{\circ }=0$. This allows us to
represent the left equivalence classes $|\kappa )_{x}$ by the columns $%
k^{\cdot }=\left[ k^{\mu }\right] $ with $k^{\mp }=\kappa ^{\mp }$ and $%
k^{\circ }=|\kappa ^{\circ }\rangle $ in the Euclidean component $\mathcal{K}%
_{x}^{\circ }\subset \mathcal{K}_{x}^{\cdot }$ as the subspace of the left
equivalence classes $|\kappa ^{\circ }\rangle =|\kappa _{\circ })$ of the
elements $\kappa _{\circ }=\left( \kappa _{b}-\delta _{u,b}\kappa
_{-}\right) _{b\in \mathfrak{b}}\in \mathfrak{A}$ such that $\kappa _{\circ
}^{\star }=\kappa ^{\circ }$. These columns are pseudo-adjoint to the rows $%
k_{\cdot }=(k_{-},k_{\circ },k_{+})$ as the right equivalence classes $%
(\kappa :=|\kappa )^{\dagger }\in \mathfrak{B}/\mathfrak{R}_{x}$ with $%
\,k_{\pm }=\kappa _{\pm }$ and $k_{\circ }=(\kappa _{\circ }$ defining the
indefinite product (\ref{one i}) in terms of the canonical pairing 
\begin{equation*}
k_{\cdot }k^{\cdot }=k_{-}k^{-}+\left\langle k_{\circ },k^{\circ
}\right\rangle +k_{+}k^{+}=\left( k^{\cdot }\mid k_{\cdot }^{\dagger
}\right) ,
\end{equation*}%
where $k^{\circ }=k_{\circ }^{\ast }\in \mathcal{K}_{x}^{\circ }$, $k^{\pm
}=k_{\mp }^{\ast }\in \mathbb{C}$ with respect to the Euclidean scalar
product $\left\langle k_{\circ },k^{\circ }\right\rangle =\left\langle
k_{\circ }^{\ast }\mid k^{\circ }\right\rangle $ of the Euclidean space $%
\mathcal{K}_{x}^{\circ }=\left\{ k^{\circ }=|\kappa ^{\circ }\rangle :\kappa
_{\circ }\in \mathfrak{A}\right\} $, and%
\begin{equation*}
\kappa _{+}^{\ast }=\sum_{b\in \mathfrak{b}}l_{x}(b^{\star })\kappa
^{b}=\kappa _{x}^{-},\;\;\kappa _{-}^{\ast }=\sum_{b\in \mathfrak{b}}\kappa
_{b}^{\ast }=\kappa ^{+}.
\end{equation*}

We notice that the representation $\delta _{\cdot }:\mathfrak{b}\ni b\mapsto
\delta _{b}$ is Hermitian: 
\begin{equation*}
(\kappa \cdot \delta _{b}\mid \kappa )=\sum_{b\in \mathfrak{b}}l(b)(\kappa
\cdot \delta _{b}\star \kappa )_{b}=(\kappa \mid \kappa \cdot \delta
_{b^{\star }}),
\end{equation*}%
and that it is well defined as right representation on $\mathfrak{B}/%
\mathfrak{R}_{x}$ (or left representation on $\mathfrak{B}/\mathfrak{R}%
_{x}^{\star }$) since $\kappa \cdot \delta _{b}\in \mathfrak{R}$ if $\kappa
\in \mathfrak{R}_{x}$: 
\begin{equation*}
(\kappa \mid \kappa )=0,\,\forall \kappa \in \mathfrak{B}\Rightarrow (\kappa
\cdot \delta _{b}\mid \kappa )=(\kappa \mid \kappa \star \delta
_{b})=0,\,\forall \kappa \in \mathfrak{B}.
\end{equation*}%
This allows us to define for each $b\in \mathfrak{b}$ an operator $(\kappa $ 
$\mathbf{j}(b)=(\kappa \cdot \delta _{b}$ such that $\mathbf{j}(b^{\star })=%
\mathbf{j}(b)^{\dagger }$ with the componentwise action 
\begin{align*}
(\kappa \cdot \delta _{b})_{-}& =\kappa _{-},\;\;(\kappa \cdot \delta
_{b})_{\circ }=\kappa _{-}(\delta _{b}-\delta _{u})+\kappa _{\circ }\cdot
\delta _{b}, \\
(\kappa \cdot \delta _{b})_{+}& =\kappa _{-}l(b)+(\kappa _{\circ }\mid
\delta _{b^{\star }}-\delta _{u})+\kappa _{+},
\end{align*}%
given as the right multiplications $\boldsymbol{k}\mapsto \boldsymbol{k}%
\mathbf{B}$, $\boldsymbol{k}\mapsto \boldsymbol{k}\mathbf{B}^{\dagger }$ of
the triangular matrices 
\begin{equation*}
\mathbf{B}=%
\begin{bmatrix}
1 & j_{\circ }^{-}(b) & j_{+}^{-}(b) \\ 
0 & j_{\circ }^{\circ }(b) & j_{+}^{\circ }(b) \\ 
0 & 0 & 1%
\end{bmatrix}%
\equiv j_{\cdot }^{\cdot }\left( b\right) ,\;j_{\cdot }^{\cdot }\left(
b^{\star }\right) =\left[ 
\begin{array}{ccc}
1 & j_{+}^{\circ }(b)^{\ast } & j_{+}^{-}(b)^{\ast } \\ 
0 & j_{\circ }^{\circ }(b)^{\ast } & j_{\circ }^{-}(b)^{\ast } \\ 
0 & 0 & 1%
\end{array}%
\right] \equiv \mathbf{B}^{\dagger }
\end{equation*}%
by the rows $\boldsymbol{k}=(k_{-},k_{\circ },k_{+})\in \mathcal{K}_{\cdot }$
(or as the left multiplications $\mathbf{Bk}$, $\mathbf{B}^{\dagger }\mathbf{%
k}$ by columns $\mathbf{k}\in \mathcal{K}_{x}^{\cdot }$). Here 
\begin{eqnarray*}
j_{+}^{-}(x,b) &=&l_{x}(b),\,\;\;\;\;(\kappa _{\circ }j_{\circ }^{\circ
}(x,b)=(\kappa _{\circ }\cdot \delta _{b}=(\kappa _{\circ }j_{x}\left(
b\right) , \\
\,\;\;j_{+}^{\circ }(x,b^{\star }) &=&\delta _{b}^{\star }\rangle _{x}=%
\mathrm{k}_{x}\left( b^{\star }\right) =\mathrm{k}_{x}^{\star }\left(
b\right) ^{\ast }=\langle \delta _{b}^{\star }|^{\ast }=j_{\circ
}^{-}(x,b)^{\ast },
\end{eqnarray*}%
where $\delta _{b}^{\star }\rangle _{x}=|\delta _{b}-\delta _{u})$ and $%
\,B_{-\nu }^{\dagger \mu }=B_{-\mu }^{\nu \ast }$ is pseudo-Euclidean
conjugation of the triangular matrix $\mathbf{B}=[B_{\nu }^{\mu }]$
corresponding to the map $\boldsymbol{k}\mapsto \boldsymbol{k}^{\dagger }$
into the adjoint columns $\mathbf{k}=\left[ k^{\mu }\right] $ with the
components $\,k^{\mu }=k_{-\mu }^{\ast }$ given by the pseudo-metric tensor $%
g^{\mu \nu }=\delta _{-\nu }^{\mu }=g_{\mu \nu }$: 
\begin{equation*}
\begin{bmatrix}
b_{-}^{-} & b_{\circ }^{-} & b_{+}^{-} \\ 
0 & b_{\circ }^{\circ } & b_{+}^{\circ } \\ 
0 & 0 & b_{+}^{+}%
\end{bmatrix}%
^{\dagger }=%
\begin{bmatrix}
0 & 0 & 1 \\ 
0 & I & 0 \\ 
1 & 0 & 0%
\end{bmatrix}%
\begin{bmatrix}
b_{-}^{-} & b_{\circ }^{-} & b_{+}^{-} \\ 
0 & b_{\circ }^{\circ } & b_{+}^{\circ } \\ 
0 & 0 & b_{+}^{+}%
\end{bmatrix}%
^{\ast }%
\begin{bmatrix}
0 & 0 & 1 \\ 
0 & I & 0 \\ 
1 & 0 & 0%
\end{bmatrix}%
=%
\begin{bmatrix}
b_{+}^{+\ast } & b_{+}^{\circ \ast } & b_{+}^{-\ast } \\ 
0 & b_{\circ }^{\circ \ast } & b_{\circ }^{-\ast } \\ 
0 & 0 & b_{\circ }^{-\ast }%
\end{bmatrix}%
\end{equation*}%
Thus we can write the constructed canonical $\dagger $-representation $%
\mathbf{j}(b)=[j_{\nu }^{\mu }(b)]$ of the monoid $\mathfrak{b}$ in the
pseudo-Euclidean space $\mathbb{K}_{x}=\mathcal{K}_{x}^{\cdot }$ of columns $%
\mathbf{k}=\left[ k^{\mu }\right] $ in terms of the usual matrix
multiplication 
\begin{eqnarray*}
&&%
\begin{bmatrix}
1 & \mathrm{k}^{\star }(a) & l(a) \\ 
0 & j(a) & \mathrm{k}(a) \\ 
0 & 0 & 1%
\end{bmatrix}%
\begin{bmatrix}
1 & \mathrm{k}^{\star }(c) & l(c) \\ 
0 & j(c) & \mathrm{k}(c) \\ 
0 & 0 & 1%
\end{bmatrix}
\\
&=&%
\begin{bmatrix}
1 & \mathrm{k}^{\star }(c)+\mathrm{k}^{\star }(a)j(c), & l(c)+\mathrm{k}%
^{\star }(a)\mathrm{k}(c)+l(a) \\ 
0 & j(a)j(c), & j(a)\mathrm{k}(c)+\mathrm{k}(a) \\ 
0 & 0 & 1%
\end{bmatrix}%
\end{eqnarray*}%
This realizes a conditionally positive function $l(b)$ as the value of the
vector form (\ref{one h}) on the column $\mathbf{e}=\left[ \delta _{+}^{\mu }%
\right] =\boldsymbol{e}^{\dagger }$ as adjoint to row $\boldsymbol{e}=\left(
1,0,0\right) $ of zero pseudonorm $\mathbf{e}^{\dagger }\mathbf{e}=e_{\mu
}e^{\mu }=0$ for each $x$ as%
\begin{equation*}
\mathbf{e}^{\dagger }\mathbf{j}(b)\mathbf{e}=e_{\mu }j_{\nu }^{\mu
}(b)e^{\nu }=j_{+}^{-}(b)=l(b).
\end{equation*}%
The proof is complete.
\end{proof}

\begin{remark}
Any indefinite-metric representation $\left( \mathcal{E}^{\cdot },j_{\cdot
}^{\cdot },c^{\cdot }\right) $ of a conditionally positive function $l$,
written in the form $l(b)=e_{\mu }j_{\nu }^{\mu }(b)e^{\nu }$ with respect
to a triangular $\star $-representation $j_{\cdot }^{\cdot }=\left[ j_{\nu
}^{\mu }\right] $ of $\mathfrak{b}$ in a pseudo-Hilbert space $\mathcal{E}%
^{\cdot }=\mathbb{C}\oplus \mathcal{E}^{\circ }\oplus \mathbb{C}$ with 
\textup{(\ref{one i})} and a zero-vector $e^{\cdot }=(e_{-},e_{\circ
},e_{+})^{\dagger }$ normalized as $e_{-}=1$, $\left\Vert e_{\circ
}\right\Vert ^{2}=-2\func{Re}e_{+}$ can be reduced to the canonical form $%
\left( \mathbb{K},\mathbf{j},\mathbf{e}\right) $ corresponding to 
\begin{equation*}
j_{+}^{-}=l,\;\;j_{+}^{\circ }=k,\;\;j_{\cdot }^{-}=k^{\star },\;\;j_{\circ
}^{\circ }=j
\end{equation*}%
with respect to the vector $\mathbf{e}=(1,0,0)^{\dagger }$ by a triangular
pseudo-isometry $\mathbf{S}:\mathbb{K}\rightarrow \mathcal{E}^{\cdot }$. In
particular, if $(\mathcal{E}^{\cdot },j_{\cdot }^{\cdot },e^{\cdot })$ is a
minimal closed representation in the sense that the vector $e^{\cdot }$ is
cyclic such that $\mathcal{E}^{\circ }$ is minimal poly-Hilbert space
generated by the action on $e^{\cdot }$ of the linear hull of operators $%
j_{\cdot }^{\circ }(\mathfrak{b})$, then it is equivalent to the closed
canonical representation on $\mathbb{K}=\mathbb{C}\oplus \mathrm{K}\oplus 
\mathbb{C}$ with the constructed minimal $\mathcal{K}^{\cdot }=\mathrm{K}$.
\end{remark}

Indeed, taking an arbitrary isometry $U:\mathcal{K}^{\circ }\rightarrow 
\mathcal{E}^{\circ }$ of a minimal space $\mathcal{K}^{\circ }$ we can
define the pseudo-isometry $\mathbf{S}$ in the form 
\begin{equation}
\mathbf{S}=%
\begin{bmatrix}
1, & e_{\circ }U, & e_{+}^{\ast } \\ 
0, & -U, & e_{\circ }^{\ast } \\ 
0, & 0, & 1%
\end{bmatrix}%
,\quad \mathbf{S}^{\dagger }\mathbf{S=I},\;\,\mathbf{S}^{\dagger }=%
\begin{bmatrix}
1, & e_{\circ }, & e_{+} \\ 
0, & -U^{\ast }, & U^{\ast }e_{\circ }^{\ast } \\ 
0, & 0, & 1%
\end{bmatrix}%
,  \label{one j}
\end{equation}%
converting the matrix $j_{\cdot }^{\cdot }(b)$ and the column $e^{\cdot }\in 
\mathcal{E}^{\cdot }$ into the canonical form 
\begin{equation*}
\mathbf{j}(b)=%
\begin{bmatrix}
1 & \mathrm{k}^{\star }(b) & l(b) \\ 
0 & j(b) & \mathrm{k}(b) \\ 
0 & 0 & 1%
\end{bmatrix}%
=\mathbf{S}^{\dagger }j_{\cdot }^{\cdot }(b)\mathbf{S},\;\;\;\;\,\mathbf{e}=%
\begin{bmatrix}
0 \\ 
0 \\ 
1%
\end{bmatrix}%
=\mathbf{S}^{\dagger }e^{\cdot },
\end{equation*}%
since $e_{\cdot }S_{\cdot }^{\cdot }=S_{\cdot }^{-}+e_{\circ }S_{\cdot
}^{\circ }+e_{+}S_{\cdot }^{+}=(1,0,0)$ if $S_{\cdot }^{-}=(1,e_{\circ
}U,e_{+})$, $S_{\cdot }^{\circ }=(0,-U,e_{\circ }^{\ast })$, $S_{\cdot
}^{+}=(0,0,1)$ for $e_{\cdot }=\left( 1,e_{\circ },e_{+}\right) $ with $%
e_{\circ }e_{\circ }^{\ast }\equiv \left\langle e_{\circ }\mid e_{\circ
}\right\rangle =e_{+}+e_{+}^{\ast }$ corresponding to $l(u)=(e_{\cdot }\mid
e_{\cdot })=0$. If the Euclidean space $\mathcal{E}^{\circ }$ is minimal
containing $\{j_{\cdot }^{\circ }(b)e^{\cdot }:b\in \mathfrak{b}\}$ \textup{(%
}or minimal closed with respect to the seminorms $\Vert k_{\circ }\Vert
^{c}=\Vert k_{\circ }j_{\circ }^{\circ }(c)\Vert ,c\in \mathfrak{b}$\textup{)%
} then, defining the operator $U$ by the isometricity condition 
\begin{gather*}
e_{\cdot }j_{\cdot }^{\cdot }(b)S_{\circ }^{\cdot }=(e_{\circ }-e_{\cdot
}j_{\circ }^{\cdot }(b))U=\mathrm{k}(b)^{\ast }, \\
(e_{\circ }-e_{\cdot }j_{\circ }^{\cdot }(a)\mid \,e_{\circ }-e_{\cdot
}j_{\circ }^{\cdot }(c))=\mathrm{k}(a)^{\ast }\mathrm{k}(c),
\end{gather*}%
we obtain a pseudo-unitary equivalence of the \textup{(}closed\textup{)}
representation $(\mathcal{E}^{\cdot },j_{\cdot }^{\cdot },e^{\cdot })$ and
the \textup{(}closed\textup{)} canonical representation $(\mathbb{K},\mathbf{%
j},\mathbf{e})$ constructed in the proof of the implication \textup{(i)} $%
\Rightarrow $ \textup{(iv)} of Theorem~\textup{\ref{T1}}.

\section{The Fock and pseudo-Fock representation of infinitely divisible
states}

We shall now describe an exponential indefinite-metric representation of the 
$\star $-monoid $\mathfrak{m}$ associated with the conditionally
positive-definite functional $\lambda \left( g\right) =\int l_{x}(g(x))%
\mathrm{d}x$ and its connection with the generalized Araki-Woods
construction \cite{4} corresponding to the chaotic infinitely divisible
state $\varphi (g)=e^{\lambda \left( g\right) }$. Unlike the Fock
representation of the Araki-Woods construction, the exponential
representation, which we will construct in a pseudo-Fock space, has the
property of decomposability in finite tensor representations, which can be
used \cite{15} to construct explicit solutions of quantum stochastic
equations even in the case of non-adapted locally integrable generators.

We call pre-Fock space $\mathcal{F}^{\circ }$ over a pre-Hilbert space $%
\mathcal{K}^{\circ }$ the linear hull $\{\mathrm{f}=\Sigma \lambda _{i}\exp
\left\{ \mathrm{k}_{i}\right\} :\lambda _{i}\in \mathbb{C},\mathrm{k}_{i}\in 
\mathcal{K}^{\circ }\}$ of \emph{exponential} vectors $\exp \left\{ \mathrm{k%
}\right\} :=\oplus _{n=0}^{\infty }\frac{1}{n!}\mathrm{k}^{\otimes n}$ as
direct weighted sums of finite tensor powers of vectors $\mathrm{k}\in 
\mathcal{K}^{\circ }$, with $\mathrm{k}^{\otimes 0}=1$ and $\mathrm{k}%
^{\otimes 1}\equiv \mathrm{k}$ such that 
\begin{equation*}
\left\langle \exp \left\{ \mathrm{k}\right\} \mid \exp \left\{ \mathrm{k}%
\right\} \right\rangle =\sum_{n=0}^{\infty }\frac{1}{n!}\left\langle \mathrm{%
k}\mid \mathrm{k}^{\prime }\right\rangle ^{n}=e^{\left\langle \mathrm{k}\mid 
\mathrm{k}^{\prime }\right\rangle }.
\end{equation*}%
This positive-definite exponential kernel describes the scalar product 
\begin{equation*}
\left\langle \mathrm{f}\mid \mathrm{f}^{\prime }\right\rangle =\Sigma
_{ij}\lambda _{i}^{\ast }\left\langle \exp \left\{ \mathrm{k}\right\} \mid
\exp \left\{ \mathrm{k}^{\prime }\right\} \right\rangle \lambda _{j}
\end{equation*}
in $\mathcal{F}^{\circ }$, and the usual Fock space is defined as a
completion of $\mathcal{F}^{\circ }$ with respect to the norm $\left\Vert 
\mathrm{f}\right\Vert =\left\langle \mathrm{f}\mid \mathrm{f}\right\rangle
^{1/2}$. Below as such $\mathcal{K}^{\circ }$ we will take the poly-Hilbert
space $\mathrm{K}=\int^{\oplus }\mathrm{K}_{x}\mathrm{d}x$ associated with
the constructed canonical representation of a conditionally positive
functional $\lambda $ on the $\star $-monoid $\mathfrak{m}$ of simple
functions $g:X\rightarrow \mathfrak{b}$, denoting by $\mathrm{K}_{\ast
}\subseteq \mathcal{K}_{\circ }$ the dual subspace of those functionals $%
k_{\circ }\in \mathcal{K}_{\circ }$ which are represented as $k_{\circ }$%
\textrm{$k$}$=\left\langle \mathrm{k}\mid \mathrm{k}\right\rangle \equiv 
\mathrm{k}^{\ast }$\textrm{$k$} on \textrm{$k$}$\in \mathrm{K}$ (such $%
k_{\circ }=\mathrm{k}^{\ast }$ are continuous with respect to all seminorms
of $\mathrm{K}$).

Thanks to the fact that the measure $\mathrm{d}x$ is atomless, the space $%
\oplus _{n=0}^{\infty }\frac{1}{n!}\mathrm{K}^{\left( n\right) }$, with $%
\mathrm{K}^{\left( n\right) }\subseteq \mathrm{K}^{\otimes n}$ consisting of
only \emph{symmetric} tensor-functions $\mathrm{f}^{\left( n\right)
}:X^{n}\rightarrow \mathrm{K}^{\otimes n}$, can be identified with the space 
$\Gamma (\mathrm{K})=\oplus _{n=0}^{\infty }\Gamma _{n}(\mathrm{K})$ of
square-integrable functions $\mathrm{f}\left( \omega \right) $ with \emph{%
arbitrary} tensor values $\mathrm{f}(\omega )\in \mathrm{K}^{\otimes }\left(
\omega \right) $ \ in \emph{full} tensor products $\mathrm{K}^{\otimes
}\left( \omega \right) =\otimes _{x\in \omega }\mathrm{K}_{x}$ of the
component spaces \thinspace $\mathrm{K}_{x}$, with the indexing sets $\omega
:$ $\left\vert \omega \right\vert =n<\infty $ taken as \emph{any} finite
subsets $\omega \subseteq X$, $\left\vert \omega \right\vert =n<\infty $.
The integrability of such tensor-valued functions $\mathrm{f}(\omega )$,
defined on the space $\Omega =\sum_{n=0}^{\infty }\Omega _{n}$ of all finite
subsets $\omega \subset X$ including empty subset $\omega =\emptyset $ with $%
\mathrm{K}^{\otimes }\left( \emptyset \right) =\mathbb{C}$ is understood
with respect to the Lebesgue measure $\mathrm{d}\omega =\prod_{x\in \omega }%
\mathrm{d}x$ with $\mathrm{d}\emptyset =1$, and the isometry of the
components $\frac{1}{n!}\mathrm{K}^{\left( n\right) }$ and $\int_{\Omega
_{n}}^{\oplus }\mathrm{K}^{\otimes }\left( \omega \right) \mathrm{d}\omega $
is given by the symmetric extension $f^{\left( n\right) }\left( x_{1},\ldots
,x_{n}\right) =f\left( \left\{ x_{1},\ldots ,x_{n}\right\} \right) $
defining $f^{\left( n\right) }$ almost everywhere on $X^{n}$ with respect to 
$\mathrm{d}x_{1}\cdots \mathrm{d}x_{n}$ such that 
\begin{equation*}
\int_{X}\cdots \int_{X}\Vert \mathrm{f}^{\left( n\right) }(x_{1},\ldots
,x_{n})\Vert ^{2}\mathrm{d}x_{1}\cdots \mathrm{d}x_{n}=n!\int_{\Omega
_{n}}\Vert \mathrm{f}(\omega _{n})\Vert ^{2}\mathrm{d}\omega _{n}\text{.}
\end{equation*}%
Denoting the series $\sum_{n=0}^{\infty }\int_{\Omega _{n}}\cdot \mathrm{d}%
\omega _{n}$ of integrals over the $n$-point subsets $\omega _{n}=\left\{
x_{1},\ldots ,x_{n}\right\} \subset X$ as the single integral on $\Omega $
with respect to the measure $\mathrm{d}\omega =\sum $ $\mathrm{d}\omega _{n}$%
, the scalar product in $\Gamma (\mathrm{K})=\int_{\Omega }^{\oplus }\mathrm{%
K}^{\otimes }\left( \omega \right) \mathrm{d}\omega $ can be written as the
single Lebesgue integral on $\Omega ,$%
\begin{equation*}
\left\langle \mathrm{f}\mid \mathrm{f}\right\rangle :=\sum_{n=0}^{\infty
}\int_{\Omega _{n}}\left\langle \mathrm{f}(\omega _{n})\mid \mathrm{f}%
(\omega _{n})\right\rangle \mathrm{d}\omega _{n}\equiv \int_{\Omega
}\left\langle \mathrm{f}(\omega )\mid \mathrm{f}(\omega )\right\rangle 
\mathrm{d}\omega
\end{equation*}%
called the \emph{multiple} Lebesgue integral on $X$. Obviously the symmetric
extension from $\Omega $ onto $\sum_{n=0}^{\infty }\frac{1}{n!}X^{n}$ of the
tensor-product functions \textrm{$f$}$\left( \omega \right) =\otimes _{x\in
\omega }\mathrm{k}\left( x\right) \equiv \mathrm{k}^{\otimes }\left( \omega
\right) $ defines almost everywhere the generating exponential functions $%
\exp \left\{ \mathrm{k}\right\} $, and 
\begin{equation*}
\left\langle \mathrm{k}^{\otimes }\mid \mathrm{k}^{\prime \otimes
}\right\rangle =e^{\left\langle \mathrm{k}\mid \mathrm{k}^{\prime
}\right\rangle },\;\;\;\;\;\left\langle \mathrm{k}\mid \mathrm{k}^{\prime
}\right\rangle =\int \left\langle \mathrm{k}(x)|\mathrm{k}^{\prime
}(x)\right\rangle _{x}\mathrm{d}x.
\end{equation*}%
since $\int_{\Omega _{n}}\Vert \mathrm{k}^{\otimes }(\omega )\Vert ^{2}%
\mathrm{d}\omega _{n}=\frac{1}{n!}\left( \int_{X}\left\Vert \mathrm{k}\left(
x\right) \right\Vert ^{2}\mathrm{d}x\right) ^{n}$ due to $\Vert \mathrm{k}%
^{\otimes }(\omega )\Vert ^{2}=\prod_{x\in \infty }\Vert \mathrm{k}(x)\Vert
_{x}^{2}$. In future we will refer to the Hilbert integral $\int_{\Omega }%
\mathrm{K}^{\otimes }\left( \omega \right) \mathrm{d}\omega $ as to the Fock
space, denoting the exponential domain in it by $\Gamma \left( \mathrm{K}%
\right) \ni \mathrm{k}^{\otimes }$.

We define decomposable operators $j(g)^{\otimes }=\oplus _{n=0}^{\infty
}j(g)^{\otimes n}$ on $\Gamma (\mathrm{K})$ by a unital $\ast $%
-representation $j:\mathfrak{m}\rightarrow \mathcal{L}(\mathrm{K})$ on $%
\mathrm{K}$ associated with the conditionally positive function $\lambda
\left( g\right) $ by means of the linear continuation $j(g)^{\otimes }%
\mathrm{f}=\Sigma \lambda _{i}j(g)^{\otimes }\mathrm{k}_{i}^{\otimes }$ of
the tensor-product operators $j(g)^{\otimes }\mathrm{k}^{\otimes }=(j(g)%
\mathrm{k})^{\otimes }$. Obviously the correspondence $g\mapsto j\left(
g\right) ^{\otimes }$ is, like $j$ itself, a unital $\ast $-representation 
\begin{equation*}
j^{\otimes }(g^{\star })=j^{\otimes }(g)^{\ast },\;\;\;\,j^{\otimes }(f\cdot
h)=j^{\otimes }(f)j^{\otimes }(h),\;\;j^{\otimes }(e)=I^{\otimes }\quad
\forall f,g,h\in \mathfrak{m}
\end{equation*}%
on the pre-Hilbert space $\Gamma \left( \mathrm{K}\right) $. In general the
operators $j(g)^{\otimes }$ are unbounded and cannot be extended onto the
complete Fock space over $\mathrm{K}$ (if only $\star $-monoid is not a
group with $g^{\star }=g^{-1}$), however we can extend them to a closed $%
\ast $-representation $j^{\otimes }\left( g\right) =j(g)^{\otimes }$ on the
completion $\mathrm{F}$ of the pre-Hilbert space $\Gamma (\mathrm{K})$ by
fundamental sequences $\mathrm{f}_{n}\in \Gamma \left( \mathrm{K}\right) $
converging in $\mathrm{F}$ with respect to all seminorms 
\begin{equation*}
\Vert \mathrm{f}\Vert ^{h}=\left( \int \Vert j^{\otimes }(\omega ,h)\mathrm{f%
}(\omega )\Vert ^{2}\mathrm{d}\omega \right) ^{1/2},\quad h\in \mathfrak{m}.
\end{equation*}%
Note that the operators $j^{\otimes }\left( g\right) $ belong to the
operator algebra $\mathcal{L}\left( \mathrm{F}\right) $ of all linear
continuous, together with their adjoints, operators $\mathrm{F}\rightarrow 
\mathrm{F}$ due to $\Vert j^{\otimes }\left( g\right) \mathrm{f}\Vert
^{h}=\Vert \mathrm{f}\Vert ^{g\cdot h}$, and in the case if they all are
bounded, $\mathrm{F}$ is usual Fock space and $\mathcal{L}\left( \mathrm{F}%
\right) =\mathcal{B}\left( \mathrm{F}\right) $. All linear functionals $%
f_{\circ }\in \mathcal{F}_{\circ }$ of the form $f_{\circ }=$\textrm{$f$}$%
^{\ast }\in \mathrm{F}_{\ast }$ are also continuous on $\mathrm{F}$ since $%
f_{\circ }\mathrm{f}=\left\langle \mathrm{f}\mid \mathrm{f}\right\rangle $
converges for any sequence converging in all seminorms $\Vert \mathrm{f}%
\Vert ^{h}$, $h\in \mathfrak{m}$.

Unfortunately, the representation $j^{\otimes }$ describes a dilation of an
infinitely divisible state $\varphi $ as a vector representation on $\mathrm{%
F}$ in the sense of the existence of an $\mathrm{f}\in \mathrm{F}$ such that 
$\varphi (g)=\left\langle \mathrm{f}\mid j^{\otimes }(g)\mathrm{f}%
\right\rangle $ for all $g\in \mathfrak{m}$ only under special `vector'
choice $\lambda \left( g\right) =\left\langle \mathrm{k}\mid (j(g)-I)\mathrm{%
k}\right\rangle $ of the logarithmic function $\lambda \left( g\right) =\ln
\varphi (g)$. If such a vector $\mathrm{k}\in \mathrm{K}$ exists, then one
can obviously take $\mathrm{f}=\exp \{-\left\langle \mathrm{k}\mid \mathrm{k}%
\right\rangle \}\mathrm{k}^{\otimes }$: 
\begin{equation*}
\left\langle \mathrm{f}\mid j^{\otimes }(g)\mathrm{f}\right\rangle =\exp
\{-\left\langle \mathrm{k}\mid \mathrm{k}\right\rangle \}\left\langle 
\mathrm{k}^{\otimes }\mid j^{\otimes }(g)\mathrm{k}^{\otimes }\right\rangle
=\exp \{\left\langle \mathrm{k}\mid (j(g)-I)\mathrm{k}\right\rangle \}.
\end{equation*}

Exploiting a similar exponential construction not for the pre-Hilbert $%
\mathrm{K}$ but for a pseudo-Hilbert extension $\mathbb{K}$ of the complex
Euclidean space $\mathrm{K}$, we can obtain a pseudo-Fock vector
representation also for a general conditionally positive form $\lambda
\left( g\right) $.

For we consider a vector-function space $\mathbb{K}=L^{1}(X)\oplus \mathrm{K}%
\oplus L^{\infty }(X)$ of the triples $\mathbf{k}=k^{-}\oplus k^{\circ
}\oplus k^{+}$, where $k^{\circ }\in \mathrm{K}$ are square integrable
vector-functions $k^{\circ }(x)\in \mathrm{K}_{x}$ from the poly-Hilbert
space $\mathrm{K}=\{\Vert k^{\circ }\Vert ^{h}<\infty :h\in \mathfrak{m}\}$,
with $k^{-}\in L^{1}(X)$ and $k^{+}\in L^{\infty }(X)$ taken as respectively
absolutely integrable and essentially bounded complex functions:%
\begin{equation*}
\Vert k^{-}\Vert _{1}=\int |k^{-}(x)|\mathrm{d}x<\infty ,\;\;\;\Vert
k^{+}\Vert _{\infty }=\mathrm{ess}\sup |k^{+}(x)|<\infty .
\end{equation*}%
We equip this complex poly-Banach space with a pseudo-Hilbert scalar product%
\begin{equation}
\left( \mathbf{k}\mid \mathbf{k}\right) =\left\langle k^{-}\mid
k^{+}\right\rangle +\left\langle k^{\circ }\mid k^{\circ }\right\rangle
+\left\langle k^{+}\mid k^{-}\right\rangle \equiv \left\langle k_{\mu
},k^{\mu }\right\rangle ,  \label{two a}
\end{equation}%
where $k_{-\nu }^{\ast }(x)=k^{\mu }(x)$ such that $\left\langle \,k_{\mu
},k^{\mu }\right\rangle =\int k_{\mu }(x)k^{\mu }(x)\mathrm{d}x$ is the
integral product corresponding to (\ref{one i}) for the column-function $%
\mathbf{k}=k_{\cdot }^{\dagger }\equiv k^{\cdot }$ adjoint to $k_{\cdot
}\left( x\right) =\left( k_{-},k_{\circ },k_{+}\right) \left( x\right) $
with the column \textbf{$k$}$\left( x\right) =\left[ k^{\mu }\left( x\right) %
\right] $ such that $\,k_{-}^{\ast }=k^{+}$, $\,k_{\circ }^{\ast }=k^{\circ
},k_{+}^{\ast }=k^{-}(x)$. Note that the products of the components $k_{\mu
} $ and $k^{\mu }$ with the same $\mu =-,\circ ,+$ are absolutely integrable
for each $\mu $, and thus all three integrals in \ref{two a} converge making 
$\mathbb{K}$ a generalized Krein space.

We define in $\mathbb{K}$ a closed decomposable $\star $-representation $(%
\mathbf{j}(g)$\textbf{$k$}$)(x)=\mathbf{j}(x,g)$\textbf{$k$}$(x)$ of $%
\mathfrak{b}$-valued functions $g(x)$ by triangular-operator functions $%
\mathbf{j}(x,g)=[j_{\nu }^{\mu }(x,g(x))]$ of the canonical form such that 
\begin{equation}
\mathbf{j}(x,g^{\star })=%
\begin{bmatrix}
1, & \mathrm{k}(x,g)^{\ast }, & l(x,g^{\star }) \\ 
0, & j(x,g)^{\ast }, & \mathrm{k}(x,g^{\star }) \\ 
0, & 0, & 1%
\end{bmatrix}%
=\mathbf{j}(x,g)^{\dagger },  \label{two b}
\end{equation}%
where the functions $l(g)\in L^{1}(X)$, $\mathrm{k}(g)\in \mathrm{K}$, $j(g):%
\mathrm{K}\rightarrow \mathrm{K}$ have been described in Theorem~\ref{T1}.

The operators $\mathbf{j}(g)$ are continuous on the whole of $\mathbb{K}$,
together with the adjoint operators $\mathbf{j}(g)^{\dagger }$, with respect
to the Hermitian form (\ref{two a}) by virtue of the inequalities 
\begin{gather*}
\Vert (\mathbf{j}(g)\mathbf{f})^{-}\Vert _{1}\leq \Vert f^{-}\Vert
_{1}+\Vert \mathrm{k}(g)\Vert \cdot \Vert f^{\circ }\Vert +\Vert l(g)\Vert
_{1}\Vert f^{+}\Vert _{\infty }<\infty , \\
\Vert (\mathbf{j}(g)\mathbf{k})^{\circ }\Vert ^{h}\leq \Vert f^{\circ }\Vert
^{g\cdot h}+\Vert \mathrm{k}(g)\Vert ^{h}\cdot \Vert e^{+}\Vert _{\infty
},\quad \Vert (\mathbf{j}(g)\mathbf{k})^{+}\Vert _{\infty }=\Vert f^{+}\Vert
_{\infty }
\end{gather*}%
for any $\mathbf{f}\in \mathbb{K}$, and satisfy conditions (\ref{one e}), (%
\ref{one f}) in the form 
\begin{equation*}
\mathbf{j}(g^{\star })=\mathbf{j}(g)^{\dagger },\;\;\mathbf{j}(f\cdot h)=%
\mathbf{j}(f)\mathbf{j}(h),\quad \mathbf{j}(e)=\mathbf{I},\quad \forall
f,g,h\in \mathfrak{m},
\end{equation*}%
where $\mathbf{I}=[\delta _{\nu }^{\mu }]$ is the unit operator in $\mathbb{K%
}$.

We consider the space $\Gamma (\mathbb{K})$ generated by `exponential'
vectors $\boldsymbol{\mathbf{k}}^{\otimes }=\oplus _{n=1}^{\infty }\mathbf{k}%
^{\otimes n}$ with a non-degenerate pseudo-Hilbert scalar product that
extends to $\Gamma (\mathbb{K})$ the Hermitian form 
\begin{equation}
(\mathbf{k}^{\otimes }\mid \mathbf{k}^{\otimes })=\exp \left\{ \int k_{\mu
}(x)k^{\mu }(x)\mathrm{d}x\right\} =e^{(\mathbf{k}\mid \mathbf{e})}
\label{two c}
\end{equation}%
Owing to the defining algebraic property 
\begin{equation*}
\Gamma (L^{1}(X)\oplus \mathrm{K}\oplus L^{\infty }(X))=\Gamma
(L^{1}(X))\otimes \Gamma (\mathrm{K})\otimes \Gamma (L^{\infty }(X))
\end{equation*}%
of the exponential functor $\Gamma $, we shall write this scalar product as
the triple integral over $\Omega $, 
\begin{equation*}
\left( \mathbf{h}\mid \mathbf{h}\right) =\iiint h(\omega ^{-},\omega ^{\circ
},\omega ^{+})\mathrm{h}(\omega ^{-},\omega ^{\circ },\omega ^{+})\mathrm{d}%
\omega ^{-}\mathrm{d}\omega ^{\circ }\mathrm{d}\omega ^{+}\equiv
\left\langle \boldsymbol{h},\mathbf{h}\right\rangle ,
\end{equation*}%
representing $\mathbf{h}\in \Gamma \left( \mathbb{K}\right) $ by the
ket-function $\mathrm{h}=\left[ \mathrm{h}(\omega ^{-},\omega ^{\circ
},\omega ^{+})\right] $, $\omega ^{\mu }\in \Omega $, and $\boldsymbol{h}=$%
\textbf{$h$}$^{\dagger }$ by the pseudo-adjoint bra-function $h(\omega
_{-},\omega _{\circ },\omega _{+})=$\textrm{$h$}$\left( \omega _{+},\omega
_{\circ },\omega _{-}\right) ^{\ast }$ with values in $\mathrm{K}_{\ast
}^{\otimes }\left( \omega _{\circ }\right) $ such that $h^{\dagger }\left(
\omega ^{-},\omega ^{\circ },\omega ^{+}\right) =$\textrm{$h$}$\left( \omega
^{+},\omega ^{\circ },\omega ^{-}\right) \in \mathrm{K}^{\otimes }\left(
\omega ^{\circ }\right) $. The exponential correspondence $k_{\cdot
}^{\otimes }\mapsto h\left( \omega _{\cdot }\right) $ for each $k_{\cdot
}=\left( k_{-},k_{\circ },k_{+}\right) $ with $k_{\circ }\in \mathrm{K}%
_{\ast }$ is given by 
\begin{equation*}
k_{\cdot }^{\otimes }(\omega _{-},\omega _{\circ },\omega
_{+})=k_{-}^{\otimes }(\omega _{-})k_{\circ }^{\otimes }(\omega _{\circ
})k_{+}^{\otimes }(\omega _{+}),\;\;\;\;\;k_{\circ }^{\otimes }(\omega
)=\otimes _{x\in \omega }k_{\circ }(x),
\end{equation*}%
with $k_{\mp }^{\otimes }(\omega )$ simply described as product functions $%
\prod_{x\in \omega }k_{\mp }(x)$ such that indeed%
\begin{equation*}
\left\langle k_{\cdot }^{\otimes },k_{\cdot }^{\dagger \otimes
}\right\rangle =\exp \left\{ \left\langle k_{-},k^{-}\right\rangle
+\left\langle k_{\circ },k^{\circ }\right\rangle +\left\langle
k_{+},k^{+}\right\rangle \right\} =e^{\left( \mathbf{k}\mid \mathbf{k}%
\right) }.
\end{equation*}%
\ \ 

The Banach space $\mathcal{F}_{\cdot }$ of such tensor-functions $h\left(
\omega _{-},\omega _{\circ },\omega _{+}\right) $ with respect to the norm 
\begin{equation*}
\Vert \mathbf{h}\Vert =\int \mathrm{d}\omega ^{-}\left( \int \mathrm{d}%
\omega ^{\circ }\underset{\omega ^{+}}{\mathrm{ess}\sup }\Vert \mathrm{h}%
(\omega ^{-},\omega ^{\circ },\omega ^{+})\Vert ^{2}\right) ^{1/2}<\infty ,
\end{equation*}%
equipped with the indefinite scalar product (\ref{two c}), will be called a 
\emph{pseudo-Fock space}. It can be easily verified that $\mathcal{F}_{\cdot
}$ contains the exponential vector $\boldsymbol{h}=k_{\cdot }^{\otimes }$ if
and only if $\left\Vert k_{-}\right\Vert _{\infty }\leq 1$, in which case $%
\left\Vert k_{\cdot }^{\otimes }\right\Vert =\exp \left\{ \left\Vert
k_{+}\right\Vert _{1}+\left\Vert k_{\circ }\right\Vert _{2}\right\} $. The
set $\mathcal{K}_{-}^{1}=\left\{ k_{\cdot }\in \mathcal{K}_{\cdot
}:\left\Vert k_{-}\right\Vert _{\infty }\leq 1\right\} $, where $\mathcal{K}%
_{\cdot }^{\dagger }=\mathbb{K}$, contains the canonical vector $e_{\cdot
}^{\otimes }$ given by the constant vector-function $e_{\cdot }(x)=\left(
1,0,0\right) $, and it is invariant under the action $k_{\cdot }\mapsto
k_{\cdot }\mathbf{j}\left( g\right) $ of $\mathfrak{m}$ since $\left(
k_{\cdot }\mathbf{j}\left( g\right) \right) _{-}=k_{-}$ for any $g\in 
\mathfrak{m}$. Therefore the completion of the linear space $\Gamma \left( 
\mathcal{K}_{-}^{1}\right) $ with respect to the Banach poly-norm $\Vert 
\boldsymbol{h}\Vert ^{f}=\left\{ \Vert \boldsymbol{h}\mathbf{j}^{\otimes
}\left( f\right) \Vert :f\in \mathfrak{m}\right\} $ is a dense subspace of $%
\mathcal{F}^{\cdot }$ which will be denoted by $\mathbb{F}_{\star }$, with $%
\mathbb{F}$ for $\mathbb{F}_{\star }^{\dagger }=\left\{ \boldsymbol{h}%
^{\dagger }:\boldsymbol{h}\in \mathbb{F}_{\star }\right\} $. ($\mathbb{F}%
_{\star }$ coincides with $\mathcal{F}_{\cdot }$ if all $j\left( g\right) $
are bounded).

The canonical exponential vector is obviously state vector for the
infinitely divisible state $\varphi (g)$ which is represented as 
\begin{equation*}
\varphi (g)=(\boldsymbol{e}^{\otimes }\mathbf{j}^{\otimes }(g)\mid 
\boldsymbol{e}^{\otimes })=\exp \{(\boldsymbol{e}\mathbf{j}(g)\mid 
\boldsymbol{e})=\mathrm{e}^{\lambda \left( g\right) }.
\end{equation*}%
What is more, as the next theorem shows, the representation $\mathbf{j}%
^{\otimes }$, compressed to the Fock subspace $\mathrm{F}\subset \mathbb{F}$
by means of a pseudo-conditional expectation 
\begin{equation*}
\epsilon \lbrack \mathbf{j}^{\otimes }(g)]:=J^{\dagger }\mathbf{j}^{\otimes
}(g)J\equiv \pi (g),
\end{equation*}%
remains multiplicative, with $J^{\dagger }\mathbf{e}^{\otimes }=1_{\emptyset
}$ defined as the vacuum state the unital $\ast $-representation $\pi $
associated with $\varphi (g)=\exp \lambda \left( g\right) $. Here $\delta
_{\emptyset }^{\omega }=1$ for $\omega =\emptyset $, $\delta _{\emptyset
}^{\omega }=0$ for $\omega \neq \emptyset $, and 
\begin{equation}
(J\mathrm{h})(\omega ^{-},\omega ^{\circ },\omega ^{+})=1_{\emptyset
}(\omega ^{-})\mathrm{h}(\omega ^{\circ }),\quad \mathrm{h}\in \mathrm{F},
\label{two d}
\end{equation}%
is a pseudo-isometry $\mathrm{F}\rightarrow \mathbb{F}$, $\;(J\mathrm{h}\mid
J\mathrm{h})=\left\langle \mathrm{h}\mid \mathrm{h}\right\rangle $ for all $%
\mathrm{h}\in \mathrm{F}$.

To obtain this result we note that any decomposable operator $\mathbf{K}%
=1\oplus \mathbf{G}\oplus \mathbf{G}^{\otimes 2}\oplus \ldots $ in $\Gamma (%
\mathbb{K})$, obtained by exponentiation $\mathbf{G}^{\otimes }$ of the
triangular operator $\mathbf{G}=\mathbf{j}(g)$, can be written in the form 
\begin{equation}
\lbrack \mathbf{Kh}](\omega ^{-},\omega ^{\circ },\omega
^{+})=\sum_{\bigsqcup\limits_{\nu \geq \mu }\omega _{\nu }^{\mu }=\omega
^{\mu }}^{\mu =-,\circ ,+}K(\mathbf{\omega })\mathrm{h}(\omega
_{-}^{-},\omega _{\circ }^{-}\sqcup \omega _{\circ }^{\circ },\omega
_{+}^{-}\sqcup \omega _{+}^{\circ }\sqcup \omega _{+}^{+}),  \label{two e}
\end{equation}%
where $\omega =\sqcup \omega _{\nu }$ denotes direct sum of pairwise
disjoint $\omega _{\nu }$ defining a decomposition $\omega _{\cdot }=\left(
\omega _{\nu }\right) $ of $\omega $. Here $K(\mathbf{\omega })$ is a
function of the table $\mathbf{\omega }=(\omega _{\nu }^{\mu })_{\nu =\circ
,+}^{\mu =-,\circ }$ of four subsets $\omega _{\nu }^{\mu }\in \Omega $ with
values in linear continuous operators 
\begin{gather}
K%
\begin{pmatrix}
\omega _{+}^{-} & \omega _{\circ }^{-} \\ 
\omega _{+}^{\circ } & \omega _{\circ }^{\circ }%
\end{pmatrix}%
:\mathrm{K}^{\otimes }(\omega _{\circ }^{-})\otimes \mathrm{K}^{\otimes
}(\omega _{\circ }^{\circ })\rightarrow \mathrm{K}^{\otimes }(\omega _{\circ
}^{\circ })\otimes \mathrm{K}^{\otimes }(\omega _{+}^{\circ }),  \notag \\
K(\mathbf{\omega })=l^{\otimes }(\omega _{+}^{-},g)\mathrm{k}^{\otimes
}(\omega _{+}^{\circ },g)j^{\otimes }(\omega _{\circ }^{\circ },g)\mathrm{k}%
^{\star \otimes }(\omega _{\circ }^{-},g),\quad \mathrm{k}^{\star }(g)=%
\mathrm{k}(g^{\star })^{\ast },  \label{two f}
\end{gather}%
where $\mathrm{K}^{\otimes }(\omega )=\bigotimes_{x\in \omega }\mathrm{K}%
_{x},$%
\begin{equation*}
l^{\otimes }(\omega )=\prod_{x\in \omega }l(x),\;\;\mathrm{k}^{\otimes
}(\omega )=\bigotimes_{x\in \omega }\mathrm{k}(x),\;\;\mathrm{k}^{\star
\otimes }(\omega )=\bigotimes_{x\in \omega }\mathrm{k}^{\star
}(x),\;\;j^{\otimes }(\omega )=\bigotimes_{x\in \omega }j(x).
\end{equation*}

\begin{theorem}
\label{T2} Let $\mathbf{K}=\bigoplus_{n=0}^{\infty }\mathbf{K}^{(n)}$ be a
decomposable operator in the pseudo-Fock space $\mathbb{F}$ defined in 
\textup{(\ref{two e})} by a linear combination of the kernels of the form 
\textup{(\ref{two f})}. Then the operators $\epsilon (\mathbf{K})=J^{\dagger
}\mathbf{K}J$, defined by pseudo-projection $J^{\dagger }:\mathbb{F}%
\rightarrow \mathrm{F}$, 
\begin{equation}
(J^{\dagger }\mathbf{h})(\omega )=\int \mathrm{h}(\omega ^{-},\omega
,\emptyset )\mathrm{d}\omega ^{-},\;\mathbf{h}\in \mathbb{F},  \label{two g}
\end{equation}%
as the adjoint to \textup{(\ref{two d}),} can be extended to a continuous
operator 
\begin{equation}
\lbrack \epsilon (\mathbf{K})\mathrm{h}](\omega )=\sum_{\upsilon \subseteq
\omega }\int K(\omega \setminus \upsilon ,\upsilon ,\vartheta _{\circ })%
\mathrm{h}(\upsilon \sqcup \vartheta _{\circ })\mathrm{d}\vartheta _{\circ },
\label{two h}
\end{equation}%
where $K(\vartheta ^{\circ },\upsilon ,\vartheta _{\circ })=\int K%
\begin{pmatrix}
\vartheta & \vartheta _{\circ } \\ 
\vartheta ^{\circ } & \upsilon%
\end{pmatrix}%
\mathrm{d}\vartheta $ defined on the completion $\mathrm{F}$ of the
pre-Hilbert space $\Gamma (\mathrm{K})$ with respect to the family of the
seminorms $\Vert \mathrm{h}\Vert ^{f}=\Vert \pi (f)\mathrm{h}\Vert $, $f\in 
\mathfrak{m}$. The map $\epsilon :\mathbf{K}\mapsto \epsilon (\mathbf{K})$
defines a Fock $\ast $-representation $\epsilon (\lambda \mathbf{K}^{\dagger
}+\lambda ^{\ast }\mathbf{K})=\lambda \epsilon (\mathbf{K})^{\ast }+\lambda
^{\ast }\epsilon (\mathbf{K})$, 
\begin{equation*}
\epsilon (\mathbf{KK}^{\dagger })=\epsilon (\mathbf{K})\epsilon (\mathbf{K}%
)^{\ast },\,\epsilon (\mathbf{I}^{\otimes })=I^{\otimes }
\end{equation*}%
of a decomposable $\dagger $-algebra of operators $\mathbf{K}$ with respect
to the involution $K^{\star }(\mathbf{\omega })=K(\mathbf{\omega }^{\prime
})^{\ast }$, where $(\omega _{\nu }^{\mu })^{\prime }=(\omega _{-\mu }^{-\nu
})$, and the associative product $[K\cdot K^{\star }](\mathbf{\omega })=$ 
\begin{equation}
=\sum_{\upsilon _{\nu }^{\mu }\subseteq \omega _{\nu }^{\mu }}^{\mu <\nu
}\sum_{\sigma _{+}^{-}\cap \tau _{+}^{-}=\upsilon _{+}^{-}}^{\sigma
_{+}^{-}\cup \tau _{+}^{-}=\omega _{+}^{-}}K%
\begin{pmatrix}
\omega _{+}^{-}\setminus \sigma _{+}^{-}, & \upsilon _{\circ }^{-}\bigsqcup
\upsilon _{+}^{-} \\ 
\omega _{+}^{\circ }\setminus \upsilon _{+}^{\circ }, & \omega _{\circ
}^{\circ }\bigsqcup \upsilon _{+}^{\circ }%
\end{pmatrix}%
K^{\star }%
\begin{pmatrix}
\omega _{+}^{-}\setminus \tau _{+}^{-}, & \omega _{\circ }^{-}\setminus
\upsilon _{\circ }^{-} \\ 
\upsilon _{+}^{-}\bigsqcup \upsilon _{+}^{\circ }, & \omega _{\circ }^{\circ
}\bigsqcup \upsilon _{\circ }^{-}%
\end{pmatrix}%
.  \label{two i}
\end{equation}%
It induces the involution $K^{\ast }(\vartheta ^{\circ },\upsilon ,\vartheta
_{\circ })=K(\vartheta _{\circ },\upsilon ,\vartheta ^{\circ })^{\ast }$ and 
$[K\cdot K^{\ast }](\vartheta ^{\circ },\upsilon ,\vartheta _{\circ })=$ 
\begin{equation*}
=\sum_{\upsilon _{\circ }\subset \vartheta _{\circ }}\sum_{\upsilon ^{\circ
}\subset \vartheta ^{\circ }}\int K(\vartheta ^{\circ }\setminus \upsilon
^{\circ },\upsilon \sqcup \upsilon ^{\circ },\upsilon _{\circ }\sqcup
\vartheta )K^{\ast }(\vartheta \sqcup \upsilon ^{\circ },\upsilon \sqcup
\upsilon _{\circ },\vartheta _{\circ }\setminus \upsilon _{\circ })\mathrm{d}%
\vartheta
\end{equation*}%
for the kernels $K(\vartheta ^{\circ },\upsilon ,\vartheta _{\circ })$, and
defines a factor algebra of the $\dagger $-algebra of operators $\mathbf{K}$
with respect to the zero $\dagger $-ideal $\{\mathbf{K}:\epsilon (\mathbf{KK}%
^{\dagger })=0\}$. The compression $\pi =\epsilon \circ j^{\otimes }$ of the 
$\ast $-representation $\epsilon $ to the operators $\mathbf{K}$ of the form 
\textup{(\ref{two f})} defined by the action \textup{(\ref{two h})}: 
\begin{equation}
\lbrack \pi (g)\mathrm{k}^{\otimes }](\omega )=\exp \left\{ \int (l(x,g)+%
\mathrm{k}^{\star }(x,g)\mathrm{k}(x,g)\mathrm{d}x\right\} (\mathrm{k}%
(g)+j(g)\mathrm{k})^{\otimes }(\omega ),  \label{two j}
\end{equation}%
of the kernels $K(\omega ^{\circ },\upsilon ,\omega _{\circ })=\exp
\{\lambda \left( g\right) \}\mathrm{k}^{\otimes }(\omega ^{\circ
},g)j^{\otimes }(\omega ,g)\mathrm{k}^{\star \otimes }(\omega _{\circ },g)$
on $\mathrm{k}^{\otimes }(\omega )=\otimes _{x\in \omega }\mathrm{k}(x)$
yields the unital $\ast $-representation $\pi :\mathfrak{m}\rightarrow 
\mathcal{L}(\mathrm{K})$,%
\begin{equation*}
\;\pi (e)=I^{\otimes },\;\;\;\pi (g^{\star })=\pi (g)^{\ast },\;\;\ \pi
(f\cdot h)=\pi (f)\pi (h)\;
\end{equation*}%
for all $f,g,h\in \mathfrak{m}$, and is associated with the infinitely
divisible state $\varphi :\mathfrak{m}\rightarrow \mathbb{C}$ in the sense
that $\varphi (g)=(1_{\emptyset }\mid \pi (g)1_{\emptyset })$, where $%
1_{\emptyset }=\mathrm{k}^{\otimes }$ for $\mathrm{k}=0$.
\end{theorem}

\begin{proof}
The operator (\ref{two d}) is a pseudo-isometry:

\begin{align*}
& (J\mathrm{h}\mid J\mathrm{h})=\int (J\mathrm{h})(\omega ^{-},\omega
^{\circ },\omega ^{+})(J\mathrm{h})^{\ast }(\omega ^{+},\omega ^{\circ
},\omega ^{-})\mathrm{d}\omega ^{-}\mathrm{d}\omega ^{\circ }\mathrm{d}%
\omega ^{+} \\
& =\int 1_{\emptyset }(\omega ^{-})\mathrm{h}(\omega ^{\circ })1_{\emptyset
}(\omega ^{+})\mathrm{h}^{\ast }(\omega ^{\circ })\mathrm{d}\omega ^{-}%
\mathrm{d}\omega ^{\circ }\mathrm{d}\omega ^{+}=\int \mathrm{h}^{\ast
}(\omega ^{\circ })\mathrm{h}(\omega ^{\circ })\mathrm{d}\omega ^{\circ
}=\left\langle \mathrm{h}\mid \mathrm{h}\right\rangle ,
\end{align*}%
and consequently, the Hermitian adjoint operator (\ref{two g}) defined by
the condition $\left\langle \mathrm{h}\mid J^{\dagger }\mathbf{h}%
\right\rangle =(\mathbf{h}\mid J\mathrm{h})$ for all $\mathrm{h}\in \mathrm{F%
},\mathbf{h}\in \mathbb{F}$, 
\begin{equation*}
\left\langle \mathrm{h}\mid J^{\dagger }\mathbf{h}\right\rangle =\int 
\mathrm{h}^{\ast }(\omega )(J^{\dagger }\mathrm{h})(\omega )\mathrm{d}\omega
=\iiint \mathrm{h}(\omega ^{-},\omega ^{\circ },\omega ^{+})\mathrm{h}^{\ast
}(\omega ^{\circ })1_{\emptyset }(\omega ^{+})\mathrm{d}\omega ^{-}\mathrm{d}%
\omega ^{\circ }\mathrm{d}\omega ^{+},
\end{equation*}%
is a pseudo-projection: $J^{\dagger }J\mathrm{h}=\mathrm{h}$ for all $%
\mathrm{h}\in \mathrm{F}$, where%
\begin{equation*}
(J\mathrm{h})(\omega _{-},\omega _{\circ },\omega _{+})=1_{\emptyset
}(\omega _{-})\mathrm{h}(\omega _{\circ })1_{\emptyset }(\omega _{+})
\end{equation*}%
is the canonical embedding $\mathrm{F}\subset \mathbb{F}$. We now show that
the action in $\mathbb{F}$ of the linear combinations of the operators $%
\mathbf{G}^{\otimes }$ with triangular $\mathbf{G}=[G_{\nu }^{\mu }]_{\nu
=-,\circ ,+}^{\mu =-,\circ ,+}$, $G_{\nu }^{\mu }=0$ for all $\mu <\nu $
with unit matrix entries $G_{-}^{-}=1=G_{+}^{+}$, can be written in the form
(\ref{two e}).

For we have 
\begin{align*}
(\mathbf{G}^{\otimes }\mathbf{k}^{\otimes })(\omega ^{\cdot })=\left( 
\mathbf{Gk}\right) ^{\otimes }(\omega ^{\cdot })& =\prod_{\mu }\left(
\sum_{\nu }G_{\nu }^{\mu }k^{\nu }\right) ^{\otimes }(\omega ^{\mu }) \\
& =\prod_{\mu }\,\sum_{\bigsqcup\limits_{\nu }\omega _{\nu }^{\mu }=\omega
^{\mu }}\prod_{\nu }(G_{\nu }^{\mu })^{\otimes }(\omega _{\nu }^{\mu
})\prod_{\nu }(k^{\nu })^{\otimes }(\omega _{\nu }^{\mu }),
\end{align*}%
where the sums over the decompositions $\omega ^{\mu }=\omega _{-}^{\mu
}\sqcup \omega _{\circ }^{\mu }\sqcup \omega _{+}^{\mu }$ in fact should be
taken only over $\omega ^{\mu }=\bigsqcup_{\nu \geq \mu }\omega _{\nu }^{\mu
}$ since $G_{\nu }^{\mu }=0$ for $\nu <\mu $. If $\omega =(\omega
^{-},\omega ^{\circ },\omega ^{+})$ do not intersect, then the same is true
for $\omega _{\nu }^{\cdot }=(\omega _{-}^{\cdot },\omega _{\circ }^{\cdot
},\omega _{+}^{\cdot })$ since $\omega _{\nu }^{\mu }\subseteq \omega ^{\mu
} $. Consequently, $\prod_{\mu }\mathrm{h}(\omega _{\cdot }^{\mu })=\mathrm{k%
}\left( \bigsqcup_{\mu }\omega _{\cdot }^{\mu }\right) $for $\mathrm{h}%
(\omega _{\cdot })=\prod_{\nu }(k^{\nu })^{\otimes }(\omega _{\nu })$, which
yields 
\begin{equation*}
(\mathbf{G}^{\otimes }\mathbf{k}^{\otimes }(\omega ^{\cdot
})=\sum_{\bigsqcup \omega _{\nu }^{.}=\omega ^{.}}\prod_{\mu ,\nu }(G_{\nu
}^{\mu })^{\otimes }(\omega _{\nu }^{\mu })(\prod_{\nu }k^{\nu })^{\otimes
}(\bigsqcup_{\mu }\omega _{\nu }^{\mu }),
\end{equation*}%
where $\bigsqcup_{\mu }\omega _{\nu }^{\mu }=\bigsqcup_{\mu \leq \nu }\omega
_{\nu }^{\mu }$, since $(G_{\nu }^{\mu })^{\otimes }(\omega
)=\bigotimes_{x\in \omega }G_{\nu }^{\mu }(x)$ is equal to zero if $\omega
=\omega _{\nu }^{\mu }\neq \emptyset $ for $\mu >\nu $. Thus we obtain (\ref%
{two e}) for exponential vectors $\mathbf{b}=\mathbf{k}^{\otimes }$ with the
kernel $K(\omega )=\prod_{\mu \leq \nu }(G_{\nu }^{\mu })^{\otimes }(\omega
_{\nu }^{\mu })$ of the form (\ref{two f}). Since this formula is linear
with respect to the kernel $\mathbf{K}$, it is also valid for linear
combinations $\mathbf{K}=\Sigma \lambda _{i}\mathbf{G}_{i}^{\otimes }$ at
least on $\Gamma (\mathbf{K})$. We now define the operator $J^{\dagger }%
\mathbf{K}J$ in $\mathrm{F}$, employing the formula 
\begin{equation*}
\int \sum_{\bigsqcup \omega _{\mu }=\omega }f(\omega _{-},\omega _{\circ
},\omega _{+})\mathrm{d}\omega =\iiint f(\omega _{-},\omega _{\circ },\omega
_{+})\mathrm{d}\omega _{-}\mathrm{d}\omega _{\circ }\mathrm{d}\omega _{+}.
\end{equation*}%
Taking into account the forms (\ref{two d}), (\ref{two g}) of the operators $%
J$, $J^{\dagger }$ we obtain for $\mathrm{h}\in \Gamma (\mathrm{K})$ the
formula 
\begin{align*}
\lbrack J^{\dagger }\mathbf{K}J\mathrm{h}](\omega )& =\int (\mathbf{K}J%
\mathrm{h})(\omega ^{-},\omega ,\emptyset )\mathrm{d}\omega ^{-} \\
& =\iiint \sum_{\omega _{\circ }^{\circ }\sqcup \omega _{+}^{\circ }=\omega
}K(\mathbf{\omega })1_{\emptyset }(\omega _{-}^{-})\mathrm{h}(\omega _{\circ
}^{-}\sqcup \omega _{\circ }^{\circ })\mathrm{d}\omega _{-}^{-}\mathrm{d}%
\omega _{\circ }^{-}\mathrm{d}\omega _{+}^{-} \\
& =\sum_{\omega _{\circ }^{\circ }\sqcup \omega _{+}^{\circ }=\omega }\int
K(\omega _{+}^{\circ },\omega _{\circ }^{\circ },\omega _{\circ }^{-})%
\mathrm{h}(\omega _{\circ }^{-}\sqcup \omega _{\circ }^{\circ })\mathrm{d}%
\omega _{\circ }^{-},
\end{align*}%
which can be written in the form (\ref{two h}) in the notation $\upsilon
=\omega _{\circ }^{\circ }$ and $\omega _{+}^{\circ }=\omega \setminus
\upsilon =\omega \cap \overline{\upsilon }$.

We shall now prove that the pseudo-conditional expectation $\mathbf{K}%
\rightarrow J^{\dagger }\mathbf{K}J$ is a $\ast $-representation on $\Gamma (%
\mathrm{K})$. To this end it is sufficient to show that this map is a
homomorphism with respect to the binary operation (\ref{two i}), involution $%
\mathbf{K}\mapsto \mathbf{K}^{\dagger }$, and the unit $\mathbf{K}=\mathbf{I}%
^{\otimes }$ on the generating elements $\mathbf{G}^{\otimes }$ for which (%
\ref{two h}) yields (\ref{two j}) for $\mathrm{h}=\mathrm{k}^{\otimes }$: $%
[J^{\dagger }\mathbf{G}^{\otimes }J\mathrm{k}^{\otimes }](\omega )=$ 
\begin{align*}
& =\sum_{\omega _{\circ }^{\circ }\sqcup \omega _{+}^{\circ }=\omega }\iint
\prod_{\mu <\nu }(G_{\nu }^{\mu })^{\otimes }(\omega _{\nu }^{\mu })\mathrm{k%
}^{\otimes }(\omega _{\circ }^{-}\sqcup \omega _{\circ }^{\circ })\mathrm{d}%
\omega _{\circ }^{-}\mathrm{d}\omega _{+}^{-} \\
& =\sum_{\omega _{\circ }^{\circ }\sqcup \omega _{+}^{\circ }=\omega
}(G_{\circ }^{\circ }\mathrm{k})^{\otimes }(\omega _{\circ }^{\circ
})(G_{+}^{\circ })^{\otimes }(\omega _{+}^{\circ })\int (G_{\circ }^{-}%
\mathrm{k})^{\otimes }(\omega _{\circ }^{-})\mathrm{d}\omega _{\circ
}^{-}\int (G_{+}^{-})^{\otimes }(\omega _{+}^{-})\mathrm{d}\omega _{+}^{-} \\
& =(G_{\circ }^{\circ }\mathrm{k}+G_{+}^{\circ })^{\otimes }(\omega )\exp
\{\int (G_{\circ }^{-}\mathrm{k}+G_{+}^{-})(x)\mathrm{d}x\}.
\end{align*}%
Using this formula we find that $[J^{\dagger }\mathbf{I}^{\otimes }J\mathrm{k%
}^{\otimes }](\omega )=\mathrm{k}^{\otimes }(\omega )$, that is, $J^{\dagger
}\mathbf{I}^{\otimes }J=I^{\otimes }$, and

\begin{equation*}
\lbrack J^{\dagger }\mathbf{G}^{\dagger \otimes }J\mathrm{k}^{\otimes
}](\omega )=(G_{\circ }^{\circ }\mathrm{k}+G_{\circ }^{-\ast })^{\otimes
}(\omega )\exp \left\{ \int (G_{+}^{\circ \ast }\mathrm{k}+G_{+}^{-\ast })(x)%
\mathrm{d}x\right\} ,
\end{equation*}%
that is, $J^{\dagger }\mathbf{K}^{\dagger }J=(J^{\dagger }\mathbf{K}J)$ for $%
\mathbf{K}=\Sigma \lambda _{i}G_{i}^{\otimes }$, $\mathbf{K}^{\dagger
}=\Sigma \lambda _{i}\mathbf{G}_{i}^{\dagger \otimes }$, and 
\begin{align*}
& [J^{\dagger }(\mathbf{FG})^{\otimes }J\mathrm{k}^{\otimes }](\omega
)=(F_{\nu }^{\circ }G_{\circ }^{\nu }\mathrm{k}+F_{\nu }^{\circ }G_{+}^{\nu
})^{\otimes }(\omega )\exp \left\{ \int (F_{\nu }^{-}G_{\circ }^{\nu }%
\mathrm{k}+F_{\nu }^{-}G_{+}^{\nu })(x)\mathrm{d}x\right\} \\
& =(F_{\circ }^{\circ }G_{\circ }^{\circ }\mathrm{k}+F_{\circ }^{\circ
}G_{+}^{\circ }+F_{+}^{\circ })^{\otimes }(\omega )\exp \left\{ \int
(G_{\circ }^{-}\mathrm{k}+F_{\circ }^{-}G_{\circ }^{\circ }\mathrm{k}%
+G_{+}^{-}+F_{\circ }^{-}G_{+}^{\circ }+F_{+}^{-})(x)\mathrm{d}x\right\} \\
& =(F_{\circ }^{\circ }\mathrm{k}+F_{+}^{\circ })^{\otimes }(\omega )e^{\int
(F_{\circ }^{-}\mathrm{k}+F_{+}^{-})(x)\mathrm{d}x}(G_{\circ }^{\circ }%
\mathrm{k}+G_{+}^{\circ })^{\otimes }(\omega )e^{\int (G_{\circ }^{-}\mathrm{%
k}+G_{+}^{-})(x)\mathrm{d}x},
\end{align*}%
where we have used the rule of multiplication 
\begin{equation*}
(\mathbf{F}\mathbf{G})_{\nu }^{\mu }=\sum_{\iota }F_{\iota }^{\mu }G_{\nu
}^{\iota }\equiv F_{\iota }^{\mu }G_{\nu }^{\iota }
\end{equation*}%
of the triangular matrices $\mathbf{F}=[F_{\nu }^{\mu }]$, $\mathbf{G}%
=[G_{\nu }^{\mu }]$, $\mu ,\nu \in \{-,\circ ,+\},\,F_{\nu }^{\mu }=0=G_{\nu
}^{\mu }$ for $\mu >\nu $, with the entries $F_{-}^{-}=1=F_{+}^{+}$, $%
G_{-}^{-}=1=G_{+}^{+}$. Thus we have proved that $\epsilon (\mathbf{F}%
^{\otimes }\mathbf{G}^{\otimes })=\epsilon (\mathbf{F}^{\otimes })\epsilon (%
\mathbf{G}^{\otimes })$, where $\epsilon (\mathbf{G}^{\otimes })\mathrm{h}%
=J^{\dagger }\mathbf{G}^{\otimes }J\mathrm{h}$ for any $\mathrm{h}=\Sigma
\lambda _{i}\mathrm{k}_{i}^{\otimes }\in \Gamma (\mathrm{K})$. We complete $%
\Gamma (\mathrm{K})$ by sequences $\mathrm{h}_{n}\in \Gamma (\mathrm{K})$
that are fundamental with respect to each of the seminorms $\Vert \mathrm{h}%
\Vert ^{f}=\Vert \epsilon \left[ \mathbf{j}^{\otimes }(f)\right] \mathrm{h}%
\Vert $, $f\in \mathfrak{m}$ (among others, also with respect to $\Vert 
\mathrm{h}\Vert e=\Vert \mathrm{h}\Vert $). Since $\pi (g)=\epsilon \lbrack 
\mathbf{j}^{\otimes }(g)]$ is a $\ast $-representation of $\mathfrak{m}$ on $%
\Gamma (\mathrm{K})$: 
\begin{equation*}
\pi (g\star f)=\epsilon \left[ \left( \mathbf{j}(g)\mathbf{j}(f^{\star
})\right) ^{\otimes }\right] =\epsilon \lbrack \mathbf{j}^{\otimes
}(g)]\epsilon \lbrack \mathbf{j}^{\otimes }(f)^{\dagger }]=\pi (g)\pi
(f)^{\ast },
\end{equation*}%
any fundamental sequence remains fundamental after multiplying it by $\pi
(g):\Vert \pi (g)\mathrm{h}\Vert ^{f}=\Vert \mathrm{h}\Vert ^{g\cdot f}$.
This allows us to extend the operators $\epsilon \lbrack \mathbf{j}^{\otimes
}(g)]=J^{\dagger }\mathbf{j}^{\otimes }(g)J$ to continuous operators $\pi
(g) $ on the completion $\mathrm{F}$ of $\Gamma \left( \mathrm{K}\right) $
with respect to the convergence described above. The continuity implies that
the algebraic relations in the decomposable $\dagger $-algebra $\mathbb{B}=%
\mathbb{CM}^{\otimes }$, where $\mathbb{M}=\int_{X}^{\oplus }\mathbb{M}_{x}%
\mathrm{d}x=\mathbf{j}\left( \mathfrak{m}\right) $, become represented in
the operator algebra $\mathcal{L}\left( \mathrm{F}\right) $ of continuous
operators \thinspace $L,L^{\ast }:\mathrm{F}\rightarrow \mathrm{F}$ on the
poly-Fock space $\mathrm{F}$. Obviously, the linear hull $\Sigma \lambda
_{i}\pi (g_{i})$ defines a $\ast $-subalgebra $\mathcal{B}$ of operators $%
\epsilon (\mathbf{K})\in \mathcal{L}(\mathrm{K})$ which is a homomorphic
image of the $\dagger $-algebra $\mathbb{B}$ genrated by linear combinations 
$\mathbf{K}=\Sigma \lambda _{i}\mathbf{j}^{\otimes }(g_{i})$ of decomposable
operators $\mathbf{G}_{i}^{\otimes }=1\oplus \mathbf{G}_{i}\oplus \mathbf{G}%
_{i}^{\otimes 2}\oplus \ldots $, where $\mathbf{G}_{i}=\mathbf{j}(g_{i})$.
We recall that the elements of $\mathbb{B}$ as decomposable operators in the
pseudo-Fock space $\mathbb{F}$ are represented by triangular kernels $%
\mathbf{K}\in \mathcal{L}\left( \mathbb{F}\right) $ with $\mathbf{K}%
^{\dagger }\in \mathcal{L}\left( \mathbb{F}\right) $ described in (\ref{two
e}) by the kernels $K^{\star }(\mathbf{\omega })=K(\mathbf{\omega }^{\prime }%
\mathbf{)}^{\ast }$, where the table $\mathbf{\omega }^{\prime }$ of four
subsets differs from $\mathbf{\omega }=(\omega _{\nu }^{\mu })\in \mathbf{%
\Omega }$ only by the interchange of $\omega _{\circ }^{-}$ and $\omega
_{+}^{\circ }$, and the multiplication $\mathbf{K}^{\dagger }\mathbf{K}$ is
defined, as in any semigroup algebra, by the operation $\mathbf{KK}^{\dagger
}=\Sigma \lambda _{i^{\prime }}\lambda _{i}^{\ast }\mathbf{j}^{\otimes
}(g_{i^{\prime }}\star g_{i})$ in $\mathbb{M}$. Here $\mathbf{KK}^{\dagger }$
is defined by the kernel (\ref{two i}) which can be verified stightforward
for the generating kernels (\ref{two f}) by virtue of $\mathbf{j}^{\otimes
}(f)\mathbf{j}^{\otimes }(g)=\mathbf{j}^{\otimes }(f\cdot g)$: $l^{\otimes
}(\omega _{+}^{-},f\cdot g)\mathrm{k}^{\otimes }\left( \omega _{+}^{\circ
},f\cdot g)j^{\otimes }(\omega _{\circ }^{\circ },f\cdot g\right) \mathrm{k}%
^{\star \otimes }(\omega _{\circ }^{-},f\cdot g)=$

\begin{eqnarray*}
\lefteqn{=[j(f)j(g)]_{\omega _{\circ }^{\circ }}^{\otimes }[j(f)\mathrm{k}%
(g)+\mathrm{k}\left( f\right) ]_{\omega _{+}^{\circ }}^{\otimes }[l(f)+%
\mathrm{k}^{\star }(f)\mathrm{k}(g)+l(g)]_{\omega _{+}^{-}}^{\otimes }[%
\mathrm{k}^{\ast }(g)+\,\mathrm{k}^{\star }(f)j(g)]_{\omega _{\circ
}^{-}}^{\otimes }} \\
&=&\sum_{\sigma \sqcup \tau \sqcup \upsilon _{+}^{-}=\omega
_{+}^{-}}l(f)_{\sigma }^{\otimes }l(g)_{\tau }^{\otimes }[\mathrm{k}^{\star
}(f)\mathrm{k}(g)]_{\upsilon _{+}^{-}}^{\otimes }j(f)_{\omega _{\circ
}^{\circ }}^{\otimes }j(g)_{\omega _{\circ }^{\circ }}^{\otimes } \\
&&\hspace{0.5cm}\otimes \sum_{\upsilon _{+}^{\circ }\sqcup \sigma
_{+}^{\circ }=\omega _{+}^{\circ }}[j(f)\mathrm{k}(g)]_{\upsilon _{+}^{\circ
}}^{\otimes }\otimes \mathrm{k}_{\sigma _{+}}^{\otimes }(f)\sum_{\sigma
_{\circ }^{-}\sqcup \tau ^{-}=\omega _{\circ }^{-}}\mathrm{k}_{\tau
^{-}}^{\star }(g)^{\otimes }\otimes \lbrack \mathrm{k}^{\star
}(f)j(g)]_{\upsilon _{\circ }^{-}}^{\otimes } \\
&=&\sum_{\upsilon _{\nu }^{\mu }\subseteq \omega _{\nu }^{\mu }}^{\mu <\nu
}\,\sum_{\sigma \sqcup \tau =\omega _{+}^{-}\setminus \upsilon
_{+}^{-}}l^{\otimes }(\sigma ,f)\mathrm{k}^{\otimes }(\omega _{+}^{\circ
}\setminus \upsilon _{+}^{\circ },f)j^{\otimes }(\omega _{\circ }^{\circ
}\sqcup \upsilon _{+}^{\circ },f)\mathrm{k}^{\star \otimes }(\upsilon
_{\circ }^{-}\sqcup \upsilon _{+}^{-},f) \\
&&\hspace{0.5cm}\times l^{\otimes }(\tau ,g)\mathrm{k}^{\otimes }(\upsilon
_{+}^{-}\sqcup \upsilon _{+}^{\circ },g)j^{\otimes }(\omega _{\circ }^{\circ
}\sqcup \upsilon _{\circ }^{-},g)\mathrm{k}^{\star \otimes }(\omega _{\circ
}^{-}\setminus \upsilon _{\circ }^{-},g).
\end{eqnarray*}%
Integrating (\ref{two i}) over $\omega _{+}^{-}\in \Omega $ we obtain the
formula of multiplication of the kernels $K(\omega _{-}^{\circ },\omega
_{\circ }^{\circ },\omega _{\circ }^{+})=\int K(\mathbf{\omega })\mathrm{d}%
\omega _{+}^{-}$:$\;\;\;\;\int [K\cdot K^{\star }](\mathbf{\omega })\mathrm{d%
}\omega _{+}^{-}=$ 
\begin{eqnarray*}
\lefteqn{=\iiint \mathrm{d}\sigma \mathrm{d}\tau \mathrm{d}\upsilon
_{+}^{-}\sum_{\upsilon _{\nu -}^{\mu }\omega _{\nu }^{\mu }}^{\mu <\nu
}K\left( \begin{array}{cc}\sigma ,\upsilon _{\circ }^{-}\sqcup \upsilon
_{+}^{-}\newline
^{\circ }\setminus \upsilon _{+}^{\circ },\omega _{\circ }^{\circ }\sqcup
\upsilon _{+}^{\circ }\end{array}\right) K^{\star }\left( \begin{array}{cc}%
\tau ,\omega _{\circ }^{-}\setminus \upsilon _{\circ }^{-}\newline
^{-}\sqcup \upsilon _{+}^{\circ },\omega _{\circ }^{\circ }\sqcup \upsilon
_{\circ }^{-}\end{array}\right) } \\
&=&\sum_{\upsilon _{\circ }^{-}\subseteq \omega _{\circ
}^{-}}\,\sum_{\upsilon _{+}^{\circ }\subseteq \omega _{+}^{\circ }}\int
K(\omega _{+}^{\circ }\setminus \upsilon _{+}^{\circ },\omega _{\circ
}^{\circ }\sqcup \upsilon _{n}^{\circ },\upsilon _{\circ }^{-}\sqcup
\upsilon _{+}^{-})K^{\ast }(\upsilon _{+}^{-}\sqcup \upsilon _{+}^{\circ
},\omega _{\circ }^{\circ }\sqcup \upsilon _{\circ }^{-},\omega _{\circ
}^{-}\setminus \upsilon _{\circ }^{-})\mathrm{d}\upsilon _{+}^{-},
\end{eqnarray*}%
where $K^{\ast }(\omega _{-}^{\circ },\omega _{\circ }^{\circ },\omega
_{\circ }^{+})=\int K^{\star }(\mathbf{\omega })\mathrm{d}\omega _{\circ
}^{-}=K(\omega _{\circ }^{+},\omega _{\circ }^{\circ },\omega _{-}^{\circ
})^{\ast }$. Thus we obtained a $\ast $-algebraic structure for
three-argument kernels connected with the Maassen-Meyer kernels $M(\vartheta
^{\circ },\vartheta ,\vartheta _{\circ })$ \cite{38}, \cite{40} by a
one-to-one M\"{o}bius transformation 
\begin{equation*}
K(\vartheta ^{\circ },\omega ,\vartheta _{\circ })=\sum_{\vartheta \subseteq
\omega }M(\vartheta ^{\circ },\vartheta ,\vartheta _{\circ })\otimes
I^{\otimes }(\omega \setminus \vartheta ).
\end{equation*}

We finally consider a $\star $-invariant subspace of the $\star $-algebra of
kernels $K_{0}(\mathbf{\omega })$ defined by the condition $\int K_{0}(%
\mathbf{\omega })\mathrm{d}\omega _{+}^{-}=0$. This is a zero ideal of the
homomorphism $\{K(\mathbf{\omega })\}\mapsto \{K(\omega _{+}^{\circ },\omega
_{\circ }^{\circ },\omega _{\circ }^{-})\}$ for which we take a convention
that it transforms the star conjugation $\star $ into $\ast $. Consequently,
it is a two-sided ideal ( $\dagger $-ideal in terms of $\mathbf{K}$, or $%
\star $-ideal in terms of the kernels $K$): 
\begin{equation*}
\int (\mathbf{KK}_{0})(\mathbf{\omega })\,\mathrm{d}\omega _{+}^{-}=0=\int (%
\mathbf{K}_{0}\mathbf{K})(\mathbf{\omega })\,\mathrm{d}\omega _{+}^{-},\quad
\forall \mathbf{K},
\end{equation*}%
which is contained in the zero ideal of the representation $\epsilon
:\epsilon (\mathbf{K}_{0})=0$ if $K_{0}(\omega _{+}^{\circ },\omega _{\circ
}^{\circ },\omega _{\circ }^{-})=0$. One can show that, owing to the fact
that the measure $\mathrm{d}x$ on $X$ is atomless, the zero ideal of the
representation $\epsilon $ is exhausted in this way. This follows from the
uniqueness of the stochastic representation (\ref{two h}), proved in terms
of the Maassen-Meyer kernels in \cite{38}, \cite{40}. Consequently, the
integral $K(\omega _{+}^{\circ },\omega _{\circ }^{\circ },\omega _{\circ
}^{-})=\int K(\mathbf{\omega })\mathrm{d}\omega _{+}^{-}$ is a homomorphism
of the factorization of the $\star $-algebra of kernels $K(\mathbf{\omega })$
also by the zero ideal of the representation $\epsilon $. The proof is
complete.
\end{proof}

\begin{remark}
We introduce four types $G_{\mu }^{\nu },\nu \neq -,\mu \neq +$, of
elementary triangular decomposable operators in $\mathrm{K}$ described by
matrices of the form 
\begin{eqnarray*}
\mathbf{G}_{\circ }^{+}(x) &=&\left[ 
\begin{array}{ccc}
1 & 0 & 0 \\ 
0 & I(x) & g_{+}^{\circ }(x) \\ 
0 & 0 & 1%
\end{array}%
\right] ,\quad \mathbf{G}_{-}^{\circ }(x)=\left[ 
\begin{array}{ccc}
1 & g_{\circ }^{-}\left( x\right) & 0 \\ 
0 & I(x) & 0 \\ 
0 & 0 & 1%
\end{array}%
\right] , \\
\mathbf{G}_{-}^{+}(x) &=&\left[ 
\begin{array}{ccc}
1 & 0 & g_{+}^{-}(x) \\ 
0 & I(x) & 0 \\ 
0 & 0 & 1%
\end{array}%
\right] ,\quad \mathbf{G}_{\circ }^{\circ }(x)=\left[ 
\begin{array}{ccc}
1 & 0 & 0 \\ 
0 & G(x) & 0 \\ 
0 & 0 & 1%
\end{array}%
\right] ,
\end{eqnarray*}%
and we write%
\begin{equation*}
G^{N}=\epsilon \lbrack (\mathbf{G}_{\circ }^{\circ })^{\otimes }]\equiv
G^{\otimes },\,e^{\lambda \left( g\right) }=\epsilon \lbrack (\mathbf{G}%
_{+}^{-})^{\otimes }],e^{\mathrm{k}^{\star }\left( g\right) A_{-}^{\circ
}}=\epsilon \lbrack (\mathbf{G}_{-}^{\circ })^{\otimes }],e^{A_{\circ }^{+}%
\mathrm{k}\left( g\right) }=\epsilon \lbrack (\mathbf{G}_{\circ
}^{+})^{\otimes }],
\end{equation*}
where $\epsilon $ is the map \textup{(\ref{two h})} for $\mathbf{K(\omega
)=\otimes }_{x\in \omega }\mathbf{G}_{\mu }^{\nu }(x)$.

Then the representation $\mathfrak{m}\ni g\mapsto \pi (g)$, associated with
the infinitely divisible state $\varphi (g)=e^{\lambda \left( g\right) }$
with respect to the vacuum vector $1_{\emptyset }\in \mathrm{K}$, can be
written as a `normally-ordered' product 
\begin{equation*}
\pi (g)=e^{\lambda \left( g\right) }e^{A_{\circ }^{+}\mathrm{k}\left(
g\right) }G^{N}e^{\mathrm{k}\left( g\right) A_{-}^{\circ }}
\end{equation*}
for all $G\in \mathfrak{m}$, defined by the functions $G(x)=j(x,g)$, $%
g_{+}^{-}(x)=l(x,g)$, $g_{\circ }^{-}(x)=\mathrm{k}^{\star }(x,g)$, $%
g_{+}^{\circ }(x)=\mathrm{k}(x,g)$.
\end{remark}

In fact, an arbitrary triangular operator $\mathbf{G}$ in $\mathrm{K}$ with
the entries $G_{-}^{-}=1=G_{+}^{+}$ can be decomposed into a
`normally-ordered' product of elementary matrices: 
\begin{equation*}
\left[ 
\begin{array}{ccc}
1 & g_{\circ }^{-} & g_{+}^{-} \\ 
0 & G & g_{+}^{\circ } \\ 
0 & 0 & 1%
\end{array}%
\right] =\left[ 
\begin{array}{ccc}
1 & 0 & g_{+}^{-} \\ 
0 & I & 0 \\ 
0 & 0 & 1%
\end{array}%
\right] \left[ 
\begin{array}{ccc}
1 & 0 & 0 \\ 
0 & I & g_{+}^{\circ } \\ 
0 & 0 & 1%
\end{array}%
\right] \left[ 
\begin{array}{ccc}
1 & 0 & 0 \\ 
0 & G & 0 \\ 
0 & 0 & 1%
\end{array}%
\right] \left[ 
\begin{array}{ccc}
1 & g_{\circ }^{-} & 0 \\ 
0 & I & 0 \\ 
0 & 0 & 1%
\end{array}%
\right] .
\end{equation*}%
Since the maps $\mathbf{G\mapsto G}^{\otimes }$ and $\mathbf{K}\mapsto
\epsilon (\mathbf{K})$ are multiplicative, we hence obtain for $\mathbf{K=G}%
^{\otimes }$ the equality 
\begin{equation*}
\epsilon \lbrack \mathbf{G}^{\otimes }]=\epsilon \lbrack (\mathbf{G}%
_{-}^{+})^{\otimes }]\epsilon \lbrack (\mathbf{G}_{-}^{+})^{\otimes
}]\epsilon \lbrack (\mathbf{G}_{\circ }^{\circ })^{\otimes }]\epsilon
\lbrack (\mathbf{G}_{-}^{\circ })^{\otimes }],
\end{equation*}%
which gives for $\mathbf{G=j}(g)$ the corresponding representation for $\pi
(g)=\epsilon \lbrack \mathbf{j}^{\otimes }(g)]$.

\section{The pseudo-Poisson structure of chaotic sates on quantum It\^{o}
algebras}

In this section we assume that the $\star $-semigroup $\mathfrak{b}$ is also
an additive group with the same neutral element $0\equiv u$, such that $%
\star $-monoid $\mathfrak{m}$ with respect to the multiplication, denoted
now by bold dot, $f\bullet h$, has also the structure of an additive group
with respect to the pointwise operations 
\begin{equation*}
(-g)(x)=-g(x),\,(f+h)(x)=f(x)+h(x),\,e(x)=0.
\end{equation*}%
In this terms $f\sqcup h$, whenever it is defined, can be written as $f+h$.
We shall also assume that $\left( f+h\right) ^{\star }=f^{\star }+h^{\star }$
and that the conditionally positive form (\ref{one b}) is an additive
homomorphism $\mathfrak{m}\mapsto \mathbb{C}$ which will be denoted as $%
\lambda \left( g\right) =\left\langle g\right\rangle $: 
\begin{equation*}
\left\langle -g\right\rangle =-\left\langle g\right\rangle ,\quad
\left\langle f+h\right\rangle =\left\langle f\right\rangle +\left\langle
h\right\rangle ,\quad \left\langle 0\right\rangle =0.
\end{equation*}%
Condition (\ref{one d}) of infinite divisibility of the state $\varphi
_{\Delta }(b)=e^{\lambda _{\Delta }(b)}$ for any integrable $\Delta
\subseteq X$ can now be written in the form of positive definiteness 
\begin{equation}
\sum_{a,c\in \mathfrak{b}}\kappa _{a}\left\langle \left( a\cdot c\right)
_{\Delta }\right\rangle \kappa _{c}^{\star }\geq 0,\quad \forall \kappa
_{b}\in \mathbb{C}:\left\vert \mathrm{supp\,}\kappa \right\vert <\infty
\label{three a}
\end{equation}%
of the function $\lambda _{\Delta }(b)=\left\langle b_{\Delta }\right\rangle 
$, where $\kappa ^{\star }=\kappa _{g^{\star }}^{\ast }$ and $b_{\Delta }$
is elementary function $b_{\Delta }(x)=b$ for $x\in \Delta $ and $b_{\Delta
}(x)=0$ otherwise, with respect to the new product $a\cdot c=a\bullet c-a-c$%
. This positive definiteness follows from the additivity of the form $%
\left\langle g\right\rangle $, which yields 
\begin{equation*}
\sum_{f,h\in \mathfrak{m}}\kappa _{f}\left\langle f\bullet h\right\rangle
\kappa _{h}^{\star }=\sum_{f,h\in \mathfrak{m}}\kappa _{f}(\left\langle
f\cdot h\right\rangle +\left\langle f\right\rangle +\left\langle
h\right\rangle )\kappa _{h}^{\star }=\sum_{f,h\in \mathfrak{m}}\kappa
_{f}\left\langle f\cdot h\right\rangle \kappa _{h}^{\star }
\end{equation*}%
for any function $g\mapsto \kappa _{g}\in \mathbb{C}$ with $\left\vert 
\mathrm{supp}\,\kappa \right\vert <\infty $ and such that $\Sigma \kappa
_{g}=0$, where on the right-hand side we can arbitrarily change the value of 
$\kappa _{e}$ since%
\begin{equation*}
0\cdot b=0\bullet a-0-a=-0=b\cdot 0
\end{equation*}
and therefore $\left\langle 0\bullet g\right\rangle =0=\left\langle g\bullet
0\right\rangle $.

We shall now assume that the additive $\star $-group $\mathfrak{b}$ has a
ring structure with respect to the new product:%
\begin{equation*}
a\cdot \left( b+c\right) =a\cdot b+a\cdot c,\;a\cdot \left( b\cdot c\right)
=\left( a\cdot b\right) \cdot c.\;\;
\end{equation*}%
Note that the associativity of this product simply follows from its
distributivity which is equivalent to the relation 
\begin{equation*}
a\cdot b+c\cdot b=b+\left( a+c\right) \cdot b.
\end{equation*}%
This is particularly easy to see if the ring $\mathfrak{b}$ has identity $1$
such that $1b=b=b1$, i.e.%
\begin{equation*}
1\cdot b-b=1+b=b\cdot 1-b,
\end{equation*}%
by virtue of the relation 
\begin{equation*}
1+a\cdot c=(1+a)\cdot (1+c),\quad \forall a,c\in \mathfrak{b}.
\end{equation*}

The function $\varphi (g)=\mathrm{e}^{\left\langle g\right\rangle }$
corresponding to the additive and positive in the above sense form $%
\left\langle g\right\rangle $ is chaotic: 
\begin{equation*}
\varphi (f\sqcup h)=e^{\left\langle f+h\right\rangle }=\varphi (f)\varphi
(h),\quad \forall f,\,h\in \mathfrak{m}:fh=0,
\end{equation*}%
and will be called the \textit{pseudo-Poisson state} over the $\star $-ring
(or $\star $-algebra) $\mathfrak{m}$ with respect to the operations $+,\cdot 
$. In other words, a pseudo-Poisson state is described by an exponential
generating functional $\varphi (f+h)=\varphi (f)\varphi (h)$ on $\mathfrak{m}%
\ni f,h$ which is positive definite in the sense of (\ref{one a}) with
respect to the operation $f\bigstar h=f+h^{\star }+f\star h$ and $e=0$,
defined pointwisely by means of the ring $\mathfrak{m}$.

By this distributivity, the canonical maps%
\begin{equation*}
\mathrm{k}:\mathfrak{m}\ni g\mapsto \mathrm{k}\left( g\right) \equiv
g\rangle \text{ \ \ and \ \ \ }\mathrm{k}^{\star }:\mathfrak{m}\ni g\mapsto 
\mathrm{k}^{\star }\left( g\right) \equiv \langle g
\end{equation*}%
defining the minimal decomposition (\ref{one e}) of the additive (linear)
positive form $\left\langle g\right\rangle $ are additive (linear), and the $%
\star $-map $i:\mathfrak{m}\ni g\mapsto j(g)-I$ satisfying in accordance
with (\ref{one f}) the conditions 
\begin{eqnarray*}
i(g)h\rangle &=&g\bullet h\rangle -g\rangle -h\rangle =g\cdot h\rangle
,\quad \forall g,h\in \mathfrak{m}, \\
\langle fi(g) &=&\langle f\bullet g-\langle g-\langle f=\langle f\cdot
g,\quad \forall g,f\in \mathfrak{m},
\end{eqnarray*}%
is also additive (linear): 
\begin{equation*}
i(f+h)=i(f)+i(h),\quad i(0)=0,\,\quad i(\lambda g)=\lambda i(g).
\end{equation*}%
Moreover, the maps $i_{x}(b)=j_{x}(b)-I_{x}$ are $\star $-representations of
the ring (algebra) $\mathfrak{b}$ in the operator $\star $-algebras $%
\mathcal{L}(\mathrm{K}_{x})=\{L:\mathrm{K}_{x}\rightarrow \mathrm{K}%
_{x},L^{\ast }\mathrm{K}_{x}\subseteq \mathrm{K}_{x}\}$ of the poly-Hilbert
spaces $\mathrm{K}_{x}=\{\mathrm{k}:\left\Vert j_{x}(a)\mathrm{k}\right\Vert
<\infty ,\forall a\in \mathfrak{b}\}$: 
\begin{eqnarray*}
i_{x}(a\cdot c) &=&i_{x}(a\bullet c)-i_{x}(a)-i_{x}(c)=j_{x}(a\bullet
c)-1-i_{x}(a)-i_{x}(c) \\
&=&(i_{x}(a)+1)(i_{x}(c)+1)-1-i_{x}(a)-i_{x}(c)=i_{x}(a)i_{x}(c).
\end{eqnarray*}%
Combining these relations and taking into account the fact that by
additivity (linearity) of $l_{x}(b)$ in the integral (\ref{one b}) we have 
\begin{equation*}
l_{x}(a\cdot c)=l_{x}(a\bullet c)-l_{x}(a)-l_{x}(c)=\mathrm{k}_{x}^{\star
}(a)\mathrm{k}_{x}(c)
\end{equation*}%
almost everywhere on $X$, we obtain decomposable $\star $-representation $%
\boldsymbol{i}(x,b)=\boldsymbol{i}_{x}(g(x))$ with four-component%
\begin{equation}
\boldsymbol{i}_{x}(b)=\left( 
\begin{array}{ll}
l_{x}(b) & \mathrm{k}_{x}^{\star }(b) \\ 
\mathrm{k}_{x}(b) & i_{x}(b)%
\end{array}%
\right) ,\quad \boldsymbol{i}_{x}(b^{\star })=\left( 
\begin{array}{ll}
l_{x}(b) & \mathrm{k}_{x}^{\star }(b) \\ 
\mathrm{k}_{x}(b) & i_{x}(b)%
\end{array}%
\right) ^{\star }  \label{three b}
\end{equation}%
of the $\star $-ring $\mathfrak{b}$ with the usual matrix Hermitian
conjugation $\boldsymbol{i}_{x}(b^{\star })=\boldsymbol{i}_{x}(b)^{\star }$
and non-usual multiplication given by the Hudson-Parthasarathy table \ \cite%
{26} 
\begin{equation}
\boldsymbol{i}_{x}(a\cdot c)=\left( 
\begin{array}{ll}
\mathrm{k}_{x}^{\star }(a)\mathrm{k}_{x}(c), & \mathrm{k}_{x}^{\star
}(a)i_{x}(c) \\ 
i_{x}(a)\mathrm{k}_{x}(c), & i_{x}(a)i_{x}(c)%
\end{array}%
\right) ,\quad \forall a,c\in \mathfrak{b}.  \label{three c}
\end{equation}%
It has a natural realization $\mathbf{i}(x,g)=\mathbf{j}(x,g)-\mathbf{I}_{x}$
given in the pseudo-Euclidean poly-Banach space $\mathbb{K}=L^{1}(X)\oplus 
\mathrm{K}\oplus L^{\infty }(X)$ by the canonical triangular representation $%
\mathbf{j}(x,g)=\mathbf{j}_{x}(g(x))$ of the $\star $-monoid $\mathfrak{m}$
with the usual matrix multiplication and non-usual pseudo-Hermitian
conjugation (\ref{two b}): 
\begin{equation}
\mathbf{i}(x,g^{\star })=\left[ 
\begin{array}{ccc}
0 & \mathrm{k}(x,g)^{\ast } & l(x,g)^{\ast } \\ 
0 & i(x,g^{\star }) & \mathrm{k}(x,g^{\star }) \\ 
0 & 0 & 0%
\end{array}%
\right] =\mathbf{i}(x,g)^{\dagger }.  \label{three d}
\end{equation}%
All that has been said means that the factor in $\mathfrak{m}/\mathfrak{i}$,
with zero $\star $-ideal $\mathfrak{i}=\{g\in \mathfrak{m}:\boldsymbol{i}%
(g)=0\}\equiv \boldsymbol{i}_{x}^{-1}(0)$, of step functions with values $%
g(x)\in \mathfrak{i}_{x}$, where $\boldsymbol{i}_{x}^{-1}(0)=\{b\in 
\mathfrak{b}:\boldsymbol{i}_{x}(b)=0\}$, can be described like $\boldsymbol{i%
}(\mathfrak{m})$ by four-component functions $\boldsymbol{g}=(g_{\nu }^{\mu
})_{\nu =\circ ,+}^{\mu =\circ ,-}$, for example of the form $\boldsymbol{g}%
(x)=\boldsymbol{i}_{x}(g(x))$ with $g_{\circ }^{\circ }=i(g)$, $g_{+}^{\circ
}=\mathrm{k}(g)$, $g_{\circ }^{-}=\mathrm{k}^{\star }(g)$, $g_{+}^{-}=l(g)$.
These form a $\star $-ring with respect to the Hermitian conjugation $\left( 
\boldsymbol{g}^{\star }(x)\right) _{\nu }^{\mu }=g_{-\mu }^{-\nu }(x)^{\ast
} $ and the table of componentwise multiplication (\ref{three c}). This
allows us to represent additive integral Hermitian forms 
\begin{equation}
\mu (g)=\int m(x,g)\mathrm{d}x,\;\;\;\;\,m(x,g)=m_{x}(g(x))  \label{three e}
\end{equation}%
on the $\star $-ring $\mathfrak{m}$ by four-component functions 
\begin{equation*}
m(x)=\left( 
\begin{array}{ll}
\mu & m_{\circ }^{+} \\ 
m_{-}^{\circ } & m_{\circ }^{\circ }%
\end{array}%
\right) (x),\,%
\begin{array}{l}
\mu (x,g)=\mu (x)g_{+}^{-}(x),\,m_{\circ }^{+}(x,g)=m_{\circ
}^{+}(x)g_{+}^{\circ }(x), \\ 
m_{-}^{\circ }(x,g)=g_{\circ }^{-}(x)m_{-}^{\circ }(x),\,m_{\circ }^{\circ
}(x,g)=\left\langle m_{\circ }^{\circ }(x),g_{\circ }^{\circ
}(x)\right\rangle ,%
\end{array}%
\end{equation*}%
in the form 
\begin{equation}
m(x,g)=\left\langle m_{\circ }^{\circ }(x),i(x,g)\right\rangle +\mathrm{m}%
^{\ast }(x)\mathrm{k}(x,g)+\mathrm{k}^{\star }(x,g)\mathrm{m}(x)+\mu
(x)l(x,g).  \label{three f}
\end{equation}%
Here $\mu (x)\in \mathbb{R}$ \textup{(}for almost all $x$\textup{)}, $%
m_{-}^{\circ }(x):\mathcal{K}_{x}\rightarrow \mathbb{C}$ is a vector linear
form $m_{-}^{\circ }=\mathrm{m}$ on the pre-Hilbert space $\mathcal{K}_{x}=\{%
\mathrm{k}_{x}^{\star }(b):b\in \mathfrak{b}\}=\mathrm{K}_{x}^{\ast }$,
adjoint to the form $m_{\circ }^{+}(x):\mathrm{K}_{x}\ni \mathrm{k}\mapsto 
\mathrm{m}^{\ast }(x)\mathrm{k}\in \mathbb{C}$, and $m_{\circ }^{\circ }(x):%
\mathcal{B}_{x}\ni B\mapsto \left\langle m_{\circ }^{\circ
}(x),B\right\rangle \in \mathbb{C}$ is an operator linear form on the $\ast $%
-subalgebra $\mathcal{B}_{x}=\{i_{x}(b):b\in \mathfrak{b}\}$ of operators $%
B,B^{\ast }:\mathrm{K}_{x}\rightarrow \mathrm{K}_{x}$. As the next theorem
shows, we thus essentially exhaust all linear positive logarithmic forms $%
\mu :\mathfrak{m}\rightarrow \mathbb{C}$ of infinitely divisible states $%
\psi (g)=e^{\mu (g)}$ on $\star $-algebras $\mathfrak{m}$, absolutely
continuous with respect to the Poisson state $\varphi (g)=e^{\left\langle
g\right\rangle }$ in the sense that $\mathfrak{i}\subseteq \mathfrak{i}^{\mu
}$. Here $\mathfrak{i}^{\mu }$ is the $\star $-ideal of step functions $%
g:X\ni x\mapsto g(x)\in \mathfrak{i}_{x}^{\mu }$ with values in two-sided
ideals 
\begin{equation}
\mathfrak{i}_{x}^{\mu }=\{b\in \mathfrak{b}:m_{x}(b)=0,m_{x}(ab)=0,%
\,m_{x}(bc)=0,\,m_{x}(abc)=0,\quad \forall a,c\in \mathfrak{b}\}.
\label{three g}
\end{equation}

\begin{theorem}
Suppose that $\mathfrak{b}$ is a $\star $-algebra over $\mathbb{C}$ and
suppose that the linear positive form \textup{(\ref{one b})} on the $\star $%
-algebra $\mathfrak{m}$ satisfies the condition 
\begin{equation*}
\forall g\in \mathfrak{m}\mathcal{\;}\exists c<\infty :\left\langle h\cdot
\left( g\star g\right) \cdot h^{\star }\right\rangle \leq c\left\langle
h\star h\right\rangle ,\quad \forall h\in \mathfrak{m}\text{,}
\end{equation*}%
of boundedness $\left\Vert i(g)\right\Vert \leq c$ of the associated
operator representation $i(g)=j(g)-I$. We equip $\mathfrak{m}$ with the
inductive convergence $g_{n}\rightarrow 0$ if $\left\Vert g_{n}\right\Vert
_{p}^{\Delta }\rightarrow 0$ for all $p=1,2,\infty $ and for some integrable 
$\Delta \in \mathfrak{F}$, where $g_{n}\in \mathfrak{m}_{\Delta }$ for all $%
n,\left\Vert g\right\Vert _{\infty }^{\Delta }=\left\Vert i(g)\right\Vert $
for $\{x\in X:g(x)\neq 0\}\subseteq \Delta $, and 
\begin{equation*}
\left\Vert g\right\Vert _{2}^{\Delta }=\left( \int_{\Delta }\left\Vert 
\mathrm{k}(x,g)\right\Vert ^{2}\,\mathrm{d}x\right) ^{1/2},\quad \left\Vert
g\right\Vert _{1}^{\Delta }=\int_{\Delta }\left\vert l(x,g)\right\vert \,%
\mathrm{d}x
\end{equation*}%
Then the following conditions are equivalent:

\begin{enumerate}
\item[\textup{(i)}] The functional $\psi (g)=e^{\mu (g)}$, continuous with
respect to the inductive convergence on $\mathfrak{m}$, is a pseudo-Poisson
state described by an absolutely continuous function $\mu _{\Delta }\left(
b\right) =\mu (b_{\Delta })$ in the sense that $\mu _{\Delta }(b)=0$ for all 
$b\in \mathfrak{b}$ if $\Delta \in \mathfrak{F}$ and $\mu _{\Delta
}=\int_{\Delta }\mathrm{d}x=0$.

\item[\textup{(ii)}] The functional $\mu :\mathfrak{m}\rightarrow \mathbb{C}$
has the integral form \textup{(\ref{three e})}, where $m_{x}:\mathfrak{b}%
\rightarrow \mathbb{C}$ is the linear function \textup{(\ref{three f})}
defined almost everywhere on $X$ by a positive numerical function $\mu
(x)\geq 0,\mathrm{ess}\sup_{x\in \Delta }\mu (x)<\infty $ for all $\Delta
\in \mathfrak{F}$ with $\mu _{\Delta }=\int_{\Delta }\mathrm{d}x<\infty $, a
vector-function $\mathrm{m}$ on $X$ with values $\mathrm{m}(x)\in \mathrm{K}%
_{x}$ defined by the values 
\begin{equation*}
\mathrm{m}^{\ast }(x)\in \mathrm{K}_{x}^{\ast },\,\int_{\Delta }\left\Vert 
\mathrm{m}(x)\right\Vert _{x}^{2}\mathrm{d}x<\infty ,\quad \forall \Delta
\in \mathfrak{F}:\mu _{\Delta }=\int_{\Delta }\mathrm{d}x<\infty ,
\end{equation*}%
of continuous \textup{(}for almost all $x\in X$\textup{)} forms $\mathrm{m}%
^{\ast }(x)\mathrm{k}=\left\langle \mathrm{m}(x)\mid \mathrm{k}\right\rangle 
$ on the Hilbert spaces $\mathrm{K}_{x}=\mathcal{K}_{x}^{\ast }$, and the
function $m_{\circ }^{\circ }$ on $X$ with values 
\begin{equation*}
m_{\circ }^{\circ }(x)\in \mathcal{B}_{x}^{\ast },\,\int_{\Delta }\underset{%
0\leq B\leq 1_{x}}{\sup }\left\langle m_{\circ }^{\circ }(x),B\right\rangle 
\mathrm{d}x<\infty ,\,\forall \Delta \in \mathfrak{F}:\mu _{\Delta
}=\int_{\Delta }\mathrm{d}x<\infty ,
\end{equation*}%
in positive forms on $C^{\ast }$-algebras $\mathcal{B}_{x}$ satisfying
almost everywhere the inequality 
\begin{equation}
\mu (x)\left\langle m_{\circ }^{\circ }(x),B^{\ast }B\right\rangle \geq
\left\Vert B\mathrm{m}(x)\right\Vert ^{2},\quad \forall B\in \mathcal{B}_{x}.
\label{three i}
\end{equation}

\item[\textup{(iii)}] There is a triangular representation 
\begin{equation*}
g\in \mathfrak{m}\mapsto \mathbf{g}(x)=\{g_{\nu }^{\mu
}(x)\},\,\;\;\;g_{-}^{\mu }=0=g_{\nu }^{+},\quad \forall \mu ,\nu \in
\{-,\circ ,+\},
\end{equation*}%
of the $\star $-algebra $\mathfrak{m}$ in the Banach space $\mathbb{K}%
=L^{1}(X)\oplus \mathcal{K}\oplus L^{\infty }(X)$ with indefinite metric 
\textup{(\ref{two a})} defined by the scalar product $\left\langle k^{\circ
}\mid k^{\circ }\right\rangle =\int \left\Vert k^{\circ }(x)\right\Vert
_{x}^{2}\mathrm{d}x$ of the Hilbert space $\mathcal{K}=\int^{\oplus }%
\mathcal{K}_{x}\,\mathrm{d}x$. This representation is locally
pseudo-unitarily equivalent to the canonical representation \textup{(\ref%
{three d})} in the sense that $\mathbf{g}(x)=\mathbf{S}^{\dagger }(x)\mathbf{%
i}(x,g)\mathbf{S}(x)$ for decomposable operators $\mathbf{S}(x)$ in $\mathbb{%
C}\oplus \mathcal{K}_{x}\oplus \mathbb{C}$ of the form \textup{(\ref{one j})}%
, and is such that 
\begin{equation}
\mu (g)=\int (\mu (x)g_{+}^{-}(x)+\left\langle M(x),g_{\circ }^{\circ
}(x)\right\rangle )\,\mathrm{d}x,\quad \forall g\in \mathfrak{m}\text{,}
\label{three j}
\end{equation}%
where $\mu \geq 0$ is a locally bounded measurable function and $M\geq 0$ is
a locally integrable function with positive values $M(x)\in \mathcal{B}%
_{x}^{\ast }$.
\end{enumerate}
\end{theorem}

\begin{proof}
First of all we notice that if the decomposable operator-functions $i_{x}(b)$
are locally bounded, then the space $\mathrm{K}$ of the canonical
representation $j(g)=I+i(g)$ of the $\star $-monoid $\mathfrak{m}$ of step
functions $g:X\rightarrow \mathfrak{b}$, complete with respect to the family
of seminorms (\ref{one g}), is a Hilbert space. This follows from the
inequality 
\begin{equation*}
\left\Vert \mathrm{k}\right\Vert ^{h}=\left\Vert j(h)\mathrm{k}\right\Vert
<\left\Vert \mathrm{k}\right\Vert +\left\Vert i(h)\mathrm{k}\right\Vert \leq
(1+\left\Vert h\right\Vert )\left\Vert \mathrm{k}\right\Vert ,
\end{equation*}%
where $\left\Vert f\right\Vert =\max_{i}\left\Vert b_{i}\right\Vert _{\Delta
(i)}<\infty $ according to (\ref{three f}) for any step integrable function $%
f(x)=b_{i},x\in \Delta (i)$, given by a finite partition $\Delta =\Sigma
\Delta (i)$ of its support $\Delta =\{x\in X:f(x)\neq 0\}$.

We shall first prove the simple implications (iii) $\Rightarrow $ (ii) $%
\Rightarrow $ (i), and next we shall construct the representation (\ref%
{three i}) of (iii) drawing on the conditions formulated in (i).

(iii) $\Rightarrow $ (ii). Suppose that $\mathbf{S}(x)$ is the triangular
transformation of the form (\ref{one j}) given by an essentially measurable
function $U(x)\in \mathcal{L}(\mathrm{K}_{x})$ with unitary values, a
function $e_{\circ }:X\ni x\mapsto e_{\circ }(x)\in \mathrm{K}_{x}^{\ast }$
given by the values $e_{\circ }^{\ast }(x)\in \mathrm{K}_{x}$ of a
vector-function $e_{\circ }^{\ast }$ with $\int_{\Delta }\left\Vert e_{\circ
}^{\ast }(x)\right\Vert _{x}^{2}\mathrm{d}x<\infty $ for all $\Delta $ such
that $\mu _{\Delta }=\int_{\Delta }\mathrm{d}x<\infty $, and a scalar
locally integrable function $e_{+}$ such that $e_{+}(x)+e_{+}^{\ast
}(x)=-\left\Vert e_{\circ }^{\ast }(x)\right\Vert _{x}^{2}$. Then $g_{\circ
}^{\circ }(x)=U^{\ast }(x)i(x,g)U(x)$, 
\begin{eqnarray*}
g_{+}^{-}(x)=e_{\circ }(x)U(x)\mathrm{k}(x,g) &&+e_{\circ
}(x)U(x)i(x,g)U^{\ast }(x)e_{\circ }^{\ast }(x) \\
&&+\mathrm{k}^{\star }(x,g)U^{\ast }(x)e_{\circ }^{\ast }(x),
\end{eqnarray*}%
and (\ref{three j}) takes the form (\ref{three e}) (\ref{three f}), where $%
\mathrm{m}^{\ast }(x)=e_{\circ }(x)U(x),\mathrm{m}(x)=U^{\ast }(x)e_{\circ
}^{\ast }(x)$ is a locally square integrable function: $\int_{\Delta
}\left\Vert \mathrm{m}(x)\right\Vert _{x}^{2}\mathrm{d}x<\infty $, and 
\begin{equation*}
\left\langle m_{\circ }^{\circ }(x),B\right\rangle =\left\langle
M(x),U^{\ast }(x)BU(x)\,\right\rangle +\mu (x)\mathrm{m}^{\ast }(x)B\mathrm{m%
}(x)
\end{equation*}%
is a positive locally integrable function: $\int_{\Delta }\left\langle
m_{\circ }^{\circ }(x),B_{x}\right\rangle \mathrm{d}x<\infty $, satisfying (%
\ref{three i}) by virtue of the positivity $\left\langle M(x),B^{\ast
}B\right\rangle \geq 0,\,\mu (x)\geq 0$ for all $x\in X$.

(ii) $\Rightarrow $ (i). If $\mu $ is the integral (\ref{three e}) of the
linear form (\ref{three f}) and (\ref{three i}) is fulfilled, then $\mu
(g^{\star }g)\geq 0$ for all $g\in \mathfrak{m}$, since$\;\ $ 
\begin{eqnarray*}
\lefteqn{0\leq m(g^{\star }g)=\left\langle m_{\circ }^{\circ },i(g^{\star
}g)\right\rangle +2\func{Re}\mathrm{m}^{\ast }\mathrm{k}(g^{\star }g)+\mu
l(g^{\star }g)} \\
&=&\left\langle m_{\circ }^{\circ },i(g)^{\ast }i(g)\right\rangle +2\func{Re}%
\mathrm{m}^{\ast }i(g)^{\ast }\mathrm{k}(g)+\mu (x)\mathrm{k}(g)^{\ast }%
\mathrm{k}(g) \\
&=&\left\langle m_{\circ }^{\circ },i(g)^{\ast }i(g)\right\rangle -\tfrac{1}{%
\mu }\left\Vert i(g)\mathrm{m}\right\Vert ^{2}+\mu ||\mathrm{k}(g)+i(g)%
\mathrm{m}/\mu ||^{2}.
\end{eqnarray*}%
By linearity of $\mu $ this is equivalent to positive definiteness 
\begin{equation*}
\sum_{a,c}\kappa _{a}^{\ast }\mu _{\Delta }(a^{\star }c)\kappa _{c}=\mu
_{\Delta }(\sum_{a,c}\kappa _{a}^{\ast }a^{\star }c\kappa _{c})=\mu _{\Delta
}(b^{\star }b)\geq 0
\end{equation*}%
of the form $\mu _{\Delta }(b)=\mu (b_{\Delta })$. Hence it follows that $%
\varphi _{\mu }(g)=e^{\mu (g)}$ is a pseudo-Poisson state given by an
absolutely continuous complex measure $\mu _{\Delta }\left( b\right)
=\int_{\Delta }m_{x}(b)\mathrm{d}x$ with density $m_{x}(b)=m(x,b_{\Delta })$%
. It is continuous with respect to the inductive convergence with respect to
the seminorms $\left\Vert g\right\Vert _{p}^{\Delta },p=1,2,\infty $,
because $\mu $ is locally bounded, $\mathrm{m}$ is locally $L^{2}$%
-integrable, and $m_{\circ }^{\circ }$ is locally $L^{1}$-integrable.

(i) $\Rightarrow $ (iii). If the function $\mu _{\Delta }(b)=\ln \psi (b)$
for a pseudo-Poisson state $\psi _{\Delta }(b)=\psi (b_{\Delta })$ on $%
\mathfrak{b}$ is absolutely continuous with respect to $\Delta \in \mathfrak{%
F}$ for every $b\in \mathfrak{b}$,\textsf{\ }then it has the form (\ref%
{three e}), where the density $m_{x}:\mathfrak{b}\rightarrow \mathbb{C}$ is
almost everywhere a linear positive functional. Since the kernel $\{g\in 
\mathfrak{m}:\,\left\Vert g\right\Vert _{p}=0,p=1,2,\infty \}$ of the
inductive convergence in $\mathfrak{m}=\cup \mathfrak{m}_{\Delta }$
coincides with the kernel $\mathfrak{i}$ of the canonical representation $%
\mathbf{i}(g)=\mathbf{j}(g)=\mathbf{I}$ in $\mathbb{K}$, which is equal,
according to its construction, to step functions $g:x\mapsto g(x)\in 
\mathfrak{i}_{x}$, where 
\begin{equation*}
\mathfrak{i}_{x}=\{b\in \mathfrak{b}:l_{x}\,(b)=0,l_{x}(ab)=0,%
\,l_{x}(bc)=0,l_{x}(abc)=0,\quad \forall a,b,c\in \mathfrak{b}\},
\end{equation*}%
the $\star $-ideal $\mathfrak{i}^{\mu }$ of functions $g\in \mathfrak{m}$
with values $g(x)$ in (\ref{three g}), corresponding to the form (\ref{three
e}) continuous in the sense that $g_{n}\rightarrow 0\Rightarrow \mu
(g_{n})\rightarrow 0$, necessarily contains $\mathfrak{i}$. This means that
a linear functional $m_{x}(b)$ that vanishes on $\mathfrak{i}_{x}$ for
almost all $x$ can be written, by this continuity, in the form (\ref{three f}%
) of a linear Hermitian functional $m_{x}(b)=m(x,b_{\Delta }),x\in \Delta $,
on the factor algebra $\mathfrak{b}/i_{x}^{-1}\left( 0\right) $ isomorphic
to the $\star $-subalgebra $i_{x}(\mathfrak{b})$ of quadruples (\ref{three b}%
) with the multiplication table (\ref{three c}). In addition, by the
Hahn-Banach theorem and the duality between $L^{p}(\Delta )$ and $%
L^{q}(\Delta )$ for $1/p+1/q=1$ we can assume that $\mu $ is locally
bounded, $\mathrm{m}$ is locally $L^{2}$-integrable, and $m_{\circ }^{\circ
} $ is locally $L^{1}$-integrable on $X$.

For every $x\in X$ we define a triangular pseudo-unitary transform of $%
\mathbf{S}(x)$ into $\mathbb{K}_{x}=\mathbb{C}\oplus \mathrm{K}_{x}\oplus 
\mathbb{C}$ of the form (\ref{one j}), where $U=-I_{x},e_{\circ }^{\ast }(x)=%
\mathrm{m}(x)$, and $e_{+}^{\ast }(x)=-\left\Vert \mathrm{m}(x)\right\Vert
_{x}^{2}/2$. Denoting $g_{\nu }^{\mu }(x)=(\mathbf{S}^{\dagger }\mathbf{i}%
(x,g)\mathbf{S}(x))_{\nu }^{\mu }$, where $\mathbf{i}(x)$ is the triangular
matrix representation (\ref{three d}) of the quadruple (\ref{three b}) for $%
b=g(x)$, we obtain 
\begin{equation*}
m(g)=\left\langle g_{\circ }^{\circ },m_{\circ }^{\circ }\right\rangle -%
\mathrm{m}^{\ast }g_{\circ }^{\circ }\mathrm{m}/\mu +\mu g_{+}^{-},
\end{equation*}%
where we have taken into account the fact that $g_{\circ }^{\circ
}(x)=i_{x}(g(x))$ and 
\begin{equation*}
\mu g_{+}^{-}=\mu l(g)+\mathrm{k}^{\star }(g)\mathrm{m}+\mathrm{m}^{\ast }%
\mathrm{k}(g)+\mathrm{m}^{\ast }i(g)\mathrm{m}/\mu .
\end{equation*}%
In this representation the positivity condition $m(x,g^{\star }g)\geq 0$
takes the form 
\begin{equation*}
\left\langle g_{\circ }^{\circ \ast }g_{\circ }^{\circ },M\right\rangle +\mu
g_{+}^{\circ \ast }g_{+}^{\circ }\geq 0,\quad \forall g\in \mathfrak{m},
\end{equation*}%
where $\left\langle B,M\right\rangle =\left\langle B,m_{\circ }^{\circ
}\right\rangle -\mathrm{m}^{\ast }B\mathrm{m}/\mu $, $B\in \mathcal{B}_{x}$,
and $g_{+}^{\circ }=\mathrm{k}_{x}(g)+i_{x}(g)\mathrm{m}$. The resulting
inequality proves that $M(x)$ is positive for $g_{+}^{\circ }(x)=0$ and $\mu
(x)\geq 0$ if $g_{\circ }^{\circ }(x)=0$. This proves the existence of
locally bounded measurable functions $\mu \geq 0$ and positive locally
integrable functions $M$ with values $M(x)\in \mathcal{B}_{x}^{\ast }$
defining the function $\mu (g)$ in the form (\ref{three j}). The proof is
complete.
\end{proof}

\begin{remark}
We consider an additive subgroup $\mathfrak{b}\subseteq \mathbb{C}\times
H\times \mathcal{L}(K)$ of the triples $b=(\beta ,\eta ,B)$ with the
involution $b^{\star }=(\beta ^{\ast },\eta ^{\#},B^{\ast })$, where $\beta
\mapsto \beta ^{\ast }\in \mathbb{C}$ is the complex conjugation, $\eta
\mapsto \eta ^{\#}\in H$ is the involution $\eta ^{\#\#}=\eta $ in a $%
\mathbb{C}$-linear subspace $H\subseteq K$ equipped with the Hermitian form $%
\left\langle \xi \mid \zeta \right\rangle =\xi ^{\#}\cdot \zeta
=\left\langle \xi \mid \zeta \right\rangle ^{\ast }$ of a pseudo-Euclidean
space $K$, and $B\mapsto B^{\ast }\in \mathcal{L}(K)$ is the Hermitian
conjugation $\left\langle B^{\ast }\xi \mid \zeta \right\rangle
=\left\langle \xi \mid B\zeta \right\rangle $ for all $\xi ,\zeta \in K$ in
the $\ast $-subalgebra $\mathcal{L}\subseteq \mathcal{L}(K)$ of operators $%
B:\eta \mapsto B\eta \in K$ leaving $H$ invariant: $BH\subseteq H$ for all $%
B\in \mathcal{L}$.

We define in $\mathfrak{b}$ the structure of a $\star $-algebra by putting 
\begin{equation*}
\lambda b=(\lambda \beta ,\lambda \eta ,\lambda B),a^{\star }c=(\xi
^{\#}\cdot \zeta ,\xi ^{\#}C+A^{\ast }\zeta ,A^{\ast }C)
\end{equation*}%
for any $\lambda \in \mathbb{C},b\in \mathfrak{b},a=(\alpha ,\xi
,A),c=(\gamma ,\zeta ,C)$, where we use the notation $\xi ^{\#}C=(C^{\ast
}\xi )^{\#}$. It is easy to prove that this distributive algebra is
associative, $(ab)c=a(bc)$, only in the case 
\begin{equation*}
(A\eta )\cdot \zeta =\xi \cdot (\eta C),(A\eta )C=A(\eta C),\quad \forall
A,\,C\in \mathcal{L},\xi ,\eta ,\zeta \in H,
\end{equation*}%
which is possible only under the condition $(A\eta )\cdot \zeta =0=\xi \cdot
(\eta C)$. This condition leads to $(A\eta )C=A(\eta C)$ if $\xi \cdot \zeta
=(\xi ^{\#}\mid \zeta )$ is a bilinear form on $H$ non-degenerate in the
sense that $\{\xi \cdot \eta =0=\eta \cdot \zeta :\forall \xi ,\zeta \in
H\}\Rightarrow \eta =0$. A simple analysis of the positivity 
\begin{equation*}
l(b^{\star }b)=\lambda \left\langle \eta \mid \eta \right\rangle
+\left\langle B\vartheta \mid \eta \right\rangle +\left\langle \eta \mid
B\vartheta \right\rangle +\left\langle \Lambda ,B^{\ast }B\right\rangle \geq
0
\end{equation*}%
of the linear $\star $-form $l(b)=\lambda \beta +\vartheta _{-}\cdot \eta
+\eta \cdot \vartheta _{+}+\left\langle \Lambda ,B\right\rangle $, where $%
\lambda =\lambda ^{\ast }$,$\;\vartheta _{+}=\vartheta =\vartheta _{-}^{\#}\ 
$,$\;\Lambda =\Lambda ^{\ast }$, leads to the conditions $\left\langle
\Lambda ,B^{\ast }B\right\rangle \geq 0$ for all $B\in \mathcal{L}$ if $%
\lambda =0$ and 
\begin{equation*}
\lambda \left\langle \eta \mid \eta \right\rangle \geq 0,\,\left\langle
\Lambda ,B^{\ast }B\right\rangle \geq \tfrac{1}{\lambda }\left\langle
B\vartheta \mid B\vartheta \right\rangle ,\quad \forall \eta \in H,B\in 
\mathcal{L}
\end{equation*}%
if $\lambda \neq 0$. The latter is possible only if the form $\left\langle
\eta \mid \eta \right\rangle =\eta ^{\#}\cdot \eta $ is definite, that is, $%
\lambda >0$ if $\left\langle \eta \mid \eta \right\rangle \geq 0$ for all $%
\eta \in H$ and $\lambda <0$ if $\left\langle \eta \mid \eta \right\rangle
\leq 0$ for all $\eta \in H$, which is a necessary condition for the
existence of a pseudo-Poisson state on $\mathfrak{b}=\mathbb{C}\times
H\times \mathcal{L}$.

Assuming without loss of generality that $\eta ^{\#}\eta \geq 0$ for all $%
\eta $ \textup{(}otherwise we have to change the notation $b\mapsto (-\beta
,\eta ,B)$ and $\eta ^{\#}\eta \mapsto -\eta ^{\#}\eta $\textup{)} we
consider the following two cases, in which $H$ is a Hilbert space with
respect to the norm $\left\Vert \eta \right\Vert =\left( \eta ^{\#},\eta
\right) ^{1/2}$, where $\left( \xi ,\zeta \right) =\frac{1}{2}\left( \xi
\cdot \zeta +\zeta \cdot \xi \right) $.
\end{remark}

\begin{example}[Gaussian state]
Let $\mathcal{L}=\{0\}$ and $\lambda =1$, that is, $b=(\beta ,\eta )$, and
let $l(b)=\left\langle \eta ,\theta \right\rangle +\beta $, where $%
\left\langle \eta ,\theta \right\rangle =2\func{Re}(\eta \mid \theta )$ for
all $\eta =\eta ^{\#}$. The algebra $\mathfrak{b}=\mathbb{C}\times H$ is now
nilpotent: $ac=(\xi ,\zeta ,0),\,abc=(0,0)$ for all $a,b,c\in \mathfrak{b}$,
and commutative, $[a,c]=ac-ca=0$, if the involution $\#$ is isometric on $H$
in $K\supseteq H$: 
\begin{equation*}
\left\langle \xi ^{\#}\mid \zeta \right\rangle =\left\langle \zeta ^{\#}\mid
\xi \right\rangle ,\quad \forall \xi ,\zeta \in H.
\end{equation*}%
The infinitely divisible functional $\varphi _{\Delta }(b)=\exp \{\left[
\beta +\left( \eta ,\theta \right) \right] \mu _{\Delta }\}$ corresponding
to the conditionally positive $\star $-form $\lambda _{\Delta }(b)=\left[
\beta +\left( \eta ,\theta \right) \right] \mu _{\Delta }$ with respect to
the Hermitian operation 
\begin{equation*}
(\alpha ,\xi )\star (\gamma ,\zeta )=(\alpha ^{\ast }+\left\langle \xi \mid
\zeta \right\rangle +\gamma ,\xi ^{\#}+\zeta ),\;\;\;\;u=\,(0,0),
\end{equation*}%
defines a generating functional $\varphi _{\Delta }(0,\eta )=1$ of the
factorial moments of a Gaussian chaotic state over $H$ with mathematical
expectation $\left\langle b_{\Delta }\right\rangle =\left( \eta ,\theta
\right) \mu _{\Delta }$ for $b=(0,\eta )$ and finite covariance $%
\left\langle b_{\Delta }^{\star }b_{\Delta }\right\rangle =\left\langle \eta
\mid \eta \right\rangle \mu _{\Delta }\in \mathbb{R}_{+}$ for every $\Delta
\in \mathfrak{F}$ such that $\mu _{\Delta }=\int_{\Delta }\mathrm{d}x<\infty 
$. This covariance is symmetric only in the commutative \textup{(}classical%
\textup{)} case, and in the converse \textup{(}quantum\textup{)} case it
satisfies the uncertainty relation 
\begin{equation*}
\left\langle a_{\Delta }^{2}\right\rangle \left\langle c_{\Delta
}^{2}\right\rangle \geq s(\xi ,\zeta )^{2}\mu _{\Delta }^{2},\quad \forall
a=(\alpha ,\xi ),c=(\gamma ,\zeta ),\xi ,\zeta \in \mathbf{\func{Re}}H
\end{equation*}%
for the commutative Heisenberg relation $[a_{\Delta },c_{\Delta }]=(\mathrm{i%
}s(\zeta ,\xi )\mu _{\Delta }0)$ corresponding to the symplectic form $s(\xi
,\zeta )=2\func{Im}(\xi \mid \zeta )$ on $\mathbf{\func{Re}}H=\{\eta \in
H:\eta ^{\#}=\eta \}$. The canonical representation \textup{(\ref{three d})}%
, defining a $\star $-representation $\mathbf{j}(g)=\mathbf{I+i}(g)$ of the $%
\star $-monoid $\mathfrak{m}$ of step functions $g:X\rightarrow \mathbb{C}%
\times H$, and the corresponding representation $\pi (g)=\epsilon \lbrack 
\mathbf{j}^{\otimes }(g)]$ in the Fock space $\mathrm{K}$, is described by
the functions $i_{x}(b)=0$,$~~\mathrm{k}_{x}(b^{\star })=\eta ^{\#}$, $~%
\mathrm{k}_{x}(b)^{\ast }=\eta ^{\ast }$,$~l_{x}(b)=\beta +\left( \eta
,\theta \right) $.
\end{example}

\begin{example}[Poisson state]
Let $H=\{0\}$ and let $\mathfrak{b}=\mathcal{L}$ be the $\ast $-algebra of
operators in $K$ bounded by the identity $I\in \mathfrak{b}$ in the sense
that 
\begin{equation*}
\forall C=B^{\ast }B\;\;\,\exists c\in \mathbb{R}_{+}:\left\langle \Lambda
,A^{\ast }CA\right\rangle \leq c\left\langle \Lambda ,A^{\ast
}A\right\rangle ,\quad \forall A\in \mathfrak{b},
\end{equation*}%
where $\Lambda $ is a linear positive form defining $l(b)=\left\langle
\Lambda ,B\right\rangle $. Bearing in mind the Gelfand-Naimark-Segal
construction, we may assume without loss of generality that this form is a
vector one, $\left\langle \Lambda ,B\right\rangle =\left\langle
eBe\right\rangle $, represented in the Hilbert space $K$ by an element $e\in
K,\left\Vert e\right\Vert ^{2}=\left\langle \Lambda ,I\right\rangle $. In
the commutative case $\mathfrak{b}$ can be identified with a subalgebra of
essentially bounded functions $b:\omega \mapsto b(\omega )\in \mathbb{C}$ on
a measurable space $\Omega $ with finite positive measure $\mathrm{d}\lambda 
$ of the mass $\lambda =\left\langle \Lambda ,I\right\rangle $ by putting $(B%
\mathrm{k})(\omega )=b(\omega )\mathrm{k}(\omega )$ on $K=L_{\lambda
}^{2}(\Omega )$, and $e(\omega )=1$ for all $\omega \in \Omega $, so that $%
l(b)=\int b(\omega )\mathrm{d}\lambda $. The infinitely divisible functional 
$\varphi _{\Delta }(b)=e^{\left\langle \Lambda ,B\right\rangle \mu _{\Delta
}}$, corresponding to the conditionally positive $\star $-form $\lambda
_{\Delta }(b)=\left\langle \Lambda ,B\right\rangle \mu _{\Delta }$ with
respect to the Hermitian operation $A\cdot C=A^{\ast }+A^{\ast }C+C$ with
the neutral element $U=0$, defines the generating functional of factorial
moments of a Poisson chaotic state over $\mathcal{L}$ with mathematical
expectation $\left\langle b_{\Delta }\right\rangle =\left\langle \Lambda
,B\right\rangle \mu _{\Delta }$ and finite covariance $\left\langle
b_{\Delta }^{\star }b_{\Delta }\right\rangle =\left\langle \Lambda ,B^{\ast
}B\right\rangle \mu _{\Delta }\in \mathbb{R}_{+}$ for each $\Delta \in 
\mathfrak{F}$ such that $\mu _{\Delta }=\int_{\Delta }\mathrm{d}x<\infty $.
This covariance is symmetric not only in the commutative \textup{(}classical%
\textup{)} case $[A,C]=AC-CA=0$, but also in the case when $\Lambda \in 
\mathcal{L}^{\ast }$ is central. The central form $\left\langle \Lambda
,B\right\rangle $, described by the condition $\left\langle \Lambda
,[A,C]\right\rangle =0$ for all $A,C\in \mathcal{L}$, defines a $\sigma $%
-finite trace on the $\ast $-algebra $\mathfrak{m}$ of step functions $%
G:X\ni x\mapsto G(x)\in \mathcal{L}$ with the integral form $\left\langle
g\right\rangle =\int \left\langle \Lambda ,G(x)\right\rangle \mathrm{d}x$ or 
$\left\langle g\right\rangle =\iint g(x,\omega )\mathrm{d}x\mathrm{d}\lambda 
$ in the case of $\mathfrak{b}\sim L_{\gamma }^{\infty }(\Omega )$.
Otherwise, the form $\left\langle \Lambda ,B\right\rangle $ can also lead to
the uncertainty relation 
\begin{equation*}
\left\langle a_{\Delta }^{2}\right\rangle \left\langle c_{\Delta
}^{2}\right\rangle \geq \left\langle \Lambda ,\frac{1}{\mathrm{i}}%
[A,C]\right\rangle ^{2}\mu _{\Delta }^{2},\quad \forall A=A^{\ast
},C=C^{\ast }.
\end{equation*}%
The canonical representation \textup{(\ref{three c})}, defining the
indefinite representation $\mathbf{j}(g)=\mathbf{I+i}(g)$ of the $\star $%
-monoid $\mathfrak{m}$ and the corresponding representation $\pi
(g)=\epsilon \lbrack \mathbf{j}^{\otimes }(g)]$ in the Fock space $\mathrm{K}
$, is described by the functions 
\begin{equation*}
i_{x}(b)=B,\;\mathrm{k}_{x}(b^{\star })=B^{\ast }e,\,\mathrm{k}_{x}^{\star
}(b)=e^{\ast }B,\,l_{x}(b)=e^{\ast }Be,
\end{equation*}%
where $e^{\ast }Be=\left\langle eBe\right\rangle =\left\langle \Lambda
,B\right\rangle $.
\end{example}

\part{Non-commutative stochastic analysis and quantum evolution in scales}

\section{Introduction}

Non-commutative generalization of the It\^{o} stochastic calculus, developed
in \cite{1}, \cite{3}, \cite{21}, \cite{36}, \cite{41} and \cite{44} gave an
adequate mathematical tool for studying the behavior of open quantum
dynamical systems singularly interacting with a boson quantum-stochastic
field. Quantum stochastic calculus also made it possible to solve an old
problem of describing such systems with continuous observation and
constructing a quantum filtration theory which would explain a continuous
spontaneous collapse under the action of such observation \cite{8}, \cite{11}
and \cite{12}. This gave examples of stochastic non-unitary, non-stationary,
and even non-adapted evolution equations in a Hilbert space whose solution
requires a proper definition of chronologically ordered quantum stochastic
semigroups and exponents of operators by extending the notion of the
multiple stochastic integral to non-commuting objects.

Here we outline the solution to this important problem by developing a new
quantum stochastic calculus in a natural scale of Fock spaces based on an
explicit definition, introduced by us in \cite{13}, of the non-adapted
quantum stochastic integral as a non-commutative generalization of the
Skorokhod integral \cite{48} represented in the Fock space. Using the
indefinite $\star $-algebraic structure of the kernel calculus, which was
obtained in the first chapter as a general property of a natural
pseudo-Euclidean representation associated with infinitely divisible states,
we establish the fundamental formula for the stochastic differential of a
function of a certain number of non-commuting quantum processes, providing a
non-commutative and non-adapted generalization of the It\^{o} formula as the
principal formula of the classical stochastic calculus. In the adapted case
this formula coincides with the well-known Hudson-Parthasarathy formula \cite%
{26} for the product of a pair of non-commuting quantum processes. In the
commutative case it gives a non-adapted generalization of the It\^{o}
formula for classical stochastic processes which was recently obtained in a
weak form by classical stochastic methods by Nualart \cite{42} in the case
of Wiener integrals. We also note that the classical stochastic calculus and
the calculus of operators in the Fock scales was worked out by the group
Hida, Kuo, Streit and Potthoff, see \cite{25} and \cite{45}, and also by
Berezanskii and Kondrat'ev \cite{19}.

Using the notion of a normal multiple quantum stochastic integral, which is
formulated below, we construct explicit solutions of quantum stochastic
evolution equations in the adapted as well as in the non-adapted case of
operator-valued coefficients and we give a simple algebraic proof of the
fact that this evolution is unitary if the generators of these equations are
pseudounitary. In the adapted stationary case the quantum stochastic
evolution was constructed by Hudson and Parthasarathy by means of the
approximation by the It\^{o} sums of quantum-stochastic generators. However,
proving unitarity by this method turned out to be a difficult problem even
in a simple case.

Within the framework of this approach Kholevo \cite{30} constructed a
solution of an adapted quantum-stochastic differential equation also for
non-stationary generators by defining the chronological exponential as a
quantum-stochastic multiplicative integral.

We note that our approach is close in spirit to the kernel calculus of
Maassen-Lindsay-Meyer \cite{36}, \cite{41}, however the difference is that
all the main objects are constructed not in terms of kernels but in terms of
operators represented in the Fock space. In addition we employ a much more
general notion of multiple stochastic integral, non-adapted in general,
which reduces to the notion of the kernel representation of an operator only
in the case of a scalar (non-random) operator function under the integral.
The possibility of defining a non-adapted single integral in terms of the
kernel calculus was shown by Lindsay \cite{37}, but the notion of the
multiple quantum-stochastic integral has not been discussed in the
literature even in the adapted case.

\section{Non-adapted stochastic integrals and differentials in Fock scale}

Let $(X,\mathfrak{F},\mu )$ be an essentially ordered space, that is, a
measurable space $X$ with a $\sigma $-finite measure $\mu :\mathfrak{F}\ni
\Delta \mapsto \mu _{\Delta }\geq 0$ and the ordering relation $x\leq
x^{\prime }$ with the property that any $n$-tuple $x=(x_{1},\ldots ,x_{n})$
is the chain $\varkappa =\{x_{1}<\cdots <x_{x}\}$ up to any permutation
modulo the product $\prod_{i=1}^{n}\mathrm{d}x_{i}$ of the measures $\mathrm{%
d}x:=\mu _{\mathrm{d}x}$. In other words, we assume that the measurable
ordering is almost linear, that is, or any $n$ the product measure of $n$%
-tuples $x\in X^{n}$ with components that are not completely increasingly
ordered is zero. Hence, in particular, it follows that the measure $\mu $ on 
$X$ is atomless. We may assume that the essential ordering on $X$ is induced
by a measurable map $t:X\rightarrow \mathbb{R}_{+}$ with respect to which $%
\mu $ is absolutely continuous in the sense of admitting the decomposition 
\begin{equation*}
\int_{\Delta }f(t(x))\mathrm{d}x=\int_{0}^{\infty }f(t)\mu _{\Delta }(t)%
\mathrm{d}t,
\end{equation*}%
for any integrable set $\Delta \subseteq X$ and any essentially bounded
function $f:\mathbb{R}_{+}\rightarrow \mathbb{C}$, where $\mu _{\Delta }(t)$
is a positive measure on $X$ for any $t\in \mathbb{R}_{+}$ and $x_{1}<\cdots
<x_{n}$ means that $t(x_{1})<\cdots <t(x_{n})$. In any case we shall assume
that we are given a map $t$ such that the above condition holds and $%
t(x)\leq t(x^{\prime })$ if $x\leq x^{\prime }$, understanding $t(x)$ as the
time at the point $x\in X$. For example, $t(x)=t$ for $x=(\mathbf{x},t)$ if $%
X=\mathbb{R}^{d}\times \mathbb{R}_{+}$ is the $(d+1)$-dimensional space-time
with the casual ordering \cite{5} and $\mathrm{d}x=\mathrm{d}\mathbf{x}%
\mathrm{d}t$, where $\mathrm{d}\mathbf{x}$ is the standard volume on $d$%
-dimensional space $\mathbb{R}^{d}\ni \mathbf{x}$.

We shall identify the finite chains $\varkappa $ with increasingly indexed $%
n $-tuples of $x_{i}\in X$, $\boldsymbol{x}=(x_{1},\ldots ,x_{n})$, $%
x_{1}<\cdots <x_{n}$, denoting by $\mathcal{X}=\sum_{n=0}^{\infty }\Gamma
_{n}$ the set of all finite chains as the union of the sets $\Gamma _{n}=\{%
\boldsymbol{x}\in X^{n}:x_{1}<\cdots <x_{n}\}$ with one-element $\Gamma
_{0}=\{\emptyset \}$ containing the empty chain as a subset of $X$: $%
\emptyset =X^{0}$. We introduce a measure `element' $\mathrm{d}\varkappa
=\prod_{x\in \varkappa }\mathrm{d}x$ on $\mathcal{X}$ induced by the direct
sum $\sum_{n=0}^{\infty }\mu _{\Delta _{n}}^{n},\Delta _{n}\in \mathfrak{F}%
^{\otimes n}$ of product measures $\mathrm{d}\boldsymbol{x}=\prod_{i=1}^{n}%
\mathrm{d}x_{i}$ on $X^{n}$ with the unit mass $\mathrm{d}\varkappa =1$ at
the point $\varkappa =\emptyset $.

Let $\{\mathrm{K}_{x}:x\in X\}$ be a family of Hilbert spaces $\mathrm{K}%
_{x} $, let $\mathcal{P}_{0}$ be an additive semigroup of positive
essentially measurable locally bounded functions $p:X\rightarrow \mathbb{R}%
_{+}$ with zero $0\in \mathcal{P}_{0}$, and let $\mathcal{P}%
_{1}=\{1+p_{0}:p_{0}\in \mathcal{P}_{0}\}$. for example, in the case $X=%
\mathbb{R}^{d}\times \mathbb{R}_{+}$ by $\mathcal{P}_{1}$ we mean the set of
polynomials $p(x)=1+\sum_{k=0}^{m}c_{k}|x|^{k}$ with respect to the modulus $%
|\mathbf{x}|=(\Sigma x_{i}^{2})^{1/2}$ of a vector $\mathbf{x}\in \mathbb{R}%
^{d}$ with positive coefficients $c_{k}\geq 0$. We denote by $\mathrm{K}(p)$
the Hilbert space of essentially measurable vector-functions $\mathrm{k}%
:x\mapsto \mathrm{k}(x)\in \mathrm{K}_{x}$ square integrable with the weight 
$p\in \mathcal{P}_{1}$: 
\begin{equation*}
\left\Vert \mathrm{k}\right\Vert (p)=\left( \int \left\Vert \mathrm{k}%
(x)\right\Vert _{x}^{2}p(x)\mathrm{d}x\right) ^{1/2}<\infty .
\end{equation*}%
Since $p\geq 1$, any space $\mathrm{K}(p)$ can be embedded into the Hilbert
space $\mathrm{K}=\mathrm{K}(1)$, and the intersection $\cap \mathrm{K}%
(p)\subseteq \mathrm{K}$ can be identified with the projective limit $%
\mathrm{K}^{+}=\lim_{p\rightarrow \infty }\mathrm{K}(p)$. This follows from
the facts that the function $\left\Vert \mathrm{k}\right\Vert (p)$ is
increasing: $p\leq q\Rightarrow \left\Vert \mathrm{k}\right\Vert (p)\leq
\left\Vert \mathrm{k}\right\Vert (q)$, and so $\mathrm{K}(q)\subseteq 
\mathrm{K}(p)$, and the set $\mathcal{P}_{1}$ is directed in the sense that
for any $p=1+r$ and $q=1+s$, $r,s\in \mathcal{P}_{0}$, there is a function
in $\mathcal{P}_{1}$ majorizing $p$ and $q$ (we can take for example $%
p+q-1=1+r+s\in \mathcal{P}_{1}$). In the case of polynomials $p\in \mathcal{P%
}_{1}$ on $X=\mathbb{R}^{d}\times \mathbb{R}_{+}$ the decreasing family $\{%
\mathrm{K}(p)\}$, where $\mathrm{K}_{x}=\mathbb{C}$, is identical with the
integer Sobolev scale of vector fields $\mathrm{k}:\mathbb{R}^{d}\rightarrow
L^{2}(\mathbb{R}_{+})$ with values $\mathrm{k}(x)(t)=\mathrm{k}(x,t)$ in the
Hilbert space $L^{2}(\mathbb{R}_{+})$ of square integrable functions on $%
\mathbb{R}_{+}$. If we replace $\mathbb{R}^{d}$ by $\mathbb{Z}^{d}$ and if
we restrict ourselves to the positive part of the integer lattice $\mathbb{Z}%
^{d}$, then we obtain the Schwartz space in the form of vector fields $k\in 
\mathrm{K}^{+}$.

The space $\mathrm{K}_{-}$, dual to $\mathrm{K}^{+}$, of continuous
functionals 
\begin{equation*}
\left\langle \mathrm{f}\mid \mathrm{k}\right\rangle =\int \left\langle 
\mathrm{f}(x)\mid \mathrm{k}(x)\right\rangle \,\mathrm{d}x,\quad \mathrm{k}%
\in \mathrm{K}^{+},
\end{equation*}%
is defined as the inductive limit $\mathrm{K}_{-}=\lim_{p\rightarrow 0}%
\mathrm{K}(p)$ in the scale $\{\mathrm{K}(p):p\in \mathcal{P}_{-}\}$, where $%
\mathcal{P}_{-}$ is the set of functions $p:X\rightarrow (0,1]$ such that $%
1/p\in \mathcal{P}_{1}$. The space $\mathrm{K}_{-}$ of such generalized
vector-functions $\mathrm{k}:X\ni x\mapsto \mathrm{k}(x)\in \mathrm{K}_{x}$
can be considered as the union $\bigcup_{p\in \mathcal{P}_{-}}\mathrm{K}(p)$
of the inductive family of Hilbert spaces $\mathrm{K}(p),p\in \mathcal{P}%
_{-} $, with the norms $\left\Vert \mathrm{k}\right\Vert (p)$, containing as
the minimal the space $\mathrm{K}=\mathrm{K}(1)$. In the extended scale $\{%
\mathrm{K}(p):p\in \mathcal{P}\}$, where $\mathcal{P}=\mathcal{P}_{-}\cup 
\mathcal{P}_{1}$, we obtain the Gel'fand chain $\mathrm{K}^{+}\subseteq 
\mathrm{K}(p^{+})\subseteq \mathrm{K}\subseteq \mathrm{K}(p_{-})\subseteq 
\mathrm{K}_{-}$, where $p^{+}\in \mathcal{P}_{1},\,p_{-}\in \mathcal{P}_{-}$%
, and $\mathrm{K}^{+}=\mathrm{K}_{-}^{\ast }$ coincides with the space of
functionals on $\mathrm{K}_{-}$ continuous with respect to the inductive
convergence. We can similarly define a Gel'fand triple $(\mathrm{F}^{+},%
\mathrm{F},\mathrm{F}_{-})$ for the Hilbert scale $\{\mathrm{F}(p):p\in 
\mathcal{P}\}$ of Fock spaces $\mathrm{F}(p)$ over $\mathrm{K}(p)$, that is,
the spaces of functions $\mathrm{f}:\varkappa \mapsto \mathrm{f}(\varkappa
)\in \mathrm{K}^{\otimes }(\varkappa )$ square integrable with weight $%
p(\varkappa )=\prod_{x\in \varkappa }p(x)$, with values in Hilbert products $%
\mathrm{K}^{\otimes }(\varkappa )=\bigotimes_{x\in \varkappa }\mathrm{K}_{x}$%
: 
\begin{equation*}
\Vert \mathrm{f}\Vert (p)=\left( \int \Vert \mathrm{f}(\varkappa )\Vert
^{2}p(\varkappa )\,\mathrm{d}\varkappa \right) ^{1/2}<\infty .
\end{equation*}%
The integral here is over all chains $\varkappa \in \mathcal{X}$ that define
the pairing on $\mathrm{F}_{-}$ by 
\begin{equation*}
\left\langle \mathrm{f}\mid \mathrm{h}\right\rangle =\int \left\langle 
\mathrm{f}(\varkappa )\mid \mathrm{h}(\varkappa )\right\rangle \,\mathrm{d}%
\varkappa ,\quad \mathrm{h}\in \mathrm{F}^{+},
\end{equation*}%
and in more detail we can write this in the form 
\begin{equation*}
\int \Vert \mathrm{f}(\varkappa )\Vert ^{2}p(\varkappa )\mathrm{d}\varkappa
=\sum_{n=0}^{\infty }\int\limits_{0\leq t_{1}<}\cdots
\int\limits_{<t_{n}<\infty }\Vert \mathrm{f}(x_{1},\ldots ,x_{n})\Vert
^{2}\prod_{i=1}^{n}p(x_{1})\mathrm{d}x_{i},
\end{equation*}%
where the $n$-fold integrals are taken over simplex domains $\Gamma _{n}=\{%
\boldsymbol{x}\in X^{n}:t(x_{1})<\cdots <t(x_{n})\}$. In a similar way as is
done in the case $X=\mathbb{R}_{+},\,t(x)=x$, one can easily establish an
isomorphism between the space $\mathrm{F}(p)$ and the symmetric or
antisymmetric Fock space over $\mathrm{K}(p)$, the isomorphism defined by
the isometry 
\begin{equation*}
\Vert \mathrm{f}\Vert (p)=\left( \sum_{n=0}^{\infty }\frac{1}{n!}\idotsint
\Vert \mathrm{f}(x_{1},\ldots ,x_{n})\Vert ^{2}\prod_{i=1}^{n}p(x_{i})%
\mathrm{d}x_{i}\right) ^{1/2},
\end{equation*}%
where the functions $\mathrm{f}(x_{1},\ldots x_{n})$ are extended to the
whole of $X^{n}$ in a suitable way.

Let $D=(D_{\nu }^{\mu })_{\nu =\circ ,+}^{\mu =-,\circ }$ be a quadruple of
functions $D_{\nu }^{\mu }$ on $X$ with values in continuous operators 
\begin{align}
& D_{+}^{-}(x):\mathrm{F}^{+}\rightarrow \mathrm{F}_{-},\;\;\;\;\;\;\
D_{\circ }^{\circ }(x):\mathrm{F}^{+}\otimes \mathrm{K}_{x}\rightarrow 
\mathrm{F}_{-}\otimes \mathrm{K}_{x},  \notag  \label{2onea} \\
& D_{+}^{\circ }(x):\mathrm{F}^{+}\rightarrow \mathrm{F}_{-}\otimes \mathrm{K%
}_{x},\;\;\;\;\;\;D_{\circ }^{-}(x):\mathrm{F}^{+}\otimes \mathrm{K}%
_{x}\rightarrow \mathrm{F}_{-},
\end{align}%
so that there is a $p\in \mathcal{P}_{1}$ such that these operators are
bounded from $\mathrm{F}(p)\supseteq \mathrm{F}^{+}$ to $\mathrm{F}(p)^{\ast
}\subseteq \mathrm{F}_{-}$, where $\mathrm{F}(p)^{\ast }=\mathrm{F}(1/p)$.
We assume that $D_{+}^{-}(x)$ is locally integrable in the sense that 
\begin{equation*}
\exists \,p\in \mathcal{P}_{1}:\Vert D_{+}^{-}\Vert
_{p,t}^{(1)}=\int_{X^{t}}\Vert D_{+}^{-}(x)\Vert _{p}\mathrm{d}x<\infty
,\quad \forall t<\infty ,
\end{equation*}%
where $X^{t}=\{x\in X:t(x)<t\}$, and $\Vert D\Vert _{p}=\sup \{\Vert D%
\mathrm{h}\Vert (p^{-1})/\Vert \mathrm{h}\Vert (p)\}$ is the norm of the
continuous operator $D:\mathrm{F}(p)\rightarrow \mathrm{F}(p)^{\ast }$ which
defines a bounded Hermitian form $\left\langle \mathrm{f}\mid D\mathrm{h}%
\right\rangle $ on $\mathrm{F}(p)$. We also assume that $D_{\circ }^{\circ
}(x)$ is locally bounded with respect to a strictly positive function $s$
such that $1/s\in \mathcal{P}_{0}$ in the sense that 
\begin{equation*}
\exists \,p\in \mathcal{P}_{1}:\Vert D_{\circ }^{\circ }\Vert
_{p,t}^{(\infty )}(s)=\mathrm{ess}\sup_{x\in X_{t}}\{s(x)\Vert D_{\circ
}^{\circ }(x)\Vert _{p}\}<\infty ,\quad \forall t<\infty ,
\end{equation*}%
where $\Vert D\Vert _{p}$ is the norm of the operator $\mathrm{F}(p)\otimes 
\mathrm{K}_{x}\rightarrow \mathrm{F}(p)^{\ast }\otimes \mathrm{K}_{x}$.
Finally, we assume that $D_{+}^{\circ }(x)$ and $D_{\circ }^{-}(x)$ are
locally square integrable with strictly positive weight $r(x)$ such that $%
1/r\in \mathcal{P}_{0}$, in the sense that 
\begin{equation*}
\exists \,p\in \mathcal{P}_{1}:\Vert D_{+}^{\circ }\Vert
_{p,t}^{(2)}(r)<\infty ,\;\;\;\;\Vert D_{\circ }^{-}\Vert
_{p,t}^{(2)}(r)<\infty ,\quad \forall t<\infty ,
\end{equation*}%
where $\Vert D\Vert _{p,t}^{(2)}(t)=(\int_{X^{t}}\Vert D(x)\Vert _{p}^{2}r(x)%
\mathrm{d}x)^{1/2}$ and $\Vert D\Vert _{p}$ are the norms, respectively, of
the operators 
\begin{equation*}
D_{+}^{\circ }(x):\mathrm{F}(p)\rightarrow \mathrm{F}(p)^{\ast }\otimes 
\mathrm{K}_{x},\;\;\;\;\;\;\,D_{\circ }^{-}(x):\mathrm{F}(p)\otimes \mathrm{K%
}_{x}\rightarrow \mathrm{F}(p)^{\ast }.
\end{equation*}%
Then for any $t\in \mathbb{R}_{+}$ we can define a generalized quantum
stochastic integral \cite{15} 
\begin{equation}
i_{0}^{t}(\mathbf{D})=\int_{X_{t}}\Lambda (\mathbf{D},\mathrm{d}%
x),\;\;\;\;\;\Lambda (\mathbf{D},\Delta )=\sum_{\mu ,\nu }\Lambda _{\nu
}^{\mu }(D_{\nu }^{\mu },\Delta )  \label{2oneb}
\end{equation}%
as the sum of four continuous operators $\Lambda _{\mu }^{\nu }(D_{\nu
}^{\mu }):\mathrm{F}^{+}\rightarrow \mathrm{F}_{-}$ for $\Delta =X^{t}$
which are operator measures on $\mathfrak{F}\ni \Delta $ with values%
\begin{eqnarray}
\lbrack \Lambda _{-}^{+}(D_{+}^{-},\Delta )\mathrm{h}](\varkappa )
&=&\int_{\Delta }[D_{+}^{-}(x)\mathrm{h}](\varkappa )\mathrm{d}x\text{ \ \
~(preservation)},  \notag \\
\lbrack \Lambda _{\circ }^{+}(D_{+}^{\circ },\Delta )\mathrm{h}](\varkappa )
&=&\sum_{x\in \Delta \cap \varkappa }[D_{+}^{\circ }(x)\mathrm{h}](\varkappa
\setminus x)\text{ \ \ (creation)},  \notag \\
\lbrack \Lambda _{-}^{\circ }(D_{\circ }^{-},\Delta )\mathrm{h}](\varkappa )
&=&\int_{\Delta }[D_{\circ }^{-}(x)\mathrm{\dot{h}}(x)](\varkappa )\mathrm{d}%
x\text{ (annihilation)},  \notag \\
\lbrack \Lambda _{\circ }^{\circ }(D_{\circ }^{\circ },\Delta )h])\varkappa
) &=&\sum_{x\in \Delta \cap \varkappa }[D_{\circ }^{\circ }(x)\mathrm{\dot{h}%
}(x)](\varkappa \setminus x)\text{ (exchange)}.  \label{2onec}
\end{eqnarray}%
Here $\mathrm{h}\in \mathrm{F}^{+},\varkappa \setminus x=\{x^{\prime }\in
\varkappa :x^{\prime }\neq x\}$ denotes the chain $\varkappa \in \mathcal{X}$
where the point $x\in \varkappa $ has been annihilated, and $\mathrm{\dot{h}}%
(x)\in \mathrm{K}_{x}\otimes \mathrm{F}^{+}$ is the pointwise derivative
defined for $\mathrm{h}\in \mathrm{F}^{+}$ almost everywhere (for $x\notin
\varkappa \in \mathcal{X}$) on $\mathcal{X}$ as the function $\mathrm{\dot{h}%
}(x,\varkappa )=\mathrm{h}(\varkappa \sqcup x)\equiv \lbrack a(x)\mathrm{h}%
](\varkappa )$, where the operation $\varkappa \sqcup \vartheta $ denotes
the union $\omega =\varkappa \cup \vartheta $ of disjoint chains $\varkappa
\cap \vartheta =\emptyset $ with pairwise comparable elements. The operator $%
a(\vartheta )\mathrm{h}(\varkappa )=\mathrm{\dot{h}}(\vartheta ,\varkappa
\setminus \vartheta )$ which annihilates the points $\vartheta
=\{x_{1},\ldots x_{n}\}\subseteq \varkappa $ in the chain $\varkappa \in 
\mathcal{X}$ defines almost everywhere $(\varkappa \cap \vartheta =\emptyset
)$ on $\mathcal{X}$ the $n$-point derivative $\mathrm{\dot{h}}(\vartheta
,\varkappa )=\mathrm{h}(\varkappa \sqcup \vartheta )$ as the Fock
representation of the $n$-th order Malliavin derivative \cite{39} at these
points. The continuity properties of this operator, which defines an
isometric map $a:\mathrm{F}\left( \frac{1}{r}+p\right) \rightarrow \mathrm{F}%
\left( \frac{1}{r}\right) \otimes \mathrm{F}(p)$, are given by the following
lemma.

\begin{lemma}
The operators $[a(\vartheta )\mathrm{h}](\varkappa )=\mathrm{h}(\varkappa
\sqcup \upsilon ),\,x\in \mathcal{X}$, defines a projection-continuous map $%
a $ of the scale $\mathrm{F}(p),p\in \mathcal{P}$ to $\mathrm{F}\left( \frac{%
1}{r}\right) \otimes \mathrm{F}(p)$, where $r^{-1}\in \mathcal{P}_{0}$, such
that $\Vert a\mathrm{h}\Vert \left( \frac{1}{r},p\right) =\Vert \mathrm{h}%
\Vert \left( \frac{1}{r}+p\right) $, and it is formally adjoint to the
creation operator%
\begin{equation*}
\lbrack a^{\ast }\mathrm{f}](\omega )=\sum_{\vartheta \subseteq \omega }%
\mathrm{f}(\vartheta ,\omega \setminus \vartheta ),\;\,\mathrm{f}\in \mathrm{%
F}(r)\otimes \mathrm{F}\left( \frac{1}{p}\right) ,
\end{equation*}%
which is a contraction in $\mathrm{F}\left( \frac{1}{q}\right) $ for $q\geq 
\frac{1}{r}+p$.
\end{lemma}

\begin{proof}
We first of all establish the principal formula of the multiple integration 
\begin{equation}
\int \sum_{\vartheta \subseteq \omega }f(\vartheta ,\omega \setminus
\vartheta )\mathrm{d}\omega =\iint f(\vartheta ,\upsilon )\mathrm{d}%
\vartheta \mathrm{d}\upsilon ,\quad \forall f\in L^{1}(\mathcal{X}^{2}),
\label{2oned}
\end{equation}%
which allows us to define the adjoint operator $a^{\ast }$. Let $f(\vartheta
,\upsilon )=g(\vartheta )h(\upsilon )$ be the product of integrable complex
functions on $\mathcal{X}$ of the form $g(\vartheta )=\prod_{x\in \vartheta
}g(x)$, $h(\upsilon )=\prod_{x\in \upsilon }h(x)$ for any $\vartheta $, $%
\upsilon \in \mathcal{X}$. Employing the binomial formula 
\begin{equation*}
\sum_{\vartheta \subseteq \omega }g(\vartheta )h(\omega \setminus \vartheta
)=\sum_{\vartheta \sqcup \upsilon =\omega }\prod_{x\in \vartheta
}g(x)\prod_{x\in \upsilon }h(x)=\prod_{x\in \omega }(g(x)+h(x)),
\end{equation*}%
and also the equality $\int f(\vartheta )\mathrm{d}\vartheta =\exp \{\int
f(x)\mathrm{d}x\}$ for $f(\vartheta )=\prod_{x\in \vartheta }f(x)$, we
obtain the formula 
\begin{equation*}
\int \sum_{\vartheta \subseteq \omega }g(\vartheta )h(\omega \setminus
\vartheta )\mathrm{d}\omega =\exp \left\{ \int (g(x)+h(x))\mathrm{d}%
x\right\} =\iint g(\vartheta )h(\upsilon )\mathrm{d}\vartheta (\mathrm{d}%
\upsilon ),
\end{equation*}%
which proves~(\ref{2oned}) on a set of product-functions $f$ dense in $L^{1}(%
\mathcal{X}\times \mathcal{X})$.

Applying this formula to the scalar product $\left\langle \mathrm{f}%
(\vartheta ,\upsilon )|\mathrm{h}(\vartheta ,\upsilon )\right\rangle \in
L^{1}(\mathcal{X}\times \mathcal{X})$, we obtain 
\begin{equation*}
\int \sum_{\vartheta \subseteq \omega }\left\langle \mathrm{f}(\vartheta
,\omega \setminus \vartheta )\mid \mathrm{h}(\omega )\right\rangle \mathrm{d}%
\omega =\iint \left\langle \mathrm{f}(\vartheta ,\upsilon )\mid \mathrm{h}%
(\vartheta \sqcup \upsilon )\right\rangle \mathrm{d}\vartheta \mathrm{d}%
\upsilon ,
\end{equation*}%
that is, $\left\langle a^{\ast }\mathrm{f}\mid \mathrm{h}\right\rangle
=\left\langle \mathrm{f}\mid a\mathrm{h}\right\rangle $, where $[a\mathrm{h}%
](\vartheta ,\upsilon )=\mathrm{h}(\upsilon \sqcup \vartheta )\equiv \mathrm{%
\dot{h}}(\vartheta ,\upsilon )$. Choosing arbitrary $\mathrm{f}\in \mathrm{F}%
(r)\otimes \mathrm{F}\left( \frac{1}{p}\right) $, we find that the
annihilation operator $a(\vartheta )\mathrm{h}=[a\mathrm{h}](\vartheta
,\cdot )$ defines an isometry $\mathrm{F}\left( \frac{1}{r}+p\right)
\rightarrow \mathrm{F}\left( \frac{1}{r}\right) \otimes \mathrm{F}(p)$ as
the operator adjoint to $a^{\ast }:\mathrm{F}(r)\otimes \mathrm{F}\left( 
\frac{1}{p}\right) \rightarrow \mathrm{F}\left( \frac{1}{q}\right) ,\,q=%
\frac{1}{r}+p$, with respect to the standard pairing of conjugate spaces $%
\mathrm{F}(p)$ and $\mathrm{F}\left( \frac{1}{p}\right) $: 
\begin{eqnarray*}
&&\iint \Vert \mathrm{\dot{h}}(\vartheta ,\upsilon \Vert
^{2}r^{-1}(\vartheta )p(\upsilon )\mathrm{d}\vartheta \mathrm{d}\upsilon \\
&=&\int \sum_{\vartheta \subseteq \omega }\Vert \mathrm{h}(\omega )\Vert
^{2}r^{-1}(\vartheta )p(\omega \setminus \vartheta )\mathrm{d}\omega =\int
\Vert \mathrm{h}(\omega )\Vert ^{2}\sum_{\vartheta \sqcup \upsilon =\omega
}r^{-1}(\vartheta )p(\upsilon )\mathrm{d}\omega \\
&=&\int \Vert \mathrm{h}(\omega )\Vert ^{2}(r^{-1}+p)(\omega )\mathrm{d}%
\omega .
\end{eqnarray*}%
Hence it follows that $a$ is projective continuous from $\mathrm{F}^{+}$ to $%
\mathrm{F}_{0}^{+}\times \mathrm{F}^{+}$, where $\mathrm{F}%
_{0}^{+}=\bigcap_{p\in \mathcal{P}_{0}}\mathrm{F}(p)$, and, in particular,
so is the one-point derivative $\mathrm{\dot{f}}(x,\upsilon )=\mathrm{f}%
(x\sqcup \upsilon )$ from $\mathrm{F}^{+}$ to $\mathrm{K}^{+}\times \mathrm{F%
}^{+}$ as a contracting map $\mathrm{F}\left( \frac{1}{r}+p\right)
\rightarrow \mathrm{F}\left( \frac{1}{p}\right) \otimes \mathrm{F}(p)$ for
all $r^{-1}\in \mathcal{P}_{0},p\in \mathcal{P}$. The lemma is proved.
\end{proof}

We are now ready to prove the inductive continuity of the integral (\ref%
{2oneb}) with respect to $\mathbf{D}=(D_{\nu }^{\mu }]$ using the inequality 
\begin{equation*}
\Vert (i_{0}^{t}(\mathbf{D})\mathrm{h})\Vert \left( \frac{1}{q}\right) \leq
\Vert D\Vert _{p,t}^{s}(r)\Vert \mathrm{h}\Vert (q),\quad \forall q\leq
r^{-1}+p+s^{-1},
\end{equation*}%
where $\Vert D\Vert _{p,t}^{s}(r)=\Vert D_{+}^{-}\Vert _{p,t}^{(1)}+\Vert
D_{+}^{\circ }\Vert _{p,t}^{(2)}(r)+\Vert D_{\circ }^{-}\Vert
_{p,t}^{(2)}(r)+\Vert D_{\circ }^{\circ }\Vert _{p,t}^{(\infty )}(s)$. We
will establish this as the single-integrant case for the multiple
generalized integral \cite{17} 
\begin{equation}
\lbrack \iota _{0}^{t}(B)\mathrm{h}](\varkappa )\Vert =\sum_{\varkappa
_{\circ }^{\circ }\sqcup \varkappa _{+}^{\circ }\subseteq \varkappa
^{t}}\int_{\mathcal{X}^{t}}\int_{\mathcal{X}^{t}}[B(\mathbf{\vartheta })%
\mathrm{\dot{h}}(\vartheta _{\circ }^{-}\sqcup \vartheta _{\circ }^{\circ
})](\varkappa _{-}^{\circ })\mathrm{d}\vartheta _{+}^{-}\mathrm{d}\vartheta
_{\circ }^{-}.  \label{2onee}
\end{equation}%
where $\varkappa ^{t}=\varkappa \cap X^{t},\mathcal{X}^{t}=\{\varkappa \in 
\mathcal{X}:\varkappa \subset X^{t}\}$ and the sum is taken over all
decompositions $\varkappa =\varkappa _{-}^{\circ }\sqcup \vartheta _{\circ
}^{\circ }\sqcup \vartheta _{+}^{\circ }$ such that $\vartheta _{\circ
}^{\circ }\in \mathcal{X}^{t}$ and $\vartheta _{+}^{\circ }\in \mathcal{X}%
^{t}$. The multi-integrant $B\left( \mathbf{\vartheta }\right) $ in general
is an operator-function of the quadruple $\mathbf{\vartheta }=(\vartheta
_{\nu }^{\mu })_{\nu =\circ ,+}^{\mu =-,\circ }$ of chains $\vartheta _{\nu
}^{\mu }\in \mathcal{X}$, defined almost everywhere by its values as the
operators 
\begin{equation*}
B%
\begin{pmatrix}
\vartheta _{+}^{-},\vartheta _{\circ }^{-} \\ 
\vartheta _{+}^{\circ },\vartheta _{\circ }^{\circ }%
\end{pmatrix}%
:\mathrm{F}^{+}\otimes \mathrm{K}^{\otimes }(\vartheta _{\circ }^{-})\otimes 
\mathrm{K}^{\otimes }(\vartheta _{\circ }^{\circ })\rightarrow \mathrm{F}%
_{-}\otimes \mathrm{K}^{\otimes }(\vartheta _{\circ }^{\circ })\otimes 
\mathrm{K}^{\otimes }(\vartheta _{+}^{\circ }).
\end{equation*}%
We will assume that these operators are bounded from $\mathrm{F}(p)$ to $%
\mathrm{F}\left( \frac{1}{p}\right) $ for some $p\in \mathcal{P}_{1}$ and
that there exist strictly positive functions $r>0,r^{-1}\in \mathcal{P}_{0}$%
, and $s>0,s^{-1}\in \mathcal{P}_{0}$, such that 
\begin{equation}
\Vert B\Vert _{p,t}^{s}(r)=\int_{\mathcal{X}^{t}}\Vert B_{+}^{-}(\vartheta
)\Vert _{p,t}^{s}\mathrm{d}\vartheta <\infty ,\quad \forall t<\infty ,
\label{2onef}
\end{equation}%
where 
\begin{equation*}
\Vert B_{+}^{-}(\vartheta _{+}^{-})\Vert _{p,t}^{s}(r)=\left( \int_{\mathcal{%
X}^{t}}\int_{\mathcal{X}^{t}}\mathrm{ess}\sup_{\vartheta _{\circ }^{\circ
}\in \mathcal{X}^{t}}(s(\vartheta _{\circ }^{\circ })\Vert B(\mathbf{%
\vartheta })\Vert _{p})^{2}r(\vartheta _{+}^{\circ }\sqcup \vartheta _{\circ
}^{-})\mathrm{d}\vartheta _{+}^{\circ }\mathrm{d}\vartheta _{\circ
}^{-})^{1/2}\right) ,
\end{equation*}%
and $s(\vartheta )=\prod_{x\in \vartheta }s(x)$, $r(\vartheta )=\prod_{x\in
\vartheta }r(x)$. We mention that the single integral~(\ref{2oneb})
corresponds to the case 
\begin{equation*}
B(\mathbf{x}_{\nu }^{\mu })=D_{\nu }^{\mu }(x),\quad B(\mathbf{\vartheta }%
)=0,\quad \forall \mathbf{\vartheta }:\sum_{\mu ,\nu }|\vartheta _{\nu
}^{\mu }|\neq 1,
\end{equation*}%
where $\mathbf{x}_{\nu }^{\mu }$ denotes one of the atomic tables 
\begin{equation}
\mathbf{x}_{+}^{-}=%
\begin{pmatrix}
x, & \emptyset \\ 
\emptyset & \emptyset%
\end{pmatrix}%
,\mathbf{x}_{+}^{\circ }=%
\begin{pmatrix}
\emptyset , & \emptyset \\ 
x & \emptyset%
\end{pmatrix}%
,\mathbf{x}_{\circ }^{-}=%
\begin{pmatrix}
\emptyset , & x \\ 
\emptyset & \emptyset%
\end{pmatrix}%
,\mathbf{x}_{\circ }^{\circ }=%
\begin{pmatrix}
\emptyset , & \emptyset \\ 
\emptyset & x%
\end{pmatrix}%
,  \label{2oneg}
\end{equation}%
determined by $x\in X$. It follows from the next theorem that the function $%
B(\mathbf{\vartheta })$ in~(\ref{2onee}) can be defined up to equivalence,
whose kernel is $B\approx 0\Leftrightarrow \Vert B\Vert _{p,t}^{s}(r)=0$ for
all $t\in \mathbb{R}_{+}$ and for some $p,r,s$. In particular, $B$ can be
defined almost everywhere only for the tables \textbf{$\vartheta $}$%
=(\vartheta _{\nu }^{\mu })$ that give disjoint decompositions $\varkappa
=\cup _{\mu ,\nu }\vartheta _{\nu }^{\mu }$ of the chains $\varkappa \in 
\mathcal{X}$, that is, are representable in the form $\mathbf{\vartheta }%
=\bigsqcup_{x\in \varkappa }\mathbf{x}$, where $\mathbf{x}$ is one of the
atomic tables~(\ref{2oneg}) with indices $\mu ,\nu $ for $x\in \vartheta
_{\nu }^{\mu }$.

\begin{theorem}
Suppose that $B(\mathbf{\vartheta })$ is a function locally integrable in
the sense of \textup{(\ref{2onef})} for some $p,r,s>0$. Then its integral 
\textup{(\ref{2onee})} is a continuous operator $T_{t}=\iota _{0}^{t}(B)$
from $\mathrm{F}^{+}$ to $\mathrm{F}_{-}$ satisfying the estimate 
\begin{equation}
\Vert T_{t}\Vert _{q}=\sup_{\mathrm{h}\in \mathrm{F}(q)}\left\{ \Vert T_{t}%
\mathrm{h}\Vert \left( \frac{1}{q}\right) /\Vert \mathrm{h}\Vert (q)\right\}
\leq \Vert B\Vert _{p,t}^{s}(r)  \label{2oneh}
\end{equation}%
for any $q\geq r^{-1}+p+s^{-1}$. The operator $T_{t}^{\ast }$, formally
adjoint to $T_{t}$ in $\mathrm{F}$, is the integral 
\begin{equation}
\iota _{0}^{t}(B)^{\ast }=\iota _{0}^{t}(B^{\star }),\quad B^{\star }%
\begin{pmatrix}
\vartheta _{+}^{-}, & \vartheta _{\circ }^{-} \\ 
\vartheta _{+}^{\circ }, & \vartheta _{\circ }^{\circ }%
\end{pmatrix}%
=B%
\begin{pmatrix}
\vartheta _{+}^{-}, & \vartheta _{+}^{\circ } \\ 
\vartheta _{\circ }^{-}, & \vartheta _{\circ }^{\circ }%
\end{pmatrix}%
^{\ast },  \label{2onei}
\end{equation}%
which is continuous from $\mathrm{F}^{+}$ to $\mathrm{F}_{-}$ with $\Vert
B^{\star }\Vert _{p}^{s,t}(r)=\Vert B\Vert _{p}^{s,t}(r)$. Moreover, the
operator-valued function $t\mapsto T_{t}$ has the quantum-stochastic
differential $\mathrm{d}T_{t}=\mathrm{d}\iota _{0}^{t}(\mathbf{D})$ in the
sense that 
\begin{equation}
i_{0}^{t}(B)=B(\emptyset )+i_{0}^{t}(\mathbf{D}),\quad D_{\nu }^{\mu
}(x)=i_{0}^{t(x)}(\dot{B}(\mathbf{x}_{\nu }^{\mu })),  \label{2onej}
\end{equation}%
defined by the quantum-stochastic derivatives $\mathbf{D}=(D_{\nu }^{\mu })$
with values \textup{(\ref{2onea})} acting from $\mathrm{F}(q)$ to $\mathrm{F}%
\left( \frac{1}{q}\right) $ and bounded almost everywhere: 
\begin{equation*}
\Vert D_{+}^{-}\Vert _{q,t}^{(1)}\leq \Vert B\Vert _{p,t}^{s}(r),\quad \Vert
D\Vert _{p,t}^{(2)}(r)\leq \Vert B\Vert _{p,t}^{s}(r),\quad \Vert D_{\circ
}^{\circ }\Vert _{q,t}^{(\infty )}(s)\leq \Vert B\Vert _{p,t}^{s}(r)
\end{equation*}%
for $D=D_{\circ }^{-}$ and $D=D_{+}^{\circ }$, $q\geq r^{-r}+p+s^{-1}$. This
differential is defined in the form of the multiple integrals \textup{(\ref%
{2onee})} with respect to $\mathbf{\vartheta }$ of pointwise derivatives $%
\dot{B}(\mathbf{x},\mathbf{\vartheta })=B(\mathbf{\vartheta }\sqcup \mathbf{x%
})$, where $\mathbf{x}$ is one of four atomic tables \textup{(\ref{2oneg})}
at a fixed point $x\in X$.
\end{theorem}

\begin{proof}
Using property (\ref{2oned}) in the form 
\begin{equation*}
\int \sum_{\sqcup \vartheta _{\nu }^{\circ }=\vartheta }f(\vartheta
_{-}^{\circ },\vartheta _{\circ }^{\circ },\vartheta _{+}^{\circ })\mathrm{d}%
\vartheta =\iiint f(\vartheta _{-}^{\circ },\vartheta _{\circ }^{\circ
},\vartheta _{+}^{\circ })\prod_{\nu }\mathrm{d}\vartheta _{\nu }^{\circ },
\end{equation*}%
it is easy to find that from the definition (\ref{2onee}) for $\mathrm{f},%
\mathrm{h}\in \mathrm{F}^{+}$ we have $\int \left\langle \mathrm{f}%
(\varkappa )\mid \lbrack T_{t}\mathrm{h}](\varkappa )\right\rangle \mathrm{d}%
\varkappa =$ 
\begin{align*}
\lefteqn{=\int_{\mathcal{X}^{t}}\mathrm{d}\vartheta _{+}^{-}\int_{\mathcal{X}%
^{t}}\mathrm{d}\vartheta _{+}^{\circ }\int_{\mathcal{X}^{t}}\mathrm{d}%
\vartheta _{\circ }^{-}\int_{\mathcal{X}^{t}}\mathrm{d}\vartheta _{\circ
}^{\circ }\left\langle \mathrm{\dot{f}}(\vartheta _{\circ }^{\circ }\sqcup
\vartheta _{+}^{\circ })\mid B(\mathbf{\vartheta })\mathrm{\dot{h}}%
(\vartheta _{\circ }^{-}\sqcup \vartheta _{\circ }^{\circ })\right\rangle }
\\
& =\int_{\mathcal{X}^{t}}\mathrm{d}\vartheta _{+}^{-}\int_{\mathcal{X}^{t}}%
\mathrm{d}\vartheta _{+}^{\circ }\int_{\mathcal{X}^{t}}\mathrm{d}\vartheta
_{\circ }^{-}\int_{\mathcal{X}^{t}}\mathrm{d}\vartheta _{\circ }^{\circ
}\left\langle B(\mathbf{\vartheta })^{\ast }\mathrm{\dot{f}}(\vartheta
_{\circ }^{\circ }\sqcup \vartheta _{+}^{\circ })\mid \mathrm{\dot{h}}%
(\vartheta _{\circ }^{-}\sqcup \vartheta _{\circ }^{\circ })\right\rangle \\
& =\int \left\langle [T_{t}^{\ast }\mathrm{f}](\varkappa )\mid \mathrm{h}%
(\varkappa )\right\rangle \mathrm{d}\varkappa ,
\end{align*}%
that is, $T_{t}^{\ast }$ acts as $\iota _{0}^{t}(B^{\star })$ in (\ref{2onee}%
) with $B^{\star }(\mathbf{\vartheta })=B(\mathbf{\vartheta }^{\prime
})^{\ast }$, where $(\vartheta _{\nu }^{\mu })^{\prime }=(\vartheta _{-\mu
}^{-\nu })$ with respect to the inverstion $-:(-,\circ ,+)\mapsto (+,\circ
,-)$. More precisely, this yields $\Vert \iota _{0}^{t}(B)\Vert _{q}=\Vert
\iota _{0}^{t}(B^{\star })\Vert _{q}$, since $\Vert T\Vert _{q}=\Vert
T^{\ast }\Vert _{q}$ by the definition (\ref{2oneh}) of $q$-norm and by 
\begin{equation*}
\sup \left\{ |\left\langle (\mathrm{f}\mid T\mathrm{h}\right\rangle |/\Vert 
\mathrm{f}\Vert (q)\Vert \mathrm{h}\Vert (q)\right\} =\sup \{|\left\langle
T^{\ast }\mathrm{f}\mid \mathrm{h}\right\rangle |/\Vert \mathrm{f}\Vert
(q)\Vert \mathrm{h}\Vert (q)\}.
\end{equation*}%
We estimate the integral $\left\langle \mathrm{f}\mid T_{t}\mathrm{h}%
\right\rangle $ using the Schwartz inequality 
\begin{equation*}
\int \Vert \mathrm{\dot{f}}(\vartheta )\Vert (p)\Vert \mathrm{\dot{h}}%
(\vartheta )\Vert (p)s^{-1}(\vartheta )\mathrm{d}\vartheta \leq \Vert 
\mathrm{\dot{f}}\Vert (s^{-1},p)\Vert \mathrm{\dot{h}}\Vert (s^{-1},p)
\end{equation*}%
and the property (\ref{2oned}) of the multiple integral according to which $%
\Vert \mathrm{\dot{f}}\Vert (s^{-1},p)=\Vert \mathrm{f}\Vert (p+s^{-1})$, $%
\Vert \mathrm{\dot{h}}\Vert (s^{-1},p)=\Vert \mathrm{h}\Vert (s^{-1}+p)$,$%
\;\left\vert (\mathrm{f}\mid T_{t}\mathrm{h})\right\vert \leq $ 
\begin{align*}
\lefteqn{\;\;\;\;\;\;\quad \leq \int_{\mathcal{X}^{t}}\mathrm{d}\vartheta
_{\circ }^{\circ }\int_{\mathcal{X}^{t}}\int_{\mathcal{X}^{t}}\Vert \mathrm{%
\dot{f}}(\vartheta _{\circ }^{\circ }\sqcup \vartheta _{+}^{\circ })\Vert
(p)\left( \int_{\mathcal{X}^{t}}\Vert B(\mathbf{\vartheta })\Vert _{p}%
\mathrm{d}\vartheta _{+}^{-}\right) \Vert \mathrm{\dot{h}}(\vartheta _{\circ
}^{-}\sqcup \vartheta _{\circ }^{\circ })\Vert (p)\mathrm{d}\vartheta
_{\circ }^{-}\mathrm{d}\vartheta _{+}^{\circ }} \\
& \quad \leq \int_{\mathcal{X}^{t}}\mathrm{d}\vartheta \left( \int_{\mathcal{%
X}^{t}}\Vert \mathrm{\dot{f}}(\vartheta \sqcup \vartheta _{+}^{\circ })\Vert
^{2}(p)\frac{\mathrm{d}\vartheta _{+}^{\circ }}{r(\vartheta _{+}^{\circ })}%
\int_{\mathcal{X}^{t}}\Vert \mathrm{\dot{h}}(\vartheta \sqcup \vartheta
_{\circ }^{-})\Vert ^{2}(p)\frac{\mathrm{d}\vartheta _{\circ }^{-}}{%
r(\vartheta _{\circ }^{-})}\right) ^{\frac{1}{2}}\Vert B_{\circ }^{\circ
}(\vartheta )\Vert _{p,t}(r) \\
& \quad \leq \int_{\mathcal{X}^{t}}\mathrm{d}\vartheta \Vert \mathrm{\dot{f}}%
(\vartheta )\Vert (r^{-1}+p)\Vert B_{\circ }^{\circ }(\vartheta )\Vert
_{p,t}(r)\Vert \mathrm{\dot{h}}(\vartheta )\Vert (r^{-1}+p) \\
& \quad \leq \underset{\vartheta \in \mathcal{X}^{t}}{\mathrm{ess}\sup }%
\{s(\vartheta )\Vert B_{\circ }^{\circ }(\vartheta )\Vert _{p,t}(r)\}\Vert 
\mathrm{f}\Vert (r^{-1}+p+s^{-1})\Vert \mathrm{h}\Vert (r^{-1}+p+s^{-1}),
\end{align*}%
where $\Vert B_{\circ }^{\circ }(\vartheta _{\circ }^{\circ })\Vert
_{p,t}(r)=(\int_{\mathcal{X}^{t}}\int_{\mathcal{X}^{t}}(\int_{\mathcal{X}%
^{t}}\Vert B(\mathbf{\vartheta })\Vert _{p}\mathrm{d}\vartheta
_{+}^{-})^{2}r(\vartheta _{\circ }^{-}\sqcup \vartheta _{+}^{\circ })\mathrm{%
d}\vartheta _{\circ }^{-}\mathrm{d}\vartheta _{+}^{\circ })^{1/2}$. Thus $%
\Vert T_{t}\Vert _{q}\leq \Vert B\Vert _{p,t}(r,s)$, where $q\geq
r^{-1}+p+s^{-1}$ and 
\begin{equation*}
\Vert B\Vert _{p,t}(r,s):=\underset{\vartheta \in \mathcal{X}^{t}}{\mathrm{%
ess}\sup }\{s(\vartheta )\Vert B_{\circ }^{\circ }(\vartheta )\Vert
_{p,t}(r)\}\leq \Vert B\Vert _{p,t}^{s}(r).
\end{equation*}%
Using the definition (\ref{2onee}) and the property 
\begin{equation*}
\int_{\mathcal{X}^{t}}\mathrm{f}(\varkappa )\mathrm{d}\varkappa =\mathrm{f}%
(\emptyset )+\int_{X^{t}}\mathrm{d}x\int_{\mathcal{X}^{t(x)}}\mathrm{\dot{f}}%
(x,\varkappa )\mathrm{d}\varkappa ,
\end{equation*}%
where $\mathrm{\dot{f}}(x,\varkappa )=\mathrm{f}(\varkappa \sqcup x)$, it is
easy to see that $[(T_{t}-T_{0})\mathrm{h}](\varkappa )=[(\iota
_{0}^{t}(B)-B(\emptyset ))\mathrm{h}](\varkappa )=$ 
\begin{align*}
& =\int_{X^{t}}\mathrm{d}x\sum_{\vartheta _{\circ }^{\circ }\sqcup \vartheta
_{+}^{\circ }\subseteq \varkappa }^{t(\vartheta _{\nu }^{\circ
})<t(x)}\{\int_{\mathcal{X}^{t(x)}}\mathrm{d}\vartheta _{+}^{-}[\int_{%
\mathcal{X}^{t(s)}}\mathrm{d}\vartheta _{\circ }^{-}(\dot{B}(\mathbf{x}%
_{+}^{-},\mathbf{\vartheta })\mathrm{\dot{h}}(\vartheta _{\circ }^{-}\sqcup
\vartheta _{\circ }^{\circ }) \\
& \quad +\dot{B}(\mathbf{x}_{\circ }^{-},\mathbf{\vartheta })\mathrm{\dot{h}}%
(x\sqcup \vartheta _{\circ }^{-}\sqcup \vartheta _{\circ }^{\circ
}))]\}(\varkappa \backslash \vartheta _{\circ }^{\circ }\backslash \vartheta
_{+}^{\circ }) \\
& +\sum_{x\in \mathcal{X}^{t}}\int_{\vartheta _{\circ }^{\circ }\sqcup
\vartheta _{+}^{\circ }\subseteq \varkappa }^{t(\vartheta _{\nu }^{\circ
})<t(x)}\{\int_{\mathcal{X}^{t(x)}}\mathrm{d}\vartheta _{+}^{-}[\int_{%
\mathcal{X}^{t(x)}}\mathrm{d}\vartheta _{\circ }^{-}(\dot{B}(\mathbf{x}%
_{+}^{\circ },\mathbf{\vartheta })\mathrm{\dot{h}}(\vartheta _{\circ
}^{-}\sqcup \vartheta _{\circ }^{\circ }) \\
& \quad +\dot{B}(\mathbf{x}_{\circ }^{\circ },\mathbf{\vartheta })\mathrm{%
\dot{h}}(x\sqcup \vartheta _{\circ }^{-}\sqcup \vartheta _{\circ }^{\circ
}))]\}((\varkappa \backslash \vartheta _{\circ }^{\circ }\backslash
\vartheta _{+}^{\circ }) \\
& =\int_{X^{t}}\mathrm{d}x[D_{+}^{-}(x)\mathrm{h}+D_{\circ }^{-}(x)\mathrm{%
\dot{h}}(x)](\varkappa )+\sum_{x\in \mathcal{X}^{t}}[D_{+}^{\circ }(x)%
\mathrm{h}+D_{\circ }^{\circ }(x)\mathrm{\dot{h}}(x)](\varkappa \setminus x).
\end{align*}%
Consequently, $T_{t}-T_{0}=\sum \Lambda _{\mu }^{\nu }(D_{\nu }^{\mu
},X^{t}) $, where $\Lambda _{\nu }^{\mu }(D,\Delta )$ are defined in (\ref%
{2onec}) as operator-valued measures on $X$ of operator-functions 
\begin{align*}
\lbrack D_{+}^{\mu }(x)\mathrm{h}](\varkappa )& =\sum_{\vartheta _{\circ
}^{\circ }\sqcup \vartheta _{+}^{\circ }\subseteq \varkappa }^{t(\vartheta
_{\nu }^{\circ })<t(x)}\int_{\mathcal{X}^{t(x)}}\mathrm{d}\vartheta
_{+}^{-}\int_{\mathcal{X}^{t(x)}}\mathrm{d}\vartheta _{\circ }^{-}[\dot{B}(%
\mathbf{x}_{+}^{\mu },\mathbf{\vartheta })\mathrm{h}(\vartheta _{\circ
}^{-}\sqcup \vartheta _{\circ }^{\circ })](\varkappa _{-}^{\circ }), \\
\lbrack D_{\circ }^{\mu }(x)\mathrm{\dot{h}}](\varkappa )& =\sum_{\vartheta
_{\circ }^{\circ }\sqcup \vartheta _{+}^{\circ }\subseteq \varkappa
}^{t(\vartheta _{\nu }^{\circ })<t(x)}\int_{\mathcal{X}^{t(x)}}\mathrm{d}%
\vartheta _{+}^{-}\int_{\mathcal{X}^{t(x)}}\mathrm{d}\vartheta _{\circ }^{-}[%
\dot{B}(\mathbf{x}_{\circ }^{\mu },\mathbf{\vartheta })\mathrm{\dot{h}}%
(\vartheta _{\circ }^{-}\sqcup \vartheta _{\circ }^{\circ })](\varkappa
_{-}^{\circ }),
\end{align*}%
acting on $\mathrm{h}\in \mathrm{F}^{+}$ and $\mathrm{\dot{h}}\in \mathrm{K}%
_{x}\otimes \mathrm{F}^{+}$, where $\varkappa _{-}^{\circ }=\varkappa \cap 
\overline{\left( \vartheta _{\circ }^{\circ }\sqcup \vartheta _{+}^{\circ
}\right) }=\varkappa \backslash \vartheta _{\circ }^{\circ }\backslash
\vartheta _{+}^{\circ }$. This can be written in terms of (\ref{2onee}) as $%
D_{\nu }^{\mu }(x)=\iota _{0}^{t}(\dot{B}(\mathbf{x}_{\nu }^{\mu }))$.
Because of the inequality $\Vert T_{t}\Vert _{q}\leq \Vert B\Vert
_{p,t}^{s}(r)$ for all $q\geq r^{-1}+p+s^{-1}$ we obtain $\Vert
D_{+}^{-}\Vert _{q,t}^{(1)}\leq \Vert B\Vert _{p,t}^{s}(r)$, since $\Vert
D_{+}^{-}(x)\Vert _{q}\leq \Vert \dot{B}(\mathbf{x}_{+}^{-})\Vert
_{p,t(x)}^{s}(r)$: 
\begin{eqnarray*}
\lefteqn{\int_{X^{t}}\Vert D_{+}^{-}(x)\Vert _{q}\mathrm{d}x\leq
\int_{X^{t}}\Vert \dot{B}(\mathbf{x}_{+}^{-})\Vert _{p,t(x)}^{s}(r)\mathrm{d}%
x} \\
={\int_{X^{t}}}\mathrm{d}x{\int_{\mathcal{X}^{t(x)}}}\Vert B_{+}^{-}(x\sqcup
\vartheta )\Vert _{p,t(x)}^{s}(r)\mathrm{d}\vartheta &=&\int_{\mathcal{X}%
^{t}}\Vert B_{+}^{-}(\vartheta )\Vert _{p,t}^{s}(r)\mathrm{d}\vartheta
-\Vert B_{+}^{-}(\emptyset )\Vert _{p,t}^{s}(r) \\
&=&\Vert B\Vert _{p,t}^{s}(r)-\Vert B_{+}^{-}(\emptyset )\Vert _{p,t}^{s}(r),
\end{eqnarray*}%
where $B_{+}^{-}(\vartheta ,\mathbf{\vartheta })=B(\mathbf{\vartheta \sqcup
\vartheta }_{+}^{-})\delta _{\emptyset }(\vartheta _{+}^{-})$ for $\mathbf{%
\vartheta }_{+}^{-}=%
\begin{pmatrix}
\vartheta , & \emptyset \\ 
\emptyset , & \emptyset%
\end{pmatrix}%
,$ $\mathbf{\vartheta }=%
\begin{pmatrix}
\vartheta _{+}^{-} & \vartheta _{\circ }^{-} \\ 
\vartheta _{+}^{\circ } & \vartheta _{\circ }^{\circ }%
\end{pmatrix}%
$.

In a similar way one can obtain 
\begin{align*}
\Vert D_{+}^{\circ }\Vert _{q,t}^{(2)}(r)& \leq \left( \int_{X^{t}}(\Vert 
\dot{B}(\mathbf{x}_{+}^{\circ })\Vert _{p,t(x)}^{s}(r))^{2}r(x)\mathrm{d}%
x\right) ^{1/2} \\
& \leq \int_{\mathcal{X}^{t}}\mathrm{d}\vartheta ^{-}\left( \int_{\mathcal{X}%
^{t}}(\Vert B_{+}(\vartheta ^{-},\vartheta ^{\circ })\Vert
_{p,t}^{s}(r))^{2}r(\vartheta ^{\circ })\mathrm{d}\vartheta ^{\circ }\right)
^{1/2}\leq \Vert B\Vert _{p,t}^{s}(r),
\end{align*}%
where $B_{+}\left( \vartheta ^{-},\vartheta ^{\circ },\mathbf{\vartheta }%
\right) =B(\mathbf{\vartheta }\sqcup \mathbf{\vartheta }_{+})\delta
_{\emptyset }(\vartheta _{+}^{-}\sqcup \vartheta _{+}^{\circ })$ for $%
\mathbf{\vartheta }_{+}=%
\begin{pmatrix}
\vartheta ^{-} & \emptyset \\ 
\vartheta ^{\circ } & \emptyset%
\end{pmatrix}%
$, 
\begin{align*}
\Vert D_{\circ }^{-}\Vert _{q,t}^{(2)}(r)& \leq \Big(\int_{X^{t}}(\Vert \dot{%
B}(\mathbf{x}_{\circ }^{-})\Vert _{p,t(x)}^{s}(r))^{2}r(x)\mathrm{d}x\Big)%
^{1/2} \\
& \leq \int_{\mathcal{X}^{t}}\mathrm{d}\vartheta _{+}\Big(\int_{\mathcal{X}%
^{t}}(\Vert B^{-}(\vartheta _{+},\vartheta _{\circ })\Vert
_{p,t}^{s}(r))^{2}r(\vartheta _{\circ })\mathrm{d}\vartheta _{\circ }\Big)%
^{1/2}\leq \Vert B\Vert _{p,t}^{s}(r),
\end{align*}%
where $B^{-}(\vartheta _{+},\vartheta _{\circ },\mathbf{\vartheta })=B(%
\mathbf{\vartheta }\sqcup \mathbf{\vartheta }^{-})\delta _{\emptyset
}(\vartheta _{+}^{-}\sqcup \vartheta _{\circ }^{-})$, $\mathbf{\vartheta }%
^{-}=%
\begin{pmatrix}
\vartheta _{+} & \vartheta _{\circ } \\ 
\emptyset & \emptyset%
\end{pmatrix}%
$.

Finally, from $\Vert D_{\circ }^{\circ }(x)\Vert _{q}\leq \Vert \dot{B}(%
\mathbf{x}_{\circ }^{\circ }\Vert _{p,t(x)}^{s}(r)$ we similarly obtain 
\begin{equation*}
\Vert D_{\circ }^{\circ }\Vert _{q,t}^{(\infty )}(s)\leq \mathrm{ess}%
\sup_{x\in X^{t}}\{s(x)\Vert \dot{B}(\mathbf{x}_{\circ }^{\circ })\Vert
_{p,t(x)}^{s}(r)\}\leq \Vert B\Vert _{p,t}^{s}(r)
\end{equation*}%
if $q\geq r^{-1}+p+s^{-s}$, which concludes the proof.
\end{proof}

\begin{remark}
The quantum-stochastic integral \textup{(\ref{2onee})} constructed in 
\textup{\cite{17}}, as well as its single variations \textup{(\ref{2oneb})}
introduced in \textup{\cite{13}}, are defined explicitly and do not require
that the functions $B$ and $\mathbf{D}$ under the integral be adapted. By
virtue of the continuity we have proved above, they can be approximated in
the inductive convergence by the sequence of integral sums $\iota
_{0}^{t}(B_{n})$, $i_{0}^{t}(\mathbf{D}_{n})$ corresponding to step
measurable operator functions $B_{n}$ and $\mathbf{D}_{n}$ if the latter
converge inductively to $B$ and $\mathbf{D}$ in the poly-norm \textup{(\ref%
{2onef})}.
\end{remark}

In fact, if there exist functions $r$, $s$ with $r^{-1}$, $s^{-1}\in 
\mathcal{P}_{0}$ and $p\in \mathcal{P}_{1}$ such that $\Vert B_{n}-B\Vert
_{p,t}^{s}(r)\rightarrow 0$, then there also exists a function $q\in 
\mathcal{P}_{1}$ such that $\Vert \iota _{0}^{t}(B_{n}-B)\Vert
_{q}\rightarrow 0$, and we have $q\geq r^{-1}+p+s^{-1}$ by the inequality 
\textup{(\ref{2oneh})}, which implies the inductive convergence $\iota
_{0}^{t}(B_{n})\rightarrow \iota _{0}^{t}(B)$ as a result of the linearity
of $\iota _{0}^{t}$. Suppose that $\mathbf{D}(x)$ is adapted in the sense
that $D_{\nu }^{\mu }(x)(\mathrm{h}^{t(x)}\otimes \mathrm{h}_{[t(x)})=%
\mathrm{f}^{t(x)}\otimes \mathrm{h}_{[t(x)}$ or 
\begin{equation*}
\lbrack D_{\nu }^{\mu }(x)\mathrm{h}](\varkappa )=[D_{\nu }^{\mu }(x)\mathrm{%
\dot{h}}(\varkappa _{\lbrack t(x)})](\varkappa ^{t(x)}),\quad \forall x\in X,
\end{equation*}%
where $\mathrm{\dot{h}}(\varkappa _{\lbrack t},\varkappa ^{t})=\mathrm{h}%
(\varkappa ^{t}\sqcup \varkappa _{\lbrack t})$ and $\varkappa ^{t}\sqcup
\varkappa _{\lbrack t}$ is the decomposition of the chain $\varkappa \in 
\mathcal{X}$ into $\varkappa ^{t}=\{x\in \varkappa :t(x)<t\}$ and $\varkappa
_{\lbrack t}=\{x\in \varkappa :t(x)\geq t\}$. In this case the above
approximation in the class of adapted step functions leads to the definition
of the quantum-stochastic integral $i_{0}^{t}(\mathbf{D})$ in the It\^{o}
sense given by Hudson and Parthasarathy for the case $X=\mathbb{R}%
_{+},t(x)=x $ as the weak limit of integrals sums 
\begin{equation*}
i_{0}^{t}(\mathbf{D}_{n})=\int_{0}^{t}\Lambda (\mathbf{D}_{n},\mathrm{d}%
x)=\sum_{i=1}^{n}D_{\nu }^{\mu }(x_{i})\Lambda _{\mu }^{\nu }(\Delta _{i}).
\end{equation*}%
Here $\mathbf{D}_{n}(x)=\mathbf{D}(x_{j})$ for $x\in \lbrack x_{j},x_{j+1})$
is an adapted approximation corresponding to the decomposition $\mathbb{R}%
_{+}=\sum_{j=1}^{n}\Delta _{i}$ into the intervals $\Delta
_{j}=[x_{j},x_{j+1})$ given by the chain $x_{0}=0<x_{1}<\cdots
<x_{n-1}<x_{n}=\infty $, and $D_{\nu }^{\mu }(x)\Lambda _{\mu }^{\nu
}(\Delta )$ is the sum of the operators \textup{(4.3)} with functions $%
D_{\nu }^{\mu }(x)$ constant on $\Delta $ which can therefore be pulled out
in front of the integrals $\Lambda _{\mu }^{\nu }$. In particular, for $%
D_{+}^{-}=0=D_{\circ }^{\circ }$ and $D_{\circ }^{-}=\widehat{1}\otimes
g=D_{+}^{\circ }$, where $\widehat{1}$ is the unit operator in $\mathrm{F}$
and $g(x)$ is a scalar locally square integrable function corresponding to
the case $\mathrm{K}_{x}=\mathbb{C}$, we obtain the It\^{o} definition of
the Wiener integral 
\begin{equation*}
\mathfrak{i}_{0}^{t}(g)=\int_{0}^{t}g(x)w(\mathrm{d}x),\quad \int_{0}^{t}g(x)%
\widehat{w}(\mathrm{d}x)=i_{0}^{t}(\mathbf{D})
\end{equation*}%
with respect to the stochastic measure $w(\Delta )$, $\Delta \in \mathfrak{F}
$ on $\mathbb{R}_{+}$, represented in $\mathrm{F}$ by the operators $%
\widehat{w}(\Delta )=\Lambda _{\circ }^{+}(\Delta )+\Lambda _{-}^{\circ
}(\Delta )$. We also note that the multiple integral \textup{(\ref{2onee})}
in the scalar case $B(\mathbf{\vartheta })=\widehat{1}\otimes b(\mathbf{%
\vartheta })$ defines the Fock representation of the generalized
Maassen-Meyer kernels \textup{\cite{21}, \cite{41}} and in the case 
\begin{equation*}
b(\mathbf{\vartheta })=\mathrm{f}({\normalsize \vartheta }_{\circ
}^{-}\sqcup {\normalsize \vartheta }_{+}^{\circ })\delta _{\emptyset }(%
{\normalsize \vartheta }_{+}^{-})\delta _{\emptyset }({\normalsize \vartheta 
}_{\circ }^{\circ }),\quad \delta _{\emptyset }({\normalsize \vartheta })=%
\begin{cases}
1, & {\normalsize \vartheta }=\emptyset , \\ 
0, & {\normalsize \vartheta }\neq \emptyset%
\end{cases}%
\end{equation*}%
it leads to the multiple stochastic integrals $\iota _{0}^{t}(B)=\widehat{I}%
_{0}^{t}(\mathrm{f})$, 
\begin{equation*}
I_{0}^{t}(\mathrm{f})=\sum_{n=0}^{\infty }\quad \idotsint\limits_{0\leq
t_{1}<\cdots <t_{n}<t}\mathrm{f}(x_{1},\ldots ,x_{n})w(\mathrm{d}%
x_{1})\ldots w(\mathrm{d}x_{n})
\end{equation*}%
of the generalized functions $\mathrm{f}\in \bigcup_{r^{-1}\in \mathcal{P}%
_{0}}\mathrm{F}(r)$, that is, to the Hida distributions \textup{\cite{23}, 
\cite{45}} of the Wiener measure $w(\Delta )$ represented as $\widehat{w}%
(\Delta )$. Thus we have constructed a general non-commutative analogue of
Hida distributions whose properties are described in the following corollary.

\begin{corollary}
\label{C 1} Suppose that the operator-function $B(\mathbf{\vartheta })=%
\widehat{1}\otimes M(\mathbf{\vartheta })$ is defined by the kernel $M$ such
that $\Vert M\Vert _{t}^{s}(r)<\infty $, 
\begin{equation*}
M%
\begin{pmatrix}
\vartheta _{+}^{-}, & \vartheta _{\circ }^{-} \\ 
\vartheta _{+}^{\circ }, & \vartheta _{\circ }^{\circ }%
\end{pmatrix}%
:\mathrm{K}^{\otimes }(\vartheta _{\circ }^{-}\sqcup \vartheta _{\circ
}^{\circ })\rightarrow \mathrm{K}^{\otimes }(\vartheta _{\circ }^{\circ
}\sqcup \vartheta _{+}^{\circ }),
\end{equation*}%
where 
\begin{equation*}
\Vert M\Vert _{t}^{s}(r)=\int_{\mathcal{X}^{t}}\mathrm{d}\vartheta _{+}^{-}%
\Big(\int_{\mathcal{X}^{t}}\mathrm{d}\vartheta _{+}^{\circ }\int_{\mathcal{X}%
^{t}}\mathrm{d}\vartheta _{\circ }^{-}\mathrm{ess}\sup_{\vartheta _{\circ
}^{\circ }\in \mathcal{X}^{t}}\{s(\vartheta _{\circ }^{\circ })\Vert
M(\vartheta )\Vert \}^{2}r(\vartheta _{+}^{\circ }\sqcup \vartheta _{\circ
}^{-})\Big)^{1/2}
\end{equation*}%
for all $t\in \mathbb{R}_{+}$ and for some $r(\vartheta )=\prod_{x\in
\vartheta }r(x)$, $s(\vartheta )=\prod_{x\in \vartheta }s(x)$; $%
r^{-1},s^{-1}\in \mathcal{P}_{0}$. Then the integral \textup{(\ref{2onee})}
defines an adapted family $T_{t}$, $t\in \mathbb{R}_{+}$, of $q$-bounded
operators $T_{t}=\iota _{0}^{t}(\widehat{1}\otimes M)$, $\Vert T_{t}\Vert
_{q}\leq \Vert M\Vert _{t}^{s}(r)$ for $q\geq r^{-1}+1+s^{-1}$, with adapted 
$p$-bounded quantum-stochastic derivatives $D_{\nu }^{\mu }(x)=\iota
_{0}^{t(x)}(\widehat{1}\otimes \dot{M}(\mathbf{x}_{\nu }^{\mu }))$.
\end{corollary}

\section{Generalized It\^{o} formula of unified quantum stochastic calculus}

Let $\mathrm{H}$ be a Hilbert space. We consider a Hilbert scale $\mathrm{G}%
(p)=\mathrm{H}\otimes \mathrm{F}(p)$, $p\in \mathcal{P}$, of complete tensor
products of $\mathrm{H}$ and the Fock spaces over $\mathrm{K}(p)$, and we
put $\mathrm{G}^{+}=\cap \mathrm{G}(p),\mathrm{G}=\mathrm{G}(1)$, and $%
\mathrm{G}_{-}=\cup \mathrm{G}(p)$ which constitute the corresponding
Gel'fand triple $\mathrm{G}^{+}\subseteq \mathrm{G}\subseteq \mathrm{G}_{-}$%
. We consider operators $T=\epsilon (K)$, not necessarily bounded, in the
Hilbert space $\mathrm{G}=\mathrm{H}\otimes \mathrm{F}$ as the $\ast $%
-representation $\epsilon $ of operator-valued kernels 
\begin{equation}
K%
\begin{pmatrix}
\omega _{+}^{-} & \omega _{\circ }^{-} \\ 
\omega _{+}^{\circ } & \omega _{\circ }^{\circ }%
\end{pmatrix}%
:\mathrm{H}\otimes \mathrm{K}^{\otimes }(\omega _{\circ }^{-}\sqcup \omega
_{\circ }^{\circ })\rightarrow \mathrm{H}\otimes \mathrm{K}^{\otimes
}(\omega _{\circ }^{\circ }\sqcup \omega _{+}^{\circ }),  \label{two1}
\end{equation}%
satisfying the integrability condition $\Vert K\Vert _{p}(r)<\infty $ for
some $r^{-1}\in \mathcal{P}_{0}$ and $p\in \mathcal{P}_{1}$, where 
\begin{equation*}
\Vert K\Vert _{p}(r)=\int \mathrm{d}\omega _{+}^{-}\left( \iint \mathrm{ess}%
\sup_{\omega _{\circ }^{\circ }}\{\frac{\Vert K(\mathbf{\omega })\Vert }{%
p(\omega _{\circ }^{\circ })}\}^{2}r(\omega _{+}^{\circ }\sqcup \omega
_{\circ }^{-})\mathrm{d}\omega _{+}^{\circ }\mathrm{d}\omega _{\circ
}^{-}\right) ^{1/2}.
\end{equation*}%
This representation $\epsilon $ is defined for $\mathrm{h}\ \in \mathrm{H}%
\otimes \mathrm{F}$ by 
\begin{equation}
\lbrack \epsilon (K)\mathrm{h}](\varkappa )=\sum_{\omega _{\circ }^{\circ
}\sqcup \omega _{+}^{\circ }=\varkappa }\iint K%
\begin{pmatrix}
\omega _{+}^{-}, & \omega _{\circ }^{-} \\ 
\omega _{+}^{\circ }, & \omega _{\circ }^{\circ }%
\end{pmatrix}%
\mathrm{h}(\omega _{\circ }^{\circ }\sqcup \omega _{\circ }^{-})\mathrm{d}%
\omega _{\circ }^{-}\mathrm{d}\omega _{+}^{-}  \label{two2}
\end{equation}%
as the vacuum-adapted operator-valued multiple integral (\ref{2onee}) with $%
t=\infty $ of the function $B(\mathbf{\vartheta })=\widehat{\delta }%
_{\emptyset }\otimes K(\mathbf{\vartheta })$, where%
\begin{equation*}
\lbrack \widehat{\delta }_{\emptyset }\mathrm{f}](\varkappa )=\mathrm{f}%
(\emptyset )\delta _{\emptyset }(\varkappa ),\;\;\delta _{\emptyset }\left(
\varkappa \right) =\{%
\begin{array}{cc}
1 & \varkappa =\emptyset \\ 
0 & \varkappa \neq \emptyset%
\end{array}%
\end{equation*}%
is the vacuum projection on $\mathrm{F}$ such that $[B(\mathbf{\vartheta })%
\mathrm{\dot{h}}(\vartheta _{\circ }^{\circ }\sqcup \vartheta _{\circ
}^{-})](\varkappa _{-}^{\circ })=0$ if $\varkappa _{-}^{\circ }=\varkappa
\backslash \vartheta _{\circ }^{\circ }\backslash \vartheta _{+}^{\circ
}\neq \emptyset $. The operator $\epsilon (K)$ can be represented
equivalently as the adapted (i.e. identity-adapted) integral (\ref{2onee})
with $t=\infty $ of a scalar-valued integrant as the function $B(\mathbf{%
\vartheta })=\widehat{1}\otimes M(\mathbf{\vartheta })$, where $\widehat{1}$
is the unit operator on $\mathrm{F}$ and $M(\mathbf{\vartheta })$ is the
generalized Maassen-Meyer kernel-integrant. This follows from 
\begin{equation*}
\lbrack B(\mathbf{\vartheta })\mathrm{\dot{h}}(\vartheta _{\circ }^{-}\sqcup
\vartheta _{\circ }^{\circ })](\varkappa _{-}^{\circ })=M(\mathbf{\vartheta }%
)\mathrm{h}(\vartheta _{\circ }^{-}\sqcup \vartheta _{\circ }^{\circ }\sqcup
\varkappa _{-}^{\circ })
\end{equation*}%
such that $\left[ \iota _{0}^{\infty }(\widehat{1}\otimes M)\mathrm{h}\right]
\left( \varkappa \right) =[\epsilon (K)\mathrm{h}](\varkappa )$ for the
kernel 
\begin{equation*}
K%
\begin{pmatrix}
\omega _{+}^{-}, & \omega _{\circ }^{-} \\ 
\omega _{+}^{\circ }, & \omega _{\circ }^{\circ }%
\end{pmatrix}%
=\sum_{\vartheta \subseteq \omega _{\circ }^{\circ }}M%
\begin{pmatrix}
\omega _{+}^{-}, & \omega _{\circ }^{-} \\ 
\omega _{+}^{\circ }, & \vartheta%
\end{pmatrix}%
\otimes I^{\otimes }(\omega _{\circ }^{\circ }\setminus \vartheta ),
\end{equation*}%
which is connected with $M$ by a one-to-one relation 
\begin{equation*}
M%
\begin{pmatrix}
\vartheta _{+}^{-}, & \vartheta _{\circ }^{-} \\ 
\vartheta _{+}^{\circ }, & \vartheta _{\circ }^{\circ }%
\end{pmatrix}%
=\sum_{\omega \subseteq \vartheta _{\circ }^{\circ }}K%
\begin{pmatrix}
\vartheta _{+}^{-}, & \vartheta _{\circ }^{-} \\ 
\vartheta _{+}^{\circ }, & \omega%
\end{pmatrix}%
\otimes (-I)^{\otimes }(\vartheta _{\circ }^{\circ }\setminus \omega ),
\end{equation*}%
where $I^{\otimes }(\upsilon )=\bigotimes_{x\in \upsilon }I_{x}$ is the unit
operator in $\mathrm{K}^{\otimes }(\upsilon )$.

According to Corollary~\ref{C 1}, $\Vert T\Vert _{q}\leq \Vert M\Vert
_{\infty }^{s}(r)$ for $q\geq r^{-1}+1+s^{-1}$. However, using the
equivalent representation (\ref{two2}) in the form of the non-adapted
integral (\ref{2onee}) of $B(\mathbf{\vartheta })=\widehat{\delta }%
_{\emptyset }\otimes K(\mathbf{\vartheta })$ and taking into account the
fact that $\Vert \widehat{\delta }_{\emptyset }\Vert _{p}=1$ for
sufficiently small $p>0$, we obtain as $p\rightarrow 0$ a more precise
estimate $\Vert T\Vert _{q}\leq \Vert K\Vert _{s^{-1}}(r)$ for $q\geq
r^{-1}+s^{-1}=\lim_{p_{0}\searrow 0}(r^{-1}+p_{0}+s^{-1})$. From this
estimate the previous one follows, since 
\begin{equation*}
\Vert \sum_{\vartheta \subseteq \omega _{\circ }^{\circ }}M(\mathbf{%
\vartheta })\otimes I^{\otimes }(\omega _{\circ }^{\circ }\setminus
\vartheta )\Vert \leq \sum_{\vartheta \subseteq \omega _{\circ }^{\circ
}}\Vert M(\mathbf{\vartheta })\Vert \leq (1+s^{-1})(\omega _{\circ }^{\circ
})\Vert M\Vert _{\infty }^{s},
\end{equation*}%
where $\Vert M\Vert _{\infty }^{s}=\mathrm{ess}\sup_{\vartheta \in \mathcal{X%
}}\{s(\vartheta )\Vert M(\vartheta )\Vert \}$, 
\begin{equation*}
s(\vartheta )=\prod_{x\in \vartheta }s(x),\;\;\;(1+s^{-1})(\omega _{\circ
}^{\circ })=\sum_{\vartheta \subseteq \omega _{\circ }^{\circ
}}s^{-1}(\vartheta )=\prod_{x\in \omega _{\circ }^{\circ }}(1+s^{-1}(x))
\end{equation*}%
and consequently $\Vert K\Vert _{p}(r)\leq \Vert M\Vert _{\infty }^{s}(r)$
for $p\geq 1+1/s$. Hence in particular there follows the existence of the
adjoint operator $T^{\ast }$ bounded in norm $\Vert T^{\ast }\Vert _{q}\leq
\Vert K^{\star }\Vert _{p}(r)=\Vert K\Vert _{p}(r)$ as the representation 
\begin{equation}
\epsilon (K)^{\ast }=\epsilon (K^{\star }),\;\,K^{\star }%
\begin{pmatrix}
\omega _{+}^{-}, & \omega _{\circ }^{-} \\ 
\omega _{+}^{\circ }, & \omega _{\circ }^{\circ }%
\end{pmatrix}%
=K%
\begin{pmatrix}
\omega _{+}^{-}, & \omega _{+}^{\circ } \\ 
\omega _{\circ }^{-}, & \omega _{\circ }^{\circ }%
\end{pmatrix}%
^{\ast }  \label{two3}
\end{equation}%
of the $\star $-adjoint kernel $K^{\star }(\mathbf{\omega })=K(\mathbf{%
\omega }^{\prime })^{\ast }$.

In the next theorem we prove that the $\star $-map $\epsilon :K\mapsto
\epsilon (K)$ is an operator representation of the $\star $-algebra of
kernels $K(\mathbf{\omega })$ satisfying the boundedness condition 
\begin{equation}
\Vert K\Vert _{\mathbf{\alpha }}=\underset{\omega =(\omega _{\nu }^{\mu })}{%
\mathrm{ess}\sup }\{\Vert K(\mathbf{\omega })\Vert /\prod_{\mu \leq \nu
}\alpha _{\nu }^{\mu }(\omega _{\nu }^{\mu })\}<\infty  \label{two4}
\end{equation}%
relative to the product of the quadruple $\mathbf{\alpha }=(\alpha _{\nu
}^{\mu })_{\nu =\circ ,+}^{\mu =-,\circ }$ of positive essentially
measurable product functions $\alpha _{\nu }^{\mu }(\omega )=\prod_{x\in
\omega }\alpha _{\nu }^{\mu }(x)$, $\omega \in \mathcal{X}$. These are
defined by an integrable function $\alpha _{+}^{-}:X\rightarrow \mathbb{R}%
_{+}$, by functions $\alpha _{+}^{\circ },\alpha _{\circ }^{-}:X\rightarrow 
\mathbb{R}_{+}$, square integrable with a certain weight $r>0$, $r^{-1}\in 
\mathcal{P}_{0}$, and by a function $\alpha _{\circ }^{\circ }:X\rightarrow 
\mathbb{R}_{+}$, essentially bounded by unity relative to some $p\in 
\mathcal{P}$: 
\begin{align}
& \Vert \alpha _{+}^{-}\Vert ^{(1)}<\infty ,\,\Vert \alpha _{+}^{\circ
}\Vert ^{(2)}(r)<\infty ,\,\Vert \alpha _{\circ }^{-}\Vert ^{(2)}(r)<\infty
,\,\Vert \alpha _{\circ }^{\circ }\Vert _{p}^{(\infty )}\leq 1,  \notag \\
& \Vert \alpha \Vert ^{(1)}=\int |\alpha (x)|\mathrm{d}x,\,\Vert \alpha
\Vert ^{(2)}(r)=\Big(\int \alpha (x)^{2}r(x)\mathrm{d}x\Big)^{1/2},  \notag
\\
& \hspace{4cm}\Vert \alpha \Vert _{p}^{(\infty )}=\underset{x}{\mathrm{ess}%
\sup }\frac{|\alpha (x)|}{p(x)}.  \label{two5}
\end{align}%
The conditional boundedness (\ref{two4}) ensures the projective boundedness
of $K$ by the inequality $\Vert K\Vert _{p}(r)\leq $ 
\begin{align}
\lefteqn{\leq \int \mathrm{d}\omega _{+}^{-}\Big(\iint \mathrm{ess}\underset{%
\omega _{\circ }^{\circ }}{\sup }\{\Vert K\Vert _{\alpha }\prod \alpha _{\nu
}^{\mu }(\omega _{\nu }^{\mu })/p(\omega _{\circ }^{\circ })\}^{2}r(\omega
_{+}^{\circ }\sqcup \omega _{\circ }^{-})\mathrm{d}\omega _{+}^{\circ }\text{%
\textrm{d}}\omega _{\circ }^{-}\Big)^{1/2},}  \label{two6} \\
& \int \alpha _{+}^{-}(\omega )\mathrm{d}\omega \Big(\int \alpha _{+}^{\circ
}(\omega )^{2}r(\omega )\mathrm{d}\omega \int \alpha _{\circ }^{-}(\omega
)^{2}r(\omega )\mathrm{d}\omega \Big)^{1/2}\mathrm{ess}\sup \frac{\alpha
_{\circ }^{\circ }(\omega )}{p(\omega )}\Vert K\Vert _{\alpha } \\
& \quad \leq \Vert K\Vert _{\alpha }\exp \left\{ \int (\alpha
_{+}^{-}(x)+r(x)(\alpha _{+}^{\circ }(x)^{2}+\alpha _{\circ }^{-}(x)^{2})/2)%
\mathrm{d}x\right\} ,  \notag
\end{align}%
where we have taken account of the fact that $\int \alpha (\omega )\mathrm{d}%
\omega =\exp \int \alpha (x)\mathrm{d}x$ for $\alpha (\omega )=\prod_{x\in
\omega }\alpha (x)$ and%
\begin{equation*}
\mathrm{ess}\sup_{\omega }\{\alpha _{\circ }^{\circ }(\omega )/p(\omega
)\}=\sup_{n}\mathrm{ess}\sup_{x\in X^{n}}\prod_{i=1}^{n}\{\alpha _{\circ
}^{\circ }(x_{i})/p(x_{i})\}=1\text{ if }\alpha _{\circ }^{\circ }\leq p.
\end{equation*}

Before we formulate the theorem we establish the following lemma.

\begin{lemma}
\label{2L 2} Suppose that the multiple quantum-stochastic integral $%
T_{t}=\iota _{0}^{t}(B)$ is defined in \textup{(\ref{2onee})} by a kernel
operator-function $B(\mathbf{\vartheta })=\epsilon (M(\mathbf{\vartheta }))$
with values in the operators of the form \textup{(\ref{two2})}, where 
\begin{equation*}
K%
\begin{pmatrix}
\upsilon _{+}^{-}, & \upsilon _{\circ }^{-} \\ 
\upsilon _{+}^{\circ }, & \upsilon _{\circ }^{\circ }%
\end{pmatrix}%
=M%
\begin{pmatrix}
\vartheta _{+}^{-}, & \vartheta _{\circ }^{-}, & \upsilon _{+}^{-}, & 
\upsilon _{\circ }^{-} \\ 
\vartheta _{+}^{\circ }, & \vartheta _{\circ }^{\circ }, & \upsilon
_{+}^{\circ }, & \upsilon _{\circ }^{\circ }%
\end{pmatrix}%
,\vartheta _{\nu }^{\mu }\in \mathcal{X},
\end{equation*}%
and $M(\mathbf{\vartheta }):\mathbf{\upsilon }\mapsto M\left( \mathbf{%
\vartheta ,\upsilon }\right) $ is a kernel-valued integrant 
\begin{equation*}
M(\mathbf{\vartheta ,\upsilon }):\mathrm{H}\otimes \mathrm{K}^{\otimes
}(\upsilon _{\circ }^{-}\sqcup \vartheta _{\circ }^{-})\otimes \mathrm{K}%
^{\otimes }(\upsilon _{\circ }^{\circ }\sqcup \vartheta _{\circ }^{\circ
})\rightarrow \mathrm{H}\otimes \mathrm{K}^{\otimes }(\upsilon _{\circ
}^{\circ }\sqcup \vartheta _{\circ }^{\circ })\otimes \mathrm{K}^{\otimes
}(\upsilon _{+}^{\circ }\sqcup \vartheta _{+}^{\circ }).
\end{equation*}%
Then $T_{t}=\epsilon (K_{t})$ for the kernel $K_{t}(\mathbf{\omega })=\nu
_{0}^{t}(\mathbf{\omega },M)$ given by the multiple counting integral on the
kernel-integrants $M$, that is, $\iota _{0}^{t}\circ \epsilon =\epsilon
\circ \nu _{0}^{t}$, where 
\begin{equation}
\nu _{0}^{t}(\mathbf{\omega },M)=\sum_{\mathbf{\vartheta }\subseteq \mathbf{%
\omega }^{t}}M(\mathbf{\vartheta },\mathbf{\omega }\setminus \mathbf{%
\vartheta }),\quad \mathbf{\omega }^{t}=(X^{t}\cap \omega _{\nu }^{\mu
})_{\nu =\circ ,+}^{\mu =-,\circ }  \label{two7}
\end{equation}%
\textup{(}the sum is taken over all possible $\vartheta _{\nu }^{\mu
}\subseteq X^{t}\cap \omega _{\nu }^{\mu }$, $\mu =-,\circ ,\nu =\circ ,+$ 
\textup{).} If $M(\mathbf{\vartheta })$ is relatively bounded in $\upsilon
_{\nu }^{\mu }\in \mathcal{X}$ for each $\mathbf{\vartheta }=\left(
\vartheta _{\nu }^{\mu }\right) $ such that 
\begin{equation*}
\Vert M(\mathbf{\vartheta })\Vert _{\mathbf{\gamma }}\leq c\prod_{\mu ,\nu
}\beta _{\nu }^{\mu }(\vartheta _{\nu }^{\mu }),\quad \beta _{\nu }^{\mu
}(\vartheta )=\prod_{x\in \vartheta }\beta _{\nu }^{\mu }(x)
\end{equation*}%
for a pair of quadruples $\mathbf{\beta }=(\beta _{\nu }^{\mu })$ , $\beta
_{\nu }^{\mu }\geq 0$ and $\mathbf{\gamma }=(\gamma _{n}^{\mu })$, $\gamma
_{\nu }^{\mu }\geq 0$ satisfying \textup{(\ref{two5})}, then the kernel $K$
is relatively bounded: $\Vert \nu _{0}^{t}(M)\Vert _{\mathbf{\alpha }}\leq c$
if $\alpha _{\nu }^{\mu }(x)\geq \beta _{\nu }^{\mu }(x)+\gamma _{\nu }^{\mu
}(x)$ for $t(x)<t$ and $\alpha _{\nu }^{\mu }(x)\geq \gamma _{\nu }^{\mu
}(x) $ for all $\mu $, $\nu $ when $t(x)\geq t$. In particular, the
generalized single integral $i_{0}^{t}(\mathbf{D})$ of the triangular
operator-integrant $\mathbf{D}\left( x\right) =\left[ \delta _{+}^{\mu
}\delta _{\nu }^{-}D_{\nu }^{\mu }\left( x\right) \right] $ with $D_{\nu
}^{\mu }(x)=\epsilon (C_{\nu }^{\mu }(x))$ is a representation $%
i_{0}^{t}\circ \epsilon =\epsilon \circ n_{0}^{t}$ of the single counting
integral 
\begin{equation*}
n_{0}^{t}(\mathbf{\omega },\mathbf{C})=\sum_{\mathbf{x}\in \mathbf{\omega }%
^{t}}C(\mathbf{x},\mathbf{\omega }\setminus \mathbf{x}),\quad C(\mathbf{x}%
_{\nu }^{\mu },\mathbf{\upsilon })=C_{\nu }^{\mu }(x,\mathbf{\upsilon }),
\end{equation*}%
of the triangular kernel-integrant $\mathbf{C}\left( x,\mathbf{\upsilon }%
\right) =\left[ \delta _{+}^{\mu }\delta _{\nu }^{-}C_{\nu }^{\mu }\left( x,%
\mathbf{\upsilon }\right) \right] $, where the sum is taken over all
possible $x\in \omega _{\nu }^{\mu }\cap X^{t}$ for $\mu =-,\circ ,\nu
=\circ ,+$, and $\mathbf{x}=\mathbf{x}_{\nu (x)}^{\mu (x)}$ is one of the
atomic matrices \textup{(\ref{2oneg})} with indices $\mu (\mathbf{x})=\mu $, 
$\nu (\mathbf{x})=\nu $, defined almost everywhere by the condition $x\in
\omega _{\nu }^{\mu }$. Moreover, we have the estimate 
\begin{equation*}
\Vert n_{0}^{t}(\mathbf{C})\Vert _{p}(r)\leq c\exp \left\{
\int_{X^{t}}(\gamma _{+}^{-}(x)+{\tfrac{1}{2}}(\gamma _{+}^{\circ
}(x)^{2}+\gamma _{\circ }^{-}(x)^{2})r(x))\mathrm{d}x\right\} .
\end{equation*}%
given the kernel-valued quadruple-integrant $C_{\nu }^{\mu }(x,\mathbf{%
\gamma })$ is relatively bounded for each $x\in X$ in $\,\mathbf{\upsilon }%
=(\upsilon _{\nu }^{\mu })$ in the sense that there exist such $\mathbf{%
\gamma }=\left( \gamma _{\nu }^{\mu }\right) $ that 
\begin{align*}
& c=\Vert C_{+}^{-}\Vert _{\mathbf{\gamma },t}^{(1)}+\Vert C_{+}^{\circ
}\Vert _{\mathbf{\gamma },t}^{(2)}(r)+\Vert C_{\circ }^{-}\Vert _{\mathbf{%
\gamma },t}^{(2)}(r)+\Vert C_{\circ }^{\circ }\Vert _{\mathbf{\gamma }%
,t}^{(\infty )}(1/p)<\infty , \\
& \Vert C_{+}^{-}\Vert _{\mathbf{\gamma },t}^{(1)}=\int_{X^{t}}\Vert
C_{+}^{-}(x)\Vert _{\mathbf{\gamma }}\mathrm{d}x,\,\Vert C\Vert _{\mathbf{%
\gamma },t}^{(2)}=\Big(\int_{X^{t}}\Vert C(x)\Vert _{\mathbf{\gamma }%
}^{2}r(x)\mathrm{d}x\Big)^{1/2}, \\
& \hspace{5cm}\Vert C_{\circ }^{\circ }\Vert _{\mathbf{\gamma },t}^{(\infty
)}\left( \frac{1}{r}\right) =\sup_{x\in X^{t}}\left\{ \frac{\Vert C_{\circ
}^{\circ }(x)\Vert _{\mathbf{\gamma }}}{p(x)}\right\} ,
\end{align*}
\end{lemma}

\begin{proof}
If $M(\mathbf{\vartheta },\mathbf{\upsilon })$ is an operator-valued
integrant-kernel that is bounded, $\Vert M\Vert _{\mathbf{\beta ,\gamma }%
}\leq c$, relative to the pair $(\mathbf{\beta ,\gamma })$, then the
relatively bounded operator $T_{t}=\epsilon (K_{t})$ is well-defined for $%
K_{t}=\nu _{0}^{t}(M)$, since 
\begin{align*}
\Vert K_{t}(\mathbf{\omega })\Vert & \leq c\sum_{\vartheta _{+}^{-}\subseteq
\omega _{+}^{-}}^{t(\vartheta _{+}^{-})<t}\sum_{\vartheta _{+}^{\circ
}\subseteq \omega _{+}^{\circ }}^{t(\vartheta _{+}^{\circ
})<t}\sum_{\vartheta _{\circ }^{-}\subseteq \omega _{\circ
}^{-}}^{t(\vartheta _{\circ }^{-})<t}\sum_{\vartheta _{\circ }^{\circ
}\subseteq \omega _{\circ }^{\circ }}^{t(\vartheta _{\circ }^{\circ
})<t}\Vert M(\mathbf{\vartheta },\mathbf{\omega }\setminus \mathbf{\vartheta 
}\Vert \\
& \leq c\prod_{\nu =\circ ,+}^{\mu =-,\circ }\sum_{\vartheta _{\nu }^{\mu
}\subseteq \omega _{\nu }^{\mu }}^{t(\vartheta _{\nu }^{\mu })<t}\beta _{\nu
}^{\mu }(\vartheta _{\nu }^{\mu })\gamma _{\nu }^{\mu }(\omega _{\nu }^{\mu
}\setminus \vartheta _{\nu }^{\mu })=c\prod_{\nu =\circ ,+}^{\mu =-,\circ
}\alpha _{\nu }^{\mu }(\omega _{\nu }^{\mu }),
\end{align*}%
where $\alpha _{\nu }^{\mu }(\omega )=\prod_{x\in \omega }^{t(x)<t}[\beta
_{\nu }^{\mu }(x)+\gamma _{\nu }^{\mu }(x)]\cdot \prod_{x\in \omega
}^{t(x)\geq t}\gamma _{\nu }^{\mu }(x)$ for $\beta _{\nu }^{\mu }(\vartheta
)=\prod_{x\in \vartheta }\beta _{\nu }^{\mu }(x)$ and $\gamma _{\nu }^{\mu
}(\upsilon )=\prod_{x\in \upsilon }\gamma _{\nu }^{\mu }(x)$. Applying the
representation (\ref{two2}) to $K_{t}(\mathbf{\omega })=\nu _{0}^{t}(\mathbf{%
\omega },M)$ it is easy to obtain the representation of the operator $%
\epsilon (K_{t})$ in the form of the generalized multiple integral (\ref%
{2onee}) of $B(\mathbf{\vartheta })=\epsilon (M(\mathbf{\vartheta }))$.
Indeed, $[T_{t}\mathrm{h}](\varkappa )=$ 
\begin{align*}
& =\sum_{\omega _{\circ }^{\circ }\sqcup \omega _{+}^{\circ }=\varkappa
}\iint \sum_{\mathbf{\vartheta }\subseteq \mathbf{\omega }^{t}}M(\mathbf{%
\vartheta },\mathbf{\omega }\setminus \mathbf{\vartheta })\mathrm{h}(\omega
_{\circ }^{\circ }\sqcup \omega _{\circ }^{-})\mathrm{d}\omega _{\circ }^{-}%
\mathrm{d}\omega _{+}^{-} \\
& =\sum_{\vartheta _{\circ }^{\circ }\sqcup \vartheta _{+}^{\circ }\subseteq
\varkappa ^{t}}\int_{\mathcal{X}^{t}}\mathrm{d}\vartheta _{\circ }^{-}\int_{%
\mathcal{X}^{t}}\mathrm{d}\vartheta _{+}^{-}\sum_{\upsilon _{\circ }^{\circ
}\sqcup \upsilon _{+}^{\circ }=\varkappa _{-}^{\circ }}\iint M(\mathbf{%
\vartheta },\mathbf{\upsilon })\mathrm{\dot{h}}(\vartheta _{\circ }^{\circ
}\sqcup \vartheta _{\circ }^{-},\upsilon _{\circ }^{\circ }\sqcup \upsilon
_{\circ }^{-})\mathrm{d}\upsilon _{\circ }^{-}\mathrm{d}\upsilon _{+}^{-},
\end{align*}%
where $\varkappa _{-}^{\circ }=\varkappa \setminus (\vartheta _{\circ
}^{\circ }\sqcup \vartheta _{+}^{\circ }),\,\mathrm{\dot{h}}(\vartheta
,\upsilon )=\mathrm{h}(\upsilon \sqcup \vartheta )$. Consequently, $%
T_{t}=\iota _{0}^{t}(B)$, where 
\begin{equation*}
\lbrack B(\mathbf{\vartheta })\mathrm{\dot{h}}(\vartheta _{\circ }^{\circ
}\sqcup \vartheta _{\circ }^{-})](\varkappa )=\sum_{\upsilon _{\circ
}^{\circ }\sqcup \upsilon _{+}^{\circ }=\varkappa }\iint M(\mathbf{\vartheta 
},\mathbf{\upsilon })\mathrm{\dot{h}}(\vartheta _{\circ }^{\circ }\sqcup
\vartheta _{\circ }^{-},\,\upsilon _{\circ }^{\circ }\sqcup \upsilon _{\circ
}^{-})\mathrm{d}\upsilon _{\circ }^{-}\mathrm{d}\upsilon _{+}^{-},
\end{equation*}%
that is, we have proved that $\epsilon \circ \nu _{0}^{t}=\iota
_{0}^{t}\circ \epsilon $. In particular, if $M(\mathbf{\vartheta },\mathbf{%
\upsilon })=0$ for $\sum |\vartheta _{\nu }^{\mu }|\neq 1$, then, obviously 
\begin{equation*}
\nu _{0}^{t}(\mathbf{\omega },M)=n_{0}^{t}(\mathbf{\omega },\mathbf{C}%
),\quad \iota _{0}^{t}(B)=i_{0}^{t}(\mathbf{D}),
\end{equation*}%
where $C_{\nu }^{\mu }(x,\mathbf{\upsilon })=M(\mathbf{x}_{\nu }^{\mu },%
\mathbf{\upsilon })\equiv C\left( \mathbf{x}_{\nu }^{\mu },\mathbf{\upsilon }%
\right) $ and $B(\mathbf{\vartheta })=0$ for $\sum |\vartheta _{\nu }^{\mu
}|\neq 1$, $D_{\nu }^{\mu }(x)=B(\mathbf{x}_{\nu }^{\mu })$. This yields the
representation $\epsilon \circ n_{0}^{t}=i_{0}^{t}\circ \epsilon $ for the
single generalized non-adapted integral (\ref{2oneb}) for $\Delta =X^{t}$ in
the form of the sum 
\begin{equation*}
\sum_{\mu ,\nu }\Lambda _{\nu }^{\mu }(\epsilon (C_{\nu }^{\mu }),\Delta
)=\epsilon \Big(\sum_{\mu ,\nu }(N_{\mu }^{\nu }(C_{\nu }^{\mu },\Delta )%
\Big),\quad N_{\mu }^{\nu }(\mathbf{\omega },C,\Delta )=\sum_{x\in \omega
_{\nu }^{\mu }\cap \Delta }C(x,\mathbf{\omega }\setminus \mathbf{x}_{\nu
}^{\mu })
\end{equation*}%
of representations of four kernel measures $N_{\nu }^{\mu }(\mathbf{\omega }%
,C_{\nu }^{\mu },\Delta )$ for that define kernel representations $\epsilon
\circ N(\Delta )=\Lambda (\Delta )\circ \epsilon $ of the canonical measures
(4.3) with $D_{\nu }^{\mu }(x)=\epsilon (C_{\nu }^{\mu }(x))$.
\end{proof}

In the following theorem, which generalizes It\^{o} formula to
noncommutative and nonadapted quantum stochastic processes $T_{t}=\epsilon
\left( K_{t}\right) $ given by an operator-valued kernel $K_{t}\left( 
\mathbf{\omega }\right) $, we use the following triangular-matrix notation 
\begin{equation*}
\mathbf{T}\left( x\right) =\left[ T\left( \mathbf{x}_{\nu }^{\mu }\right) %
\right] ,\;\;T\left( \mathbf{x}\right) =\nabla _{\mathbf{x}%
}T_{t}|_{t=t\left( x\right) }
\end{equation*}%
for the \emph{quantum stochastic germs} $\nabla _{\mathbf{x}}T=\epsilon
\left( \dot{K}\left( \mathbf{x}\right) \right) $ given by the point
derivatives of the kernel $\dot{K}\left( \mathbf{x},\mathbf{\upsilon }%
\right) =K\left( \mathbf{\upsilon }\sqcup \mathbf{x}\right) $,\ with $T_{\nu
}^{\mu }\left( x\right) =T\left( \mathbf{x}_{\nu }^{\mu }\right) $ equal
zero for $\mu =+$ or $\nu =-$ and $T_{-}^{-}\left( x\right) =T_{t\left(
x\right) }=T_{+}^{+}\left( x\right) $. We notice that if $K_{t}\left( 
\mathbf{\omega }\right) =K_{0}\left( \mathbf{\omega }\right)
+n_{0}^{t}\left( \mathbf{C}\left( \mathbf{\omega }\right) \right) $,
corresponding to the single-integral representation $T_{t}-T_{0}=i_{0}^{t}%
\left( \mathbf{D}\right) $\ with $\mathbf{D}\left( x\right) =\epsilon \left( 
\mathbf{C}\left( x\right) \right) $, then $\dot{K}_{t}(\mathbf{x},\mathbf{%
\upsilon })=K_{t}(\mathbf{\upsilon }\sqcup \mathbf{x})$ is given by%
\begin{equation*}
\dot{K}_{t}(\mathbf{x},\mathbf{\upsilon })=\dot{K}_{t\wedge t\left( x\right)
}\left( \mathbf{x},\mathbf{\upsilon }\right) +\sum_{\mathbf{z}\in \mathbf{%
\upsilon }}^{t\left( x\right) \leq t\left( z\right) <t}C\left( \mathbf{z},%
\mathbf{\upsilon }\backslash \mathbf{z}\sqcup \mathbf{x}\right) .
\end{equation*}%
This proves that $\dot{K}_{t}(\mathbf{x},\mathbf{\upsilon })$ doesn't depend
on $t\in (t\left( x\right) ,t^{+}(x)]$, where $t^{+}\left( x\right) =\min
\left\{ t\left( x^{\prime }\right) >t\left( x\right) :x^{\prime }\in \sqcup
\upsilon _{\nu }^{\mu }\right\} $, and therefore the right limit%
\begin{equation*}
\dot{K}_{t\left( x\right) ]}(\mathbf{x},\mathbf{\upsilon }):=\lim_{t\searrow
t\left( x\right) }\dot{K}_{t}(\mathbf{x},\mathbf{\upsilon })=\dot{K}%
_{t\left( x\right) }\left( \mathbf{x},\mathbf{\upsilon }\right) +C\left( 
\mathbf{x},\mathbf{\upsilon }\right)
\end{equation*}%
trivially exists for each $\mathbf{x}\in \left\{ \mathbf{x}_{\nu }^{\mu
}\right\} $ and $\mathbf{\upsilon }$ with $\dot{K}_{t\left( x\right) ]}(%
\mathbf{x}_{-}^{-},\mathbf{\upsilon })=K_{t\left( x\right) }(\mathbf{%
\upsilon })=\dot{K}_{t\left( x\right) ]}(\mathbf{x}_{-}^{-},\mathbf{\upsilon 
})$ for $\dot{K}_{t}(\mathbf{x}_{-}^{-},\mathbf{\upsilon })=K_{t}(\mathbf{%
\upsilon })=\dot{K}_{t}(\mathbf{x}_{+}^{+},\mathbf{\upsilon })$ due to the
independency of $K\left( \mathbf{\omega }\right) $ on $\omega _{-}^{-}$ and $%
\omega _{+}^{+}$. We may assume that the germs $\nabla _{\mathbf{x}%
}T_{t}=\epsilon \left( \dot{K}_{t}\left( \mathbf{x}\right) \right) $ also
converge from the right to $G\left( \mathbf{x}\right) =T\left( \mathbf{x}%
\right) +D\left( \mathbf{x}\right) $ with $D\left( \mathbf{x}\right)
=\epsilon \left( C\left( \mathbf{x}\right) \right) $ at $t\searrow t\left(
x\right) $ for $x\in X$ corresponding to each atomic table $\mathbf{x}$ in (%
\ref{2oneg}) as they have limits $\epsilon \left( \dot{K}_{t\left( x\right)
}\left( \mathbf{x}_{-}^{-}\right) \right) =T_{t\left( x\right) }=\epsilon
\left( \dot{K}_{t\left( x\right) }\left( \mathbf{x}_{+}^{+}\right) \right) $
for $\mathbf{x}\in \left\{ \mathbf{x}_{-}^{-},\mathbf{x}_{+}^{+}\right\} $.
As it is proved in the following theorem, these germ-limits $G\left( \mathbf{%
x}\right) $ are given by the matrix elements $D\left( \mathbf{x}_{\nu }^{\mu
}\right) $ of the QS-derivatives $\mathbf{D}=\left[ D_{\nu }^{\mu }\left( 
\mathbf{x}\right) \right] $ at least in the case $K_{t}=\nu _{0}^{t}(M)$ 
\textup{(\ref{two7}) }corresponding to the multiple integral representation $%
T_{t}=\iota _{0}^{t}(B)$ \textup{(}see \textup{(\ref{2onee}))} with $B(%
\mathbf{\vartheta })=\epsilon (M(\mathbf{\vartheta }))$.

\begin{theorem}
\label{2T 2} If kernel $K(\mathbf{\omega })$ is relatively bounded, then the
same is true for the kernel $K^{\star }(\mathbf{\omega }):\Vert K^{\star
}\Vert _{\mathbf{\gamma }}=\Vert K\Vert _{\mathbf{\gamma }^{\prime }}$,
where $%
\begin{pmatrix}
\gamma _{+}^{-} & \gamma _{\circ }^{-} \\ 
\gamma _{+}^{\circ } & \gamma _{\circ }^{\circ }%
\end{pmatrix}%
^{\prime }=%
\begin{pmatrix}
\gamma _{+}^{-} & \gamma _{+}^{\circ } \\ 
\gamma _{\circ }^{-} & \gamma _{\circ }^{\circ }%
\end{pmatrix}%
$, and the operator $T^{\ast }=\epsilon (K^{\star })$, as well as the
operator $T=\epsilon (K)$, is $q$-bounded by the estimate \textup{(\ref{two6}%
)} for $q\geq p+1/r$. For any such kernels $K(\mathbf{\vartheta })$ and $%
K^{\star }(\mathbf{\vartheta })$, bounded relative to the quadruples $%
\mathbf{\alpha }=(\alpha _{\nu }^{\mu })$ and $\mathbf{\gamma }=(\gamma
_{\nu }^{\mu })$ of functions $\alpha _{\nu }^{\mu }(x),\,\gamma _{\nu
}^{\mu }(x)$ satisfying \textup{(\ref{two5})} the operator 
\begin{equation*}
\epsilon (K)\epsilon (K)^{\ast }=\epsilon (K\cdot K^{\star }),\quad \epsilon
(I^{\otimes })=I
\end{equation*}%
is well-defined as a $\ast $-representation of kernel product \textup{(2.9)}
of \textup{Chapter I} with the estimate $\Vert K\cdot K^{\star }\Vert _{%
\mathbf{\beta }}\leq \Vert K\Vert _{\mathbf{\alpha }}\Vert K^{\star }\Vert _{%
\mathbf{\gamma }}$ if $\beta _{\nu }^{\mu }\geq (\mathbf{\alpha \cdot \gamma 
})_{\nu }^{\mu }$, where $(\mathbf{\alpha \cdot \gamma })_{\nu }^{\mu
}(x)=\sum \alpha _{\lambda }^{\mu }(x)\gamma _{\nu }^{\lambda }(x)$ is
defined by the product of triangular matrices 
\begin{equation*}
\left[ 
\begin{array}{ccc}
1 & \alpha _{\circ }^{-} & \alpha _{+}^{-} \\ 
0 & \alpha _{\circ }^{\circ } & \alpha _{+}^{\circ } \\ 
0 & 0 & 1%
\end{array}%
\right] \left[ 
\begin{array}{ccc}
1 & \gamma _{\circ }^{-} & \gamma _{+}^{-} \\ 
0 & \gamma _{\circ }^{\circ } & \gamma _{+}^{\circ } \\ 
0 & 0 & 1%
\end{array}%
\right] =\left[ 
\begin{array}{ccc}
1 & \alpha _{\circ }^{-}\gamma _{\circ }^{\circ }+\gamma _{\circ }^{-}, & 
\gamma _{+}^{-}+\alpha _{\circ }^{-}\gamma _{+}^{\circ }+\alpha _{+}^{-} \\ 
0, & \alpha _{\circ }^{\circ }\gamma _{\circ }^{\circ }, & \alpha _{\circ
}^{\circ }\gamma _{+}^{\circ }+\gamma _{+}^{\circ } \\ 
0, & 0, & 1%
\end{array}%
\right] .
\end{equation*}%
Let $T_{t}=\epsilon (K_{t})$ with $\dot{K}_{t}\left( \mathbf{x,\upsilon }%
\right) $ defining the right limit $\nabla _{\mathbf{x}}T_{t}|_{t=t\left(
x\right) ]}=\epsilon (\dot{K}_{t(x)]}(\mathbf{x}))$ of $\nabla _{\mathbf{x}%
}T_{t}$ at $t\searrow t\left( x\right) $. Let $\mathbf{T}(x)=[T_{\nu }^{\mu
}(x)]$ and $\mathbf{G}(x)=[G_{\nu }^{\mu }(x)]$ denote the triangular
matrices of germs $T\left( \mathbf{x}\right) =\nabla _{\mathbf{x}%
}T_{t}|_{t=t\left( x\right) }$ and $G(\mathbf{x})=\nabla _{\mathbf{x}%
}T_{t}|_{t=t\left( x\right) ]}$ as operator-valued matrix elements 
\begin{equation}
T_{\nu }^{\mu }\left( x\right) =\epsilon (\dot{K}_{t(x)}(\mathbf{x}_{\nu
}^{\mu })),\quad G_{\nu }^{\mu }(x)=\epsilon (\dot{K}_{t(x)]}(\mathbf{x}%
_{\nu }^{\mu }))  \label{two8}
\end{equation}%
corresponding to point-derivatives $\dot{K}_{t}\left( \mathbf{x}_{\nu }^{\mu
}\right) $ at $t=t\left( x\right) $ and their right limits at $t=t\left(
x\right) ]$ respectively. Then the operator-functions $D_{\nu }^{\mu
}(x)=G_{\nu }^{\mu }(x)-T_{\nu }^{\mu }(x)$ are quantum-stochastic
derivatives of the function $t\mapsto T_{t}$ which define the QS
differential $\mathrm{d}T_{t}=\mathrm{d}i_{0}^{t}(\mathbf{D})$ in the
difference form so that $T_{t}-T_{0}=i_{0}^{t}(\mathbf{G}-\mathbf{T})$.
Moreover, $T_{t}^{\ast }-T_{0}^{\ast }=i_{0}^{t}(\mathbf{G}^{\dagger }-%
\mathbf{T}^{\dagger })$, and we have the generalized non-adapted It\^{o}
formula 
\begin{equation}
T_{t}T_{t}^{\ast }-T_{0}T_{0}^{\ast }=i_{0}^{t}(\mathbf{TD}^{\dagger }+%
\mathbf{DT}^{\dagger }+\mathbf{DD}^{\dagger })=i_{0}^{t}(\mathbf{GG}%
^{\dagger }-\mathbf{TT}^{\dagger }),  \label{two9}
\end{equation}%
where $\mathbf{D}\mapsto \mathbf{D}^{\dagger }$ is the pseudo-Euclidean
conjugation $[D_{\nu }^{\mu }(x)]^{\dagger }=[D_{-\mu }^{-\nu }(x)]^{\ast }$
of the triangular operators 
\begin{equation*}
\mathbf{T}=\left[ 
\begin{array}{ccc}
T & T_{\circ }^{-} & T_{+}^{-} \\ 
0 & T_{\circ }^{\circ } & T_{+}^{\circ } \\ 
0 & 0 & T%
\end{array}%
\right] ,\quad \mathbf{D}=\left[ 
\begin{array}{ccc}
0 & D_{\circ }^{-} & D_{+}^{\circ } \\ 
0 & D_{\circ }^{\circ } & D_{+}^{\circ } \\ 
0 & 0 & 0%
\end{array}%
\right] ,\quad \mathbf{G}=\left[ 
\begin{array}{ccc}
T & G_{\circ }^{-} & G_{+}^{-} \\ 
0 & G_{\circ }^{\circ } & G_{+}^{\circ } \\ 
0 & 0 & T%
\end{array}%
\right]
\end{equation*}%
with the standard block-matrix multiplication $(\mathbf{TG})_{\nu }^{\mu
}=\Sigma T_{\lambda }^{\mu }G_{\nu }^{\lambda }$.
\end{theorem}

\begin{proof}
The adjoint operators $\epsilon (K)$ and $\epsilon (K^{\star })$, which
define the $\ast $-representation (\ref{two2}) with respect to the kernels $%
K $ bounded in the sense of (\ref{two4}) and (\ref{two5}), are $q$-bounded
for $q\geq p+1/r$ by the estimate $\Vert \epsilon (K)\Vert _{q}\leq \Vert
K\Vert _{p}(r)$ and inequality (\ref{two6}), which leads to the exponential
estimate 
\begin{equation*}
\Vert \epsilon (K)\Vert _{q}\leq \Vert K\Vert _{\mathbf{\alpha }}\exp
\{\Vert \alpha _{+}^{-}\Vert ^{(1)}+{\tfrac{1}{2}}(\Vert \alpha _{+}^{\circ
}\Vert ^{(2)}(r)^{2}+\Vert \alpha _{\circ }^{-}\Vert ^{(2)}(r)^{2})\}.
\end{equation*}%
The formula for the kernel multiplication $K^{\star }\cdot K$, which
corresponds to the operator product $\epsilon (K^{\star })\epsilon (K)$, has
already been found for scalar $\mathrm{H}=\mathbb{C}$ in the case of linear
combinations of exponential kernels 
\begin{equation*}
\mathbf{f}^{\otimes }(\mathbf{\vartheta })=f_{+}^{-}(\vartheta
_{+}^{-})f_{+}^{\circ }(\vartheta _{+}^{\circ })\otimes f_{\circ }^{\circ
}(\vartheta _{\circ }^{\circ })\otimes f_{\circ }^{-}(\vartheta _{\circ
}^{-}),
\end{equation*}%
where $f_{\nu }^{\mu }(\vartheta )=\bigotimes_{x\in \vartheta
}f(x)(f_{+}^{-}(\vartheta )=\prod_{x\in \vartheta }f_{+}^{-}(x))$. We shall
now verify this formula for operator-valued kernels $K(\omega )$ and $%
K^{\star }(\omega )$, noticing that their product is $\mathbf{\beta }$%
-bounded for $\mathbf{\beta }=\mathbf{\alpha \cdot \gamma }$, since $%
\left\Vert K^{\star }\cdot K\right\Vert (\omega )\leq $ 
\begin{eqnarray*}
&\leq &\sum \left\Vert K\left( 
\begin{array}{ll}
\omega _{+}^{-}\setminus \sigma _{+}^{-}, & \upsilon _{\circ }^{-}\sqcup
\upsilon _{+}^{-} \\ 
\omega _{+}^{\circ }\setminus \upsilon _{+}^{\circ }, & \omega _{\circ
}^{\circ }\sqcup \upsilon _{+}^{\circ }%
\end{array}%
\right) \right\Vert \cdot \left\Vert K^{\star }\left( 
\begin{array}{ll}
\omega _{+}^{-}\setminus \tau _{+}^{-}, & \omega _{\circ }^{-}\setminus
\upsilon _{\circ }^{-} \\ 
\upsilon _{+}^{-}\sqcup \upsilon _{+}^{\circ }, & \omega _{\circ }^{\circ
}\sqcup \upsilon _{\circ }^{-}%
\end{array}%
\right) \right\Vert \\
&\leq &\left\Vert K\right\Vert _{\mathbf{\alpha }}\left\Vert K^{\star
}\right\Vert _{\mathbf{\gamma }}\sum \mathbf{\alpha }^{\otimes }\left( 
\begin{array}{ll}
\omega _{+}^{-}\setminus \sigma _{+}^{-}, & \upsilon _{\circ }^{-}\sqcup
\upsilon _{+}^{-} \\ 
\omega _{+}^{\circ }\setminus \upsilon _{+}^{\circ }, & \omega _{\circ
}^{\circ }\sqcup \upsilon _{+}^{\circ }%
\end{array}%
\right) \mathbf{\gamma }^{\otimes }\left( 
\begin{array}{ll}
\omega _{+}^{-}\setminus \tau _{+}^{-}, & \omega _{\circ }^{-}\setminus
\upsilon _{\circ }^{-} \\ 
\upsilon _{+}^{-}\sqcup \upsilon _{+}^{\circ }, & \omega _{\circ }^{\circ
}\sqcup \upsilon _{\circ }^{-}%
\end{array}%
\right) \\
&=&\left\Vert K^{\star }\right\Vert _{\mathbf{\gamma }}\left\Vert
K\right\Vert _{\mathbf{\alpha }}(\mathbf{\alpha \cdot \gamma })^{\otimes
}(\omega );(\mathbf{\alpha \cdot \gamma })_{\nu }^{\mu }=\sum_{\mu \leq
\lambda \leq \nu }\alpha _{\lambda }^{\mu }\gamma _{\nu }^{\lambda },
\end{eqnarray*}%
where we have employed the multiplication formula $\mathbf{\alpha }^{\otimes
}\cdot \mathbf{\gamma }^{\otimes }=(\mathbf{\alpha }\cdot \mathbf{\gamma }%
)^{\otimes }$ for scalar exponential kernels 
\begin{equation*}
\mathbf{\beta }^{\otimes }(\mathbf{\omega })=\prod \beta _{\nu }^{\mu
}(\omega _{\nu }^{\mu });\,\beta _{\nu }^{\mu }(\omega )=\prod_{x\in \omega
}\beta _{\nu }^{\mu }(x):(\mathbf{\alpha \cdot \gamma })_{\nu }^{\mu
}(x)=\sum \gamma _{\lambda }^{\mu }(x)\alpha _{\nu }^{\lambda }(x).
\end{equation*}%
Using the main formula (\ref{2oned}) of the scalar integration we write the
scalar square of the action (\ref{two2}) in the form $\left\Vert \epsilon (K)%
\mathrm{h}\right\Vert ^{2}=$ 
\begin{eqnarray*}
&=&\int \left\Vert \sum_{\omega _{\circ }^{\circ }\sqcup \omega _{+}^{\circ
}=\varkappa }\iint K^{\star }(\boldsymbol{\omega })\mathrm{h}(\omega _{\circ
}^{\circ }\sqcup \omega _{\circ }^{\circ })\mathrm{d}\omega _{+}^{-}\mathrm{d%
}\omega _{\circ }^{-}\right\Vert ^{2}\mathrm{d}\varkappa \\
&=&\int \sum_{\sigma _{\circ }^{\circ }\sqcup \sigma _{+}^{\circ }=\varkappa
}\sum_{\tau _{\circ }^{\circ }\sqcup \tau _{+}^{\circ }=\varkappa
}\left\langle K^{\star }(\boldsymbol{\sigma })\mathrm{h}(\sigma _{\circ
}^{-}\sqcup \sigma _{\circ }^{\circ })\mid K^{\star }(\boldsymbol{\tau })%
\mathrm{h}(\tau _{\circ }^{-}\sqcup \tau _{\circ }^{\circ })\right\rangle 
\mathrm{d}\varkappa \\
&=&\iiint \langle K^{\star }\left( 
\begin{array}{cc}
\sigma _{+}^{-}, & \sigma _{+}^{\circ } \\ 
\upsilon _{\circ }^{-}\sqcup \upsilon _{+}^{-}, & \upsilon _{\circ }^{\circ
}\sqcup \upsilon _{+}^{\circ }%
\end{array}%
\right) \mathrm{h}(\upsilon _{\circ }^{\circ }\sqcup \upsilon _{+}^{\circ
}\sqcup \sigma _{+}^{\circ })\mid \\
&&\quad K^{\star }\left( 
\begin{array}{cc}
\tau _{+}^{-}, & \tau _{\circ }^{-} \\ 
\upsilon _{+}^{\circ }\sqcup \upsilon _{+}^{-}, & \upsilon _{\circ }^{\circ
}\sqcup \upsilon _{\circ }^{-}%
\end{array}%
\right) \mathrm{h}(\upsilon _{\circ }^{\circ }\sqcup \upsilon _{\circ
}^{-}\sqcup \tau _{\circ }^{-})\rangle \mathrm{d}\sigma \,\mathrm{d}x\,%
\mathrm{d}\upsilon \\
&=&\iiint \langle \mathrm{h}(\upsilon _{\circ }^{\circ }\sqcup \upsilon
_{+}^{\circ }\sqcup \sigma _{+}^{\circ })\mid K\left( 
\begin{array}{cc}
\sigma _{+}^{-}, & \upsilon _{\circ }^{-}\sqcup \upsilon _{+}^{-} \\ 
\sigma _{+}^{\circ }, & \upsilon _{\circ }^{\circ }\sqcup \upsilon
_{+}^{\circ }%
\end{array}%
\right) \\
&&\quad \times K^{\star }\left( 
\begin{array}{cc}
\tau _{+}^{-}, & \tau _{\circ }^{-} \\ 
\upsilon _{+}^{\circ }\sqcup \upsilon _{+}^{-}, & \upsilon _{\circ }^{\circ
}\sqcup \upsilon _{\circ }^{-}%
\end{array}%
\right) \mathrm{h}(\upsilon _{\circ }^{\circ }\sqcup \upsilon _{\circ
}^{-}\sqcup \tau _{\circ }^{-})\rangle \mathrm{d}\sigma \,\mathrm{d}x\,%
\mathrm{d}\upsilon \\
&=&\int (\mathrm{h}(\varkappa )\mid \sum_{\omega _{\circ }^{\circ }\sqcup
\omega _{+}^{\circ }=\vartheta }\iint (K\cdot K^{\star })(\omega )\mathrm{h}%
(\omega _{\circ }^{-}\sqcup \omega _{\circ }^{\circ })\mathrm{d}\omega
_{+}^{-}\mathrm{d}\omega _{\circ }^{-})\mathrm{d}\vartheta ,
\end{eqnarray*}%
where $\upsilon _{\circ }^{\circ }=\sigma _{\circ }^{\circ }\cap \tau
_{\circ }^{\circ }$, $\upsilon _{+}^{\circ }=\sigma _{\circ }^{\circ }\cap
\tau _{+}^{\circ }$, $\upsilon _{\circ }^{-}=\tau _{\circ }^{\circ }\cap
\sigma _{+}^{\circ }$, and $\upsilon _{+}^{-}=\sigma _{+}^{\circ }\cap \tau
_{+}^{\circ }$. Since $\mathrm{h}\in \mathrm{H}\otimes \mathrm{F}(q)$ is
arbitrary, this proves the formula for the kernel multiplication for $%
K^{\star }$ and $K$, which extends to any relatively bounded kernels $M$ and 
$K$ because of the polarization formula for the Hermitian function $K^{\star
}\cdot K$.

We shall now consider the stochastic differential $\mathrm{d}T_{t}$ of the
multiple integral $T_{t}=\iota _{0}^{t}(B)$ of the operator function $B(%
\mathbf{\vartheta })=\epsilon (M(\mathbf{\vartheta }))$ defined by the
quantum-stochastic derivatives 
\begin{equation*}
D_{\nu }^{\mu }(x)=\iota _{0}^{t(x)}(\dot{B}(\mathbf{x}_{\nu }^{\mu
}))=\epsilon (C_{\nu }^{\mu }(x)),
\end{equation*}%
representing the differences of the kernels 
\begin{equation*}
C_{\nu }^{\mu }(x,\mathbf{\upsilon })=\nu _{0}^{t(x)}(\mathbf{\upsilon },%
\dot{M}(\mathbf{x}_{\nu }^{\mu }))=\dot{K}_{t(x)]}(\mathbf{x}_{\nu }^{\mu },%
\mathbf{\upsilon })-\dot{K}_{t(x)}(\mathbf{x}_{\nu }^{\mu },\mathbf{\upsilon 
}).
\end{equation*}%
Here $\nu _{0}^{t}(\mathbf{\upsilon },\dot{M}(\mathbf{x}))=\sum_{\mathbf{%
\vartheta }\subseteq \mathbf{\upsilon }^{t}}M(\mathbf{\vartheta }\sqcup 
\mathbf{x},\mathbf{\upsilon }\setminus \mathbf{\vartheta })$, $\mathbf{x}$
is one of the atomic matrices (\ref{2oneg}), and 
\begin{eqnarray*}
\dot{K}_{t(x)}(\mathbf{x},\mathbf{\upsilon }) &=&\sum_{\mathbf{\vartheta }%
\subseteq \upsilon ^{t(x)}}M(\mathbf{\vartheta },(\mathbf{\upsilon }\sqcup 
\mathbf{x})\setminus \mathbf{\vartheta })=K_{t(x)}(\mathbf{\upsilon }\sqcup 
\mathbf{x}), \\
\dot{K}_{t(x)]}(\mathbf{x},\mathbf{\upsilon }) &=&\sum_{\mathbf{\vartheta }%
\subseteq \upsilon ^{t(x)}\sqcup \mathbf{x}}M(\mathbf{\vartheta },(\mathbf{%
\upsilon }\sqcup \mathbf{x})\setminus \mathbf{\vartheta }) \\
&=&K_{t(x)}(\mathbf{\upsilon }\sqcup \mathbf{x})+\sum_{\mathbf{\vartheta }%
\subseteq \mathbf{\upsilon }^{t(x)}}M(\mathbf{\vartheta }\sqcup \mathbf{x},%
\mathbf{\upsilon }\setminus \mathbf{\vartheta }) \\
&=&\dot{K}_{t(x)}(\mathbf{x},\mathbf{\upsilon })+\nu _{0}^{t(x)}(\mathbf{%
\upsilon },\dot{M}(\mathbf{x})).
\end{eqnarray*}%
We note that $K_{t]}(\mathbf{\omega })=\sum_{\mathbf{\vartheta }\subseteq 
\mathbf{\omega }^{t]}}M(\mathbf{\vartheta },\mathbf{\omega }\setminus 
\mathbf{\vartheta })=K_{t_{+}}(\mathbf{\omega })$, where $t_{+}=\min \{t(%
\mathbf{x})>t:\mathbf{x}\in \mathbf{\omega }\}$, $\mathbf{\omega }^{t]}=\{%
\mathbf{x}\in \mathbf{\omega }:t(x)\leq t\}$, so that $\dot{K}_{t(x)]}(%
\mathbf{x,\upsilon })=\dot{K}_{t}(\mathbf{x,\upsilon })$ for any $t\in (t(%
\mathbf{x}),t_{+}(\mathbf{x})]$. Thus the derivatives $D_{\nu }^{\mu
}(x),x\in X^{t}$, defining the increment $T_{t}-T_{0}=i_{0}^{t}(\mathbf{D})$%
, can be written in the form of the differences 
\begin{equation*}
D_{\nu }^{\mu }(x)=\epsilon \lbrack \dot{K}_{t(x)]}(\mathbf{x}_{\nu }^{\mu
})]-\epsilon \lbrack \dot{K}_{t(x)}(\mathbf{x}_{\nu }^{\mu })]
\end{equation*}%
of the operators (\ref{two8}). If we consider $\dot{K}_{t}(\mathbf{x})$ as
one of the four entries $\dot{K}_{t}(\mathbf{x}_{\nu }^{\mu })=K_{t}(x)_{\nu
}^{\mu }$ in the triangular operator kernel $\mathbf{K}_{t}(x)$ with $%
K_{t}(x)_{-}^{-}=K_{t(x)}=K_{t}(x)_{+}^{+}$, we define the triangular
functions 
\begin{equation*}
\mathbf{T}(x)=\epsilon (\mathbf{K}_{t(x)}(x)),\;\;\;\;\,\mathbf{G}%
(x)=\epsilon (\mathbf{K}_{t(x)]}(x)).
\end{equation*}%
This allows us to obtain the quantum non-adapted It\^{o} formula in the form 
\begin{equation*}
T_{t}T_{t}^{\ast }-T_{0}T_{0}^{\ast }=i_{0}^{t}(\mathbf{TD}^{\dagger }+%
\mathbf{DT}^{\dagger }+\mathbf{DD}^{\dagger }),
\end{equation*}%
where $\mathbf{D}(x)=\mathbf{G}(x)-\mathbf{T}(x)$. This is a consequence of
the fact that the map (\ref{two2}) is a $\star $-homomorphism, $%
T_{t}T_{t}^{\ast }=\epsilon (K\cdot K^{\star })$, and the formula (\ref{two
i}) of Chapter I for the product of the operator kernels $K_{t}$ and $%
K_{t}^{\star }$, which can be written in the form 
\begin{equation*}
(K_{t}\cdot K_{t}^{\star })(\mathbf{\omega }\sqcup \mathbf{x}_{\nu }^{\mu
})=\sum_{\lambda =\mu }^{\nu }[K_{t}(x)_{\lambda }^{\mu }\cdot K_{t}^{\star
}(x)_{\nu }^{\lambda }](\mathbf{\omega })=[\mathbf{K}_{t}(x)\mathbf{K}%
_{t}^{\dagger }(x)]_{\nu }^{\mu }(\mathbf{\omega }),
\end{equation*}%
where the right-hand side is computed as an entry in the product of
triangular matrices $\mathbf{K}(x)=[K_{\nu }^{\mu }(x)]$ which defines the
multiplication of the entries as operator-valued kernels $K_{t}(x,\mathbf{%
\omega })_{-}^{-}=K_{t}(\mathbf{\omega })=K_{t}(x,\mathbf{\omega })_{+}^{+}$%
, $\dot{K}(\mathbf{x},\mathbf{\omega })=K(\mathbf{\omega }\sqcup \mathbf{x})$%
. For from (\ref{two i}) of Chapter I we obtain 
\begin{eqnarray*}
\lbrack K\cdot K^{\star }](\mathbf{\omega }\sqcup \mathbf{x}_{\circ }^{\circ
}) &=&[\dot{K}(\mathbf{x}_{\circ }^{\circ })\cdot \dot{K}^{\star }(\mathbf{x}%
_{\circ }^{\circ })](\mathbf{\omega }), \\
\lbrack K\cdot K^{\star }](\mathbf{\omega }\sqcup \mathbf{x}_{+}^{\circ })
&=&[K\cdot \dot{K}^{\star }(\mathbf{x}_{\circ }^{-})+\dot{K}(\mathbf{x}%
_{\circ }^{-})\dot{K}^{\star }(\mathbf{x}_{\circ }^{\circ })](\mathbf{\omega 
}), \\
\lbrack K\cdot K^{\star }](\mathbf{\omega }\sqcup \mathbf{x}_{\circ }^{-})
&=&[\dot{K}(\mathbf{x}_{\circ }^{\circ })\dot{K}^{\star }(\mathbf{x}%
_{+}^{\circ })+\dot{K}(\mathbf{x}_{+}^{\circ })\cdot K^{\star }](\mathbf{%
\omega }), \\
\lbrack K\cdot K^{\star }](\mathbf{\omega }\sqcup \mathbf{x}_{+}^{-})
&=&[K\cdot \dot{K}^{\star }(\mathbf{x}_{+}^{-})+\dot{K}(\mathbf{x}_{\circ
}^{-})\cdot \dot{K}^{\star }(\mathbf{x}_{+}^{\circ })+\dot{K}(\mathbf{x}%
_{+}^{-})\cdot K^{\star }](\mathbf{\omega }),
\end{eqnarray*}%
which are the matrix elements of%
\begin{equation*}
\lbrack K\cdot K^{\star }](\mathbf{\omega }\sqcup \mathbf{x})=\left[ \dot{K}%
(x)_{\lambda }^{\mu }\cdot \dot{K}_{t}^{\star }(x)_{\nu }^{\lambda }\right]
\left( \mathbf{\omega }\right) =\left( \mathbf{K}\cdot \mathbf{K}^{\dagger
}\right) \left( x,\mathbf{\omega }\right) .
\end{equation*}%
This allows us to write $\epsilon \lbrack (\dot{K}_{t}\cdot \dot{K}%
_{t}^{\star })($\textbf{$x$}$_{\nu }^{\mu })]=\sum_{\lambda =\mu }^{\nu
}\epsilon (\dot{K}_{t}(x)_{\lambda }^{\mu }\dot{K}_{t}^{\star }(x)_{\nu
}^{\lambda })$ in the form of the triangular operator 
\begin{equation*}
\epsilon (\mathbf{K}_{t}(x)\mathbf{K}_{t}^{\dagger }(x))=\epsilon (\mathbf{K}%
_{t}(x))\epsilon (\mathbf{K}_{t}(x))^{\ast }
\end{equation*}%
which is the product of the triangular matrices $\mathbf{T}_{t}(x)$ and $%
\mathbf{T}_{t}^{\dagger }(x)$ with operator product of the entries. We put $%
t=t(x)$ and $t=t_{+}(x)$ in the formula, and we obtain 
\begin{equation*}
\epsilon \lbrack (\mathbf{K}_{t(x)]}\cdot \mathbf{K}_{t(x)]}^{\dagger })(x)-(%
\mathbf{K}_{t(x)}\cdot \mathbf{K}_{t(x)}^{\dagger })(x)]=\mathbf{G}(x)%
\mathbf{G}^{\dagger }(x)-\mathbf{T}(x)\mathbf{T}^{\dagger }(x),
\end{equation*}%
which allows us to write the stochastic derivative of the quantum
non-adapted process $T_{t}T_{t}^{\ast }$ in the form 
\begin{equation*}
\mathrm{d}(T_{t}T_{t}^{\ast })=\mathrm{d}i_{0}^{t}(\mathbf{GG}^{\dagger }-%
\mathbf{TT}^{\dagger }),
\end{equation*}%
corresponding to (\ref{two9}). The theorem has been proved.
\end{proof}

\begin{remark}
Using the non-adapted table of stochastic multiplication 
\begin{eqnarray*}
\mathbf{G}^{\dagger }\mathbf{G}-\mathbf{T}^{\dagger }\mathbf{T} &=&\mathbf{D}%
^{\dagger }\mathbf{T}+\mathbf{T}^{\dagger }\mathbf{D}+\mathbf{D}^{\dagger }%
\mathbf{D} \\
&=&\left[ 
\begin{array}{ccc}
0, & T^{\ast }D_{\circ }^{-}, & T^{\ast }D_{+}^{-}+D_{+}^{-\ast }T \\ 
0, & 0, & D_{\circ }^{-\ast }T \\ 
0, & 0, & 0%
\end{array}%
\right] +\left[ 
\begin{array}{ccc}
0, & D_{+}^{\circ \ast }D_{\circ }^{\circ }, & D_{+}^{\circ \ast
}D_{+}^{\circ } \\ 
0, & D_{\circ }^{\circ \ast }D_{\circ }^{\circ }, & D_{\circ }^{\circ \ast
}D_{+}^{\circ } \\ 
0, & 0, & 0%
\end{array}%
\right] \\
&&+\left[ 
\begin{array}{ccc}
0, & D_{+}^{\circ \ast }T_{\circ }^{\circ }+T_{+}^{\circ \ast }D_{\circ
}^{\circ }, & D_{+}^{\circ \ast }T_{+}^{\circ }+T_{+}^{\circ \ast
}D_{+}^{\circ } \\ 
0, & D_{\circ }^{\circ \ast }T_{\circ }^{\circ }+T_{\circ }^{\circ \ast
}D_{\circ }^{\circ }, & D_{\circ }^{\circ \ast }T_{+}^{\circ }+T_{\circ
}^{\circ \ast }D_{+}^{\circ } \\ 
0, & 0, & 0%
\end{array}%
\right]
\end{eqnarray*}%
we can write \textup{(\ref{two9})} in a weak form 
\begin{eqnarray}
&&\left\Vert T_{t}\mathrm{h}\right\Vert ^{2}-\left\Vert T_{0}\mathrm{h}%
\right\Vert ^{2}=\int_{X^{t}}2\mathbf{\func{Re}}\left\langle T_{t(x)}\mathrm{%
h}\mid D_{+}^{-}(x)\mathrm{h}+D_{\circ }^{-}(x)\mathrm{\dot{h}}%
(x)\right\rangle \mathrm{d}x  \label{two10} \\
&&+\int_{X^{t}}\left[ \left\Vert D_{+}^{\circ }(x)\mathrm{\dot{h}}+D_{\circ
}^{\circ }(x)\mathrm{\dot{h}}(x)\right\Vert ^{2}+2\mathbf{\func{Re}}%
\left\langle a(x)T_{t(x)}\mathrm{h}\mid D_{+}^{\circ }(x)\mathrm{\dot{h}}%
+D_{\circ }^{\circ }(x)\mathrm{\dot{h}}(x)\right\rangle \right] \mathrm{d}x,
\notag
\end{eqnarray}%
where $a(x)T_{t(x)}\mathrm{h}=T_{+}^{\circ }(x)\mathrm{h}+T_{\circ }^{\circ
}(x)\mathrm{\dot{h}}(x)$. This formula is valid for any non-adapted single
integrals $T_{t}=T_{0}+i_{0}^{t}(\mathbf{D})$ with square integrable values $%
T_{t}\mathrm{h}$ for all $\mathrm{h}\in \mathrm{G}^{+}$ if we understand by $%
a(x)$ the annihilation operator $[a(x)T_{t(x)}\mathrm{h}](\varkappa
)=[T_{t(x)}\mathrm{h}](x\cup \varkappa )$ at the point $x\in X$.
\end{remark}

Indeed, taking into account that 
\begin{equation*}
\left\langle \mathrm{f}\mid i_{0}^{t}(\mathbf{D})\mathrm{h}\right\rangle
=\int_{X^{t}}[\left\langle \mathrm{f}\mid D_{+}^{-}(x)\mathrm{h}+D_{\circ
}^{-}\mathrm{\dot{h}}(x)\right\rangle +\left\langle \mathrm{\dot{f}}(x)\mid
D_{+}^{\circ }(x)\mathrm{h}+D_{\circ }^{\circ }(x)\mathrm{\dot{h}}%
(x)\right\rangle ]\mathrm{d}x,
\end{equation*}%
we readily obtain the weak form of the non-adapted It\^{o} formula if we
substitute $\mathbf{D}^{\dagger }\mathbf{T}+\mathbf{D}^{\dagger }\mathbf{D}+%
\mathbf{T}^{\dagger }\mathbf{D}$ in place of $\mathbf{D}$. This formula can
also be obtained by a direct computation 
\begin{equation*}
\left\Vert i_{0}^{t}(\mathbf{D})\mathrm{h}\right\Vert ^{2}+2\mathbf{\func{Re}%
}\left\langle i_{0}^{t}(\mathbf{D})\mathrm{h}\mid T_{0}\mathrm{h}%
\right\rangle =\left\Vert T_{t}\mathrm{h}\right\Vert ^{2}-\left\Vert T_{0}%
\mathrm{h}\right\Vert ^{2}
\end{equation*}%
without assuming that the family $T_{t}$ is defined by the kernels \textup{(%
\ref{two7})} which represent it in the form of the multiple stochastic
integral \textup{(\ref{2onee})} of $B=\epsilon (M)$. For we compute the
square of the norm of the full single integral 
\begin{eqnarray*}
\lbrack i_{0}^{t}(\mathbf{D})\mathrm{h}](\varkappa )
&=&\int_{X^{t}}[D_{+}^{-}(x)\mathrm{h}+D_{\circ }^{-}(x)\mathrm{\dot{h}}%
(x)](\varkappa )\mathrm{d}x \\
&&+\sum_{x\in \varkappa ^{t}}[D_{+}^{\circ }(x)\mathrm{h}+D_{\circ }^{\circ
}(x)\mathrm{\dot{h}}(x)](\varkappa \setminus x),
\end{eqnarray*}%
and we obtain $\Vert i_{0}^{t}(\mathbf{D}\mathrm{h}\Vert ^{2}=\Vert \int
\Vert ^{2}+2\,\mathbf{\func{Re}}\left\langle \sum \mid \int \right\rangle
+\Vert \sum \Vert ^{2}$, where 
\begin{eqnarray*}
\left\Vert \int \right\Vert ^{2} &=&\int_{X^{t}}\int_{X^{t}}\left\langle
D_{+}^{-}(z)\mathrm{h}+D_{+}^{-}(z)\mathrm{\dot{h}}(z)\mid D_{+}^{-}(x)%
\mathrm{h}+D_{\circ }^{-}(x)\mathrm{\dot{h}}(x)\right\rangle \mathrm{d}x\,%
\mathrm{d}z \\
&=&\int_{X^{t}}2\mathbf{\func{Re}}\left\langle \int_{X^{t(x)}}\left[
D_{+}^{-}(z)\mathrm{h}+D_{\circ }^{-}(z)\mathrm{\dot{h}}(z)\mathrm{d}z\right]
\mid D_{+}^{-}(x)\mathrm{h}+D_{\circ }^{-}(x)\mathrm{\dot{h}}%
(x)\right\rangle \mathrm{d}x,
\end{eqnarray*}%
the cross-term can be written as $2\,\mathbf{\func{Re}}\left\langle \sum
\mid \int \right\rangle =$ 
\begin{align*}
& 2\,\mathbf{\func{Re}}\int \left\langle \sum_{z\in \varkappa
^{t}}[D_{+}^{\circ }(z)\mathrm{h}+D_{\circ }^{\circ }(z)\mathrm{\dot{h}}%
(z)](\varkappa \setminus z)\mid \int_{X^{t}}[D_{+}^{-}(x)\mathrm{h}+D_{\circ
}^{-}(x)\mathrm{\dot{h}}(x)](\varkappa )\mathrm{d}x\right\rangle \mathrm{d}%
\varkappa \\
& =\int_{X^{t}}2\,\mathbf{\func{Re}}\int \left\langle \sum_{z\in \varkappa
^{t(x)}}[D_{+}^{\circ }(z)\mathrm{h}+D_{\circ }^{\circ }(z)\mathrm{\dot{h}}%
(z)](\varkappa )\mid D_{+}^{-}(x)\mathrm{h}+D_{\circ }^{-}(x)\mathrm{\dot{h}}%
(x)\right\rangle \mathrm{d}\varkappa \,\mathrm{d}x \\
& +\int_{X^{t}}2\,\mathbf{\func{Re}}\left\langle
a(x)\int_{X^{t(x)}}[D_{+}^{-}(z)\mathrm{h}+D_{\circ }^{-}(z)\mathrm{\dot{h}}%
(z)]\mathrm{d}s\mid D_{+}^{\circ }(x)\mathrm{h}+D_{\circ }^{\circ }(x)%
\mathrm{\dot{h}}(x)\right\rangle \mathrm{d}x
\end{align*}%
and $\left\Vert \sum \right\Vert ^{2}-\int \sum_{x\in \varkappa
^{t}}\left\Vert [D_{+}^{\circ }(x)\mathrm{h}+D_{\circ }^{\circ }(x)\mathrm{%
\dot{h}}(x)](\varkappa \setminus x)\right\Vert ^{2}\mathrm{d}\varkappa =$ 
\begin{align*}
& =\left\Vert \sum \right\Vert ^{2}-\int \sum_{x\in \varkappa
^{t}}\left\Vert [D_{+}^{\circ }(x)\mathrm{h}+D_{\circ }^{\circ }(x)\mathrm{%
\dot{h}}(x)](\varkappa \setminus x)\right\Vert ^{2}\mathrm{d}\varkappa \\
& =\int \sum_{x,z\in \varkappa ^{t}}^{x\neq z}\left\langle [D_{+}^{\circ }(z)%
\mathrm{h}+D_{\circ }^{\circ }(z)\mathrm{\dot{h}}(z)](\varkappa \setminus
z)\mid \lbrack D_{+}^{\circ }(x)\mathrm{h}+D_{\circ }^{\circ }(x)\mathrm{%
\dot{h}}(x)](\varkappa \setminus x)\right\rangle \mathrm{d}\varkappa \\
& =\int_{X^{t}}2\func{Re}\int \left\langle a(x)\sum_{z\in \varkappa ^{t(x)}}%
\mathrm{f}\left( z,\varkappa \backslash z\right) \mid \mathrm{f}\left(
x,\varkappa \right) \right\rangle \mathrm{d}\varkappa \,\mathrm{d}x,
\end{align*}%
where $\mathrm{f}\left( x,\varkappa \right) =[D_{+}^{\circ }(x)\mathrm{h}%
+D_{\circ }^{\circ }(x)\mathrm{\dot{h}}(x)]\left( \varkappa \right) $. Here
we have used \textup{(\ref{2oned})} in the form 
\begin{equation*}
\int \sum_{x\in \varkappa ^{t}}\left\langle \mathrm{f}(x,\varkappa )\mid 
\mathrm{h}(x,\varkappa \setminus x)\right\rangle \mathrm{d}\varkappa
=\int_{X^{t}}\int \left\langle \mathrm{f}(x,\varkappa \sqcup x)\mid \mathrm{h%
}(x,\varkappa )\right\rangle \mathrm{d}\varkappa \,\mathrm{d}x,
\end{equation*}%
which gives the It\^{o} term of the Hudson-Parthasarathy formula for the
adapted integrals of the form 
\begin{equation*}
\int \sum_{x\in \varkappa ^{t}}\left\Vert [D_{+}^{\circ }(x)\mathrm{h}%
+D_{\circ }^{\circ }(x)\mathrm{\dot{h}}(x)](\varkappa \setminus
x)\right\Vert ^{2}\mathrm{d}\varkappa =\int_{X^{t}}\left\Vert D_{+}^{\circ
}(x)\mathrm{h}+D_{\circ }^{\circ }(x)\mathrm{\dot{h}}(x)\right\Vert ^{2}%
\mathrm{d}x,
\end{equation*}%
and $[a(x)\mathrm{f}(x)](\varkappa )=\mathrm{f}(x,\varkappa \sqcup x)$ is
the annihilation operator at $x\in X$. Adding up all three integrals, we
obtain 
\begin{align*}
\lefteqn{\left\Vert i_{0}^{t}(\mathbf{D})\mathrm{h}\right\Vert
^{2}=\int_{X^{t}}2\func{Re}\left\langle (i_{0}^{t(x)}(\mathbf{D})\mathrm{h}%
\mid D_{+}^{-}(x)\mathrm{h}+D_{\circ }^{-}(x)\mathrm{\dot{h}}%
(x)\right\rangle \mathrm{d}x} \\
& +\int_{X^{t}}\left\Vert D_{+}^{\circ }(x)\mathrm{h}(x)+D_{\circ }^{\circ
}(x)\mathrm{\dot{h}}(x)\right\Vert ^{2}\mathrm{d}x \\
& +\int_{X^{t}}2\,\func{Re}\left\langle a(x)i_{0}^{t(x)}(\mathbf{D})\mathrm{h%
}\mid D_{+}^{\circ }(x)\mathrm{h}(x)+D_{\circ }^{\circ }(x)\mathrm{\dot{h}}%
(x)\right\rangle ^{2}\mathrm{d}x,
\end{align*}%
which leads to the weak form \textup{(\ref{two10})} of the non-adapted
generalization of the quantum It\^{o} formula for $T_{t}=T_{0}+i_{0}^{t}(%
\mathbf{D})$.

If $T_{t}=\epsilon (K_{t})$ is the representation \textup{(\ref{two2})} of
the kernel \textup{(\ref{two6})}, then obviously 
\begin{equation*}
\lbrack \epsilon (K_{t})\mathrm{h}](\varkappa \sqcup x)=[\epsilon (\dot{K}%
_{t}(x_{+}^{\circ }))\mathrm{h}+\epsilon (\dot{K}(x_{\circ }^{\circ }))%
\mathrm{\dot{h}}(x)](\varkappa ),
\end{equation*}%
and therefore $a(x)T_{t(x)}\mathrm{h}=T_{+}^{\circ }(x)\mathrm{h}+T_{\circ
}^{\circ }\mathrm{\dot{h}}(x)$.

In particular, in the scalar case $\mathrm{K}_{x}=\mathbb{C}$ for $%
D_{+}^{-}=0=D_{\circ }^{\circ }$, $D_{\circ }^{-}(x)=D(x)=D_{+}^{\circ }(x)$
and $T_{\circ }^{\circ }(x)=T_{t(x)}$, $T_{\circ }^{-}(x)=T_{+}^{\circ
}(x)\equiv \partial \left( x\right) T$ we obtain 
\begin{multline*}
\left\Vert T_{t}\mathrm{h}\right\Vert ^{2}-\left\Vert T_{0}\mathrm{h}%
\right\Vert ^{2}=\int_{X^{t}}2\func{Re}\left\langle T_{t(x)}\mathrm{h}\mid 
\mathrm{d}T_{t(x)}\mathrm{h}\right\rangle \\
\quad +\int_{X^{t}}[\left\Vert D(x)\mathrm{h}\right\Vert ^{2}+2\,\func{Re}%
\left\langle \partial \left( x\right) T\mathrm{h}\mid D(x)\mathrm{h}%
\right\rangle ]\mathrm{d}x,
\end{multline*}%
where $\partial \left( x\right) T\mathrm{h}=a(x)T_{t(x)}\mathrm{h}-T_{t(x)}%
\mathrm{\dot{h}}(x)\equiv \lbrack a(x),\,T_{t(x)}]\mathrm{h}$. This gives
the It\^{o} formula for the normally-ordered non-adapted integral 
\begin{equation*}
T_{t}-T_{0}=\int_{X^{t}}(\Lambda _{\circ }^{+}(\mathrm{d}x)D(x)+D(x)\Lambda
_{-}^{\circ }(\mathrm{d}x))=\int_{X^{t}}\mathrm{d}T_{t(x)}
\end{equation*}%
with respect to the Wiener stochastic measure $w(\Delta )$, $\Delta \in 
\mathfrak{F}$, which is represented in $\mathrm{F}$ by commuting operators $%
\widehat{w}(\Delta )=\Lambda _{\circ }^{+}(\Delta )+\Lambda _{-}^{\circ
}(\Delta )$. Consider a particular case when the operators $T_{0},\,D(x)$,
and consequently $T_{t}$ are anticipating functions $T_{0}(w),\,D(x,w)$, and 
$T_{t}(w)$ of $w$, that is, $T_{0}=T_{0}(\widehat{w}),\,D(x)=D(x,\widehat{w}%
) $, and $T_{t}=T_{t}(\widehat{w})$. Then the operators $T(x)=[a(x),%
\,T_{t(x)}]=\epsilon (\dot{K}_{t(x)}(x))$ are defined by the Malliavin
derivative $\partial _{x}T_{t}(w)|_{t=t\left( x\right) }$ as the Wiener
representation of the pointwise derivative $\dot{K}_{t(x)}(x,\varkappa
)=K_{t(x)}(x\cup \varkappa )$ of operator-valued kernels in the multiple
stochastic integral $T_{t}(w)=\int K_{t}(\varkappa )w(\mathrm{d}\varkappa
)=I(K_{t})$. In this particular case \textup{(\ref{two10})} was recently
obtained by Nualart in \textup{\cite{42}}.

We note that in the adapted case we always have $T_{\circ }^{\circ
}(x)=T_{t(x)}\otimes I(x)$ and $T_{\nu }^{\mu }(x)=0$ for $\mu \neq \nu $
except, possibly, $T_{+}^{-}(x)=\epsilon (K_{+}^{-}(x))$. Hence we readily
obtain the following result.

\begin{corollary}
\label{2C 2} The quantum stochastic process $T_{t}=\epsilon (K_{t})$ is
adapted if and only if the kernel process $K_{t}$ is adapted in the sense
that 
\begin{equation*}
K_{t}(\sigma ,\upsilon ,\tau )=\int K_{t}\left( 
\begin{array}{ll}
\omega , & \tau \\ 
\sigma , & \upsilon%
\end{array}%
\right) \mathrm{d}\omega =\delta _{\emptyset }(\sigma _{\lbrack
t})I^{\otimes }(\upsilon _{\lbrack t})\delta _{\emptyset }(\tau _{\lbrack
t})\otimes K_{t}(\sigma ^{t},\upsilon ^{\tau },\tau ^{t}),
\end{equation*}%
where $\delta _{\emptyset }(\varkappa )=1$ if $\varkappa =\emptyset $, $%
\delta _{\emptyset }(\varkappa )=0$ if $\varkappa \neq \emptyset $, $%
I^{\otimes }(\varkappa )=\bigotimes_{x\in \varkappa }I(x)$, $\varkappa
^{t}=\varkappa \cap X^{t}$, $\varkappa _{\lbrack t}=\{x\in \varkappa
:t(x)\geq t\}$. The quantum-stochastic It\^{o} formula \textup{(\ref{two9})}
for such processes can be written in the strong form 
\begin{eqnarray*}
T_{t}^{\ast }T_{t}-T_{0}^{\ast }T_{0} &=&\int_{X^{t}}(T_{t(x)}^{\ast }%
\mathrm{d}T(x)+\mathrm{d}T^{\ast }(x)T_{t(x)}+\mathrm{d}T^{\ast }(x)\mathrm{d%
}T(x)) \\
&=&i_{0}^{t}(\mathbf{G}^{\dagger }\mathbf{G}-T^{\ast }T\otimes 1),
\end{eqnarray*}%
where 
\begin{gather*}
\mathrm{d}T(x)=\Lambda (\mathbf{D},\mathrm{d}x),\,\mathrm{d}T^{\ast
}(x)=\Lambda (\mathbf{D}^{\dagger },\mathrm{d}x), \\
\mathrm{d}T^{\ast }(x)\mathrm{d}T(x)=\Lambda (\mathbf{D}^{\dagger }\mathbf{D}%
,\mathrm{d}x),\quad 1(x)=\left[ 
\begin{array}{lll}
1 & 0 & 0 \\ 
0 & I(x) & 0 \\ 
0 & 0 & 1%
\end{array}%
\right] ,
\end{gather*}%
and in the weak form as \textup{(\ref{two10})}, where $a(x)T_{t(x)}\mathrm{h}%
=[T_{t(x)}\otimes I(x)]\mathrm{\dot{h}}(x)$,
\end{corollary}

\section{Non-stationary quantum evolutions and chronological products}

We have proved the continuity of the $\ast $-representation $\epsilon $ of
an inductive $\star $-algebra $\mathfrak{B}$ of relatively bounded
operator-valued kernels $K(\mathbf{\omega })$ in the operator $\ast $%
-algebra $\mathfrak{B}(\mathrm{G}^{+})$ of the inductive limit $\mathrm{G}%
^{+}=\cap \mathrm{G}(p)$, and this property allows us to construct the
quantum-stochastic functional calculus. Namely, if $K=f(Q_{1},\ldots ,Q_{m})$
is an analytic function of the kernels $Q_{i}\in \mathfrak{B}$ as the limit
of polynomials $K_{n}$ with fixed ordering of non-commuting $Q_{1}\ldots
,Q_{n}$, the limit taken in the sense of $\Vert K_{n}-K\Vert _{\mathbf{%
\alpha }}\rightarrow 0$ for $(p,q)$-admissible quadruple $\mathbf{\alpha }%
=(\alpha _{\nu }^{\mu })$ of positive functions $\alpha _{\nu }^{\mu }(x)>0$%
, then $T=\epsilon (K)$ is an ordered function $f(Z_{1},\ldots ,Z_{m})$ of
operators $Z_{i}=\epsilon (Q_{i})$ as the limit $\Vert T_{n}-T\Vert
_{q}\rightarrow 0$ for $q\geq p+1/r$ of the corresponding polynomials $%
T_{n}=\epsilon (K_{n})$. The function $T^{\ast }=f^{\star }(Z_{1}^{\ast
},\ldots ,Z_{m}^{\ast })$ with transposed order of action of the operators $%
Z_{i}^{\ast }=\epsilon (Q_{i}^{\star })$ is also defined as a $q$-bounded
operator $T^{\ast }=\epsilon (K^{\star })$ in the scale $\{\mathrm{G}(p)\}$
for $K^{\star }=f^{\star }(Q_{1}^{\star },\ldots ,Q_{m}^{\star })$.

The differential form of this unified QS calculus is given by the
non-commutative and non-adapted generalization of the function It\^{o}
formula 
\begin{equation}
\mathrm{d}Z_{t}=\mathrm{d}i_{0}^{t}(\mathbf{A})\;\;\;\Rightarrow \;\;\;%
\mathrm{d}f(Z_{t})=\mathrm{d}i_{0}^{t}(f(\mathbf{Z}+\mathbf{A})-f(\mathbf{Z}%
)),  \label{three1}
\end{equation}%
defined for any analytic function $T_{t}=f(Z_{t})$ of an operator-valued
quantum stochastic curve $Z_{t}=\epsilon (Q_{t})$ as the generalized
QS-differential of $\epsilon (K_{t})$ for $K_{t}=f(Q_{t})$ as soon as this
function is well-defined also on the germs $\mathbf{Y}\left( x\right) =$%
\textbf{$Z$}$\left( x\right) +\mathbf{A}\left( x\right) $, \textbf{$Z$}$%
\left( x\right) =\left[ Z\left( \mathbf{x}_{\nu }^{\mu }\right) \right] $ of 
$Z_{t}$ as the triangular matrix-functions with the elements $Y\left( 
\mathbf{x}\right) =\mathbf{\nabla }_{\mathbf{x}}Z_{t\left( x\right) ]}$, $%
Z\left( \mathbf{x}\right) =\nabla _{\mathbf{x}}Z_{t\left( x\right) }$ for $%
\mathbf{x}\in \left\{ \mathbf{x}_{\nu }^{\mu }\right\} $. Here 
\begin{equation*}
T_{\nu }^{\mu }(x)=f(\mathbf{Z})_{\nu }^{\mu }(x),\quad G_{\nu }^{\mu }(x)=f(%
\mathbf{Z}+\mathbf{A})_{\nu }^{\mu }(x),
\end{equation*}%
where $f(\mathbf{Z})(x)=f(\mathbf{Z}(x))$ is a triangular matrix which as an
analytic function of the triangular matrix 
\begin{equation*}
\mathbf{Z}(x)=\left[ \epsilon \left( \dot{K}\left( \mathbf{x}_{\nu }^{\mu
}\right) \right) \right] =\epsilon \left( \mathbf{K}\left( x\right) \right)
,\;\;\dot{K}\left( \mathbf{x}_{\nu }^{\mu },\mathbf{\upsilon }\right)
=K\left( \mathbf{\upsilon }\sqcup \mathbf{x}_{\nu }^{\mu }\right)
\end{equation*}%
with the elements representing $\dot{Q}_{t(x)}(\mathbf{x})$ and $\dot{Q}%
_{t(x)]}(\mathbf{x})$, respectively, as 
\begin{eqnarray*}
\mathbf{Z}(x) &=&\left[ \epsilon (\dot{Q}_{t(x)}(\mathbf{x}_{\nu }^{\mu }))%
\right] ,\;\;\mathbf{Y\left( x\right) =}\left[ Z_{\nu }^{\mu }(x)+A_{\nu
}^{\mu }(x)\right] ,\quad \\
A_{\nu }^{\mu }(x) &=&\epsilon \left[ \dot{Q}_{t(x)]}(\mathbf{x}_{\nu }^{\mu
})-\dot{Q}_{t(x)}(\mathbf{x}_{\nu }^{\mu })\right] .
\end{eqnarray*}%
For an ordered function $T_{t}=f(Z_{1t},\ldots ,Z_{mt})$ this can be written
in terms of $Z_{it}$, with the differential $\mathrm{d}Z_{it}=\mathrm{d}%
i_{0}^{t}(A_{i})$, and $\mathbf{Y}_{i}=\mathbf{Z}_{i}+\mathbf{A}_{i}$ as 
\begin{equation*}
\mathrm{d}T_{t}=\mathrm{d}i_{0}^{t}(f(\mathbf{Y}_{1},\ldots ,\mathbf{Y}%
_{m})-f(\mathbf{Z}_{1},\ldots ,\mathbf{Z}_{m})).
\end{equation*}%
In particular, if all triangular operator-matrices $\{\mathbf{Y}_{i},\mathbf{%
Z}_{i}\}$ commute, then we can obtain the exponential function $T_{t}=\exp
\{Z_{t}\}$ for $Z_{t}=\sum_{i=1}^{m}Z_{it}$ as a solution of the following
quantum-stochastic non-adapted differential equation: 
\begin{equation}
\mathrm{d}T_{t}=\mathrm{d}i_{0}^{t}[\mathbf{T}(\mathbf{S}-\widehat{\mathbf{1}%
})],\,T_{0}=I,  \label{three2}
\end{equation}%
where $\mathbf{S}=\exp \big\{\sum_{i=1}^{m}\mathbf{A}_{i}\big\}$, $\mathbf{T}%
=\exp \big\{\sum_{i=1}^{m}\mathbf{Z}_{i}\big\}$ corresponding to%
\begin{equation*}
\mathbf{TS}=\exp \big\{\sum_{i=1}^{m}\mathbf{Y}_{i}\big\}.
\end{equation*}

We shall now deal with the problem of solving a general quantum-stochastic
equation 
\begin{equation}
T_{t}=T_{0}^{t}+i_{0}^{t}(\mathbf{TA}^{t})  \label{three3}
\end{equation}%
of type (\ref{three2}) corresponding to the integral equation $%
T_{t}=I+i_{0}^{t}(\mathbf{TA})$ with $T_{0}^{t}=I$ and $\mathbf{A}^{t}(x)=%
\mathbf{S}(x)-\hat{1}(x)$ independent of $t$. In general $T_{0}^{t}$ is
given as a nonadapted function of $t\in \mathbb{R}_{+}$ with values in
continuous operators $\mathrm{G}^{+}\rightarrow \mathrm{G}_{-}$, and $%
\mathbf{A}^{t}(x)=[A^{t}(x)_{\nu }^{\mu }]$ is a triangular matrix-function
of $x\in X$, where $A^{t}(x)_{\nu }^{\mu }=0$ for $\mu =+$ or $\nu =-$ and
the non-zero values are continuous operators 
\begin{align*}
& A_{+}^{-}(x):\mathrm{G}^{+}\rightarrow \mathrm{G}_{-},\quad A_{\circ
}^{\circ }(x):\mathrm{G}^{+}\otimes \mathrm{K}_{x}\rightarrow \mathrm{G}%
_{-}\otimes \mathrm{K}_{x}, \\
& A_{+}^{\circ }(x):\mathrm{G}^{+}\rightarrow \mathrm{G}_{-}\otimes \mathrm{K%
}_{x},\quad A_{\circ }^{-}(x):\mathrm{G}^{+}\otimes \mathrm{K}%
_{x}\rightarrow \mathrm{G}_{-},
\end{align*}%
for example $T_{0}^{t}=T_{0}U_{0}^{t}$, $\mathbf{A}^{t}(x)=\mathbf{A}%
(x)(U_{t(x)}^{t}\otimes \mathbf{I}(x))$, where $\{U_{s}^{t}:t>s\in \mathbb{R}%
_{+}\}$ is a given two-parameter family of evolution operators on $\mathrm{G}%
^{+}$. First of all we prove the following lemma.

\begin{lemma}
\label{2L 3} Suppose that the operator-functions 
\begin{equation*}
T_{0}^{t}=\epsilon (K_{0}^{t}),\quad A^{t}(x)_{\nu }^{\mu }=\epsilon (L^{t}(%
\mathbf{x}_{\nu }^{\mu }))_{\nu =\circ ,+}^{\mu =-,\circ }
\end{equation*}%
are the representations \textup{(\ref{two2})} of the kernel functions $%
K_{0}^{t}($\textbf{$\omega $}$)$, $L^{t}(\mathbf{x}_{\nu }^{\mu },\mathbf{%
\upsilon })$, where $\mathbf{\omega }=(\omega _{\nu }^{\mu }),\omega _{\nu
}^{\mu }\in \mathcal{X},\mathbf{\upsilon }=(\upsilon _{\nu }^{\mu
}),\upsilon _{\nu }^{\mu }\in \mathcal{X}$, and $\mathbf{x}_{\nu }^{\mu }$
are the atomic tables \textup{(\ref{2oneg})}. Then the integral equation 
\textup{(\ref{three3})} is the operator representation $T_{t}=\epsilon
(K_{t})$ of a triangular system of recurrence equations 
\begin{equation}
K_{t}(\mathbf{\omega })=K_{0}^{t}(\mathbf{\omega })+\sum_{x\in \mathbf{%
\omega }^{t}}[K_{t(x)}\cdot L_{x}^{t}](\mathbf{\omega }),  \label{three4}
\end{equation}%
where the kernel-operators $L_{x}^{t}(\mathbf{\omega })$ are defined almost
everywhere (for pairwise disjoint $\left( \omega _{\nu }^{\mu }\right) _{\nu
=+,\circ }^{\mu =-,\circ }$) as%
\begin{equation*}
L_{x}^{t}(\mathbf{\upsilon }\sqcup \mathbf{x}_{\nu }^{\mu }):=L^{t}\left( 
\mathbf{x}_{\nu }^{\mu },\mathbf{\upsilon }\right) \equiv L^{t}\left( x,%
\mathbf{\upsilon }\right) _{\nu }^{\mu }
\end{equation*}%
by the matrix elements of $\mathbf{L}^{t}(x,\mathbf{\upsilon })$, with $%
L_{x}^{t}(\mathbf{\omega })=0$ if $x\notin \sqcup \,\omega _{\nu }^{\mu }$,
and $K_{t(x)}\cdot L_{x}^{t}$ is the kernel product. The solution of \textup{%
(\ref{three4})} is uniquely defined almost everywhere (if $t(x)\neq
t(x^{\prime })$ for all $x\neq x^{\prime }\in \sqcup \omega _{\nu }^{\mu }$)
as the sum 
\begin{equation*}
K_{t}(\mathbf{\omega })=\sum_{\mathbf{\vartheta }\subseteq \mathbf{\omega }%
^{t}}M_{t}(\mathbf{\vartheta },\mathbf{\omega }\setminus \mathbf{\vartheta }%
)=\nu _{0}^{t}(\mathbf{\omega },M_{t})
\end{equation*}%
of chronological kernel products 
\begin{equation}
M_{t}(\mathbf{\vartheta },\mathbf{\upsilon })=[K_{0}^{t(x_{1})}\cdot
L_{x_{1}}^{t(x_{2})}\cdot \ldots \cdot L_{x_{m-1}}^{t(x_{m})}\cdot
L_{x_{m}}^{t}](\mathbf{\vartheta }\sqcup \mathbf{\upsilon })  \label{three5}
\end{equation}%
over all decompositions \textbf{$\vartheta $}$=\mathbf{x}_{1}\sqcup \ldots
\sqcup \mathbf{x}_{m}$ of the tables \textbf{$\vartheta $}$=(\vartheta _{\nu
}^{\mu })$ into atomic tables $\mathbf{x}_{i}$, each of the form \textup{(%
\ref{2oneg})}, with the correspondence $x_{i}\in \vartheta _{\nu }^{\mu
}\Leftrightarrow \mathbf{x}_{i}=\mathbf{x}_{i\nu }^{\mu }$. It gives the
unique solution \textup{(\ref{three3})} in the form of the generalized
multiple integral 
\begin{equation*}
T_{t}=\iota _{0}^{t}(B_{t}),\,\;\;\;\;B_{t}(\mathbf{\vartheta })=\epsilon
(M_{t}(\mathbf{\vartheta })),
\end{equation*}%
if Fock representation $B_{t}(\mathbf{\vartheta })$ of the products \textup{(%
\ref{three5})} satisfies the condition $\Vert B_{t}\Vert _{p}^{s}(r)<\infty $
for some admissible $p\in \mathcal{P}_{1}$ and $r^{-1}$, $s^{-1}\in \mathcal{%
P}_{0}$.
\end{lemma}

\begin{proof}
We substitute $T_{0}^{t}=\epsilon (K_{0}^{t})$, $\mathbf{A}^{t}(x)=\epsilon (%
\mathbf{L}^{t}(x))$, and $T_{t}=\epsilon (K_{t})$ in (\ref{three3}) and we
take into account the fact that%
\begin{equation*}
\mathbf{T}(x)\mathbf{A}^{t}(x)=\epsilon (\mathbf{K}_{t(x)}(x)\cdot \mathbf{L}%
^{t}(x)),
\end{equation*}%
where $\mathbf{K}_{t}(x)=[K_{t}(x)_{\nu }^{\mu }]$ is a triangular matrix, $%
K_{t}(x)_{\nu }^{\mu }=0$ if $\mu >\nu $, with the non-zero kernel entries 
\begin{equation*}
K_{t}(x,\mathbf{\upsilon })_{-}^{-}=K_{t}(\mathbf{\upsilon })=K(x,\mathbf{%
\upsilon })_{+}^{+},\;\;K_{t}(x,\mathbf{\upsilon })_{\nu }^{\mu }=K_{t}(%
\mathbf{\upsilon }\sqcup \mathbf{x}_{\nu }^{\mu }),
\end{equation*}%
and $\mathbf{L}^{t}(x)=[L^{t}(x)_{\nu }^{\mu }]$, where $L^{t}(x)_{\nu
}^{\mu }=0$ if $\mu >\nu $ and the entries 
\begin{equation*}
L^{t}(x,\mathbf{\upsilon })_{-}^{-}=0=L^{t}(x,\mathbf{\upsilon }%
)_{+}^{+},\;\;\;x\notin \sqcup \upsilon _{\nu }^{\mu },\;\;L^{t}(x)_{\nu
}^{\mu }=\dot{L}_{x}^{t}(\mathbf{x}_{\nu }^{\mu })
\end{equation*}%
are defined by the kernels $L_{x}^{t}(\mathbf{\omega })=L^{t}(\mathbf{x},%
\mathbf{\omega }\setminus \mathbf{x})$, with $L_{x}^{t}(\mathbf{\omega })=0$
if $x\notin \omega _{\nu }^{\mu }$ for all $\mu \neq +,\nu \neq -$, in the
same way as the entries $K_{t}(x,\mathbf{\upsilon })_{\nu }^{\mu }$ are
defined by the kernels $K_{t}(\mathbf{\omega })$. As a result we found that (%
\ref{three3}) is satisfied if 
\begin{multline*}
K_{t}(\mathbf{\omega })=K_{0}^{t}(\mathbf{\omega })+\sum_{\mathbf{x}\in 
\mathbf{\omega }^{t}}[\mathbf{K}_{t(x)}(x)\mathbf{L}^{t}(x)]_{\nu (\mathbf{x}%
)}^{\mu (\mathbf{x})}(\mathbf{\omega }\setminus \mathbf{x}) \\
=K_{0}^{t}(\mathbf{\omega })+\sum_{\mu <+}^{\nu >-}\,\,\sum_{x\in \omega
_{\nu }^{\mu }}^{t(x)<t}[K_{t(x)}\cdot L_{x}^{t}](\mathbf{x}_{\nu }^{\mu
}\sqcup \mathbf{\omega }\setminus \mathbf{x}_{\nu }^{\mu }),
\end{multline*}%
which corresponds to (\ref{three4}). The solution of this equation for any
table $\mathbf{\omega }=(\omega _{\nu }^{\mu })_{\nu =\circ ,+}^{\mu
=-,\circ }$ with chronologically ordered entries is represented as the sum (%
\ref{two6}) of the chronological products (\ref{three5}) of the
operator-valued kernels $M_{t}(\mathbf{\emptyset },\mathbf{\omega }%
)=K_{0}^{t}(\mathbf{\omega })$ and $L_{x}^{t}(\mathbf{\omega })$, since 
\begin{align*}
K_{t}(\mathbf{\omega })& =\sum_{\mathbf{\vartheta }\subseteq \mathbf{\omega }%
^{t}}M_{t}(\mathbf{\vartheta },\mathbf{\omega }\setminus \mathbf{\vartheta }%
)=M_{t}(\emptyset ,\mathbf{\omega })+\sum_{|\mathbf{\vartheta }|\geq 1}^{%
\mathbf{\vartheta }\subseteq \mathbf{\omega }^{t}}M_{t}(\mathbf{\vartheta },%
\mathbf{\omega }\setminus \mathbf{\vartheta }) \\
& =M_{t}(\emptyset ,\mathbf{\omega })+\sum_{\mathbf{x}\in \mathbf{\omega }%
^{t}}\sum_{\mathbf{\vartheta }\in \mathbf{\omega }^{t(x)}}M_{t}(\mathbf{%
\vartheta }\sqcup \mathbf{x},\mathbf{\omega }\setminus (\mathbf{\vartheta }%
\sqcup \mathbf{x})) \\
& =K_{0}^{t}(\mathbf{\omega })+\sum_{\mathbf{x}\in \mathbf{\omega }%
^{t}}\sum_{\mathbf{\vartheta }\subseteq \mathbf{\omega }^{t(x)}}[M_{t(x)}%
\cdot L_{x}^{t}](\mathbf{\omega })=K_{0}^{t}(\mathbf{\omega })+\sum_{\mathbf{%
x}\in \mathbf{\omega }^{t}}[K_{t(x)}\cdot L_{x}^{t}](\mathbf{\omega }),
\end{align*}%
where we have used the representation (\ref{three5}) in the recurrent form 
\begin{equation*}
M_{t}(\mathbf{\vartheta }\sqcup \mathbf{x},\mathbf{\upsilon })=[M_{t(x)}(%
\mathbf{x})\cdot L^{t}(\mathbf{x})]_{\nu (\mathbf{x})}^{\mu (\mathbf{x})}(%
\mathbf{\vartheta }\sqcup \mathbf{\upsilon })=[M_{t(x)}\cdot L_{x}^{t}](%
\mathbf{x}\sqcup \mathbf{\vartheta }\sqcup \mathbf{\upsilon }).
\end{equation*}%
This defines the representation of the solution $T_{t}=\epsilon (K_{t})$ in
the form of the non-adapted quantum-stochastic integral (\ref{2onee}) of $%
B_{t}=\epsilon (M_{t})$, since by Lemma~\ref{2L 2}, $\epsilon \circ \nu
_{0}^{t}=\iota _{0}^{t}\circ \epsilon $ if the integrability condition (\ref%
{2onef}) is fulfilled.
\end{proof}

\begin{theorem}
\label{2T 3} Suppose that $U_{s}^{t}=\epsilon (V_{s}^{t})$ is the
representation on $\mathrm{G}$ of the evolution family $\{V_{s}^{t}:t\geq
s\in \mathbb{R}_{+}\}$ of relatively bounded operator-valued kernels 
\begin{equation*}
V_{s}^{t}%
\begin{pmatrix}
\omega _{+}^{-}, & \omega _{\circ }^{-} \\ 
\omega _{+}^{\circ }, & \omega _{\circ }^{\circ }%
\end{pmatrix}%
:\mathrm{H}\otimes \mathrm{K}^{\otimes }(\omega _{\circ }^{-})\otimes 
\mathrm{K}^{\otimes }(\omega _{\circ }^{\circ })\rightarrow \mathrm{H}%
\otimes \mathrm{K}^{\otimes }(\omega _{\circ }^{\circ })\otimes \mathrm{K}%
^{\otimes }(\omega _{+}^{\circ }),
\end{equation*}%
satisfying the condition $V_{r}^{s}\cdot V_{s}^{t}=V_{r}^{t}$ for all $r<s<t$%
, where the representation is considered with respect to the kernel product 
\textup{(\ref{two9})} of \textup{Chapter I} with the unit $V_{t}^{t}(\mathbf{%
\omega })=I\otimes \mathbf{I}^{\otimes }(\mathbf{\omega })$. Suppose that 
\begin{equation*}
K_{0}^{t}(\mathbf{\omega })=[K_{0}^{s}\cdot V_{s}^{t}](\mathbf{\omega }%
),\,\;\;\;\;L_{x}^{t}(\mathbf{\omega })=[L_{x}^{s}\cdot V_{s}^{t}](\mathbf{%
\omega }),\quad \forall t>s
\end{equation*}%
are the kernel products defining the representation \textup{(\ref{three4})}
of equation \textup{(\ref{three3})}. Then the kernels chronological product 
\begin{equation}
K_{t}(\mathbf{\omega })=[K_{0}^{t(x_{1})}\cdot F_{x_{1}}^{t(x_{2})}\cdot
\ldots \cdot F_{x_{n-1}}^{t(x_{n})}\cdot F_{t(x_{n})}^{t}](\mathbf{\omega })
\label{three6}
\end{equation}%
for $F_{x}^{t}(\mathbf{\omega })=L_{x}^{t}(\mathbf{\omega })+V_{t(x)}^{t}(%
\mathbf{\omega })$, $\mathbf{\omega }^{t}=\mathbf{x}_{1}\sqcup \ldots \sqcup 
\mathbf{x}_{n}$, is a unique solution of the system \textup{(\ref{three4})}
for almost all $\mathbf{\omega }=(\omega _{\nu }^{\mu })$ \textup{(}if $%
t(x)\neq t(x^{\prime })$ for all $x\neq x^{\prime }\in \sqcup \omega _{\nu
}^{\mu }$\textup{)}. This yields the representation of the solution of 
\textup{(\ref{three3})} in the form $T_{t}=\epsilon (K_{t})$ defined on $%
\mathrm{G}$ for each $t$ as a relatively bounded operator if the product 
\textup{(\ref{three6})} satisfies the condition $\Vert K_{t}\Vert _{\mathbf{%
\alpha }}<\infty $ with respect to the norm \textup{(\ref{two4})} for the
quadruple $\mathbf{\alpha }=(\alpha _{\nu }^{\mu })$ of functions admissible
in the sense of \textup{(\ref{two5})} and equal to zero for $t(x)>t$. The
operators $T_{t}$ are isometric, that is, $T_{t}^{\ast }T_{t}=\hat{I}$ 
\textup{(}unitary: $T_{t}^{\ast }=T_{t}^{-1}$\textup{)} if and only if the
operators $T_{0}$ and $U_{s}^{t},t\geq s\geq 0$, are isometric \textup{(}%
unitary\textup{)}. Consequently, for all $t$ we have $%
T_{0}^{t}=T_{0}U_{0}^{t}$ and the triangular operator-matrices $\mathbf{S}%
(x)=[S_{\nu }^{\mu }(x)]$ that define the generators of the equation \textup{%
(\ref{three3})} in the form $\mathbf{A}^{t}(x)=(\mathbf{S}(x)-\hat{\mathbf{1}%
})(U_{t(x)}^{t}\otimes \mathbf{1}(x))$ are pseudo-isometric, that is, $%
\mathbf{S}^{\dagger }(x)\mathbf{S}(x)=\hat{I}\otimes \mathbf{1}(x)$ \textup{(%
}pseudo-unitary: $\mathbf{S}^{\dagger }(x)=\mathbf{S}(x)^{-1}$\textup{)},
and such that 
\begin{align}
& S_{\circ }^{\circ }(x)^{\ast }S_{\circ }^{\circ }(x)=\hat{I}\otimes
I(x),\;\,S_{+}^{-}(x)^{\ast }+S_{+}^{\circ }(x)^{\ast }S_{+}^{\circ
}(x)+S_{+}^{-}(x)=0,  \notag  \label{three7} \\
& S_{\circ }^{-}(x)^{\ast }+S_{\circ }^{\circ }(x)^{\ast }S_{+}^{\circ
}(x)=0,\;\;\;\;\;\;\,S_{+}^{\circ }(x)^{\ast }S_{\circ }^{\circ
}(x)+S_{\circ }^{-}(x)=0
\end{align}%
\textup{(}and $S_{\circ }^{\circ }(x)$ are unitary, that is, $S_{\circ
}^{\circ }(x)^{\ast }=S_{\circ }^{\circ }(x)^{-1}$\textup{)} for almost all $%
x\in X^{t}$.
\end{theorem}

\begin{proof}
Suppose that $\mathbf{\upsilon }=\mathbf{\upsilon }_{0}\sqcup \mathbf{%
\upsilon }_{1}\sqcup \ldots \sqcup \mathbf{\upsilon }_{m}$ is a
decomposition of the table $\mathbf{\upsilon }=(\upsilon _{\nu }^{\mu })=%
\mathbf{\omega }\setminus $\textbf{$\vartheta $} into the subtables $\mathbf{%
\upsilon }_{i}=\mathbf{x}_{i}^{1}\sqcup \ldots \sqcup \mathbf{x}_{i}^{n_{i}}$
determined by the points $x_{i}\in X^{t}$ of the atomic tables $\mathbf{x}%
_{i}$ in the chronological decomposition \textbf{$\vartheta $}$=\mathbf{x}%
_{1}\sqcup \ldots \sqcup \mathbf{x}_{m}$, so that $t(x_{i})<t(x_{i}^{1})<%
\cdots <t(x_{i}^{n_{i}})<t(x_{i+1}),t(x_{0})=0$. Then 
\begin{align*}
K_{0}^{t(x_{1})}& =K_{0}^{t(x_{0}^{1})}\cdot
V_{t(x_{0}^{1})}^{t(x_{0}^{2})}\ldots V_{t({x_{0}^{n_{i}}}%
)}^{t(x_{1})},L_{x_{i}}^{t(x_{i+1})}=L_{x_{i}}^{t(x_{i}^{1})}\cdot
V_{t(x_{i}^{1})}^{t(x_{i}^{2})}\ldots V_{t(x_{i}^{n}i)}^{t(x_{i+1})}, \\
K_{t}(\mathbf{\omega })& =\sum_{\mathbf{\vartheta }\subseteq \mathbf{\omega }%
^{t}}[K_{0}^{t(x_{1})}\cdot L_{x_{1}}^{t(x_{2})}\ldots
L_{m-1}^{t(x_{m})}\cdot L_{x_{m}}^{t}](\mathbf{\vartheta }\sqcup \mathbf{%
\upsilon }) \\
& =[K_{0}^{t(z_{1})}\cdot
(V_{t(z_{1})}^{t(z_{1})}+L_{z_{1}}^{t(z_{2})})\ldots
(V_{t(z_{n})}^{t}+L_{z_{n}}^{t})](\mathbf{\omega }),
\end{align*}%
where the points $z_{1},\ldots ,z_{n}\in X^{t}$, $t(z_{1})<\ldots <t(z_{n})$
define the decomposition $\mathbf{\omega }=\sqcup \mathbf{z}_{i}$ into
atomic tables (\ref{2oneg}). Thus the chronological products (\ref{three6})
of the kernels $F_{z}^{t}=L_{z}^{t}+V_{t(z)}^{t}$ defines a unique solution
of the system (\ref{three4}), which is a pseudo-isometric (pseudo-unitary)
kernel if and only if the same is true for each factor $%
K_{0}^{t(z_{1})},F_{z_{1}}^{t(z_{2})},\ldots ,F_{z_{n}}^{t}$. If, in
addition, the kernel $K_{t}(\mathbf{\omega })$ is locally bounded for each $%
t $ relative to the quadruple $\mathbf{\alpha }=(\alpha _{\nu }^{\mu })$ of
positive functions $\alpha _{\nu }^{\mu }(x)$ locally integrable in the
sense 
\begin{equation*}
\int_{X^{t}}\alpha _{+}^{-}(x)\mathrm{d}x<\infty ,\;\;\int_{X^{t}}(\alpha
_{+}^{\circ }(x)^{2}+\alpha _{\circ }^{-}(x)^{2})r(x)\mathrm{d}x<\infty ,\;\;%
\underset{x\in X^{t}}{\mathrm{ess}\sup }\frac{\alpha _{\circ }^{\circ }(x)}{%
p(x)}<\infty ,
\end{equation*}%
then, in accordance with Theorem~\ref{2T 2}, the representation (\ref{two2})
defines the map $\epsilon :K_{t}\rightarrow T_{t}$ as $\star $-homomorphism
in the $\ast $-algebra of $q$-bounded, $q\geq p+1/r$, operators on $\mathrm{G%
}^{+}$ satisfying the exponential estimate (\ref{two6}). Moreover, $T_{t}$
is an isometry (a unitary operator) if the kernel $K_{t}$ is
pseudo-symmetric, that is, $K_{t}^{\star }\cdot K_{t}=I\otimes 1^{\otimes }$
(pseudo-unitary, that is, $K_{t}^{\star }=K_{t}^{-1}$), with respect to the
kernel product (\ref{two9}) of Chapter I and the pseudo-involution $%
K_{t}\mapsto K_{t}^{\star }$. For any chronologically ordered collection $%
\mathbf{\omega }=(\omega _{\nu }^{\mu })$ this is guaranteed by the
corresponding properties of the kernels $K_{0},V_{s}^{t},s\leq t$ and $F_{z}$
(for almost all $z\in X^{t}$) by virtue of the representation of (\ref%
{three6}) in the form of a finite product of the kernels $%
K_{0}=K_{0}^{0},V_{0}^{t(x_{1})}$, and $F_{z}=F_{z}^{t(z)},V_{t(z)}^{t},z\in
\omega ,t>t(z)$. Hence the kernel-matrices $\mathbf{F}(x)=[F_{\nu }^{\mu
}(x)]$, with the entries 
\begin{equation*}
F_{\nu }^{\mu }(x)=0,\mu >\nu ,\;\;F_{-}^{-}(x)=I=F_{+}^{+}(x),\;\;F_{\nu
}^{\mu }(x)=\dot{F}_{x}(\mathbf{x}_{\nu }^{\mu }),
\end{equation*}%
are pseudo-isometric (pseudo-unitary). This implies that the operators $%
T_{0}=\epsilon (K_{0})$ and $U_{s}^{t}=\epsilon (V_{s}^{t})$ are isometric
(unitary) and the triangular matrix $\mathbf{S}(x)=[\epsilon (F_{\nu }^{\mu
}(x))]$ (where $S_{\nu }^{\mu }(x)=0$ if $\mu >\nu $ $%
S_{-}^{-}(x)=I=S_{+}^{+}(x)$ and $S_{\nu }^{\mu }(x)=\epsilon (\dot{F}%
(x_{\nu }^{\mu }))$ if $\mu \neq +,\nu \neq -$), defining the generator $%
A(x)\equiv A^{t(x)}(x)$ as $\mathbf{S}(x)-I\otimes \mathbf{1}(x)$, is
pseudo-isometric (pseudo-unitary).

By virtue of the uniqueness of the representation $T_{0}=\epsilon
(K_{0}),U_{s}^{t}=\epsilon (V_{s}^{t})$, and $\mathbf{S}(x)=\epsilon (%
\mathbf{F}(x))$ up to the $\star $-ideal described in Section 2 of Chapter
I, the resulting conditions are necessary and sufficient for the solution $%
T_{t}=\epsilon (K_{t})$ of the non-adapted quantum-stochastic equation (\ref%
{three3}), uniquely (up to the ideal mentioned above) determined by the
pseudo-isometric (pseudo-unitary) kernels (\ref{three6}), to be isometric
(unitary). Writing the condition $\mathbf{S}^{\dagger }\mathbf{S}=\hat{I}%
\otimes \mathbf{1}$ in terms of the matrix entries $S_{\nu }^{\mu
}(x),S_{-\nu }^{\dagger \mu }=S_{-\mu }^{\nu \ast }$, we obtain the system (%
\ref{three7}): 
\begin{equation*}
\lbrack \mathbf{S}^{\dagger }\mathbf{S}](x)=%
\begin{bmatrix}
1, & S_{+}^{\circ }(x)^{\ast }, & S_{+}^{-}(x)^{\ast } \\ 
0, & S_{\circ }^{\circ }(x)^{\ast }, & S_{\circ }^{-}(x)^{\ast } \\ 
0, & 0, & 1%
\end{bmatrix}%
\begin{bmatrix}
1, & S_{\circ }^{-}(x), & S_{+}^{-}(x) \\ 
0, & S_{\circ }^{\circ }(x), & S_{+}^{\circ }(x) \\ 
0, & 0, & 1%
\end{bmatrix}%
=I\otimes 
\begin{bmatrix}
1, & 0, & 0 \\ 
0, & I(x), & 0 \\ 
0, & 0, & 1%
\end{bmatrix}%
.
\end{equation*}%
Thus Theorem~\ref{2T 3} is proved.
\end{proof}

\begin{remark}
Suppose that the evolution family $\{U_{s}^{t}\}$ is a solution of the
non-stochastic non-adapted equation 
\begin{equation}
U_{s}^{t}=\hat{I}+\int_{s\leq t(x)<t}U_{s}^{t(x)}S_{+}^{-}(x)\mathrm{d}%
x,\quad S<t,  \label{three8}
\end{equation}%
and in the dissipative case $S_{+}^{-}(x)+S_{+}^{-}(x)^{\ast }<0$ this
solution is defined as an adapted family of contractions $U_{s}^{t}:\mathrm{G%
}\rightarrow \mathrm{G},\Vert U_{s}^{t}\Vert \leq \mathbf{1}$. Then the
solution of the differential equation \textup{(\ref{three2})} can be written
in the form of a purely stochastic quantum multiple integral $T_{t}=\iota
_{0}^{t}(B^{t})$ satisfying \textup{(\ref{three3})} with $%
T_{0}^{t}=U_{0}^{t} $ and the generators $\mathbf{A}^{t}(x)=\mathbf{A}%
(x)(U_{t(x)}^{t}\otimes \mathbf{1}(x))$, where 
\begin{equation*}
A_{+}^{-}(x)=0,A_{+}^{\circ }(x)=S_{+}^{\circ }(x),A_{\circ
}^{-}(x)=S_{\circ }^{-}(x),A_{\circ }^{\circ }(x)=S_{\circ }^{\circ }(x)-%
\hat{I}\otimes I(x).
\end{equation*}%
In the case when the operator function $S_{+}^{-}(x)$ is locally absolutely
integrable in the sense that $\int_{X^{t}}\Vert S_{+}^{-}(x)\Vert \mathrm{d}%
x<\infty $ for all $t$, and if we have 
\begin{equation*}
U_{s}^{t}=\sum_{n=0}^{\infty }\quad \idotsint\limits_{s<t(x_{1})<\cdots
<t(x_{n})<t}S_{+}^{-}(x_{1})\ldots S_{+}^{-}(x_{n})\prod_{i=1}^{n}\mathrm{d}%
x_{i}=\int_{\mathcal{X}_{s}^{t}}S_{+}^{-}(\vartheta )\mathrm{d}\vartheta ,
\end{equation*}%
where $\mathcal{X}_{s}^{t}=\{\vartheta \in \mathcal{X}:\vartheta \leq
\lbrack s,t)\},S_{+}^{-}(x_{1},\ldots ,x_{n})=S_{+}^{-}(x_{1})\ldots
S_{+}^{-}(x_{n})$, then this representation can be directly obtained by the
integration with respect to $\omega _{+}^{-}\in \mathcal{X}$ of the kernel $%
K_{t}(\mathbf{\omega })=[F_{z_{1}}\ldots F_{z_{n}}](\mathbf{\omega })$,
defined for $\mathbf{\omega }^{t}=\mathbf{z}_{1}\sqcup \ldots \sqcup \mathbf{%
z}_{n}$ as the chronological product of the kernels $F_{x}(\mathbf{\omega }%
)=F_{\nu }^{\mu }(x,\mathbf{\omega }\setminus \mathbf{x}_{\nu }^{\mu })$ for
all $\mathbf{\omega }=(\omega _{\nu }^{\mu })$ if $x\in \omega _{\nu }^{\mu
} $ and $F_{x}(\mathbf{\omega })=I\otimes \mathbf{1}^{\otimes }(\mathbf{%
\omega })$ if $x\notin \sqcup \omega _{\nu }^{\mu }$, which correspond to
the representation $S_{\nu }^{\mu }(x)=\epsilon (F_{\nu }^{\mu }(x))$.
\end{remark}

For we write the solution of the equation $T_{t}=\hat{I}+i_{0}^{t}(\mathbf{T}%
(\mathbf{S}-\hat{\mathbf{I}}))$ in the form $T_{t}=\epsilon (K_{t})$, where $%
K_{t}$ is the kernel \textup{(\ref{three6})} with $K_{0}^{t}=I^{\otimes }$
and $F_{x}^{t}=F_{x}$ independent of $t$. We denote by $\{z_{1},\ldots
,z_{n}\}$ the subchain of the chain $\{x_{1},\ldots x_{m}\}$ of the
decomposition $\mathbf{\omega }^{t}=\mathbf{x}_{1}\sqcup \ldots \sqcup 
\mathbf{x}_{m}$ that corresponds to the elements $z_{i}\notin \omega
_{+}^{-} $, and we write the integral of $K_{t}(\mathbf{\omega })$ with
respect to $\omega _{+}^{-}\in \mathcal{X}$ in the form of a multiple
integral in $\vartheta _{i}\in \mathcal{X}_{t(z_{i})}^{t(z_{i+1})},i=0,1,%
\ldots ,n$, where $t(z_{0})=0,t(z_{n+1})=t$, and $z_{i}\in \mathcal{X}%
,i=1,\ldots ,n$. Then, in accordance with \textup{(\ref{two9})} of \textup{%
Chapter I} we obtain the kernel chronological product 
\begin{equation*}
K_{t}(\omega ^{\circ },\upsilon ,\omega _{\circ })=[V_{0}^{t(z_{1})}\cdot
F_{z_{1}}\cdot V_{t(z_{1})}^{t(z_{2})}\cdots F_{z_{n}}\cdot
V_{t(z_{n})}^{t}](\omega ^{\circ },\upsilon ,\omega _{\circ }).
\end{equation*}%
Here in the square brackets we have the product of integral kernels $%
F_{x}(\omega ^{\circ },\upsilon ,\omega _{\circ })=\int F_{x}%
\begin{pmatrix}
\omega & \omega _{\circ } \\ 
\omega ^{\circ } & \upsilon%
\end{pmatrix}%
\mathrm{d}\omega $ and 
\begin{equation*}
V_{s}^{t}(\omega ^{\circ },\upsilon ,\omega _{\circ })=\sum_{n=0}^{\infty
}\quad \idotsint\limits_{s\leq t(x_{1})<\cdots
<t(x_{n})<t}[F_{+}^{-}(x_{1})\ldots F_{+}^{-}(x_{n})](\omega ^{\circ
},\upsilon ,\omega _{\circ })\prod_{i=1}^{n}\mathrm{d}x_{i},
\end{equation*}%
where $[F_{+}^{-}(x)](\omega ^{\circ },\upsilon ,\omega _{\circ
})=F_{+}^{-}\left( 
\begin{array}{cc}
x, & \omega _{\circ } \\ 
\omega ^{\circ }, & \upsilon%
\end{array}%
\right) ,x\in X$. On the other hand, we can obtain the same result if we
integrate the kernel product 
\begin{equation*}
K_{t}(\mathbf{\omega })=[V_{0}^{t(z_{1})}\cdot F_{z_{1}}\cdot
V_{t(z_{1})}^{t(z_{2})}\cdots F_{z_{n}}V_{t(z_{n})}^{t}](\mathbf{\omega })
\end{equation*}%
with respect to $\omega _{+}^{-}\in \mathcal{X}$, where the kernels $%
V_{s}^{t}(\mathbf{\omega })=[F_{x_{1}}\ldots F_{x_{n}}](\mathbf{\omega })$ 
\textup{(}for $X_{s}^{t}\cap \omega _{+}^{-}=x_{1}\sqcup \ldots \sqcup x_{n}$%
\textup{)} define the representation $U_{s}^{t}=\epsilon (V_{s}^{t})$ of the
solution of \textup{(\ref{three8})} for $S_{+}^{-}(x)=\epsilon
(F_{+}^{-}(x)) $. Putting $F_{x}(\mathbf{\omega })=I\otimes 1^{\otimes }(%
\mathbf{\omega })$ for $x\in \omega _{+}^{-}\cap X^{t}$ and taking into
account the consistency condition $V_{r}^{s}\cdot V_{s}^{t}=V_{r}^{t}$, we
find the solution of \textup{(\ref{three2})} as the solution of \textup{(\ref%
{three3})} with the generators $A^{t}(x)_{\nu }^{\mu }=\epsilon (L^{t}(%
\mathbf{x}_{\nu }^{\mu })) $, where 
\begin{equation*}
L^{t}(\mathbf{x}_{\nu }^{\mu },\mathbf{\upsilon })=[(F_{x}-\mathbf{I}%
^{\otimes })\cdot V_{t(x)}^{t}](\mathbf{\upsilon }\sqcup \mathbf{x}_{\nu
}^{\mu })=0\text{ for }(\mu ,\nu )=(-,+).
\end{equation*}%
This solution can be written in the form of the quantum-stochastic multiple
non-adapted integral \textup{(\ref{2onee})} of $B_{t}(\mathbf{\vartheta }%
)=\epsilon (M_{t}(\mathbf{\vartheta }))$, where $M_{t}(\mathbf{\vartheta },%
\mathbf{\upsilon })$ is defined in \textup{(\ref{three5})} by the kernels $%
K_{0}^{t}=V_{0}^{t}$ and $L_{x}^{t}=(F_{x}-\mathbf{I}^{\otimes })\cdot
V_{t(x)}^{t}$. The operator-function $B_{t}(\mathbf{\vartheta })$ is equal
to zero if $\vartheta _{+}^{-}\neq \emptyset $, since the product \textup{(%
\ref{three5})} is zero for $x_{i}\in \vartheta _{+}^{-}$. From this we can
readily obtain the following corollary.

\begin{corollary}
\label{2C 3} Suppose that $S_{\nu }^{\mu }(x)=F_{\nu }^{\mu }(x)\otimes \hat{%
\mathbf{1}}$, where $F_{+}^{-}(x)$ are closed dissipative operators such
that there exists a consistent family $\{V_{s}^{t}\}$ of contractions in $%
\mathrm{H}$ which allows us to write the solution of \textup{(\ref{three8})}
in the form $U_{s}^{t}=V_{s}^{t}\otimes \hat{\mathbf{1}}$. \textup{(}It is
sufficient, for example, to require that $F_{+}^{-}(x)$ be locally
absolutely integrable, that is, $\int_{X_{t}}\Vert F_{+}^{-}(x)\Vert \mathrm{%
d}x<\infty $ for all $t$.\textup{)}

Suppose that the operator-functions 
\begin{equation*}
F_{+}^{\circ }(x):\mathrm{H}\rightarrow \mathrm{H}\otimes \mathrm{K}%
_{x},\quad F_{\circ }^{-}(x):\mathrm{H}\otimes \mathrm{K}_{x}\rightarrow 
\mathrm{H}
\end{equation*}%
are locally square integrable in the sense that 
\begin{equation*}
\Vert F\Vert _{t}^{(2)}(r)=\Big(\int_{X^{t}}\Vert F(x)\Vert ^{2}r(x)\mathrm{d%
}x\Big)^{1/2}<\infty
\end{equation*}%
and $\Vert F_{\circ }^{\circ }\Vert _{t,p}^{(\infty )}=\mathrm{ess}%
\sup_{x\in X^{t}}\{\Vert F_{\circ }^{\circ }(x)\Vert /p(x)\}\leq 1$ for some 
$r^{-1}\in \mathcal{P}_{0}$ and $p\in \mathcal{P}_{1}$. Then the solution $%
T_{t}=\iota _{0}^{t}(B)$, $B(\mathbf{\vartheta })=M(\mathbf{\vartheta }%
)\otimes \hat{\mathbf{1}}$ of the quantum-stochastic equation \textup{(\ref%
{three2})} is uniquely determined for each $t\geq 0$ as a relatively bounded
operator $T_{t}=\epsilon (K_{t})$ representing by means of \textup{(\ref%
{two8})} the adapted chronological products 
\begin{equation}
K_{t}(\omega ^{\circ },\upsilon ,\omega _{\circ })=V_{0}^{t(x_{1})}\odot F(%
\mathbf{x}_{1})\odot V_{t(x_{1})}^{t(x_{2})}\odot \ldots \odot F(\mathbf{x}%
_{n})\odot V_{t(x_{n})}^{t}.  \label{three9}
\end{equation}%
Here $\{x_{1},\ldots ,x_{n}\}=(\omega ^{\circ }\sqcup \upsilon \sqcup \omega
_{\circ })\cap X^{t}$ is the chronologically ordered chain $0\leq
t(x_{1})<\cdots <t(x_{n})<t$, $\,\mathbf{x}=\mathbf{x}_{+}^{\circ }$ if $%
x\in \omega ^{\circ }$, $\,\mathbf{x}=\mathbf{x}_{\circ }^{\circ }$ if $x\in
\upsilon $, $\,\mathbf{x}=\mathbf{x}_{\circ }^{-}$ if $x\in \omega _{\circ }$
are atomic tables \textup{(\ref{2oneg})}, $F(\mathbf{x}_{\nu }^{\mu
})=F_{\nu }^{\mu }(x)$ is one of the three functions $F_{+}^{\circ
},F_{\circ }^{\circ },F_{\circ }^{-}$, and $\odot $ denotes the semi-tensor
product defined recurrently by 
\begin{equation*}
K(\mathbf{\upsilon })\odot F(\mathbf{\vartheta })=(K(\mathbf{\upsilon }%
)\otimes I^{\otimes }(\vartheta _{\circ }^{\circ }\sqcup \vartheta
_{+}^{\circ }))(F(\mathbf{\vartheta })\otimes I^{\otimes }(\upsilon _{\circ
}^{-}\sqcup \upsilon _{\circ }^{\circ })
\end{equation*}%
where $\upsilon _{\circ }^{-}=\omega _{\circ }$,$\,\upsilon _{\circ }^{\circ
}=\upsilon $,$\,\upsilon _{+}^{\circ }=\omega ^{\circ }$,$\,\mathbf{%
\vartheta }=\mathbf{x}_{\circ }^{-},\mathbf{x}_{\circ }^{\circ },\mathbf{x}%
_{+}^{\circ }$, and $F(\mathbf{x}_{+}^{-})=V_{t(x)}^{t}$.

Moreover, the family $T_{t}$ is adapted, it can be written as the purely
quantum-stochastic integral \textup{(\ref{2onee})} of the Maassen-Meyer
kernels 
\begin{equation*}
M_{t}(\omega ^{\circ },\upsilon ,\omega _{\circ })=V_{\circ
}^{t(x_{1})}\odot L(x_{1})\odot V_{t(x_{1})}^{t(x_{2})}\odot \cdots \odot
L(x_{n})\odot V_{t(x_{n})}^{t},
\end{equation*}%
where $\omega ^{\circ }\sqcup \upsilon \sqcup \omega _{\circ
}=\{x_{1},\ldots ,x_{n}\}$, $L(\mathbf{x}_{\nu }^{\mu })=F(\mathbf{x}_{\nu
}^{\mu })-I\otimes \delta _{\nu }^{\mu }\equiv L_{\nu }^{\mu }(x)$, and the
following estimate holds: 
\begin{equation}
\Vert T_{t}\Vert _{p}(r)\leq \exp \big\{{\textstyle{\frac{1}{2}}}%
\int_{X^{t}}(\Vert L_{\circ }^{-}(x)\Vert ^{2}+\Vert L_{+}^{\circ }(x)\Vert
^{2})r(x)\mathrm{d}x\big\}.  \label{three10}
\end{equation}
\end{corollary}

In fact since $\Vert V_{s}^{t}\Vert \leq 1$, the kernels (\ref{three9}) are
bounded: 
\begin{equation*}
\Vert K_{t}(\omega ^{\circ },\upsilon ,\omega _{\circ })\leq \Vert
F_{+}^{\circ }(\omega ^{\circ })\Vert _{t}\Vert F_{\circ }^{\circ }(\upsilon
)\Vert \,\Vert F_{\circ }^{-}(\omega _{\circ })\Vert _{t},
\end{equation*}%
relative to $\Vert F(\omega )\Vert _{t}=\prod_{x\in \omega ^{t}}\Vert
F(x)\Vert $. To obtain (\ref{three10}) we use (\ref{two6}), where we put $%
\alpha _{+}^{\circ }(x)=\Vert L_{+}^{\circ }(x)\Vert $ and $\alpha _{\circ
}^{-}(x)=\Vert L_{\circ }^{-}(x)\Vert $ for $x\in X^{t}$, $\,\alpha _{\circ
}^{\circ }(x)=0=\alpha _{\circ }^{-}(x)$ for $x\in X^{t}$, $\alpha
_{+}^{\circ }(x)=0=\alpha _{\circ }^{-}(x)$ for $t(x)\geq t$, $\,\alpha
_{\circ }^{\circ }(x)=\Vert F_{\circ }^{\circ }(x)\Vert $ for $x\in X^{t}$, $%
\,\alpha _{\circ }^{\circ }(x)=1$ for $t(x)\geq t$, and $\alpha
_{+}^{-}(x)=0 $ for all $x\in X$, and now the estimate (\ref{three10})
corresponds to $\Vert T_{t}\Vert _{\mathbf{\alpha }}=1$.

\begin{example}
We construct the solution of \textup{(\ref{three2})} corresponding to the
pseudo-unitary operators $\mathbf{S}(x)=\mathbf{F}(x)\otimes \hat{\mathbf{1}}
$ with the triangular operators $\mathbf{F}(x)=e^{\mathrm{i}\mathbf{H}(x)}$,
where $\mathbf{H}^{\dagger }(x)=\mathbf{H}(x)$ are pseudo-selfadjoint
operators with the entries $H_{\nu }^{\mu }=0$ for $\mu =+$ or $\nu =-$, $%
H_{\circ }^{-}(x)^{\ast }=H_{+}^{\circ }(x)$, and $H_{\circ }^{\circ
}(x)^{\ast }=H_{\circ }^{\circ }(x)$. We assume that the local absolute
integrability condition $\Vert F_{+}^{-}\Vert _{t}^{(1)}=\int_{X^{t}}\Vert
F_{+}^{-}(x)\Vert \mathrm{d}x<\infty $ is satisfied, which leads, since $%
\mathbf{F}$ is pseudo-unitary, to 
\begin{equation*}
\Vert F_{+}^{\circ }\Vert _{t}^{(2)}=\Big(\int_{X^{t}}\Vert F_{+}^{\circ
}(x)\Vert ^{2}\mathrm{d}x\Big)^{1/2}<\infty ,\quad \Vert F_{\circ }^{-}\Vert
_{t}^{(2)}=\Big(\int_{X^{t}}\Vert F_{\circ }^{-}(x)\Vert ^{2}\mathrm{d}x\Big)%
^{1/2}<\infty
\end{equation*}%
and $\Vert F_{\circ }^{\circ }\Vert _{t}^{(\infty )}=\mathrm{ess}\sup_{x\in
X^{t}}\Vert F_{\circ }^{\circ }(x)\Vert =1$. We can now define the operators 
$T_{t}=\epsilon (K_{t})$ as the representations of the chronologically
ordered products $K_{t}(\mathbf{\omega })=F(\mathbf{x}_{1})\odot \cdots
\odot F(\mathbf{x}_{n})$ for $\sqcup _{i=1}^{n}\mathbf{x}_{i}=\mathbf{\omega 
}^{t}$, where $F(\mathbf{x}_{\nu }^{\mu })=F_{\nu }^{\mu }(x)$ are entries
in the exponential matrix $\exp \{i\mathbf{H}(x)\}$. We compute these
entries by induction finding the powers $\mathbf{H}^{0}=\mathbf{I},\,\mathbf{%
H}^{1}=\mathbf{H}$, 
\begin{equation*}
\mathbf{H}^{2}=%
\begin{bmatrix}
0, & H_{\circ }^{-}H_{\circ }^{\circ }, & H_{\circ }^{-}H_{+}^{\circ } \\ 
0, & H_{\circ }^{\circ }H_{\circ }^{\circ }, & H_{\circ }^{\circ
}H_{+}^{\circ } \\ 
0, & 0, & 0%
\end{bmatrix}%
,\,\mathbf{H}^{n+2}=%
\begin{bmatrix}
0, & H_{\circ }^{-}H_{\circ }^{\circ _{n-1}}, & H_{\circ }^{-}H_{\circ
}^{\circ _{n}}H_{+}^{\circ } \\ 
0, & H_{\circ }^{\circ _{n+2}}, & H_{\circ }^{\circ _{n-1}}H_{+}^{\circ } \\ 
0, & 0, & 0%
\end{bmatrix}%
.
\end{equation*}%
As a result we obtain $\mathbf{F}=\sum_{n=0}^{\infty }(\mathrm{i}\mathbf{H}%
^{n}/n!$ as the triangular matrix with 
\begin{align*}
F_{\nu }^{\mu }& =0,\mu >\nu ,\;\;~\;\;\;F_{-}^{-}=I=F_{+}^{+}, \\
F_{\circ }^{\circ }& =e^{\mathrm{i}H_{\circ }^{\circ
}},\;\;\;F_{+}^{-}=H_{\circ }^{-}[(e^{\mathrm{i}H_{\circ }^{\circ
}}-I_{\circ }^{\circ }-\mathrm{i}H_{\circ }^{\circ })/H_{\circ }^{\circ
}H_{\circ }^{\circ }]H_{+}^{\circ }+\mathrm{i}H_{+}^{-}, \\
F_{+}^{\circ }& =[(e^{\mathrm{i}H_{\circ }^{\circ }}-I_{\circ }^{\circ
})/H_{\circ }^{\circ }]H_{+}^{\circ },\,\;\;\;\;\;\ F_{\circ }^{-}=H_{\circ
}^{-}[(e^{\mathrm{i}H_{\circ }^{\circ }}-I_{\circ }^{\circ })/H_{\circ
}^{\circ }].
\end{align*}%
Substituting the adjoint operators $H_{\circ }^{-},H_{+}^{\circ }$ in the
form 
\begin{equation*}
H_{\circ }^{-}=F^{\ast }H_{\circ }^{\circ }-\mathrm{i}E^{\ast
},\,H_{+}^{\circ }=H_{\circ }^{\circ }F+\mathrm{i}E,
\end{equation*}%
where the operators $E(x)$ are uniquely determined by the conditions $%
H_{\circ }^{\circ }(x)E(x)=0$, we can obtain the following canonical
decomposition for the operators%
\begin{equation*}
L_{\nu }^{\mu }(x)=F_{\nu }^{\mu }(x)-I\otimes \delta _{\nu }^{\mu }I(x)
\end{equation*}
of the unitary quantum-stochastic evolution $T_{t}$: 
\begin{equation*}
\begin{pmatrix}
L_{+}^{-} & L_{\circ }^{-} \\ 
L_{+}^{\circ } & L_{\circ }^{\circ }%
\end{pmatrix}%
=%
\begin{pmatrix}
F^{\ast }L_{\circ }^{\circ }F, & F^{\ast }L_{\circ }^{\circ } \\ 
L_{\circ }^{\circ }F, & L_{\circ }^{\circ }%
\end{pmatrix}%
+%
\begin{pmatrix}
\frac{1}{2}E^{\ast }E, & E^{\ast } \\ 
-E, & 0%
\end{pmatrix}%
+%
\begin{pmatrix}
\mathrm{i}H, & 0 \\ 
0, & 0%
\end{pmatrix}%
,
\end{equation*}%
where $H=H_{+}^{-}-F^{\ast }H_{\circ }^{\circ }F$, $L_{\circ }^{\circ }=\exp
\{\mathrm{i}H_{\circ }^{\circ }\}-I_{\circ }^{\circ }$. Each of these three
tables $\mathbf{L}_{i},\,i=1,2,3$, corresponds to a pseudo-unitary matrix $%
\mathbf{F}_{i}=\mathbf{I}+\mathbf{L}_{i}$, these matrices commute, and we
have $\prod_{i=1}^{3}\mathbf{F}_{i}=\mathbf{I}+\sum_{i=1}^{3}\mathbf{L}_{i}=%
\mathbf{F}$ by the orthogonality of $\mathbf{L}_{i}$. The first matrix can
be diagonalized by means of the pseudo-unitary transform $\mathbf{F}%
_{0}^{\dagger }\mathbf{F}_{1}\mathbf{F}_{0}$ so that 
\begin{equation*}
\mathbf{F}_{0}=%
\begin{bmatrix}
1, & F^{\ast }, & -K \\ 
0, & I, & -F \\ 
0, & 0, & 1%
\end{bmatrix}%
,\quad \mathbf{F}_{0}^{\dagger }\mathbf{L}_{1}\mathbf{F}_{0}=%
\begin{bmatrix}
0, & 0, & 0 \\ 
0, & L_{\circ }^{\circ }, & 0 \\ 
0, & 0, & 0%
\end{bmatrix}%
,
\end{equation*}%
where $K=F^{\ast }F/2$. This defines the decomposition of the quantum
stochastic evolution into three types:

\begin{enumerate}
\item Poissonian quantum unitary evolution, which is given by the diagonal
matrix $\mathbf{F}$ corresponding to $H_{\nu }^{\mu }=0$ except $\mu ,\nu =0$%
: 
\begin{equation*}
T_{t}=\epsilon (K_{t})=F_{[0,t)}^{\rhd },\quad F_{[0,t)}^{\rhd }=:\exp \{%
\mathrm{i}\int_{X^{t}}H_{\circ }^{\circ }(x)\Lambda _{\circ }^{\circ }(%
\mathrm{d}x)\}:
\end{equation*}%
where $[F_{[0,t)}^{\rhd }\mathrm{h}](\varkappa )=F_{\circ }^{\circ
}(x_{1})\odot \cdots \odot F_{\circ }^{\circ }(x_{n})\mathrm{h}(\varkappa )$
for the chain $\varkappa ^{t}=\{x_{1},\ldots ,x_{n}\}$, $t(x_{1})<\cdots
<t(x_{n})$;

\item Brownian quantum unitary evolution corresponding to $H_{\circ }^{\circ
}=0=H_{+}^{-}$ and $\mathrm{i}H_{+}^{\circ \ast }=E=\mathrm{i}H_{\circ }^{-}$%
;

\item Lebesgue quantum unitary evolution corresponding to $H_{\nu }^{\mu }=0$
for all $(\mu ,\nu )\neq (-,+)$: 
\begin{equation*}
T_{t}=\epsilon (K_{t})=\int_{\mathcal{X}^{t}}\mathrm{i}^{|\varkappa |}\left(
\prod_{x\in \varkappa }^{\rightarrow }H_{+}^{-}(x)\right) \mathrm{d}%
\varkappa =\overrightarrow{\exp }\Big\{\mathrm{i}\int_{X^{t}}H_{+}^{-}(x)%
\mathrm{d}x\Big\}\otimes \hat{\mathbf{1}},
\end{equation*}%
where $\prod\limits_{x\in \varkappa }^{\rightarrow
}H_{+}^{-}(x)=H_{+}^{-}(x_{1})\ldots H_{+}^{-}(x_{n})$ for $\varkappa
=\{x_{1}<\cdots <x_{n}\}$.
\end{enumerate}
\end{example}

\end{document}